\documentclass[9pt]{article}

\usepackage{fullpage}

\usepackage{amsmath,amssymb,amscd,amsthm,amsfonts}
\usepackage{graphicx,subfigure}
\usepackage{wrapfig}
\usepackage[ocgcolorlinks,colorlinks=true,urlcolor=blue,linkcolor=blue]{hyperref}
\usepackage{pstricks}
\usepackage{dsfont}
\usepackage{pifont}
\usepackage{stmaryrd}
\usepackage{bbold}

\newtheorem{theorem}{Theorem}[section]
\newtheorem{lemma}[theorem]{Lemma}

\newtheorem*{theoremp}{Theorem}

\newtheorem{proposition}[theorem]{Proposition}
\newtheorem{prop}[theorem]{Proposition}

\newtheorem{corollary}[theorem]{Corollary}

\theoremstyle{definition}

\newtheorem{conjecture}[theorem]{Conjecture}

\newtheorem{oproblem}[theorem]{Open Problem}

\newtheoremstyle{named}{}{}{\itshape}{}{\bfseries}{.}{.5em}{\thmnote{#3}#1}
\theoremstyle{named}
\newtheorem*{namedtheorem}{}

\theoremstyle{remark}

\numberwithin{equation}{section}

\newcommand{\s}{\mathbb{S}}
\newcommand{\B}{\mathbb{B}}
\newcommand{\N}{\mathbb{N}}
\newcommand{\R}{\mathbb{R}}

\newcommand{\Q}{\mathbb{Q}}

\newcommand{\Z}{\mathbb{Z}}

\newcommand{\F}{\mathcal{F}}

\newcommand{\rank}{\operatorname{rank}}
\newcommand{\real}{{\mathbb R}}

\newcommand{\ve}[1]{\boldsymbol{#1}}
\newcommand{\aaa}{\boldsymbol{a}}
\newcommand{\bb}{\boldsymbol{b}}

\newcommand{\ee}{\boldsymbol{e}}

\newcommand{\uu}{\boldsymbol{u}}
\newcommand{\vv}{\boldsymbol{v}}
\newcommand{\xx}{\boldsymbol{x}}
\newcommand{\yy}{\boldsymbol{y}}
\newcommand{\zero}{\boldsymbol{0}}
\newcommand{\HH}{\mathcal{H}}
\newcommand{\K}{\mathsf{K}}
\newcommand{\T}{\mathsf{T}}

\newcommand{\myS}{\mathbb{S}}
\newcommand{\Ex}{\mathrm{E}}
\newcommand{\RC}{\mathcal{R}}
\newcommand{\PP}{\mathcal{P}}

\newcommand{\pth}[1]{\ensuremath{\left(#1\right)}}

\DeclareMathOperator{\conv}{conv}
\DeclareMathOperator{\inter}{int}

\DeclareMathOperator{\vol}{vol}

\DeclareMathOperator{\alt}{alt}
\DeclareMathOperator{\cd}{cd}
\DeclareMathOperator{\KG}{KG}
\DeclareMathOperator{\sd}{sd}

\DeclareMathOperator{\im}{im}

\usepackage{soul}

\newcommand{\csd}{\textsc{ColorfulSimp-Depth}}
\newcommand{\Halfspacedepth}{\textsc{Halfspace-Depth}}
\newcommand{\Simplicialdepth}{\textsc{Simplicial-Depth}}
\newcommand{\Rayshootingdepth}{\textsc{Rayshooting-Depth}}
\newcommand{\Ojadepth}{\textsc{Oja-Depth}}
\newcommand{\Regressiondepth}{\textsc{Regression-Depth}}

\newcommand{\closestpoint}{\texttt{c}}
\renewcommand{\vol}{\texttt{vol}}
\renewcommand{\P}{\mathcal{P}}

\newcommand{\D}{\mathcal{D}}
\renewcommand{\L}{\mathcal{L}}

\newcommand{\RR}{\R}
\newcommand{\eps}{\varepsilon}
\newcommand{\etal}{\emph{et al.}}

\begin{document}

\title{The discrete yet ubiquitous theorems of Carath\'eodory, Helly, Sperner, Tucker, and Tverberg}

\date{\today}

\author{
  Jesus A. De Loera\footnote{University of California, Department of Mathematics, Davis, CA 95616, USA.}  \and
  Xavier Goaoc\footnote{Universit\'e de Lorraine, CNRS, Inria, LORIA, F-54000 Nancy, France.} \and
  Fr\'ed\'eric Meunier\footnote{Universit\'e Paris Est, CERMICS, ENPC, F-77454, Marne-la-Vall\'ee, France.} \and
  Nabil Mustafa\footnote{Universit\'e Paris-Est, Laboratoire d'Informatique Gaspard-Monge, Equipe A3SI, ESIEE Paris, Noisy-le-Grand, France.}
  }

\maketitle

\begin{abstract}
We discuss five fundamental results of discrete mathematics: the lemmas 
of Sperner and Tucker from combinatorial topology, and the theorems of 
Carath\'eodory, Helly, and Tverberg from combinatorial geometry. We explore 
their connections and emphasize their broad impact in application areas such as 
data science, game theory, graph theory, mathematical optimization, computational 
geometry, etc.
\end{abstract}

\setcounter{tocdepth}{2}

\tableofcontents

\section{Introduction} \label{s:intro}

This article surveys the theory and applications of five elementary
theorems. Two of them, due to Sperner and Tucker, are from
combinatorial topology and are well-known for being the discrete
analogues of Brouwer's fixed point theorem and the Borsuk-Ulam
theorem. The other three, due to Carath\'eodory, Helly, and Tverberg,
are the pillars of combinatorial convexity. These theorems are between
fifty and one hundred years old, which is not very old as far as
mathematics go, but have already produced a closely-knit family of
results in combinatorial geometry and topology.  They have also found spectacular
applications in, among others places, mathematical optimization,
equilibrium theorems for games, graph theory, fair-division problems, 
the theory of geometric algorithms and data analysis.

The first goal of this paper is to introduce some of the many
reformulations and variations of our five theorems and explore how
these results fit together. It is convenient to split this
presentation into two parts. Sections~\ref{s:topology}
and~\ref{s:convexgeo}, discuss Sperner-Tucker and 
Carath\'eodory-Helly-Tverberg respectively. At a coarse level, the former deals with
combinatorial topology and the latter deals with combinatorial
geometry. In each case, we include a special section on algorithmic
aspects of these results relevant later for applications.

The second goal of this survey is to sample some of the many
applications of our five theorems.  There, we proceed by broad areas
and examine in Sections~\ref{s-games+fairness+independence}
to~\ref{s:datapoints}, examples from game theory and fair division,
from graph theory, from optimization, and from geometric data
analysis. Some of our illustrations are classical (e.g., Nash
equilibria, von Neumann's min-max theorem, linear programming), others
are more specialized (e.g., Dol'nikov's colorability defect or the
polynomial partitioning technique). We aim to show that our five
theorems provide simple proofs of them all. This led us to present
some new proofs, for instance for \emph{Meshulam's lemma} (Section~\ref{s:topology}) or for the
\emph{ham sandwich theorem} (Section~\ref{s:datapoints}).

The research topics that we discuss are vibrant and have already prompted a number of prior 
surveys~\cite{amentaetal2017helly,50tverbergsurvey,barany+soberonsurvey,Blagojevic+Zieglersurvey,DGKsurvey63,Eckhoff:1993survey,kalai1995survey,holmsen+wenger,3nziegler}, 
but other surveys were focused on a single one of the five theorems or did not cover applications. The important developments that we present here are from the past few years and 
emphasize both a global view and the value of geometric and topological ideas for modern applied and computational mathematics.  This research area abounds with open questions,  
all the more enticing because they can  often be stated without much technical apparatus. We made a particular effort to stress some of them.

\subsection{The five theorems at a glance}

Let us start with a classical rendition of \emph{Brouwer's fixed point
theorem}. If you stand in your favorite Parisian boulangerie holding a map of the city
in your hands, then crumple and squeeze it (without ripping it apart,  
\makebox[\linewidth][s]{mind you) and throw it to
  the ground, some point on the map must have landed right on top of its precise}

\begin{figure}[hptb]
  \begin{center}
    \includegraphics[page=2]{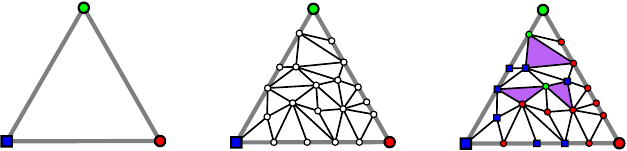}
  \end{center}
  \caption{Sperner's lemma in the plane. Start with a triangle with
    vertices colored red, green, blue (left). Subdivide it into
    smaller triangles that only meet at a common edge or a common
    vertex (center). Color every new vertex on an edge of the original
    triangle like either of the vertices of that original edge, and
    color the remaining vertices arbitrarily (right). At least one of
    the smaller triangles has vertices with pairwise distinct colors
    (cf. the shaded triangles on the right). \label{figsperner}}
\end{figure}

\smallskip\noindent
\begin{minipage}{4.1cm}
  \begin{center}
    \includegraphics[page=1,width=4 cm]{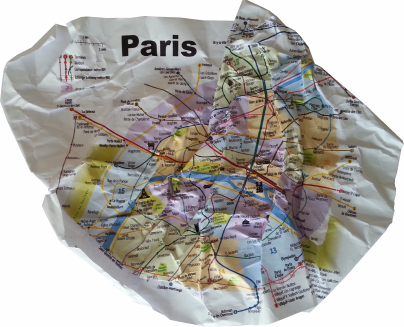}
  \end{center}
\end{minipage}
\hfill
\begin{minipage}{12cm}
  location. Brouwer's theorem follows from the classical \emph{Sperner
    lemma} on the labeling of triangulations (see
  Figure~\ref{figsperner}). Surprisingly, all that is needed to prove
  Sperner's lemma is to understand why a house with an odd number of
  openings (doors and windows) must have a room with an odd number of
  openings. This simplicity and its amazing applications attracted the
  attention of popular newspapers~\cite{NYT-sperner} and video
  sites~\cite{PBS-sperner}. Sperner's lemma is one of five theorems
  and we present it in detail in Section \ref{sec:sperner+tucker}.
\end{minipage}

\noindent
\begin{minipage}{12cm}
  ~In the game of Hex, two players take turns coloring, black and
  white, the hexagonal cells of an $11 \times 11$ diamond-shape board
  (see picture on the right); the opposite sides of the board have matching
  colors, and the player that manages to connect the two sides of his/her
  color wins (here, black wins). Since its invention by Hein in 1942,
  there has never been a draw in Hex. The fact that  there is always
  a winner happens to have a geometric explanation: for any
  triangulation of the projective plane and any two-coloring of its
  vertices, one of the color classes spans a
  \makebox[\linewidth][s]{non-contractible
    cycle~\cite{Hex-Tucker}. (To see that this implies the impossibility of a draw}
\end{minipage}
\hfill
\begin{minipage}{4cm}
  \begin{center}
    \includegraphics[page=3]{figures-final}
  \end{center}
\end{minipage}

\begin{figure}[hptb]
  \begin{center}
    \includegraphics[page=4]{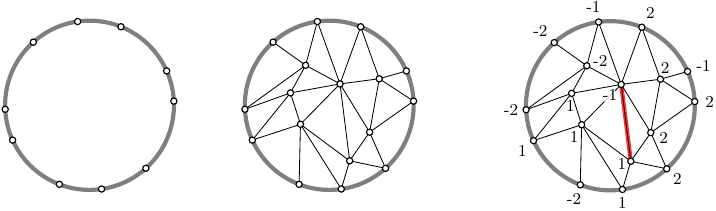}
  \end{center}
  \caption{Tucker's lemma in the plane. Start with a symmetric
    subdivision of the circle (left), and extend it into a
    triangulation of the disk (center). Label every vertex of the
    triangulation by $\{-2,-1,1,2\}$ so that antipodal points on the
    circle get opposite labels (right). There must exist an edge with
    opposite labels. \label{figtucker}}
\end{figure}

\smallskip\noindent  in Hex, take the dual of
the hexagonal cell decomposition to obtain a triangulation of the
diamond, then carefully identify the boundaries to turn that diamond
into a projective plane.)  This geometric property is equivalent to
the two-dimensional case of Tucker's lemma, whose statement is given
in the caption of Figure~\ref{figtucker}. Tucker's lemma is discussed
in Section \ref{sec:sperner+tucker}, see in particular the detailed
discussion following Proposition~\ref{tuckerlemma}.

As a matter of fact, Gale~\cite{gale+hex} proved that the game of Hex cannot end
in a draw using Brouwer's fixed point theorem and Nash~\cite{nashhex}
proved that for boards of arbitrary size, the first player has a
winning strategy. Another application of Tucker's lemma is the
\emph{ham sandwich theorem}, which says that any three finite measures
in $\R^3$ (such as a piece of bread, a slice of cheese, and a slice of
ham for an open-faced sandwich) can be simultaneously bisected by a
plane.

\noindent
\begin{minipage}{4cm}
  \begin{center}
    \includegraphics[page=5]{figures-final}
  \end{center}
\end{minipage}
\hfill
\begin{minipage}{12cm}
Let us now consider finite point sets in the plane. It turns out that
any seven points can be partitioned into three parts so that the 
triangles, segments, and points that they form have a point in common;
for example, the seven points on the left admit $\{1,7\}$, $\{2,4,6\}$,
and $\{3,5\}$ as such a partition. This is the simplest case of 
\emph{Tverberg's theorem}. Tverberg's theorem will be discuss at length in Section \ref{tv-section}. 
As the number of points grows, so does the
number of possible parts in which we can partition the points while assuring all the convex hulls
intersect: $10$ points allow four parts, $13$ points allow five parts, \ldots, and in general 
$3r-2$ points allow $r$ parts. A similar phenomenon holds in arbitrary
dimension: any set of$(r-1)(d+1)+1$ points \makebox[\linewidth][s]{can be partitioned into
$r$ parts whose convex hulls intersect. Coming back to our}
\end{minipage}

\vspace{0.1cm}\noindent  
(uncrumpled) map of Paris, consider the $302$ points that represent
the subway stations. By Tverberg's theorem, they can be partitioned
into $101$ parts, so that the corresponding 101 triangles and segments all
intersect in a common point $c$. Observe that any line passing through $c$
must leave at least $101$ subway stations on either of its (closed)
sides. The \emph{median} of  a list of real numbers separates the list in half
by sizes. The properties of point $c$ make it an acceptable 
two-dimensional generalization of the median, but now for the
set of subway stations. More generally, the \emph{centerpoint
  theorem}, which follows from Tverberg's theorem, asserts that for
any finite measure $\mu$ in the plane there is a point $c_\mu$ such
that $\mu(H) \ge \frac13 \mu(\R^2)$ for every halfplane $H$ that
contains $c_\mu$.  As we will show in Section \ref{depthsec} centerpoints are very important objects in applications
and have influenced geometry too, e.g., Tverberg himself was motivated to prove his famous theorem (which we discuss
at length later), with intention of finding an elegant proof of the centerpoint theorem. See the end of his classic paper \cite{Tverberg:1966tb}.

\bigskip

\noindent
\begin{minipage}{12cm}
\quad In the (fully supervised) \emph{classification} problem in
machine learning, one is given a data set (e.g., images), each with a
tag (e.g., indicating whether the image depicts a cat or a car), and
one is presented with new data to be tagged. A natural approach is to
map the set of data to a set of points in some geometric space and
look for a simple separation in the space (e.g., a line or a circle in the plane) that separates the
  points with different tags; for instance, perceptron neural networks --
some of the basic
\makebox[\linewidth][s]{classifiers -- look for a hyperplane best separating
the tagged point sets.}
\end{minipage}
\hfill
\begin{minipage}{4cm}
  \begin{center}
    \includegraphics[page=6]{figures-final}
  \end{center}
\end{minipage}

\vspace{0.1cm}\noindent It is easy to conclude that a separator exists
by producing an explicit hyperplane. \emph{Kirchberger's theorem}
states that it is also easy to certify when \emph{no} line separator
exists. In the plane, when no line separator exists, 
\makebox[\linewidth][s]{then there must
exist a point of one of the colors, say blue, contained in a triangle
of the opposite color, red, or}

\vspace{0.1cm}\noindent
\begin{minipage}{4cm}
  \begin{center}
    \includegraphics[page=7]{figures-final}
  \end{center}
\end{minipage}
\hfill
\begin{minipage}{12cm}
a red segment intersecting a blue segment. The set of all lines that
define a halfplane containing a given point of $\R^2$ defines a convex
cone in $\R^3$, so Kirchberger essentially reduces to an intersection
property: if a finite family of convex sets in $\R^d$ has empty
intersection, then some $d+1$ of these sets already have empty
intersection. This property is also known as \emph{Helly's theorem}
and is one of our main theorems. The curious reader may check that the
centerpoint theorem, discussed above, also follows easily from Helly's
theorem. It will be carefully studied in Section \ref{s:Helly}.
\end{minipage}

\bigskip

\noindent
\begin{minipage}{12cm}
  \quad Let us finally turn our attention to the geometry underlying
  the popular \emph{magic squares} from ancient China. A {magic
    square} is a $n \times n$ square grid of non-negative real numbers
  such that the entries along any row, column, and diagonals, all add
  up to the same value. Look at the four $3 \times 3$ examples on the
  right. It turns out that any $3 \times 3$ magic square can be
  written as a linear combination, with non-negative coefficients, of
  only three of these four magic squares! In fact, for any $n$ there
  exists a finite set $X_n$ of $n\times n$ magic squares
  \parbox[s]{\linewidth}{such that any other $n \times n$ magic square can
    be written using only $(n-1)^2-1$}
\end{minipage}
\hfill
\begin{minipage}{5cm}
  \begin{center}
\begin{tabular}{|c|c|c|}
\hline 0&2&1\\ \hline 2&1&0\\ \hline 1&0&2\\ \hline
\end{tabular}
\quad
\begin{tabular}{|c|c|c|}
\hline 2&0&1\\ \hline 0&1&2\\ \hline 1&2&0\\ \hline
\end{tabular}
\\
\medskip
\begin{tabular}{|c|c|c|}
\hline 1&2&0\\ \hline 0&1&2\\ \hline 2&0&1\\ \hline
\end{tabular}
\quad
\begin{tabular}{|c|c|c|}
\hline 1&0&2\\ \hline 2&1&0\\ \hline 0&2&1\\ \hline
\end{tabular}
  \end{center}
\end{minipage}

\vspace{0.1cm}\noindent elements of $X_n$. This last statement follows from
\emph{Carath\'eodory's theorem}, which we will study carefully in
Section \ref{caratheodory-thms} : any vector in a cone in $\R^d$ is a
non-negative linear combination of extreme rays of the cone, and only
dimension many are used in the linear combination. Indeed, the set of
$n \times n$ magic squares forms a polyhedral cone in a vector space
of dimension $(n-1)^2-1$.  It may come as a surprise that no one knows
what $X_n$ is for all $n \ge 6$. See \cite{Maya+Yo+Raymond}.

\smallskip\noindent
\begin{minipage}{4cm}
  \begin{center}
    \includegraphics[page=8]{figures-final}
  \end{center}
\end{minipage}
\hfill
\begin{minipage}{12cm}
\quad A colorful generalization of Carath\'eodory's theorem asserts that if
three polygons in the plane, one with red vertices, one with green
vertices, and one with blue vertices, all contain a given point $p \in
\R^2$, then there exists a colorful triangle, using a vertex of each
color, that also contains $p$. This implies that for the centerpoint
$c$ that we constructed earlier for the Parisian subway stations from
Tverberg's theorem, at least $\binom{101}{3}$ of the triangles spanned
by the subway stations contain $c$.  In fact, 
\makebox[\linewidth][s]{there is a
quantitatively stronger statement given by the \emph{first selection}}
\end{minipage}

\smallskip
\noindent
\emph{lemma}.  It
states that for any set of $n$ points in the plane, there exists apoint covered by at least $\frac29 \binom{n}3$ of the triangles they
span. We will see more about this topic in Section \ref{depthsec}.

\subsection{Notation and preliminaries}

We collect in this subsection notation, terminology, and general basic
background on combinatorics, geometry, and topology that will be used
in the rest of this survey.  The advanced reader may want to skip or
move quickly through this section. For a more thorough introduction to
the topics listed here, we recommend the classical books and textbooks
in combinatorial convexity~\cite{Bar2002,handbookofconvexgeo,gruber2007convex,Mbook} as well
as~\cite[$\mathsection 5.3$]{Sch03}. For topological combinatorics and combinatorial aspects
of algebraic topology see ~\cite{Longueville-book,matousek2003using,munkres1984elements}.

Given $n \in \N$, we write $[n]$ to denote the set $\{1,2,\ldots,
n\}$. If $X$ is a set and $k \in \N$, we write $\binom{X}{k}$ for the
set of $k$-element subsets of $X$. The notation $\tilde{O}(\cdot)$
denotes asymptotic notation where we ignore poly-logarithmic factors:
$f(n) = \widetilde{O}(g(n))$ if there exists $k \in \N$ such that
$f(n) = O(g(n) \log^k g(n))$.

We denote by $(\ve e_1, \ve e_2, \ldots, \ve e_d)$ the
orthonormal frame of $\R^d$. Given two vectors $\xx$ and $\yy$ in
$\R^d$, we write $\xx \le \yy$ to mean that $x_i \le y_i$ for $i =1,
2, \ldots, d$. We write $\B^d = \{ \ve x \in \R^d \colon \sum_{i=1}^d
x_i^2 \leq 1 \}$ for the unit ball in $\R^d$ and $\s^{d} = \{ \ve x
\in \R^{d+1} \colon \sum_{i=1}^{d+1} x_i^2 = 1 \}$ for the unit sphere
in $\R^{d+1}$.
 
\subsubsection{Polytopes, simplices, polyhedra, cones}

Let $A \subseteq \R^d$ be a set. The {\em convex hull} of $A$, denoted
by $\conv(A)$, is the intersection of all convex sets containing
$A$. In other words, $\conv(A)$ is the smallest convex set containing $A$. It
is well-known that
\[ \conv(A)=\left\{\sum_{i=1}^{n} \gamma_i \ve a_i\colon n\in \N,\; \ve a_i \in A,\; 
\gamma_i \geq 0, \ \text{and} \ \gamma_1+\cdots+\gamma_n=1\right\}.\]
A {\em polytope} is the convex hull of a finite set of points in
$\R^d$. Here are a few examples. The convex hull of affinely
independent points is a \emph{simplex}; the {\em standard}
$k$-dimensional simplex $\Delta_{k}$ is $\conv(\{\ve
e_1,\dots,\ve e_{k+1}\})$. With $e_i$ is a the $i$-standard unit vector. 
The convex hull of $\ve e_1,-\ve e_1,\ve
e_2,-\ve e_2,$ $\dots,\ve e_k,-\ve e_k$ is the $k$-dimensional
\emph{cross-polytope}. The convex hull of all vectors with $0,1$
entries is the $d$-dimensional \emph{hypercube}. A \emph{face} of a
polytope is its intersection with a hyperplane that avoids its
relative interior. Faces of dimension $0$ are \emph{vertices} and
inclusion-maximal faces are \emph{facets}. A face of a polytope
(resp. simplex) is also a polytope (resp. simplex). There is a face of
dimension $-1$, the empty set. 

A {\em polyhedron} is the intersection of finitely many halfspaces in
$\R^d$. In particular, any polyhedron can be represented as $\{\xx \in
\R^d\colon A \xx \le \ve b\}$ where $A$ is an $n \times d$ matrix and $\ve b
\in \R^n$. A \emph{polyhedral cone} is a polyhedron that is closed
under addition and scaling by a positive constant. In particular, any
polyhedral cone writes as $\{\xx \in \R^d: A \xx \ge \ve 0\}$.

The \emph{polar} of a point $\ve x \in \R^d\setminus\{\ve 0\}$ is the
halfspace $\{\ve y \in \R^d\colon \ve x \cdot \ve y \le 1\}$ (and
vice-versa) and the \emph{dual} of $\ve x$ is the hyperplane $\{\ve y
\in \real^d\colon \ve x \cdot \ve y = 1\}$ bounding its polar; we speak of
polarity or duality \emph{relative to $\ve o$} to represent the
polarity or duality after the coordinate system has been translated to
have its origin in $\ve o$.

The \emph{Weyl-Minkowski theorem} asserts that the polytopes are
exactly the bounded polyhedra. (This can be proven using polarity
arguments.) Polytopes can thus be represented both as convex hulls of
finitely many points, or as intersections of finitely many
halfspaces. To see that these viewpoints are complementary, we invite
the reader to prove that projections or intersections of polytopes are
polytopes. The Weyl-Minkowski theorem similarly implies that any
polyhedral cone also writes as the set of convex combinations of a finite set of rays
emanating from the origin. By a theorem of Motzkin et
al.~\cite{motzkin-double}, every polyhedron $P$ decomposes into the
\emph{Minkowski sum} of a polytope $Q$ and a polyhedral cone $C$: $P =
\{\xx + \yy: \xx \in Q, \yy \in C\}$.

A polyhedron is \emph{pointed} if the largest affine subspace it
contains is zero-dimensio\-nal. Any polyhedron can be decomposed
into the Minkowski sum of a pointed polyhedron and a vector space.

\subsubsection{Simplicial complexes}

A \emph{geometric simplicial complex} is a family of simplices with two properties:
the intersection of any two distinct simplices
is a face of both of them; and it contains all the faces of  every member of the family. 
Simplicial complexes may also be considered abstractly by only retaining which
sets of vertices span a simplex. Formally, an \emph{abstract
  simplicial complex} $\K$ is a family of finite subsets (the
\emph{faces}) of some ground set (the \emph{vertices}) that is closed
under taking subsets: if $\sigma \in\K$ and $\tau \subseteq \sigma$,
then $\tau \in\K$. We write $\sigma_n$ for the \emph{abstract
  $n$-dimensional simplex} consisting of all subsets of $[n+1]$. A
\emph{facet} of a simplicial complex is an inclusion-maximal face.

Let us stress that the meaning of the word \emph{face} (or
\emph{facet}) depends on the context and can denote a polytope (for
polytopes), a simplex (for simplices), or a set of vertices (for
abstract simplicial complexes). In particular, we consider $\Delta_n$
to be a polytope, so it has $n+1$ faces of maximal dimension $n-1$; in
contrast, $\sigma_n$ has a single face of maximal dimension $n$.

Given a geometric simplicial complex $\K$, we let $|\K|$ denote the
underlying topological space, that is $|\K| = \bigcup_{\sigma \in \K}
\sigma$. If $\mathsf{L}$ is the abstract simplicial complex obtained from
$\K$, we say that $\K$ is a \emph{geometric realization} of $\mathsf{L}$ and
put $|\mathsf{L}| = |\K|$. (The reader can check that all this is well-defined
up to homeomorphism.) A \emph{triangulation} of a topological space $X$ is a
(geometric or abstract) simplicial complex whose underlying
topological space is homeomorphic to $X$.

\subsubsection{Homology}

We will use some basic notions of homology, mostly simplicial homology
over $\Z$ or $\Z_q = \Z/q\Z$. To allow readers unacquainted with
homology to appreciate at least our simplest examples, we recall here
the basic definitions. An important idea in homology theory is that
topological spaces can be studied by associating to them some groups,
called homology groups. These groups can be defined geometrically (in
singular homology) or combinatorially (in simplicial homology) from a
triangulation of the space. In the cases that we consider, these
approaches produce isomorphic groups, and we mostly work with
simplicial homology.

Given a simplicial complex $\K$, we denote by $C_i(\K,\Z_2)$ the set of
finite formal sums of $i$-dimensional faces of~$\K$, and
$C_{\bullet}(\K,\Z_2) = \bigoplus_{i=0}^{\infty} C_i(\K,\Z_2)$ is the
\emph{chain complex} of $\K$. The map that sends every $i$-face of
$\K$ to the formal sum, with coefficients in $\Z_2$, of its
$(i-1)$-faces extends linearly, over $\Z_2$, to a map $\partial_i:
C_i(\K,\Z_2) \to C_{i-1}(\K,\Z_2)$. Notice that $C_{\bullet}(\K,\Z_2)$
has an additive group structure and that the sum of the $\partial_i$
is a morphism from $C_{\bullet}(\K,\Z_2)$ into itself; this morphism
is called the \emph{boundary map} of $C_{\bullet}(\K,\Z_2)$. Note that
with $\Z_2$-coefficients finite formal sums are simply subsets and the
boundary operator maps to proper subsets that appear an odd number of
times. It turns out that $\partial_{i-1} \circ \partial_i = 0$, so we
can define the \emph{$i$-th homology group} of $\K$ over $\Z_2$
coefficients as the quotient group $H_i(\K,\Z_2) = \ker \partial_i /
\im \partial_{i+1}$.

Intuitively, the rank of $H_i(\K,\Z_2)$ relates to the number of
independent holes of dimension $i$ in $|\K|$; for example, the rank of $H_0(\K,\Z_2)$ counts the number of connected components of $|\K|$. In
particular, if $\K$ is a single vertex then all its homology groups
are trivial except the $0$-th one. The \emph{reduced homology groups}
$\widetilde{H}_i$ modify slightly $H_0$ so that it is trivial for
connected sets; in dimension $i \ge 1$, homology groups and reduced
homology groups coincide and $\widetilde{H}_{-1}$ is defined as $0$
for nonempty complexes. Going from $\Z_2$ to other coefficient groups
involves only one technical complication: the definition of
$\partial_i$ involves some sign book-keeping, so as to ensure that
$\partial_{i-1} \circ \partial_i = 0$. See~\cite{munkres1984elements}
for details.

\section{Combinatorial topology}
\label{s:topology}

From combinatorial topology, we will focus on two results about 
labeled or colored triangulations of simplicial complexes: \emph{Sperner's lemma and Tucker's lemma}. 
The importance of these two lemmas owes much to their
special position at the crossroads of topology and combinatorics.
As we alluded in the introduction, Sperner's and Tucker's lemmas are the combinatorial equivalent 
versions to the famous topological theorems of Brouwer and Borsuk-Ulam, respectively.
Their combinatorial nature makes them particularly well-suited for computations and applications too.

Combinatorial structures have been used in algebraic topology since
Poincar\'e's foundational \emph{analysis situs}, so it is not
surprising that some topological questions may be studied by
combinatorial methods. The lemmas of Sperner and Tucker are well-known
for offering an elementary access, via labelings of combinatorial
structures, to important results in topology such as the theorems of
Brouwer and Borsuk-Ulam.

It is perhaps less obvious that some combinatorial problems may be
studied by topological methods. A seminal example of topological
methods applied to combinatorics was the use by L. Lov\'asz of the
Borsuk-Ulam theorem to settle a conjecture of Kneser on the chromatic
number of certain graphs (see Section~\ref{s-graphs}). His paper opened up the
application of topological methods in combinatorics that are now
common tools. These techniques appear in several
books~\cite{Longueville-book,matousek2003using} and
surveys~\cite{Bjorner-survey,Karasev-survey}. In many cases the
topological methods hinge on the theorems of Brouwer or Borsuk-Ulam;
as we discuss in the application sections, on several occasions the
topological machinery can be made implicit, and the combinatorial
question settled directly by the lemmas of Sperner or Tucker.

\subsection{Sperner and Tucker} \label{sec:sperner+tucker}

A \emph{labeling} of a simplicial complex $\K$ by a set $S$ is a map
from the vertices of $\K$ to $S$; the \emph{label} of a vertex is its
image. Sperner's lemma gives a sufficient condition for the existence
of a \emph{fully-labeled} facet, that is a facet whose vertices have
pairwise distinct labels. (Sometimes the labels are called
\emph{colors} and fully-labeled faces are called \emph{colorful}; we
will avoid this terminology in this paper to avoid confusion with the
colorful theorems in convex geometry that we discuss in
Sections~\ref{s:convexgeo} and~\ref{s-graphs}.)

\begin{namedtheorem}[Sperner's lemma]
  Assume that the vertices of a finite triangulation $\T$ of a
  simplex~$\Delta$ are labeled so that any vertex lying in a face of
  $\Delta$ has the same label as one of the vertices of that face. If
  the vertices of $\Delta$ are given pairwise distinct labels, then the
  number of facets of $\T$ whose vertices have pairwise distinct
  labels is odd. 
\end{namedtheorem}

\bigskip

\begin{figure}[h]
  \begin{center}
    \includegraphics[page=9]{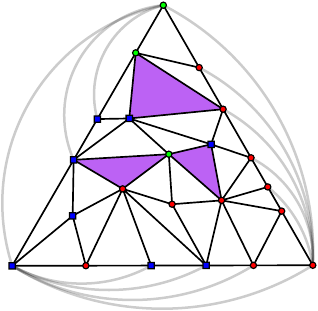}
  \end{center}
  \caption{A Sperner labeling of a triangulation of $\Delta_2$
    illustrated by colors.  The fully-labeled triangles are shown
    shaded. The gray edges augment the triangulation of $\Delta_2$
    into a triangulation of~$\s^2$. \label{figspernerproof}}
\end{figure}

\noindent
We call a labeling that satisfies the assumptions of Sperner's lemma a
\emph{Sperner labeling}. A more general version holds for
\emph{pseudomanifolds}, i.e., for pure $d$-dimensional simplicial
complexes where every face of dimension $d-1$ is contained in at most
two facets. (Recall that a simplicial complex is \emph{pure} if all
its inclusion-maximal faces have the same dimension.) In particular,
any triangulation of a $d$-dimensional manifold is a pseudomanifold of
dimension $d$. The {\em boundary} of a pseudomanifold is the
subcomplex generated by its $(d-1)$-dimensional simplices that are
faces of exactly one $d$-dimensional simplex.

\begin{proposition}
  \label{prop:sperner_even}
  Any labeling by $[d+1]$ of a $d$-dimensional pseudomanifold without boundary has an
  even number of fully-labeled facets.
\end{proposition}

\noindent
Proposition~\ref{prop:sperner_even} follows from a simple parity
argument.  Consider the graph where the nodes are the facets and where
the edges connect pairs of facets that share a $(d-1)$-dimensional
face whose vertices use every label in $[d]$. This graph has only
nodes of degree $0$, $1$, or~$2$, so it consists of vertex-disjoint
cycles and paths. The nodes of degree~$1$ are exactly the
fully-labeled facets and there are evenly many of them (twice the
number of paths). Coming back to the remark in introduction: this is
where it is useful to understand why a house with an odd number of
openings has a room with an odd number of openings.

Clearly Sperner's lemma follows from Proposition~\ref{prop:sperner_even}.
For this observe that any Sperner
labeling of a $d$-dimensional simplex $\Delta_d$ extends into a
triangulation of $\s^d$ where (i) the outer vertices of $\Delta_d$
form a fully-labeled facet, and (ii) no other added facet is
fully-labeled (as illustrated in Figure~\ref{figspernerproof}). Knowing a
vertex of degree one allows easily, by path following, to find another
one; we come back in Section~\ref{fixed-point-computing} on
algorithmic applications of this idea. Other arguments can be used 
to prove Sperner's lemma \cite{McLennan+Tourky2008}.

\bigskip

\begin{figure}
  \begin{center}
    \includegraphics[page=10]{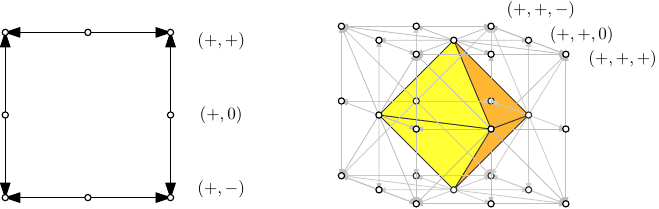}
    \caption{Geometric interpretation of the dominance graph for
      $\preceq$ in the octahedral Tucker lemma for $n=2$ (left) and $n=3$ (right).\label{f:geometry-octahedral-Tucker}}
  \end{center}
\end{figure}

Now we will state the octahedral Tucker lemma. This is a rather
streamlined version of Tucker's lemma that already suffices for all
our applications. Given vectors of signs $\ve x, \ve y
\in\{+,-,0\}^n$, we write $\ve x\preceq\ve y$ if every nonzero
coordinate of~$\ve x$ is the same as in $\ve y$. We let $\xx^+$ denote
the set of indices $i$ such that $x_i = +$, and $\xx^-$ similarly. In
particular, $\xx \preceq \yy$ if and only if $\xx^+ \subseteq \yy^+$
and $\xx^- \subseteq \yy^-$. Note that each vector of signs uniquely
identifies a coordinate (sub)-orthant, and that the order $\preceq$
indicates containment.  There is an interpretation of $(\{+,-,0\}^n,
\preceq)$ as a simplicial complex illustrated in
Figure~\ref{f:geometry-octahedral-Tucker}.  By $\pm a$ we mean a
choice of one of the two scalars $-a$ or $a$.

\begin{namedtheorem}[Octahedral Tucker lemma] \label{tuckerOlemma}
  Let $\lambda:\{+,-,0\}^n\setminus\{\zero\}\rightarrow
  \{\pm 1, \pm 2, \ldots, \pm m\}$ be such that
  $\lambda(-\ve x)=-\lambda(\ve x)$ for all $\xx$. If
  $\lambda(\xx)+\lambda(\yy)\neq 0$ for all $\xx\preceq\yy$, then $n
  \le m$.
\end{namedtheorem}

\noindent
The octahedral Tucker lemma was apparently first stated
explicitly by Ziegler~\cite[Lemma~4.1]{Z-GK}, following its implicit
use by Matou\v{s}ek~\cite{Matousek:2004hm} in his proof of the lower
bound on the chromatic number of Kneser graphs from Tucker's lemma 
(see Section~\ref{s-chromatic}).
Several classical proofs of Proposition~\ref{tuckerlemma} can be found
in Matou\v{s}ek's book~\cite{matousek2003using}. As we explain in
Section~\ref{subsec:continuousversionspernertucker}, the lemmas of Sperner and Tucker are
indeed ``topological'' in that they essentially state that certain
chain maps, namely those induced by the labeling maps, are non-trivial
in simplicial homology with coefficients over $\Z_2$.

\begin{figure}[h]
  \centering
      \includegraphics[page=11]{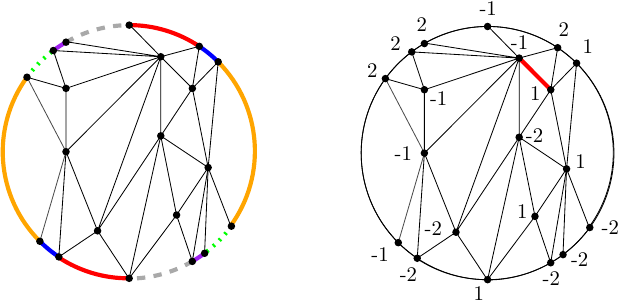}
  \caption{An illustration of Tucker's lemma: a triangulation of $\B^2$
    that induces a symmetric triangulation of $\s^1$. On the left the antipodal simplices on 
    the boundary are painted with the same color.  On the right side, a labeling of the vertices 
    that is antipodal on the boundary.}  \label{fig:anti_sym_tri}
\end{figure}

A more common version of Tucker's lemma deals with triangulations of
a ball instead of an octahedron. Tucker's lemma gives a lower
bound on the number of distinct labels used by labelings that avoid
certain local patterns. We say that a triangulation $\T$ of $\B^d$
induces a \emph{symmetric triangulation} of $\s^{d-1}$ if its boundary
$\partial\T$ forms a centrally-symmetric triangulation
of~$\s^{d-1}$. A labeling of a symmetric triangulation of $\s^d$ by
integers is \emph{antipodal} if antipodal vertices have opposite
labels.

\begin{proposition}[Tucker lemma] \label{tuckerlemma}
  Let $\T$ be a triangulation of $\B^d$ that induces a symmetric
  triangulation of $\s^{d-1}$. Let $\lambda$ be a labeling of the
  vertices of $\T$ by $\{\pm 1,\ldots,\pm d\}$. If $\lambda(-\ve
  v)=-\lambda(\ve v)$ for every vertex $\ve v$ of $\partial\T$,
  then there exists an edge $\ve u \ve v$ in $\T$ with $\lambda(\ve
  u)+\lambda(\ve v)=0$.
\end{proposition}

\noindent
(See Figure~\ref{fig:anti_sym_tri} for an illustration.) The octahedral
version  is obtained by applying
Proposition~\ref{tuckerlemma} to the barycentric subdivision $\T$ of
the $n$-dimensional \emph{cross-polytope} $\Diamond^n = \conv \{\pm
\ve e_i\}_{1 \le i \le
  n}$~\cite[Theorem~2.3.2]{matousek2003using}. Indeed, consider a map
$\lambda: \{+,-,0\}^n\setminus\{\ve 0\} \to \{\pm 1, \ldots, \pm
m\}$ such that $\lambda(-\ve x)=-\lambda(\ve x)$ for all $\xx$ and
$\lambda(\xx)+\lambda(\yy)\neq 0$ for all $\xx\preceq\yy$. Any $\ve x
\in \{+,-,0,\}^n\setminus\{\ve 0\}$ can be interpreted as the vertex in
$\T$ corresponding to the face $\conv(\{\ve e_i: i\in \ve x^+\} \cup
\{-\ve e_i: i \in \ve x^-\})$ of $\Diamond^n$. The edges in $\T$ connect
pairs $\ve x$ and $\ve y$ such that $\xx\preceq\yy$. Defining
$\lambda(\ve 0) = m+1$, we get a labeling of all vertices of $\T$
satisfying $\lambda(-\ve v)=-\lambda(\ve v)$ for every vertex $\ve v$
of $\partial\T$ and having no edge $\ve u \ve v$ in $\T$ with
$\lambda(\ve u)+\lambda(\ve v)=0$. By Proposition~\ref{tuckerlemma},
we must have $m+1 > n$.

As stated in the introduction, Tucker's lemma is equivalent to
  the fact that for any triangulation of the projective plane and any
  two-coloring of its vertices, one of the color classes spans a
  non-contractible cycle. Indeed, such a two-coloring can be seen as a
  two-coloring of a triangulation of the disk, with symmetric vertices
  of the boundary getting identical colors.  If all monochromatic
  cycles were contractible, we could easily choose a sign for each
  vertex and get a labeling that would contradict Tucker's lemma.  The
  reverse implication is also easy and left to the reader.

\bigskip

Consider a triangulation $\T$ as in Proposition~\ref{tuckerlemma} and
a labeling $\lambda$ of its vertices by $\{\pm 1, \ldots, \pm m\}$; a
$k$-dimensional face of $\T$ is \emph{alternating} if its vertices can
be indexed $\ve v_{i_0}, \ldots \ve v_{i_k}$ so that $0 < \lambda(\ve
v_{i_0}) < -\lambda(\ve v_{i_1}) < \cdots < (-1)^{k}\lambda(\ve
v_{i_{k}})$ or  if $0 > \lambda(\ve v_{i_0}) > -\lambda(\ve v_{i_1}) >
\cdots > (-1)^{k}\lambda(\ve v_{i_{k}})$. In the first case we call the simplex
\emph{positively alternating} and in the second case \emph{negatively alternating}.
A lemma due to Fan~\cite{kyfan1952} generalizes Tucker's lemma in terms of 
a parity counting of alternating simplices.

\begin{theorem}[Fan's lemma]
  \label{our-fan}
  Let $\T$ be a triangulation of $\B^d$ that induces a symmetric
  triangulation of $\s^{d-1}$. Let $\lambda$ be a labeling of the
  vertices of $\T$ by $\{\pm 1,\ldots,\pm m\}$ such that $\lambda(-\ve
  v)=-\lambda(\ve v)$ for every vertex $\ve v$ of $\partial\T$.
  If no two adjacent vertices have opposite labels, then $\T$ has an
  odd number of alternating facets.
\end{theorem}

\noindent
Fan's lemma readily implies Tucker's lemma since the existence of an
alternating $d$-dimensional face implies $m \ge d+1$. Going in the
other direction, it was only recently observed by
Alishahi~\cite{alishahi2017} that an existential version of Fan's
lemma is easily derived from Tucker's lemma.

Let us illustrate that, surprisingly, Fan's lemma can be easier to
prove than Tucker's lemma: this is one example where a stronger
hypothesis facilitates induction. We give an inductive proof for a
\emph{flag of hemispheres}, i.e., a triangulation $\T$ of $\B^d$ such
that the restriction of $\T$ on $H_i^+$ and on $H_i^-$ triangulates
each of them, where $H_i^+$ and $H_i^-$ are the $i$-dimensional
hemispheres
$$H_i^+=\{\ve x\in\s^{d-1}\colon x_{i+1}\geq 0,
x_{i+2}=\cdots=x_d=0\},$$ $$H_i^-=\{\ve x\in\s^{d-1}\colon x_{i+1}\leq
0, x_{i+2}=\cdots=x_d=0\}.$$
(Prescott and Su~\cite{prescott+su} gave another combinatorial
constructive proof for this special case.) Consider the graph whose
nodes are the facets of $\T$ and whose edges connect pairs of facets
that share a $(d-1)$-dimensional face that is positively
alternating. We augment this graph with an extra node $s$ and add
edges connecting $s$ to all facets of $\T$ that have a
$(d-1)$-dimensional positively alternating face on $\partial\T$; in a
sense, $s$ represents the ``outer facet''. Apart from $s$, all nodes
have degree $0$, $1$, or $2$. The nodes of degree $1$ are exactly the
$d$-dimensional alternating facets. The triangulation $\T$ refines
$H_{d-1}^+$, which is homeomorphic to $\B^{d-1}$. So, by induction in
dimension $d-1$, the number of $(d-1)$-dimensional alternating faces
of $\partial\T$ in $H_{d-1}^+$ is odd; the same holds for the
$(d-1)$-dimensional faces of $\partial\T$ in $H_{d-1}^-$. The symmetry
of $\partial\T$, and that of the labeling, imply that the degree of
$s$ is odd; it follows that there is an odd number of $d$-dimensional
alternating facets in $\T$.

The above elementary proof requires that the triangulation restricts
nicely to lower-dimensional spheres to allow induction.  Two proofs of
Theorem~\ref{our-fan} can be found in the literature, both for an equivalent version with
a sphere instead of a ball (later we will explain this equivalence for our own proof of Theorem \ref{our-fan}). 
On the one hand, \v{Z}ivaljevi\'c~\cite{MR2728499} observed that the labeling is
essentially a classifying map that is unique up to $\Z_2$-homotopy, so
the number of alternating facets (mod $2$) reformulates as the cap
product of a certain cohomology class with a certain homology
class. On the other hand, Musin~\cite{Musin:2012tq} builds a
simplicial $\Z_2$-map from the triangulation to a $d$-dimensional
polytope for which the following holds: a simplex is alternating if
and only if its image by this simplicial map contains $0$ in its
convex hull; a degree argument allows then to conclude. It turns out
that Alishahi's idea to derive Fan's lemma from Tucker's lemma leads
to a short proof of Theorem~\ref{our-fan}, also based on a degree
argument. Before we spell out this (original) proof let us stress that
the following question remains:

\begin{oproblem}
  Give a direct combinatorial proof of Fan's lemma (as stated in Theorem
  \ref{our-fan}) and of Tucker's lemma (Proposition~\ref{tuckerlemma})
  valid for \emph{any} centrally symmetric triangulation.
\end{oproblem}

Let us now prove Theorem~\ref{our-fan} for an arbitrary triangulation
$\T$ of $\B^d$. Let $\lambda$ be a labeling of the vertices of $\T$
satisfying the conditions of Theorem~\ref{our-fan}. We first turn $\T$
into a triangulation $\T'$ of $\s^d$ by gluing two antipodal copies of
$\B^d$, each triangulated by $\T$; notice that the number of
\emph{positively alternating} facets of $\T'$ equals the number of
\emph{alternating} facets of $\T$ (both positive and negative ones),
since the negatively alternating facets in one copy of $\T$ become
positively alternating in the other copy. We define next a labeling
$\mu$ of the vertices of $\sd\T'$ by $\{\pm 1,\ldots,\pm (d+1)\}$: a
vertex $v$ of $\sd\T'$ corresponds to a simplex $\tau_v$ of $\T'$, and
we set the absolute value of $\mu(v)$ to be the number of vertices of
the largest alternating face of $\tau_v$, and its sign according to
alternation. This sign is defined uniquely since there cannot be
maximal alternating faces of both types in $\tau$ (this can be checked
using the fact that no adjacent vertices in $\T$ can have opposite
labels). Now, a crucial observation is that if $\sigma$ is an
alternating facet (for $\lambda$), then $\sd\sigma$ contains exactly
one alternating facet (for $\mu$) of the same type; and if $\sigma$ is
not an alternating facet, then $\sd\sigma$ contains no alternating
facet. At this point, to establish Theorem~\ref{our-fan}, it suffices
thus to prove that $\sd\T'$ contains an odd number of positively
alternating facets. This fact follows from basic degree theory. The
assumptions on $\lambda$ guarantee that $\mu$ induces an antipodal
simplicial map from $\sd\T'$ to $\partial\Diamond^{d+1}$, the boundary
of the $(d+1)$-dimensional cross-polytope, whose vertices are
identified with the elements in $\{\pm 1,\ldots,\pm (d+1)\}$. Any
antipodal self-map of $\s^d$ is of odd degree~\cite[Theorem
  4.3.32]{MR1430097}. Thus, denoting by $t\in C_d(\sd\T',\Z_2)$ the
formal sum of all facets of $\sd\T'$ and by $z\in
C_d(\partial\Diamond^{d+1},\Z_2)$ the formal sum of all facets of
$\partial\Diamond^{d+1}$, we must have $\mu_{\sharp}(t)=z$. The only
simplices that are mapped to the simplex
$\{+1,-2,\ldots,(-1)^d(d+1)\}$ by $\mu_{\sharp}$ are the positively
alternating ones, so there are an odd number of them.

\subsection{Continuous versions} \label{subsec:continuousversionspernertucker} 

One of the most famous theorems about fixed points is due to the Dutch
mathematician L.\,E.\,J.\ Brouwer and states that any continuous
function from a ball into itself has a fixed point. Brouwer's theorem
is often seen as a continuous version of Sperner's lemma (without the
oddness assertion): they can be deduced easily from one another.

Let us sketch how Brouwer's theorem follows from Sperner's lemma (we discuss the other direction in
Section~\ref{s:s-t-generalizations}). Contrary to Brouwer's original proof, which says nothing about 
how to find the fixed point or a good approximation of it, this proof has computational implications (see
Section~\ref{fixed-point-computing}). Without loss of generality we
take the $d$-dimensional standard simplex $\Delta_d \subset {\mathbb
  R}^{d+1}$ as a realization of a ball (it is easy to set up a
homeomorphism). Then we triangulate the simplex $\Delta_d$ and design
a labeling of that triangulation tailored to the continuous
function~$f$ under consideration.  Specifically, we associate to a
vertex $\ve v$ of the triangulation the label $i$ if the $i$-th
barycentric coordinate of $\ve v$ is larger than the $i$-th
barycentric coordinate of its image $f(\ve v)$. (So, intuitively, if $\ve v$
is labeled $i$, then $f$ moves $\ve v$ outwards the $i$-th vertex of the
standard simplex.) Note that, unless the vertex $\ve v$ is a fixed point, there must be at least one such index.
If there are several such indices, simply make an arbitrary choice
among them. Now, the labeling we provided satisfies the assumptions of
Sperner's lemma, so we can find a fully-labeled simplex of $\T$ such that the $i$-th 
barycentric coordinate of the vertex labeled $i$ is decreased by $f$. Re-triangulate $\Delta_d$ again and again
adding more and more points in such a way that the maximum diameter of
the simplices appearing in the triangulation goes to zero in the
limit. At each step we find a fully-labeled simplex. The barycenters
of all such simplices will produce an infinite sequence of points and,
since it is a bounded sequence, it contains a convergent subsequence.
Let $\xx^*$ be the limit of this subsequence. Since the map $f$ is
continuous, the $i$-th barycentric coordinate of $\xx^*$ is larger or
equal than the $i$-th barycentric coordinate of $f(\xx^*)$ for every
$i$ and therefore $\xx^*$ is a fixed point of the map.

\bigskip

\begin{figure}
  \begin{center}
    \includegraphics[page=12]{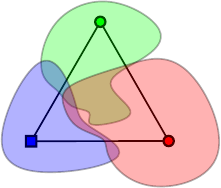}
    \caption{\small An illustration of the KKM lemma in two dimensions.}
    \label{f:KKM}
  \end{center}
\end{figure}

The Knaster-Kuratowski-Mazurkiewicz theorem, also known as the KKM
theorem, is a classical consequence of Sperner's lemma or Brouwer's
theorem. It is used for instance in game theory or for the study of
variational inequalities. Consider the $d$-dimensional simplex 
$\Delta_d =\conv \{\ve e_i \colon 1 \le i \le d+1\}$ and $d+1$ closed sets
$C_1, C_2, \ldots, C_{d+1}$ in $\R^d$. This theorem, illustrated in
Figure~\ref{f:KKM}, ensures that if for every $I \subseteq [d+1]$ the
face $\conv\{\ve e_i \colon i \in I\}$ of $\Delta_d$ is covered by $\bigcup_{i
  \in I} C_i$, then $\bigcap_{i=1}^{d+1}C_i$ is nonempty. This
statement is somehow reminiscent of Helly's theorem. 
A corollary of the KKM theorem can actually be used to prove it, see
Section~\ref{s:Helly}. This corollary states that if $d+1$ closed sets
$C_1, C_2, \ldots, C_{d+1}$ are such that each of them contains a
distinct facet of $\Delta_d$ and such that their union covers
$\Delta_d$, then their intersection is nonempty (this statement is also
called a \emph{dual KKM theorem} in \cite{asadaetal}). To see that it is a
consequence of the KKM theorem, assign number $i$ to the facet covered
by $C_i$; number the vertices of $\Delta_d$ so that vertex $i$ is on
facet $i$; the $C_i$'s satisfy then the condition of the KKM theorem.
  
Several variations of the original KKM theorem exist.
Gale~\cite{Gale-colorfulkkm} proved the following \emph{colorful}
version: given $d+1$ different KKM covers $\{C^1_i\}_{i=1}^{d+1},
\{C^2_i\}_{i=1}^{d+1}$, \ldots, $\{C^{d+1}_i\}_{i=1}^{d+1}$ of the
$d$-simplex, there exists a permutation $\pi$ of the numbers $[d+1]$
such that $\bigcap_{i=1}^{d+1} C_{\pi(i)}^i \neq \varnothing$.
Clearly choosing all the covers to be the same recovers the classical
version of KKM stated above.  Gale's colorful KKM theorem has an
intuitive interpretation first stated by Gale himself: ``if each of
three people paint a triangle red, white and blue according to the KKM
rules of covering, then there will be a point which is in the red set
of one person, the white set of another, the blue of the third''. A
recent strengthening of this colorful theorem~\cite{asadaetal} states
that, given $d$ KKM covers $\{C^1_i\}_{i=1}^{d+1},
\{C^2_i\}_{i=1}^{d+1} \ldots, \{C^{d}_i\}_{i=1}^{d+1}$ of the
$d$-simplex $\Delta_d$, then there exist a point $\xx$ in $\Delta_d$
and $d+1$ bijections $\pi_i\colon[d] \rightarrow[d+1] \setminus \{i\}$
for $i=1,\dots,d+1$, such that $\xx \in \bigcap_{j=1}^d
C_{\pi_i(j)}^j$ for every $i$. It is interesting to note the proofs of
these colorful results combine degree theory with Birkhoff's theorem
on doubly-stochastic matrices. Finally, we note \cite{Musin:2017} has
a common generalizations of Sperner, Tucker, KKM, Fan.

\bigskip

Another fascinating and very useful consequence of Brouwer's theorem
is Kakutani's 1941 fixed-point theorem. It deals, not with real-valued
functions, but with \emph{set-valued functions}, where points are
mapped to subsets.  For a suitable notion of continuity, it ensures
that for any continuous function $F$ mapping points of a ball to
convex subsets of it there is an $\ve x$ such that $\ve x\in F(\ve
x)$. Kakutani's theorem is especially useful in game theory, its most
traditional application being the Nash theorem, see
Section~\ref{s-games+fairness+independence}.
  
  \bigskip
  
Similar to Sperner's lemma, Tucker's lemma has continuous and covering
versions. The continuous version is the celebrated \emph{Borsuk-Ulam
  theorem}, which has many applications in discrete geometry,
combinatorics, and topology. It asserts that there is no continuous
function from $\s^d$ into $\s^{d-1}$ that commutes with the central
symmetry. Nice proofs of the Borsuk-Ulam theorem from Tucker's lemma,
as well as equivalent formulations and many applications, can be found
for instance in the books of
Matou\v{s}ek~\cite[Chapter~2.3]{matousek2003using} and de
Longueville~\cite[Chapter~1]{Longueville-book}. Just as KKM is the
covering version of Sperner's lemma, Tucker's lemma has a covering
version, the \emph{Lyusternik-Schnirel'mann theorem} ~\cite{LStheorem}.  It
states that, if the sphere $\s^d$ is covered by $d+1$ closed subsets,
then one of the sets must contain two antipodal points. This
theorem and some of its extensions (e.g., those due to K. Fan) have
found many applications in other areas of mathematics, for instance,
as for the KKM theorem, in the study of variational inequalities.

\subsection{Generalizations and variations}\label{s:s-t-generalizations}

A labeling $\lambda$ of a simplicial complex $\T$ by a set $S$ can be
interpreted as a map from the vertices of $\T$ to the vertices of some
abstract simplicial complex $\K$ with vertex set $S$. This viewpoint leads
to several interesting developments.

\bigskip

A first idea is to extend $\lambda$ into a \emph{linear map} $f$ from
$|\T|$ into $|\K|$. For a Sperner labeling, $f$ maps $\Delta_d$ into
itself. Composing $f$ with a suitable orthogonal transformation
ensures that any fixed point of the resulting map, which must exist by
Brouwer's theorem, lies in a fully-labeled simplex; this is the
standard proof of Sperner's lemma from Brouwer's theorem~\cite[Section
  1.1]{Longueville-book}. Using this idea, Sperner's lemma can be 
  generalized to triangulations of arbitrary polytopes by De Loera, Peterson and Su~\cite{DeLoera:2002hj}.

\begin{proposition}[Polytopal Sperner lemma] \label{p:polytopalsperner}
  Let $P \subset \R^d$ be a polytope with $n$ vertices, $\T$ a
  triangulation of $P$, and $\lambda$ a labeling of the vertices of
  $\T$ by $[n]$. If the vertices of $P$ have pairwise
  distinct labels and every vertex of~$\T$ lying on a face $F$ of the
  boundary of $P$ has same label as some vertex of $F$, then $\T$
  contains at least $n-d$ fully-labeled facets.
\end{proposition}

\begin{figure}[h]
  \begin{center}
    \includegraphics[page=13]{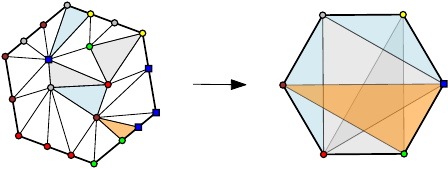}
  \end{center}
  \caption{A Sperner-type labeling of a triangulation of a polygon and
    the associated surjective map to the
    polygon.\label{f:polytopalsperner}}
\end{figure}

\noindent
The gist of the proof is that $\lambda$ extends into a piece-wise
linear map from $|\T|$ to~$P$ (as illustrated in
Figure~\ref{f:polytopalsperner}). This map can be shown to be
surjective, so its image provides a covering of $P$ by simplices
spanned by its vertices. The number of such simplices required to
cover $P$ is at least $n-d$, and each of them is the image (under the
extension of $\lambda$) of a fully-labeled facet of $\T$.  This
approach generalizes to non-convex polyhedral
manifolds~\cite{meunier-sperner} and broader classes of
manifolds~\cite{Musin:2012tq, Musin:2016}. The reader interested in recent
progress on lower bounds for the number of fully-labeled facets is
refered to the work of Asada et al.~\cite{asadaetal}.

\bigskip

Another idea is to extend $\lambda$ into a chain map $\lambda_\sharp:
C_{\bullet}(\T,R) \to C_{\bullet}(\K,R)$. Depending on the coefficients
ring $R$ and the complex $\K$, one gets different generalizations of
the classical statements. This extension goes as follows. Send every
simplex $\{\ve v_{0}, \ldots, \ve v_{k}\}$ of $\T$ to $\{\lambda(\ve v_{0}),
\ldots, \lambda(\ve v_{k})\}$ if the $\lambda(\ve v_{i})$ are pairwise
distinct and to $0$ otherwise; the linear extension of this map to
$C_{\bullet}(\T,\Z_2)$ is the map $\lambda_\sharp$, and it commutes
with $\partial$ (it is a \emph{chain map}).

\medskip

For $R=\Z_2$, we obtain a short proof by induction of Sperner's
lemma. Assume that $\lambda$ is a Sperner labeling of a triangulation
$\T$ of $\Delta_d$, so $\lambda_\sharp: C_{\bullet}(\T,\Z_2) \to
C_{\bullet}(\sigma_d,\Z_2)$. Let $t$ denote the sum of $d$-simplices
of $\T$ and $\sigma$ the unique $d$-simplex of $\sigma_d$. Observe
that $\lambda_\sharp(t) = \ell \sigma$, where $\ell$ is the number of
fully-labeled simplices in $\T$ mod $2$. A simple induction shows that
$\lambda_\sharp(\partial t) = \partial \sigma$. Thus, $\partial \lambda_\sharp(t)
= \lambda_\sharp(\partial t) \neq 0$. It comes that $\ell \neq 0$, and $\T$
has an odd number of fully-labeled simplices.

\medskip

For $R=\Z$, the same proof gives a Sperner-type result for oriented
simplices. The \emph{orientation} of a fully-labeled simplex
$\{\ve v_1,\ve v_2, \ldots, \ve v_{d+1}\}$ of $\T$, where $\ve v_i$
has label $i$, is defined as the sign ($1$ or $-1$) of the determinant
\[ \left|\begin{array}{cccc} \ve v_1 & \ve v_2 & \ldots & \ve v_{d+1} \\ 1 & 1 & \ldots & 1\end{array}\right|.\]
Specifically, this proves that in any Sperner labeling of a
triangulation of $\Delta_d$ by $[d+1]$ the orientations
of the fully-labeled simplices add up to $1$~\cite[Theorem~2]{Fan68}. 
This approach also yields a proof of (a signed version of) Fan's
Lemma~\cite[Theorem~1.10]{Longueville-book}, again for triangulations with flags of hemispheres.

\medskip

For $R = \Z_q$, this idea leads to generalizations of the lemmas of
Tucker and Fan with more than two
``signs''~\cite{HaSaScZi09,Me05}. The general insight, following a
generalization by Dold~\cite{Dold:1983wr} of the Borsuk-Ulam theorem,
is to replace $\s^d$ and the antipodality by a suitable simplicial
complex on which $\Z_q$ acts freely. The resulting $\Z_q$-Fan lemmas
actually provide combinatorial proofs of Dold's theorem. The
deduction of the $\Z_q$-Fan lemmas from their $\Z_q$-Tucker versions
works as explained above in the $\Z_2$ case~\cite{alishahi2017}.

\medskip

A variant of the above argument yields a new and elementary proof of a relative of Sperner's lemma, due to Meshulam~\cite{Mes01}. 
It is a powerful result that has found many applications in graph theory and in discrete geometry, such as the recently-found 
generalization of the colorful Carath\'eodory theorem by Holmsen and Karasev~\cite{Holmsen+Karasev}. 
We come back later to some of its applications in discrete geometry (Section \ref{caratheodory-thms}) and in graph theory (Section~\ref{subsubsec:hall}). 
Meshulam's lemma was first explicitly derived from (a special version of) Leray's acyclic cover theorem (its first use was implicit; see 
Kalai and Meshulam~\cite[Proposition~3.1]{Kalai:2005cm} for an explicit statement and proof for rational homology). 

For the sake of presentation, we
consider here homology over $\Z_2$ but the proof generalizes,
\emph{mutatis mutandis}, to an arbitrary ring. Given a simplicial
complex $\K$ and a subset $X$ of its vertices, we denote by $\K[X]$
the simplicial complex formed by the simplices of $\K$ whose vertices
are in~$X$.

\begin{prop}[Meshulam lemma]\label{p:meshulamlemma}
  Let $\lambda$ be a labeling of the vertices of a simplicial complex
  $\K$ by $[d+1]$. If $\widetilde{H}_{|I|-2}\pth{\K[\lambda^{-1}(I)]}$
  is trivial for every nonempty $I \subseteq [d+1]$, then $\K$ contains
  a fully-labeled $d$-dimensional face.
\end{prop}

\noindent
The proof goes as follows. Let $\lambda_{\sharp}:C_\bullet(\K) \to
C_\bullet(\sigma_d)$ denote the chain map induced by $\lambda$ and
recall that it maps simplices with repeated labels to $0$. We build a
chain map $f_\sharp\colon C_\bullet(\sigma_d) \to C_\bullet(\K)$ such
that $\lambda_\sharp \circ f_\sharp = \operatorname{id}_{\sigma_d}$;
the identity $(\lambda_\sharp \circ f_\sharp)([d+1]) = [d+1]$ then
ensures that $f_\sharp([d+1])$ contains an odd number of fully-labeled
simplices, proving the statement. We build $f_\sharp$ by increasing dimension. We start by setting
$f_\sharp(\{i\})$ to be some (arbitrary) vertex in
$\lambda^{-1}(\{i\})$ for every $i \in [d+1]$; this is possible
because $\widetilde{H}_{-1}(\K[\lambda^{-1}(\{i\})])$ is
trivial. Assume that $f_\sharp$ is defined over all chains up to
dimension $k$, that it maps any subset $I$ of cardinality at most $k+1$ to $C_\bullet(\K[\lambda^{-1}(I)])$, and that it commutes with the boundary operator. Now, for
any subset $I$ of cardinality $k+2$, we have
\[ \partial\pth{\sum_{i \in I}
f_\sharp(I\setminus \{i\})} = \partial f_\sharp\pth{ \sum_{i \in I} I\setminus \{i\}} = \partial f_\sharp(\partial I) = f_\sharp(\partial^2 I) = 0,\]
so the $k$-chain $\sum_{i \in I} f_\sharp(I\setminus \{i\})$ is a
cycle in $C_\bullet(\K[\lambda^{-1}(I)])$. The assumption of the lemma
ensures that it is the boundary of some chain $\gamma \in
C_{|I|-1}(\K[\lambda^{-1}(I)])$, and we set $f_\sharp(I)$ to be
$\gamma$. To see that $(\lambda_\sharp \circ f_\sharp )(I) = I$ for any
$I\subseteq [d+1]$, first note that this is straightforward for
$|I|=1$. For the general case, remark that
\[ \partial \lambda_\sharp \pth{ f_\sharp(I)} = \lambda_\sharp \pth{\partial f_\sharp(I)} = \sum_{i \in I} \lambda_\sharp \pth{f_\sharp(I\setminus \{i\})} = \sum_{i \in I} I\setminus \{i\},\]
so we can assume by induction that $\partial \lambda_\sharp \pth{
  f_\sharp(I)} = \partial I$. We conclude by observing that
$f_\sharp(I) \in C_{|I|-1}(\K[\lambda^{-1}(I)])$ means that
$\lambda_\sharp\pth{f_\sharp(I)}$ is supported only on subsets of $I$,
so it must be that $\lambda_\sharp \pth{ f_\sharp(I)} =I$.  

\bigskip

The parity argument used to prove Proposition~\ref{prop:sperner_even}
can also be found, specialized to a certain labeled pseudomanifold,
in Scarf's proofs of the integer Helly theorem~\cite{Sca1977} and his
classical lemma in game theory~\cite{game-scarf}. There is a related
unbounded polar version that will be useful in Section~\ref{s:kernels}.

\begin{corollary}[{\cite[Theorem~3]{KiPa09}}]\label{cor:KP}
  Let $P$ be a $d$-dimensional pointed polyhedron whose characteristic
  cone is generated by $d$ linearly independent vectors and whose
  facets are labeled by~$[d]$. If no facet containing the $i$-th
  extreme direction is labeled $i$, then there exists a extreme point
  incident to a facet of each label.
\end{corollary}

Another recent variation of Sperner's lemma, motivated by applications
in approximation algorithms, asks for the minimum possible number of
facets in the Sperner labeling that must be
non-uniquely-labeled. Mirzakhani and Vondrak~\cite{mirza+vondrak,MR3726616}
settled this question for certain triangulations of the simplex, for
which they provided the optimal Sperner labeling. They then proposed
two very interesting open questions.
   
\begin{oproblem}
  Is there a theorem that interpolates between the result above (lower
  bound on the number of simplices with at least two different colors)
  and the original Sperner's lemma (lower bound on the number of
  simplices with vertices of different color) by predicting how many
  simplices are colored with at least $j$ different colors? How does
  such theorem depend on the structure of the particular
  triangulation?
\end{oproblem}

It must be mentioned that Tucker-type theorems and Sperner-type
theorems are related to each other in interesting ways. For example,
it is known that the Borsuk-Ulam theorem implies the Brouwer
fixed-point theorem, but at the combinatorial level Nyman and Su
\cite{Nyman+Su2013} proved that Fan's lemma implies Sperner's theorem
too.  Other interconnections can be found in
\cite{prescott+su,spencer+su}.

\subsection{Computational considerations} \label{fixed-point-computing} 

The proof of Sperner's lemma given for
Proposition~\ref{prop:sperner_even} builds a graph where every vertex
has degree zero, one, or two, then exhibits a vertex of degree one and
argues that any \emph{other} vertex of degree one must correspond to a
fully-labeled simplex. This provides a simple algorithm for finding a
fully-labeled simplex: just follow the path!  We can combine this
simple path-following algorithm for finding fully-labeled simplices
with the proof of Brouwer's theorem, presented at the beginning of
Section \ref{subsec:continuousversionspernertucker}, and provide a
method for finding an \emph{approximate fixed-point} of the continuous
map $f$. Again assume we are given a continuous map $f: \Delta_d \to
\Delta_d$ and an $\varepsilon>0$. Our goal is to find $\xx \in
\Delta_d$ such that $\|f(\xx)-\xx\| \le \varepsilon$. For this, it
suffices to compute a triangulation of $\Delta_d$ with simplices of
diameter sufficiently small, depending on $\varepsilon$ and the
modulus of continuity of $f$, label it as in the proof of Brouwer's
theorem, and any vertex $\xx$ of a fully-labeled simplex does the job
(this fact is more easily formalized by using the $\ell_\infty$ norm
on the barycentric coordinates). This template of proof was first
presented in \cite{scarf1967} and is quite flexible, e.g., it applies
to non-contracting functions. We left out many details, for instance
the choice of the triangulation to speed up the algorithm and the
estimation of the modulus of continuity. The interested reader can
find more details on methods to compute approximate fixed points based
on these ideas in~\cite{todd1993new}.

\medskip

The theory of computational complexity is a formal way for
computer scientists to classify the inherent difficulty of
computational problems. Families of problems, called \emph{complexity
  classes}, collect problems of equivalent difficulty (a complete
introduction can be found in \cite{Arora+Barak:complexity}). Famous
complexity classes of course include the class P and the class NP,
but we briefly discuss here, and in Section \ref{compu-convgeo},
about the complexity classes that relate to computational versions of
our five central theorems.

The \emph{path-following algorithm} for computing the fully-labeled
simplex for Sperner's lemma is representative of the \emph{PPAD
  class}, a complexity class well-suited for studying computational
problems whose solution is known to exist, but finding it is not that
easy. This class was presented by Papadimitriou~\cite{ppad-original},
see also~\cite{Nisanetal-2007}. The prototypical problem of the class
PPAD (which abbreviates \emph{Polynomial Parity Argument for Directed
  graphs}) assumes an underlying directed graph where every vertex has
in- and out-degrees at most one; the graph may be implicit, and all
that is required is the existence of a function that computes the
neighborhood of a given vertex in time polynomial in the encoding of
that vertex. The PPAD problem asks, given the encoding of an
unbalanced vertex (that is, with different in- and out-degrees), to
compute the encoding of \emph{another} unbalanced vertex. (Note that
this computational problem does not easily reduce to a meaningful
associated decision problem, since the existence of this other
unbalanced vertex follows systematically from parity considerations.)
A problem is in the \emph{PPAD class} if it has a polynomial reduction
to the PPAD problem, and a problem from the PPAD class is
\emph{PPAD-complete} if every problem from the PPAD-class has a
polynomial reduction to that problem. (All reductions here are meant
in the usual sense of polynomial reductions~\cite[$\mathsection
  2$]{Nisanetal-2007}.)  PPAD-completeness results imply conditional
lower bounds in the following sense: one cannot solve a PPAD-complete
problem substantially faster than by path-following, unless there is
also a substantially better method than path-following for the PPAD
problem (and similarly for every other problem in the PPAD class). As
in the case of the P vs NP problem, failure over time to improve on
even the most streamlined of these problems supports the conjecture
that \emph{none} of these methods can be substantially improved.

The ``Sperner problem'' -- where one asks, given a Sperner labeling, for
a fully-labeled simplex -- is well-known to be PPAD-complete in any
fixed dimension. (Formalizing this problem properly requires some
care, for instance the definition of a canonical sequence of
triangulations with simplices of vanishing diameter, which we do not
discuss.)  Papadimitriou's seminal paper, which started the theory of
PPAD problems~\cite{ppad-original}, settled the three-dimensional case
and listed the two-dimensional case as an important open problem; it
was settled in the positive a decade later by Chen and
Deng~\cite{chen2009complexity}.

While Tucker's lemma can also be proved via a path-following
argument~\cite{freundtodd}, the computational problem associated to
Tucker's lemma is not known to belong to the PPAD class: contrary to
Sperner's lemma, there is no natural orientation of the edges of the
underlying graph. The suitable complexity class to use for the
``Tucker problem'' is a superclass of the PPAD class, the class {\em
  PPA}.  Here PPA abbreviates \emph{Polynomial Parity Argument for
  graphs}. This class was introduced at the same time as the PPAD
class, its definition is almost the same: instead of working with
directed graphs, one works with undirected ones. The underlying graph
defining the PPA problem has all its vertices of degree at most two
and asks, given the encoding of a vertex of degree one, to output the
encoding of another degree one vertex. PPA contains PPAD but it is a
famous problem to decide whether the two classes are actually the
same, already asked in the paper founding this topic
\cite{ppad-original}. Experts believe that these two classes are
different \cite{GoldbergPapa17}.

\begin{oproblem} 
  Are the complexity classes PPA and PPAD equal?  
\end{oproblem}

As for the Sperner problem, the Tucker problem is PPA-complete already
in dimension two (see Aisenberg et al.~\cite{aisenbergetal2015}, who
corrected an earlier, wrong, assertion of
PPADness). P{\'a}lv{\"o}lgyi~\cite{palvolgyi20092d} introduced a clean
variation of this problem -- the octahedral Tucker problem -- together
with the open question below: Given a function
$\lambda:\{+,-,0\}^n\setminus\{\zero\}\rightarrow \{\pm 1, \pm 2,
\ldots, \pm m\}$, computable in time polynomial in $n$ and such that
$n > m$ and $\lambda(-\ve x)=-\lambda(\ve x)$ for all $\xx$, compute
$\xx\preceq\yy$ such that $\lambda(\xx)+\lambda(\yy) = 0$.  Note that
contrary to the Tucker problem we just discussed, the dimension is not
part of the input. The computational complexity of the algorithmic
version of the Octahedral Tucker lemma had been an outstanding
challenging problem, but the paper \cite{dengetal2017} resolved this
problem by proving Octahedral Tucker to be PPA-complete.

\bigskip

We will also discuss in our applications, in particular
Section~\ref{s:nash-complexity}, some consequences of the lemmas of
Sperner or Tucker whose computational versions may be complete for two
related complexity classes. The first class is FIXP, introduced by
Etessami and Yannakakis~\cite{etessami+yannakakis}. It consists of the
problems whose resolution on an instance $\ell$ reduces to the
computation of a fixed point of some function $F_\ell$ that can be
expressed by the operations $\{+,*,-,/,\max,\min\}$ with rational
constants and functions and computed in time polynomial in the size of
$\ell$; this extends PPAD, which coincides with the case of linear
functions. The second class, called $\exists\R$~\cite{Schaefer2015}
(sometimes abbreviated ETR for existence of real solutions, see
\cite{gargetal2015}), studies problems that reduce to deciding the
emptiness of a general semi-algebraic set, i.e., the set of real
solutions of a system of inequalities with polynomials as
constraints. These two complexity classes are relevant in Section
\ref{s-games+fairness+independence}.

\section{Combinatorial convexity}
\label{s:convexgeo}

We now focus on three classical combinatorial theorems about convex
sets first identified in the early 20th century. These are the
theorems of Carath\'eodory, Helly and Tverberg.  The importance of
convexity in applications, and hence of these three theorems, owes
much to the computational effectiveness of convex optimization
algorithms both in practice and in theory
\cite{bertsekasbook,boyd2004convex}. This encourages applied
mathematicians to look for convexity, or for ways to approximate
complicated sets using convex sets. Surprisingly, convexity appears in
unexpected settings.

Extensive surveys were devoted to (subsets and variations on) some of
these three theorems by Danzer, Gr\"unbaum, and
Klee~\cite{DGKsurvey63}, Eckhoff \cite{Eckhoff:1993survey}, Holmsen
and Wenger~\cite{holmsen+wenger} and Amenta, De Loera, and
Sober\'on. \cite{amentaetal2017helly}.  An account of early variations
of Carath\'eodory's theorem is in the memoir by Reay~\cite{Reay-mem}.

There is an abundant literature on \emph{axiomatic convexity}, which
studies analogues of the theorems of Carath\'eodory, Helly, and
Tverberg, not over Euclidean spaces as we do here, but over purely
combinatorial abstract settings, for instance in the convexity spaces
defined by arbitrary graphs, finite geometries, matroids, greedoids,
etc. The three theorems play a significant, and interesting role,
there too, but we do not cover this topic here. We refer the
interested reader to the
references~\cite{Duchetsurvey1987,Kay:1971uf,convexityspaces-vandevel}.

\subsection{Carath\'eodory} \label{caratheodory-thms}

We will first consider Carath\'eodory-type theorems that certify
membership of a point in the convex hull of a set via linear
non-negative combinations.  The original theorem of
Carath\'eodory~\cite{originalCaratheodory} asserts that any point in
the convex hull of a finite point set in $\R^d$ is a convex
combination of some at most $d+1$ of these points.  Equivalently, if a
vector $\ve b$ belongs to the \emph{cone} of $X = \{\vv_1, \vv_2,
\ldots, \vv_n\} \subset \R^d$ (i.e., the \emph{positive hull}
of all non-negative \emph{real} linear combinations of vectors in $X$),
then $\ve b$ is a positive combination of at most $d$ vectors of
$X$. To see this, let $A = \pth{\begin{matrix} \vv_1 & \vv_2 & \cdots
    & \vv_n \end{matrix}}$ and assume that $\ve{\tilde{x}}$ is a
solution of
\begin{equation}\label{eq:Cara}
  \begin{array}{rcl}
    A \ve x & = & \ve b\\
    \ve x & \ge & \ve 0.
  \end{array}
\end{equation}
If the support of $\ve{\tilde{x}}$ has size at least $d+1$, then $A
\ve x = \ve 0$ has some nontrivial solution $\ve z$ with support
contained 

\smallskip\noindent
\begin{minipage}{12cm}
  in the support of $\ve{\tilde{x}}$. For an adequate value of
$t$, the vector $\ve{\tilde{x}}+t\ve z$ is a solution of System~\eqref{eq:Cara} with
  smaller support than $\ve{\tilde{x}}$. A closer examination of this
  argument yields that, in the plane, any point in the convex hull of
  four points lies in two of the triangles they span (as illustrated
  on the right). The following strengthening of Carath\'eodory's
  Theorem will be useful in optimization.
\end{minipage}
\hfill
\begin{minipage}{4cm}
  \centering \includegraphics[page=14]{figures-final}
\end{minipage}

\begin{proposition}\label{p:cara++}
  Any point in the convex hull of (at least) $d+2$ points in $\R^d$ lies in the
  convex hull of at least two $(d+1)$-element subsets of these points.
\end{proposition}

\noindent
A geometric proof of Proposition~\ref{p:cara++} 
with any one of the $d+2$ points, say $\ve p$, and the point $\ve x$ in their convex 

\smallskip\noindent
\begin{minipage}{4cm}
  \centering \includegraphics[page=15]{figures-final}
\end{minipage}
\hfill
\begin{minipage}{12cm}
hull.  Shoot a ray from $\ve p$ to $\ve x$ and collect the (at most)
$d$ vertices of the face of a (triangulation of the) convex hull
through which this ray exits; these $d$ points and $\ve p$ contain the
point $\ve x$ in their convex hull, and any of the $d+2$ points can be
used as origin of the ray. The figure on the left illustrates this
process in the case of a (blue) point $\ve x$ contained in the convex
hull of the (white) vertices of a cube in $\R^3$.  The green
projections to faces help us find the four white points containing
$\ve x$.
\end{minipage}

\bigskip

Several variants of Carath\'eodory's theorem have been developed. For
instance, Steinitz \cite {originalsteinitz} proved that a point in the
interior of the convex hull of a set $S$ lies in the interior of the
convex hull of some $2d$ points of $S$. A related classical result is
the \emph{Krein-Milman theorem}~\cite{originalkreinmilman}: if $C$ is
a convex set, then every point in $C$ is a convex combination of its
extreme points, i.e., those that are not convex combinations of others
in the set. Allowing more flexible representations,
Dobbins~\cite{dobbins2015point} proved that any point in an
$ab$-dimensional polytope is the barycenter of $a$ points on its
$b$-dimensional skeleton (another proof of Dobbins' theorem is shown
in \cite{blafrickzie}).

\bigskip

One of the most applicable and powerful variants is the \emph{colorful
  Carath\'eodory} theorem. We saw this in the introduction already.

\begin{namedtheorem}[Colorful Carath\'eodory theorem]
  Let $C_1, C_2, \ldots, C_{d+1}$ be point sets in $\R^d$. If a point
  $\ve p$ is in the convex hull of every $C_i$, then there exist $\ve x_1 \in
  C_1$, $\ve x_2 \in C_2$, \ldots, $\ve x_{d+1} \in C_{d+1}$ such that $\ve p$
  lies in the convex hull of $\{\ve x_1,\ve x_2, \ldots, \ve x_{d+1}\}$.
\end{namedtheorem}
  
\noindent
In other words, if the origin is contained in the convex hull of each
of $d+1$ point sets $C_1, C_2, \ldots, C_{d+1}$ (the color classes),
then it is contained in a colorful ``simplex'', i.e., one where each
of the vertices comes from a different $C_i$ (with the understanding
that this ``simplex'' may be degenerate). This theorem was discovered
by B\'ar\'any~\cite{baranys-caratheodory} who showed that a
colorful simplex that minimizes the distance to the origin must
contain it. Indeed, when the distance is still positive, it is
attained on a facet and the vertex opposite to that facet may be
changed to further decrease the distance to the origin. This approach
inspired new proofs and algorithms, by minimization, of the colorful
Carath\'eodory theorem \cite{barany1995caratheodory}  and other results such
as Tverberg's theorem~\cite{Roudneff-ejc,Tverberg1981,TverbergVrecica}
(see more on Tverberg's theorem later in Section \ref{tv-section}). 

An alternative proof of the colorful Carath\'eodory theorem applies
Meshulam's lemma (Proposition~\ref{p:meshulamlemma}) to the join of
two abstract simplicial complexes built on top of
$\bigcup_{i=1}^{d+1}C_i$: one has a simplex for any subset of points with
no repeated color, and the other has a simplex for every subset of
points not surrounding the origin. (The labeling is given by the
identification of the vertices of the join to $\bigcup_{i=1}^{d+1}C_i$.)
This approach emerged from a colorful Helly theorem of Kalai and
Meshulam~\cite{Kalai:2005cm} (see below) and later allowed a purely
combinatorial generalization of the colorful  Carath\'eodory theorem by
Holmsen~\cite{Holmsen-colorful} where the geometry and the
colorfulness are abstracted away into, respectively, an oriented
matroid and a matroid. See also the paper \cite{Holmsen+Karasev}.

The assumption of the colorful Carath\'eodory theorem ensures that not
only one, but actually many colorful simplices exist; we come back to this
question when discussing simplicial depth in Section~\ref{depthsec}. This
also underlies its connection to Tverberg's theorem which we discussed in Section \ref{tv-section}.

Many variations and strengthenings of the colorful Carath\'eodory
theorem have been explored starting with B\'ar\'any's seminal
paper~\cite{baranys-caratheodory} and other collaborations
\cite{barany1995caratheodory}.  Recent strengthenings include Deza et
al. colorful simplicial depth \cite{Dezaetal2006} (discussed in
Section \ref{depthsec}) and Frick and Zerbib's common generalization
of the colorful Carath\'eodory theorem and the KKM theorem
\cite{frick+zerbib}.  Another key variation, discovered independently
by Arocha et al.~\cite{ABBFM09} and Holmsen et al.~\cite{HPT08}, is
that the assumption that the convex hull of each $C_i$, $1 \leq i \leq
d+1$, contains the origin can be weakened to only require that the
convex hull of each pair $C_i \cup C_j$, $1 \leq i < j \leq d+1$,
contains the origin. There are examples showing that it is not
sufficient that the convex hulls of triples contain the origin, but
weaker relaxations are possible~\cite{deza+meunier:cara}.  Arocha et
al.~\cite{ABBFM09} also proved another ``very colorful Carath\'eodory
theorem''.

Via point-hyperplane duality, one can derive from the colorful Carath\'eodory theorem, a colorful theorem of Helly. 
We will not discuss this in detail, but let us explain the basics. Consider $d+1$ families  $F_1, F_2,\dots,F_{d+1}$ of convex sets inside $\R^d$ (researchers
think of these as colors).  Assume that every \emph{colorful} selection of $d+1$ of the sets, i.e., one set from each $F_i$, has a non-empty common intersection.  Then, the classical colorful Helly theorem of Lov\'asz  (see \cite{baranys-caratheodory}) says that there is at least one family $F_i$, 
whose sets have a non-empty intersection.  Here is now the dual of the \emph{very} colorful Carath\'eodory theorem of \cite{ABBFM09}: 
given a finite family of halfspaces in $\R^d$ colored with $d+1$ colors, if every colorful selection of $d+1$ halfspaces has a non-empty common intersection, 
then there exist \emph{two} colors classes all of whose members intersect. For more examples and references, see \cite{Martinez-Sandovaletal18} and 
references therein. 

The proof of the colorful Cara\-th\'eo\-dory theorem also implies that
given $d+1$ point sets $C_1, \ldots, C_{d+1}$ and a convex set $C$,
either one $C_i$ can be separated from~$C$ by a hyperplane, or there
exists a colorful simplex intersecting $C$.  Building on this, Mustafa
and Ray~\cite{MR16} showed that given $\lfloor \frac{d}{2} \rfloor +1$
sets of points in $\R^d$ and a convex object $C$, then either one of
the sets can be separated from $C$ by a constant number of
hyperplanes, or there is a $\lfloor \frac{d}{2} \rfloor$-dimensional
colorful simplex intersecting $C$.

\bigskip

The \emph{integer} Carath\'eodory problem considers a finite set
  $X \subset \Z^d$ and $\ve v \in \Z^d$ in its positive hull
  and asks whether $\ve v$ can be written as a non-negative,
  \emph{integer} linear combination of some elements of $X$, and, if
  true, how many elements are needed. The answer to the first
question is negative in general. For instance consider
\begin{equation}\label{eq:non-hilbert}
  X = \left\{ \left(\begin{array}{c} 1 \\ 0 \end{array}\right),
    \left(\begin{array}{c} 1 \\ 2 \end{array}\right)\right\} \quad
  \hbox{and} \quad \ve v = \left(\begin{array}{c} 1 \\ 1
  \end{array}\right),
\end{equation}
$\ve v$ is an integral vector of the positive hull of $X$ but is not an
\emph{integer} non-negative combination of elements of $X$.

It is thus natural to restrict one's attention to subsets $X \subset
\Z^d$ such that every integral point of the cone of $X$ can be written
as a non-negative integer combination of elements of $X$; such sets
are called \emph{Hilbert generating sets}. This restriction is
reasonable given that the integer points of any rational polyhedral
cone $C$ have a finite Hilbert generating set. However, even in this setting,
there is no version of an integer Carath\'eodory theorem with a bound on
the size of the representation depending only on the dimension. Take for example
\begin{equation}\label{eq:hilbert-unbounded}
  \begin{aligned}
X_n & = \{2^i\ve e_j + \ve e_d \colon 0 \le i \le n-1 \hbox{ and } 1 \le j
\le d-1\} \subset \Z^d \\ & \hbox{and} \quad \ve v_n = (2^n-1,2^n-1,
\ldots, 2^n-1, n(d-1))^T,
  \end{aligned}
\end{equation}
$\ve v_n$ can be written as an integer combination over $X_n$, but any
such combination requires at least $n$ summands. Notice that the
coefficients in Example~\eqref{eq:hilbert-unbounded} grow quickly with
$n$. Such growth is necessary to force larger and larger sums because a
Carath\'eodory-type theorem is possible if one wants to bound the number
of summands in terms of the dimension and the size of the
coordinates. The best upper bound in that direction was recently
obtained in \cite{sparseintcaratheo}. An earlier bound appears in \cite{eisenbrandshmonin-caratheodory}. 

\begin{theorem}\label{Bound_via_siegel}
  Let $X=\{{\ve x}_1,\ldots, {\ve x}_t\} \subseteq \Z^d\setminus \{\ve
  0\}$ be a finite set, let $\| X \|_{\infty} = \max_{{\ve x} \in X}
  \| {\ve x} \|_{\infty}$, let $W = \pth{\begin{matrix} \ve x_1 &
      \cdots & \ve x_t\end{matrix}}$, and let $\Lambda$ denote the
  sublattice of $\Z^t$ of the integer points in the row space of
  $W$. Any vector representable as non-negative integer combination
  over $X$ can be written as a combination of at most $\min\{ \rank W
  + \log \det(\Lambda), \ 2d\log(2\sqrt{d}||X||_{\infty}) \}$ terms.
\end{theorem}

\noindent
The proof of Theorem~\ref{Bound_via_siegel} starts from some non-negative
integer combination and uses some element of the kernel of $W$ to
eliminate one of the summands. This is very similar to the classical
proof of the real-valued Carath\'eodory theorem, but now the kernel
element must be an \emph{integral} vector with coordinate entries in
$\{-1,1,0\}$. That such a vector exists is no longer a rank argument,
but follows from Siegel's lemma~(see \cite{siegel,vaalerbest}).

Interestingly, a full-fledged integer Carath\'eodory theorem, depending only on the
dimension, does exist for Hilbert bases of pointed cones. Let us explain. 
First of all, a \emph{Hilbert basis} is an
inclusion-minimal Hilbert generating set. We say a cone is \emph{pointed} if
it contains no linear subspace other than the nullspace. It is known a pointed
cone has a unique Hilbert basis~ (see e.g., \cite[Corollary
  2.6.4]{alggeo4optimization}). In contrast, when the cone is not pointed, there is no
  uniqueness, for instance $\{
(x,y) : x+y=0 \}$ has two Hilbert bases $\{(1,-1),(-1,1)\}$ and $\{(2,-2),
(-1,1)\}$.  As we will see in
Sections \ref{s:hilbertbasisinIP1} and \ref{s:hilbertbasisinIP2},
Hilbert bases play an important role in optimization theory and in the
solution of integer optimization problems. Here is the best known upper bound
of the number of Hilbert basis elements necessary to write a vector. This is due to
Seb\H{o}~\cite{sebocaratheodory}:

\begin{theorem}\label{thm:sebo}
  If the pointed cone $C$ is generated by a Hilbert basis $X \subseteq \Z^d$, then
  any of its integral points can be written as a non-negative integer
  combination of at most $2d-2$ elements of $X$.
\end{theorem}

\noindent
A weaker version of Seb\H{o}'s theorem, with a constant $2d-1$, was
previously obtained by Cook, Fonlupt, and
Schrijver~\cite{cook+fonlupt+schrijver}. Note that the sets $X_n$ in
Example~\eqref{eq:hilbert-unbounded} do define pointed cones, but are
not the Hilbert bases of those cones.  Seb\"o's theorem gives an upper bound, 
but the best known lower bound on the size of the linear
combination is only $d+\lfloor \frac{d}{6} \rfloor$ for $d \ge
6$~\cite{brunsetal}, so that leaves an open important problem:

\begin{oproblem}
  Determine the best possible constant for the integer
  Carath\'eo\-dory theorem on Hilbert bases of pointed cones.
\end{oproblem}

\noindent
The answer is known to be $d$ for $d=3$~\cite{sebocaratheodory} and in
some special cases such as the cone formed by the bases of any
matroid~\cite{gijswijt+regts}.

\bigskip

We conclude with an \emph{approximate} Carath\'{e}odory theorem, 
recently recovered by Barman~\cite{Barman}, which has an interesting 
application in game theory (see Section~\ref{s:nash-complexity}). Informally, 
it says that any point in the convex hull of a point set $X \subseteq \R^d$ can be
approximated by a convex combination of few elements of $X$. The
precise relation between the quality of approximation and the size of
the convex combination is quantified as follows:

\begin{theorem}\label{thm:caratheodory-p}
  Let $p \in [2, \infty)$ and let $X \subseteq \R^d$. For any point $\ve a
  \in \conv(X)$ and any $\varepsilon >0$, there is a point $\ve b$ such
  that $(i)$ $\|\ve a-\ve b\|_p \le \varepsilon$, and $(ii)$ $\ve b$ can be expressed as a linear
  combination of at most $4p \pth{\frac{\max_{\ve x \in X} \| \ve x
    \|_p}{\varepsilon}}^2$ vectors from $X$.
\end{theorem}

\noindent
Observe that the number of points used to represent the approximation
$\ve b$ to $\ve a$ is independent of the ambient dimension $d$. The
point $\ve b$ can in fact be chosen as the barycenter of points of $X$
with non-negative integer weights, the sum of the weights being at
most $k = 4p \pth{\frac{\max_{\ve x \in X} \| \ve x
    \|_p}{\varepsilon}}^2$. Barman's nice probabilistic proof writes
$a$ as a barycentric combination of $d+1$ points of $X$ (by
Carath\'eodory's theorem) and finds the $k$ points adding to $\ve b$
by sampling those $d+1$ points using the weights as probabilities
(some special care must be taken to ensure a bound independent of the
dimension). Theorem~\ref{thm:caratheodory-p} can also be derived from
Maurey's lemma in functional analysis (see
Pisier~\cite{pisier1980remarques} and
Carl~\cite{carl1985inequalities}). See also~\cite{BHPR16} for the
derivation of a related theorem using the Perceptron
algorithm~\cite{N62}. A very recent new generalization of Theorem
\ref{thm:caratheodory-p} was presented by Adiprasito, B\'ar\'any and
Mustafa in \cite{adiprasito+barany+mustafa}.  They proved that, given
a point set $P\subset \mathbb{R}^d$ of cardinality $n$, a point $\ve a
\in \mathrm{conv} (P)$, and an integer $r\le d$, $r \le n$, then there
exist a subset $Q\subset P$ of $r$ elements such that the distance
between $\ve a$ and $\mathrm{conv} (Q)$ is less then $\mathrm{diam}
P/\sqrt {2r}$.  Here the diameter of $P$ is the largest distance
between a pair of points in $P$.

\subsection{Helly}\label{s:Helly}

\begin{figure}[h]
  \begin{center}
  \includegraphics[page=16]{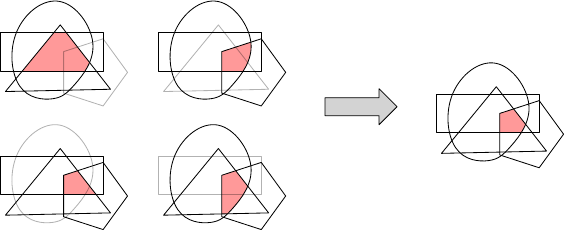}
  \caption{Helly's theorem in the plane.\label{f:Helly-pic}}
  \end{center}
\end{figure}

\emph{Helly's  theorem} asserts that for a finite family of convex
subsets of $\R^d$ with at least $d+1$ members, if every $d+1$ members intersect, 
then the whole family intersects. In the contrapositive, the empty intersection of
finitely many convex sets in $\R^d$ is always witnessed by the empty
intersection of some $d+1$ of the sets. See Figure~\ref{f:Helly-pic}

The special case of Helly's theorem where each subset is a halfspace
is of particular interest. Since a family of halfspaces not containing
the origin has empty intersection if and only if their inner normals
positively span the space, Helly's theorem for halfspaces is
equivalent to Carath\'eodory's theorem for their polars. Note that
the polar of a set of hyperplanes is a set of ray vectors. The case of
halfspaces implies the general case because, given a family of convex
sets in $\R^d$, we can replace each set by a polytope that it contains
without altering the intersection patterns. It suffices to take a
witness point in the intersection of every subfamily, and then replace each
set by the convex hull of the witness points that it contains.

Helly's original proof used the separability of compact convex sets by
hyperplanes to set up an induction on the dimension. It is spelled out
in the survey of Danzer et al.~\cite[Section~1]{DGKsurvey63} along
with references to eight other proofs. The most common proof deduces
Helly's theorem from \emph{Radon's lemma} (the case $r=2$ of Tverberg's
theorem)~\cite[Section~1.3]{Mbook}. Starting with $k\ge d+2$ convex
sets $C_1, C_2, \ldots, C_k$ where any $k-1$ intersect, one picks a
witness point $\ve w_i \in \bigcap_{j \neq i} C_j$. Partition $\ve w_1,\ve w_2,
\ldots,\ve w_k$ into two subsets with intersecting convex hulls, and
observe that this intersection point lies in every convex set. A proof of Helly's theorem, due to
Krasnosselsky~\cite{krano}, fits well with the theme of our
survey as it uses the KKM theorem. Lift the $\ve w_i$ to the vertex $\ve e_i$ of the
simplex $\Delta_{k-1}$. The
map $\ve e_i \mapsto \ve w_i$ extends into a linear map
$f:\Delta_{k-1} \to \R^d$, and setting $D_i = f^{-1}(\conv \{\ve w_j :
j \neq i\})$ produces $k$ closed subsets of $\Delta_{k-1}$. Each facet
of $\Delta_{k-1}$ is covered by a distinct $D_i$ and the $D_i$'s
cover~$\Delta_{k-1}$ by Carath\'eodory's theorem in $\R^d$. The KKM
theorem (see Section~\ref{subsec:continuousversionspernertucker}) 

\smallskip\noindent
\begin{minipage}{12cm}
yields a point $\ve w \in \bigcap_i D_i$ and the point $f(\ve w)$
is contained in every~$C_i$. Chakerian showed in~\cite{Chakerian-HellyBrouwer} that Helly's theorem
also follows from Brouwer's fixed point theorem, in a similar fashion
as the proof of the KKM theorem: the function, instead of moving every
point $\ve x$ by the vector of its distances to each $C_i$, moves
every point $\ve x$ to the barycenter of its projections on each
$C_i$. (See the right-hand picture for an illustration: the
projections of $\ve x$ to the three convex sets, red, green, and blue,
are shown by the white points; the black square is
their barycenter).
\end{minipage}
\hfill
\begin{minipage}{4cm}
\centering  \includegraphics[page=17]{figures-final}
\end{minipage}

\bigskip

\noindent
\begin{minipage}{3.5cm}
  \centering
  \includegraphics[page=18]{figures-final}
\end{minipage}
\hfill
\begin{minipage}{12cm}
\quad Helly's theorem for integral points was first established by
Doignon~\cite[Proposition~4.2]{Doi1973}, and rediscovered later by
Scarf~\cite{Sca1977} and Bell~\cite{Bell:1977tm}.
Hoffman~\cite{Hoffman:1979ix} observed that the techniques apply for
more cases than the integer lattice. The proofs of Doignon, Bell, and
Hoffman hinge on the following insight: if a polytope in $\R^k$ has
$m$ vertices, each with integer coordinates, and contains no other
integral point, then no two vertices may have all their coordinates of
the same parity (else their mid-point would yield a contradiction) and
thus $m$  is at most~$2^k$. Bell's proof starts with a family of $m$
halfspaces in $\R^k$ whose intersection contains no integer point and
such that removing any \makebox[\linewidth][s]{halfspace would enlarge the intersection to
include some integer point.} 
\end{minipage}

\smallskip\noindent
Translating each halfspace in the direction of its outer normal until every facet contains a witness
point with integer coordinates, one gets witness points that must be
distinct and form a polytope as above, so $m$ is at most $2^k$.
Hoffman's proof is more complicated but holds in a more general
axiomatic setting. Scarf's proof is algorithmic and relies on
Sperner's lemma (see also the variation by
Todd~\cite{todd1977number}). Remark that the equivalence via polarity
between Helly's and Carath\'eodory's theorems in~$\R^d$ does not carry
over to~$\Z^d$, as the bounds are respectively $2^d$ and at most
$2d-2$ (by Theorem~\ref{thm:sebo}).

\bigskip

Some of our applications will use the following version where some of
the coordinates in the intersection may be required to be integers.

\begin{namedtheorem}[Mixed Helly theorem] \label{mixedhelly}
  Let $\F$ be a finite family of convex sets in $\R^{d+k}$ of cardinality at least $(d+1)2^k$. 
  If every $(d+1)2^k$ members of $\F$ have a common point whose last $k$
  coordinates are integer numbers, then all members of $\F$ have an
  intersection point whose last $k$ coordinates are integer numbers.
\end{namedtheorem}

\noindent
The mixed Helly theorem was announced by Hoffman~\cite{Hoffman:1979ix}
as one of the outcomes of his axiomatic setting; he however deferred
details of the proof for the mixed analogue of the property of the
Doignon theorem to a forthcoming paper, which never appeared. A
complete proof came decades later and is due to Averkov and
Weismantel~\cite{AW2012}. Their proof proceeds in two steps. They
start with a family of halfspaces in $\R^{d+k}$ whose intersection is
a nonempty full-dimensional polytope $P$ containing no point whose
last $k$ coordinates are integer (the general case follows). They
project $P$ onto the last $k$ coordinates, obtaining a polytope $T(P)$
in $\R^k$ with no integer point; Doignon's theorem then ensures that
at most $2^k$ of the halfspaces supporting the facets of $T(P)$
already intersect with no integer point. By Carath\'eodory's theorem
in $\R^{d+k}$, each $k$-dimensional halfspace is the projection of the
intersection of some at most $d+1$ of the original halfspaces; the
bound follows.

\bigskip

\emph{Fractional} versions of Helly's theorem play an important role
in the study of sampling and hitting geometric set systems. Here,
\emph{fractional} means that one only assumes that some constant fraction of
the subfamilies (of a given size) intersect, and concludes that some
constant fraction of the whole family intersect.

\begin{theoremp}[Fractional Helly theorem]\label{thm:fracHelly}
  Let $0 < \alpha \le 1$ and $\F$ be a family of $n$ convex sets in
  $\R^d$. If $\alpha \binom{n}{d+1}$ of the $(d+1)$-element subsets of
  $\F$ have non-empty intersection, then some $(1-(1-\alpha)^{\frac1{d+1}})n$
  elements of $\F$ intersect.
\end{theoremp}

\noindent
The first result in this direction was proven by Katchalski and
Liu~\cite{Kat79frac}. Starting with a family of $n$ convex sets, they
assign to any subfamily the lexicographically minimum point in their
intersection. The set of points lexicographically \emph{larger} than a
given point is convex, so Helly's theorem ensures that the minimal
point of the intersection of $k \ge d$ convex sets is also the minimal
point in the intersection of some $d$ among them. A weak version of
the above theorem, where the size of the intersecting subfamily is
only guaranteed to be at least~$\frac{\alpha}{d+1}n$, then follows
from a pigeonhole argument.

There are few settings in which fractional Helly theorems are
known. On the one hand, Matou\v{s}ek~\cite{Matousek:2004cs} proved,
via a general sampling technique due to
Clarkson~\cite{clarkson-rs,clarkson-shor}, that any set system with
bounded VC dimension affords a fractional Helly theorem; his approach
holds for other measures of complexity than the VC
dimension~\cite{pinchasi-fracHelly}); we come back to the notions of
VC dimension in Section~\ref{s:geompart}. On the other hand,
B\'ar\'any and Matou\v{s}ek~\cite{baranymatousek} established a
fractional Helly theorem for lattices, including over the integers. It is surprising that
they only have to check the non-empty integral intersection of a positive fraction of
$(d+1)$-tuples, instead of the expected  $2^d$-tuples of intersections.

The bound of $(1-(1-\alpha)^{\frac1{d+1}})n$ in
Theorem~\ref{thm:fracHelly} is sharp and was obtained by
Kalai~\cite{Kalai:1984bg} and, independently, by
Eckhoff~\cite{eckhoff1985upper}, via a study of nerve complexes that
led to a more general \emph{topological} point of view. The
\emph{nerve} of a family of convex sets is the abstract simplicial
complex with a vertex for every set in the family, and a simplex for
every intersecting sub-family. Helly's theorem and its
fractional version easily translate in terms of nerves: the former
states that the nerve cannot contain the boundary of a $(\ge
d)$-dimensional simplex without containing the simplex, and the latter
asserts that if the nerve contains a positive fraction of the
$d$-dimensional faces, then it must contain a simplex of dimension a
positive fraction of $n$. Kalai's proof uses his technique of
\emph{algebraic shifting}~\cite{Kalai01algebraicshifting} to study how
the number of simplices of various dimensions behaves as the nerve is
simplified through a sequence of \emph{$d$-collapses}, a type of
filtration available to nerves of convex
sets~\cite{tancer2013intersection}. If a complex is $d$-collapsible,
then all its subcomplexes have trivial homology in dimension $d$ and
above, i.e., it is \emph{$d$-Leray}. Kalai's proof of
Theorem~\ref{thm:fracHelly} extends to $d$-Leray complexes
(see~\cite[$\mathsection 5.2$]{HT06}), and Alon et al.~\cite{AKMM}
further proved that families of subsets of $\R^d$ that are closed
under intersection and whose nerve is $d$-Leray also admit weak
$\varepsilon$-nets and $(p,q)$-theorems; examples of such families
include good covers (this follows from the Nerve theorem~\cite{nerve})
and acyclic covers~\cite{acyclic}. Topological versions of Helly theorem
have a further application in geometric group theory \cite{Farb}.

\bigskip

Fairly general \emph{topological} Helly theorems can be derived from
non-embeddability results via a construction reminiscent, again, of
the setup of the KKM theorem. Let us illustrate the basic idea with
five sets in the plane. If any four intersect and any three have
path-connected intersection, then we can draw $K_5$, the complete
graph on $5$ vertices, \emph{inside} the family by placing each vertex
in the intersections of four sets (different vertices missing
different sets) and connecting any two vertices by a path contained in
the intersection of the three sets that contain them both. By a
classical theorem of Kuratowski for planar graphs \cite{diestel-graphtheory}, 
there exist two edges that have no common vertex and that intersect. 
This intersection point must be in all five sets. An induction on the same idea 
yields that in a family of planar sets, where intersections are empty or path-connected, if
every four sets intersect, then they all must intersect.  
In higher dimension, where all graphs can be drawn without crossing, the same
approach can be combined with non-embeddability results derived from
the Borsuk-Ulam theorem, e.g., the \emph{Van Kampen-Flores theorem}
which states that $\Delta_{2k+2}^{(k)}$ does not embed in
$\R^{2k}$~\cite{vanKampen:KomplexeInEuklidischenRaeumen-1932,Flores:NichtEinbettbar-1933};
cf.~\cite[Chapter~5]{matousek2003using}. The discussion we present below on the
topological Tverberg theorem is also connected to embeddability of complexes, 
e.g., the paper \cite{Blagojevic:2014ul} proves that  the topological Radon theorem 
implies the Van Kampen-Flores theorem.

The most general result in this direction~\cite{BH-journal} is that any family $\F$ of subsets in~$\R^d$
admits a Helly-type theorem in which the constant that replaces $d+1$ in the case of convex sets is 
bounded as a function of the dimension $d$ and
\[ b = \max_{G \subseteq \F, 0 \le i \le \lceil \frac{d}2\rceil-1} \tilde{\beta}_i(\cap G),\]
where $\tilde{\beta}_i(X)$ denotes the $i$-th reduced Betti numbers,
over $\Z_2$, of a space $X$.

Non-embeddability arguments and the study of nerves offer two
different pathways to topological Helly theorems. While the former
allows more flexible assumptions, the latter offers more powerful
conclusions in the form of a sharp fractional Helly theorem. It is not
known whether the benefits of both approaches could be combined.

\begin{oproblem}
  Given $b$ and $d$, is there a fractional Helly theorem for families
  $\F$ of subsets in~$\R^d$ where $\tilde{\beta}_i(\cap G) \le b$ for
  any $G \subseteq \F$ and any $0 \le i \le \lceil\frac{d}2\rceil-1$?
\end{oproblem}

\noindent
This open question relates to a more
systematic effort to build a theory of \emph{homological VC
  dimension}~\cite[Conjectures~6 and~7]{kalai_conjectures}.
  There are some recent results for planar sets with connected intersections
(the case $d=2$ and $b=0$)~\cite{NervesMinors}.

Helly-type theorems are too many to list them all in this survey, but we wish to point out
at least one more variation.  \emph{Quantitative} Helly's theorem were introduced by B\'ar\'any, Katchalski, and Pach in \cite{baranykatchalskipach}. 
In this family of Helly-style theorems one is not content with a non-empty intersection of a family, but the intersections must have measurable or enumerable 
information in the hypothesis and the conclusion. Typical measurements include the volume, the diameter, or the number of points in a lattice. Motivated 
by applications in optimization, in the last two years several papers have been published on this subject, both for continuous \cite{Brazitikos2017,DeLoeraetalquant,Naszodi2016} 
and for discrete \cite{alievetal,Averkovetal-tightbounds,chestnutetal,Rolnick+Soberon} quantitative Helly-type theorems.  For other recent Helly-type theorems see \cite{amentaetal2017helly}. 

We next discuss the fifth remarkable theorem of our survey.

\subsection{Tverberg} \label{tv-section}

\noindent Tverberg-type theorems allow for the partition of finite
point sets so that the convex hulls of the parts intersect. In its
original form we have:

\begin{namedtheorem}[Tverberg theorem]
  \label{thm:Tverberg}
  Any set of at least $(r-1)(d+1)+1$ points in $\R^d$ can be
  partitioned into $r$ subsets whose convex hulls all have at least
  one point in common.
\end{namedtheorem}

\begin{figure}
  \centering
  \includegraphics[page=19]{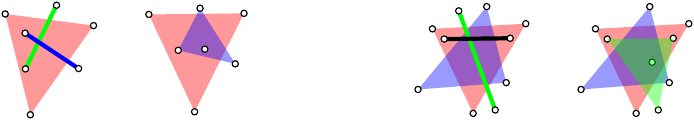}
  \caption{Tverberg's theorem in the plane: the two types of partitions for $r=3$ (left) and for $r=4$ (right).\label{f:Tverberg-partitions}}
\end{figure}

\noindent
Such a division in parts is often called a \emph{Tverberg partition};
{\blue see Figure~\ref{f:Tverberg-partitions}}. The case $r=2$ is
known as \emph{Radon's lemma}~\cite{originalRadon} and the case $d=2$,
but general $r$, was proven by Birch~\cite{Birch:1959ii}, before
Tverberg proved the general statement.

Tverberg's first proof of his theorem~\cite{Tverberg:1966tb} relies on
a deformation argument: start with a configuration with a known
Tverberg partition, and move the points continuously to the target
configuration. This process is such that, while the number of Tverberg
partitions may change, there will always be one present.  A simpler
proof consists in arguing that a partition of the point set minimizing
an adequate function must be a Tverberg partition; this idea, which
originates in B\'ar\'any's proof of the colorful Carath\'eodory
theorem, was gradually refined by Tverberg~\cite{Tverberg1981},
Tverberg and Vre\'cica~\cite{TverbergVrecica} and
Roudneff~\cite{Roudneff-ejc} (Roudneff minimizes the sum of the
squared distances between a point and the convex hulls of the
parts). Another proof, due to Sarkaria~\cite{Sarkaria:1992vt}, uses
bilinear algebra to deduce Tverberg's theorem from the colorful
Carath\'eodory theorem; the idea behind Sarkaria's proof was later
made simpler (using explicit tensors instead of number fields) and
more algorithmic by B\'ar\'any and Onn in~\cite{baranyonn-colorfulLP}.
Recently, B\'ar\'any and Sober\'on revisited these ideas and prove a
new generalization of Tverberg's theorem using affine combinations
\cite{Tverbergplusminus}.

\bigskip

Just as for Carath\'eodory's and Helly's theorems, there is an
\emph{integer Tverberg} theorem.  Its most recent version guarantees
that any set of at least $(r-1)d 2^d +1$ integer points in $\Z^d$ can
be partitioned into $r$ parts whose convex hulls have a point of
$\Z^d$ in common~\cite{integertverberg2017}.  The proof of this upper
bound goes as follows. From the integer Helly theorem, one can prove
that any finite set of integer points $S \subset {\mathbb Z}^d$ has an
integer \emph{centerpoint}: a point $\ve p$ such that for every
hyperplane $H$ containing $\ve p$, one of its half-spaces contains at
least $|S|/2^d$ points from $S$. Using this centerpoint, it is not
difficult to see that $(r-1)d 2^d +1$ points can be partitioned into
$r$ pairwise disjoint simplices all containing $\ve p$.  We will see
more about centerpoints in Sections \ref{centerpointsinoptima} and
\ref{half-space depth:subsec}.

The upper bound on the integer Tverberg number is not known to be
sharp, and the best lower bound of $2^d (r-1)+1$ is due to Doignon (this
result was communicated in \cite{Eckhoff:2000jw}). Recently in \cite{lowdimTverberg}
the authors showed that the Tverberg number in $\mathbb{Z}^2$ is exactly $4r-3$ when $r\geq 3$ and 
improved the upper bounds for the Tverberg numbers of $\mathbb{Z}^3$.

\begin{figure}
  \centering
  \includegraphics[page=20]{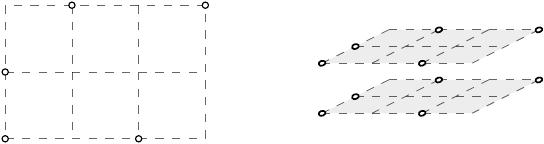}
  \caption{The lower bound for the integral Radon Theorem. Left: Five integral points in the plane with no integral Radon partition. 
  Right: A configuration of $k$ integral points in $\R^d$ with no Radon partition can be turned into a configuration of $2k$ integral points in $\R^{d+1}$ with no Radon partition.\label{f:Radon-integral}}
\end{figure}

The special case of bipartition, i.e., $r=2$, is called the
\emph{integral Radon theorem}.  A sharper upper bound of $ d2^d-d+3$
and a lower bound of $\frac54 2^{d} +1$ were established by
Onn~\cite{onn+radon} (see Figure~\ref{f:Radon-integral}). Even low
dimensional cases are hard.  The only sharp bound known, also due to
Onn, is for $r=d=2$: any six integral points in the plane have an
integral Radon partition. An upper bound of $17$ for the case $d=3$
was proven by Bezdek and Blokhuis~\cite{Bezdek2003}.

\begin{oproblem}
  Determine the exact value of the integer Tverberg numbers. In particular, is the
  integer Radon number for $d=3$ less than 17? Is it bigger than 11?
\end{oproblem}

\noindent
More generally, there is the notion of an \emph{integer quantitative
 Tverberg number}: any set of at least \linebreak $\pth{(2^d-2)\left\lceil
  \tfrac23(k+1)\right\rceil+2}(r-1)kd+k$ integer points can be
partitioned into $r$ parts whose convex hulls have $k$ integral points
in common~\cite{integertverberg2017}; similar results hold more
generally for sets that are discrete, i.e., intersect any
compact in only finitely many points, for instance the difference
between a lattice and one of its sublattices. Recent improvements on
the quantitative integer Helly theorem~\cite{Averkovetal-tightbounds,chestnutetal} 
leads to sharper upper bounds for the Tverberg version.

\begin{oproblem}
  Determine tighter lower and upper bounds on the integer quantitative Tverberg numbers. 
\end{oproblem}

\bigskip

Tverberg's theorem can be understood as stating that for any
\emph{linear} map from $\Delta_{(r-1)(d+1)}^{(d)}$, the
$d$-dimensional skeleton of the simplex with $(r-1)(d+1)+1$ vertices,
into $\R^d$ there must exist $r$ disjoint simplices whose images
intersect. This reformulation invites the question, going as far back as 1979 
(see ending of the important paper~\cite{Bajmoczy:1979bj}), whether the same
conclusion holds for all \emph{continuous} maps. In other words, is
there a topological Tverberg theorem? For $r=2$, this is the question
of non-embeddability discussed above in relation to topological Helly
theorems.

\begin{figure}
  \centering
  \includegraphics[page=21]{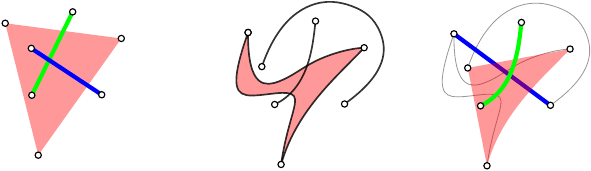}
  \caption{The topological Tverberg theorem in the plane ($r=3$). Left: A configuration of seven 
  points and its Tverberg partition into three parts. Center: The linear map used in the left is 
  deformed continuously so as to break the previous Tverberg partition, every edge not represented is kept straight. 
  Right: Nonetheless, a new Tverberg partition emerges.\label{f:Tverberg-topo}. 
  Note that the topological Tverberg conjecture remains open in dimension two when $r$ is not a prime-power. 
  See~\cite{50tverbergsurvey} and, for the two-dimensional case,~ \cite{topotverbergwinding}.}
\end{figure}

Positive answers were obtained first for $r=2$ (the \emph{topological
  Radon theorem}) by Bajm\'oczy and B\'ar\'any~\cite{Bajmoczy:1979bj},
then for $r$ prime by B\'ar\'any et al.~\cite{Barany:1981vh}, and for
$r$ a power of a prime by \"Ozaydin~\cite{Oza87} and independently,
but later, by Volovikov~\cite{Volovikov-TT} and
Sarkaria~\cite{Sarkaria-TT}.  Matou\v{s}ek~\cite[Chapters~5
  and~6]{matousek2003using} offers an accessible introduction to the
techniques behind the topological Tverberg theorem. 

For $r=2$,  the proof of the topological Radon theorem uses the notion of \emph{deleted product} of a geometric simplicial complex $\K$ with itself, defined as 
\[
{\K_{\Delta}^{2}}=\{ \sigma \times \tau : \sigma, \tau \in \K, \ \sigma \cap \tau = \varnothing \}.
\]
Now, for contradiction, suppose there exists a continuous map $f:|\K| \to \R^d$ with the property that points in distinct faces are mapped to distinct points.
It induces another continuous map $\tilde{f}:  | {\K_{\Delta}^{2}}| \to \s^{d-1}$  where $ \tilde{f}(x_1,x_2) =\frac{f(x_1)-f(x_2)}{\|f(x_1)-f(x_2)\|}.$
The map $\tilde{f}$ commutes with the central symmetries of $\s^{d-1}$ and ${\K_{\Delta}^{2}}$, where the central symmetry of ${\K_{\Delta}^{2}}$
is the exchange of the two components.  When $\K$ is the boundary of the $d$-dimensional simplex, 
${\K_{\Delta}^{2}}$ is homotopy equivalent to $\s^d$ (this is not trivial).
Since the Borsuk-Ulam theorem prevents the existence of
an antipodal map from $\s^d$ to $\s^{d-1}$, the continuous $f$ will have
two faces intersecting in its image, which gives us a contradiction.

More generally, to prove the topological Tverberg theorem when $r$ is a prime, 
one may start with a map from $\K=\Delta_{(r-1)(d+1)}^{(d)}$
to $\R^d$ with no $r$-wise intersection and use it to build another
map from an $r$-fold deleted product ${\K_{\Delta}^{r}}$, replace the antipodality by
the action of the symmetric group, and apply a generalization of the
Borsuk-Ulam theorem such as Dold's theorem~\cite{Dold:1983wr}. The
case of $r$ a prime power is technically more involved. (Note that
this outline leaves out some issues such as dimension-reduction
considerations~\cite{topotverbergwinding}.)

By the late 1990's, the widespread belief that B\'ar\'any's question
had a positive answer for every $r$ and $d$ was known as the
\emph{topological Tverberg conjecture}. It was only recently refuted
by Frick~\cite{frick15}, who completed an approach of Mabillard and
Wagner~\cite{Mabillard2014} building on \"Ozaydin's
work~\cite{Oza87}. In a nutshell, \"Ozaydin proved that an equivariant
map from the adequate $r$-fold product $\tilde{X}$ exists if and only
if $r$ is not a prime power. Mabillard and Wagner proposed an
isotopy-based approach to construct a map with no $r$-wise
intersection when such an equivariant map exists, but could only
develop it in codimension larger than what the topological Tverberg
conjecture allows. Frick overcame this codimension restriction,
producing the first series of counter-examples to the topological
Tverberg conjecture. The current state of affairs is that a
counter-example is known for every $r$ that is not a prime power and
for every $d \ge 2r$. See the survey \cite{50tverbergsurvey} for more
details.

\bigskip

To conclude our discussion of ``all things Tverberg'', let us highlight some 
natural variants inspired by Tverberg's theorem for which only very partial results are known.
One can find the most recent variants and extensions of Tverberg's theorem 
in \cite{barany+soberonsurvey, DeLoeraetal-tv-2018, Por-2018}.

A \emph{tolerant Tverberg theorem}, due to Sober\'{o}n and
Strausz~\cite{Soberon:2012er}, asserts that any set of
$(t+1)(r-1)(d+1)+1$ points can be partitioned into $r$ parts such
that, after deletion of any $t$ points, what remains is a Tverberg
partition. This bound was improved to $r(t+2)-1$ for $d=1$ and to
$2r(t+2)-1$ for $d=2$~\cite{Mulzer:2013je} (bound for $d=1$ is tight). 
Recently there have been two significant improvements,  Garc{\'i}a-Col{\'i}n 
et al. \cite{nataliaetal2017} gave an asymptotically tight bound for the tolerant 
Tverberg  theorem when the dimension and the size of the partition are fixed. 
Later, in \cite{pablo2018}, Sober\'on used the probabilistic method to give 
another asymptotic bound that is polynomial in all three parameters. 
Still we can ask for precise values.

\begin{oproblem}
  What is the smallest number $n$ such that any set of $n$ points in
  $\R^d$ has a Tverberg partition into $r$ parts that tolerates the
  deletion of $t$ points?
\end{oproblem}

\noindent
A related Carath\'eodory-type variation of Tverberg's
theorem~\cite{ABBFM09} considers $r$ linear maps $f_1,\dots,f_r$,
assumes that $f_1(e) \cap f_2(e) \cap \dots \cap f_r(e)$ is non-empty
for every $1$-dimensional edge $e$ of $\Delta_{(r-1)(d+1)}$, and
concludes the existence of disjoint faces in the simplex
$\Delta_{(r-1)(d+1)}$, $\sigma_1,\sigma_2,\dots,\sigma_r$ of
dimensions summing to $(r-1)(d+1)+1-r$ and such that $f_1(\sigma_1)
\cap f_2(\sigma_2) \cap \dots \cap f_r(\sigma_r)$ is not empty.
  
A conjectured \emph{relaxed} version of Tverberg's theorem, due to
Reay, goes as follows. Denote by $T(d,r,k)$ the minimum positive
integer number~$n$ such that any set of ~$n$ points $a_1,\ldots,a_n$
in~$\R^d$ (not necessarily distinct) admits a partition into~$r$
pairwise disjoint sets $A_1, \dots, A_r$ such that any size~$k$
subfamily of $\{\conv(A_1), \conv(A_2),$ $\ldots,\conv(A_r) \}$ has a
nonempty intersection.  Note that Tverberg's theorem says $T(d,r,r) =
(r-1)(d+1)+1$; Reay conjectured that the Tverberg constant is tight
even for smaller values of ~$k$:

\begin{conjecture}
  $T(d,r,k) =(r-1)(d+1)+1$ for all $2 \le k \le r$. 
\end{conjecture}

\noindent
Note that if all $\conv(A_i)$ intersect then they intersect $k$-wise,
so $T(d,r,k)$ is at most $T(d,r,r)$ and in particular $T(d,r,k)$ is
finite for any $k \le r$. Moreover, Helly's theorem ensures that
$T(d,r,k) = T(d,r,r)$ for every $d+1 \le k \le r$. The conjecture is
known to hold for $d + 1 \le 2k - 1$ or $k < r < \frac{d+1}{d+1-k}k$,
and some weaker bounds are known in several other cases; we refer the
interested reader to~\cite{asadaetal2} and the discussion therein.
 
Let us take this opportunity to mention another famous conjecture of
flavor similar to Reay's conjecture. A \emph{thrackle} is a graph
that can be drawn in the plane in such a way that any pair of edges
intersects precisely once, either at a common vertex or a transverse
intersection point.

\begin{conjecture}[Conway's Thrackle conjecture]
  For any thrackle, the number of edges is at most the number of vertices. 
\end{conjecture}

\noindent
The conjecture is known to hold if all edges are drawn as straight
line segments~\cite{E46} (it is akin to Reay's setup for $k=2$). We
refer the interested reader to the recent progress of Fulek and
Pach~\cite{fulek2010} and the discussion and references therein.

The next big open question was stated in 1979 by Sierksma (he offered an entire Dutch cheese as a prize for whoever
can solve this problem). In unpublished mimeographed notes he conjectured about the \emph{number} of distinct Tverberg 
partitions of a set of points guaranteed to exist for   $(r-1)(d+1)+1$ points in $\R^d$. 

\begin{conjecture}[Sierksma]
  Any set of $(r-1)(d+1)+1$ points in $\R^d$ has at least $((r-1)!)^d$
  distinct Tverberg partitions into $r$ parts.
\end{conjecture}

\noindent
We do not state the lower bounds here, as they are a bit cumbersome, but
lower bounds for the number of Tverberg partitions were first obtained when $r$ is prime by Vu{\v{c}}i{\'{c}} and {\v{Z}}ivaljevi{\'{c}} in~\cite{vucicrade}. They
used topological tools to settle this. Hell showed that these bound also hold when $r$ is a prime power  \cite{hell1}. Later Hell, without any topology, provided
better bounds for the case of the plane~\cite{Hell2008}. 

Last, but certainly not least, the \emph{colorful} Tverberg conjecture
was formulated by B\'ar\'any and Larman~\cite{Barany:1992tx} some 25
years ago, but only a few results are known today (see~\cite{barany+soberonsurvey,blamaszie,Blagojevicetal2015,colortv-ZivaljevicV92} and the references therein).

\begin{conjecture}\label{conjecture-colorful-Tverberg}
  Let $F_1, F_2, \ldots, F_{d+1} \subset \R^d$ be sets of $r$ points
  each. There exists a partition of $\bigcup_{i=1}^{d+1}F_i$ into $r$
  sets $A_1, A_2, \ldots, A_{r}$, of $d+1$ points each, such that every
  $A_i$ contains exactly one point from every $F_j$ and $\bigcap_{j=1}^r A_j \not= \varnothing$.
\end{conjecture}

\noindent
Further conjectures along the lines of the B\'ar\'any-Larman
conjecture, with colorful, discrete and quantitative flavors were
formulated by De Loera et al. ~\cite{integertverberg2017}.

\subsection{Computational considerations} \label{compu-convgeo}

We continue our discussion of computational issues begun in Section \ref{fixed-point-computing}.
We remark that the Carath\'eodory and Helly theorems,  in their classical real-valued versions from Section \ref{s:convexgeo}, 
are dual to each other and, essentially, 
if one has an algorithm to find the Carath\'eodory decomposition of a vector in terms of other vectors, one 
also has an algorithm for finding an intersection point for a family of convex sets. For Helly's theorem, one 
wishes to find a point in the intersection of convex sets. The problem of finding such an intersection point can be 
thought  of as a special case of the problem of  minimizing a convex function over convex sets. We will see 
explicit cases for this problem later in Section \ref{s:optimization}, where the convex sets are explicitly given by 
constraints (i.e., equations and inequalities), but for now all that we need to know is that: (a) a whole range of different algorithms for solving such problems 
exist,  (b) some of these algorithms are in fact efficient, and (c) depending on the type of input constraints (e.g., convex sets 
defined linear inequalities versus arbitrary constraints) one can be even more efficient  \cite{bertsekasbook,boyd2004convex}. By the convexity assumption, 
a local minimum is also a global minimum and, thanks to the Helly and Carath\'eodory theorems, 
there are nice necessary and sufficient conditions for when we have found an optimum. We discuss more about this in Section \ref{s:optimization}.

Convex optimization problems have been classified in levels of increased computational difficulty and different 
specialized algorithms are available (e.g., Least squares, linear programming,  Conic optimization, Semidefinite programming, etc). 
For example, Carath\'eodory's theorem, in the simplest real-valued form presented in Section \ref{caratheodory-thms}, can be formulated as a linear programming problem and
Helly's theorem for halfspaces is reducible  to linear programming too. The good news is linear programs have efficient algorithmic 
solutions, both in theory and in practice \cite{Sch86}. Still, even if we move to the most general version of (real-valued) Helly's theorem
of finite family of arbitrary convex sets, the challenge is to solve a \emph{convex programming} problem that has many types of algorithms. 
We cannot cover them here but we recommend  \cite{bertsekasbook,boyd2004convex} for excellent introductions to convexity algorithms and computational methods.

Compare now the good news above to the bad news involving the integer and mixed-integer versions of the Carath\'eodory and Helly theorems.  
We are now in the realm of  combinatorial, integer, and mixed-integer programming where, for the most part, the problems are not efficiently solvable. 
Solving combinatorial or mixed-integer optimization problems, that is, finding an optimal solution to such problems, can be a difficult task.  
Unlike linear programming, whose feasible region is a convex (polyhedral) set, in combinatorial problems, one must search a lattice of feasible points or,
 in the mixed-integer case, a set of disjoint half-lines or line segments to find an optimal solution. 
Thus, unlike linear programming, finding a global optimum to the problem requires us to relax, approximate, decompose the solution space,
and sometimes we are forced to enumerate all possibilities.
 
As an example of the higher computational difficulty of discrete versions of the Helly and Carath\'eodory theorems let us look at the problem of computing a Hilbert basis. 
This was a particularly simple case for the integer version of Carath\'eodory's theorem presented in Section \ref{caratheodory-thms}. Alas, it was  proved in 
\cite{durandetal} that deciding whether a given solution belongs to the Hilbert basis of a given system is coNP-complete. Thus, even in this tame case, the integral 
Carath\'eodory property is hard to realize computationally. 

We chose to highlight the colorful Carath\'eodory theorem because it
is so general and because it can be used to prove the original
Carath\'eodory theorem and many other existence theorems in
high-dimensional discrete geometry, such as Tverberg's theorem or the
centerpoint theorem (see Section \ref{s:datapoints} for details).
While the original Carath\'eodory's theorem can be cast as a linear
program and thus a solution can be implemented in polynomial time,
much less is known about the algorithmic complexity of its colorful
version. More precisely, the algorithmic colorful Carath\'eodory
problem is the computational problem of finding such a colorful choice
of elements as described in the theorem. Despite several efforts in
the past, the computational complexity of the colorful Carath\'eodory
problem in arbitrary dimension is still open.  In
\cite{PPAD-ColorfulCaratheodory}, Meunier et al. showed that the problem
lies in the complexity class PPAD. 

\begin{oproblem}
  What is the complexity of finding a colorful simplex under the 
hypotheses of the colorful Carath\'eodory theorem?
\end{oproblem}

This question was formulated for the first time by B\'ar\'any and Onn
in \cite{baranyonn-colorfulLP} and there they formulated a general family of related questions 
that come under the name \emph{colorful linear programming}. 

Meunier and Sarrabezolles~\cite{meunier-sarrabezolles2} have
shown that a closely related problem is PPAD-complete: given $d+1$ pairs of
points $P_1,\dots,P_{d+1}\in \Q^d$ and a colorful choice that contains
the origin in its convex hull, find another colorful choice of points
that contains the origin in their convex hull.

Since we have no exact combinatorial polynomial-time algorithms for
the colorful Carath\'eodory theorem, approximation iterative algorithms
are of interest. This was first considered in \cite{baranyonn-colorfulLP}, but other researchers, e.g.,  
\cite{mulzerstein2015},  have approached this problem too.

Let us now speak about computational complexity of Tverberg's theorem.
Sarka\-ria's proof of Tverberg's theorem (later simplified by B\'ar\'any
and Onn \cite{barany1995caratheodory}) gives a polynomial-time way to
compute a Tverberg partitions from a colorful Carath\'eodory choice
with the origin in its convex hull. In this way, the computational
issues about Tverberg's theorem are closely connected to computational issues
regarding the colorful Carath\'eodory theorem. Since one can calculate a Tverberg
partition from a colorful selection of Carath\'eodory, one can show
Tverberg's theorem belongs to the class PPAD.  One of the simplest, yet most
frustrating, aspects of Tverberg's result is that it is not clear how
to find a Tverberg partition. So it is natural to ask:

\begin{oproblem}
  Is there a polynomial-time algorithm to find a Tverberg partition
  when one exists? That is, given $n=(d+1)(m+1)+1$ points in $\R^d$,
  compute, in time polynomial in $n$, a partition into $m$ parts with
  intersecting convex hulls.
\end{oproblem}

Since finding an $m$-Tverberg partition is an open question,
approximate versions of Tverberg's theorem are of interest.  Mulzer et
al. \cite{mulzerstein2015} designed a deterministic algorithm that
finds a Tverberg partition into $n/4(d+1)^3$ parts in time $d^{O(\log d)}
n$.  This means that for every fixed dimension one can compute an
approximate Tverberg point (and hence also an approximate centerpoint)
in linear time. Rolnick and Sober\'on \cite{rolnick+soberon2} proposed probabilistic
algorithms for computing Tverberg partitions into $n/d^3$ parts with 
error probability $\epsilon$, and with time complexity that is weakly polynomial 
in $n,d, \log(\frac{1}{\epsilon})$.

\section{Games and fair division}
\label{s-games+fairness+independence}

Mathematics and the social sciences have had rich interactions since
Condorcet's seminal work on the analysis of voting systems. The
relevance of (combinatorial) convexity and topology to this
interdisciplinary research was first established in the 1940-50's
through the work of Nash, von Neumann, Gale, Shapley, Scarf and many
others, and it has been confirmed in the following decades in the
development of fair-division algorithms and computational social
choice (see e.g., \cite{bramsetal1,bramsetal2, collectionmathdecisionsgames,Nisanetal-2007} 
and the many references therein). This section shows how our five discrete theorems appear in these
topics too.

\subsection{Strategic games}

Game theory studies a broad range of ``games'' that model situations
where several agents collaborate or compete. Strategic games model
the situations where $N$ players interact by each choosing from
finitely many ``strategies'' to play and enjoy a payoff depending on
the strategies chosen by all players (his or her choice included).
Formally, each player is modeled by the pair $(S_i,u^i)$ where
$S_i$ is a finite set of
strategies available to him or her and a payoff function $u^i :
S_1\times\cdots\times S_N \to \R$. A central theme in the theory of
strategic games is the search for equilibria, where each player's
choice is the best response to the other players' choices.

\subsubsection{Nash equilibria}

Formally, a \emph{Nash equilibrium in pure strategies} is a choice of
strategy for each player $s_1 \in S_1, \ldots, s_N \in S_N$ such that
for $i=1,\ldots,N$ and all $g \in S_i$
\[  u^i(s_1, \ldots, s_{i-1},s_i,s_{i+1}, \ldots, s_N) \ge  u^i(s_1, \ldots, s_{i-1}, g, s_{i+1}, \ldots, s_N).\]
Let us illustrate pure Nash equilibria with the \emph{max-cut game},
where an arbitrary graph $G = (V, E)$ is fixed and each vertex $x \in
V$ represents a player. Each player $x$ chooses from two strategies
$S_x=\{1, -1\}$, and his/her payoff function is the number of neighbors of
$x$ in $G$ with a different strategy:
\[ u^x(s_1, \dots, s_{|V|}) = |\{y \colon xy \in E, \ \hbox{and} \  s_x \not= s_y\}|.\]
Any bipartition of $V$ that maximizes the number of edges between the
two parts, also called \emph{maximum cut}, is a pure Nash
equilibrium. Indeed, if a player could strictly increase his or her
payoff by switching strategy, then this switch would increase the
value of the cut.

\medskip

Now, unfortunately not every $N$-player game has a \emph{pure} Nash
equilibrium. For example, consider the \emph{matching penny
  game}. Players Alice and Bob simultaneously select heads or tails of
a coin. If the choice is the same, then Alice wins one penny and Bob
loses a penny. If they choose differently, then Bob wins a penny and
Alice loses a penny. Each player thus has a choice of two strategies,
and the payoffs for each player can be recorded in two $2 \times 2$
matrices ($A$ and $B$, for each player).
Note that the pure strategies alone offer no Nash equilibrium as this
is a winner-takes-all situation. (In the matching penny game, the payoff
matrices can be put together to show $A+B=0$; this is an example of a
\emph{zero-sum game}, a notion to which we come back later.)

\bigskip

It turns out that equilibria always exist when considering a
\emph{randomized choice of the strategies}. We now allow more freedom
to the players by making them choose  not a single strategy, but a probability
distribution over all their strategies. Once choices are made, a
random strategy is selected for each player from his or her
distribution, and the (random) payoff is determined. Each player then
wants to maximize her or his expected payoff. Formally, a \emph{mixed
  strategy} for player $i$ is a probability distribution $m_i$ on the
set of pure strategies $S_i$. For instance, in the matching penny game, 
this means that each player decides his or her move according to
a biased coin flip and is free to choose the bias. Note that the
mixed strategies include all the pure strategies as a special case too.

The set of all possible mixed strategies are the vectors that lie in the convex
polytope $M=\prod \Delta_{|S_i|-1}$. We define a product
probability measure on $S=S_1\times\cdots\times S_N$ by
$P_m(s)=\prod_{i=1}^N m_i(s_i)$ where $s=(s_1,s_2,\dots,s_N)$.
Therefore the expected payoff for this probability distribution
$P_m(s)$ for the $i$-th player is
 $$U^i(m_1, \ldots, m_N)=U^i(m)=\sum_{s=(s_1,s_2,\dots,s_N) \in S} P_m(s) u^i(s).$$
The mixed strategies $m=(m_1, \ldots, m_N)$ form a \emph{Nash
  equilibrium in mixed strategies} if for each player $i$ and for all 
probability distributions $p$ on $S_i$, modifying $m_i$ to $p$ does not increase
the expected payoff with respect to the choices of other players, that is
\[ U^i(m_1, \ldots,m_i, \ldots, m_N) \ge U^i(m_1, \ldots, m_{i-1}, p, m_{i+1},   \ldots, m_N).\]
The literature has plenty of examples of two-player
games~\cite[$\mathsection 8.1$]{matousek2007understanding}. Thinking
about three or more players is more delicate, as illustrated by Nash's
three-man poker game~\cite[p. 293]{nashpaper2}.

 The existence of Nash equilibria for mixed strategies -- the theorem,
 for which John Nash received the Nobel prize -- is one of the most
 celebrated applications of combinatorial topology, following from
 Brouwer's theorem. See \cite{nashpaper,nashpaper2}.

\begin{theorem}[Nash's theorem]\label{thm:Nash}
  Every $N$-player game with continuous payoff functions 
  has at least one Nash equilibrium in mixed strategies.
\end{theorem}

\noindent
Nash's original, very short, proof~\cite{nashpaper} makes strong use
of the combinatorial topology and convexity arguments.  It considers
the set-valued function that maps each $N$-tuple of mixed strategies
$m=(m_1, \ldots, m_N)$ to the set of $N$-tuples $(t_1, \ldots, t_N)$
where $t_i$ is a best response to $(m_1, \ldots, m_{i-1}, m_{i+1},
\ldots, m_N)$. Because the probability distributions on $S_i$ are the
points of the simplex $\Delta_{|S_i|-1}$, Kakutani's theorem (we saw
this in Section \ref{s:topology} after Brouwer) ensures this function
has a fixed point, which is the desired equilibrium.

Nash gave a second proof using Brouwer's fixed point theorem
\cite{nashpaper2}. For this, Nash constructed a continuous function $f$
from the polytope $M$, associated to the game above,
into itself.  For $m \in M$, we define $f(m)$ component-wise as follows: 
$f(m)=(f_i(m),\dots,f_N(m))$, and each entry $f_i(m)$ is equal to $(f_{i1}(m),f_{i2}(m),\dots,f_{it_i}(m))$ where
 {\small
 $$f_{ij}(m)=\frac{m_{ij} + \max(0, u^i(m_1, \ldots, m_{i-1},
     s^{(i)}_j, m_{i+1}, \ldots, m_N)- u^i(m_1, \ldots,m_i, \ldots,
     m_N) )} {1 + \sum_{k=1}^{|S_i|} \max(0, u^i(m_1, \ldots, m_{i-1},
     s^{(i)}_k, m_{i+1}, \ldots, m_N)- u^i(m_1, \ldots,m_i, \ldots,
     m_N) )},$$ } where the $s_j^{(i)}$'s are the pure strategies available for player $i$. 
 \medskip
 Nash showed this function $f$ is a
 multivariate continuous map from the polytope $M$ into itself and
 thus, by Brouwer's fixed point theorem, it must have at least one
 fixed point. Nash then went to show that any fixed point of this
 function is in fact a Nash equilibrium.

 \bigskip

 The polytope $M=\prod
  \Delta_{|S_i|-1}$ is actually a Cartesian product of simplices, sometimes called a \emph{simplotope}. The special structure of
  simplotopes has been exploited for the computation of fixed-points (see \cite{simplotopes,todd1993new} and references there) and in the algebraic solution 
  of equilibrium problems via polynomial equations \cite{mclennan2005}. There is one more reason why knowing $M$ is a polyhedron is useful. 
  Despite its geometric beauty and intricacy, there are several modeling limitations with the notion of Nash equilibria. For example, 
the assumption behind Nash mixed strategies equilibria is that the choices 
of each player are independent of those of his/her opponents, but that may not hold. Alternative mathematical models 
that adjust the definition of Nash mixed strategies to allow \emph{correlated equilibria} appear in \cite{aumann}. In other words, the expression
for the payoff function does not use anymore the easy product probability structure of a simplotope, but can have a more complicated polyhedral geometry. 
For example, it is known  that all correlated equilibria are described by a finite set of linear inequalities and that it is non-empty, independently of 
Nash theorem. See \cite{gilboa+zemel} and its references.

\subsubsection{Two-player games}\label{s:nash-complexity}

Nash equilibria have been largely investigated in the area of
algorithmic game theory, see for instance the introductory chapter of~\cite[$\mathsection 2$]{Nisanetal-2007}. 
We only discuss here some of their relations to combinatorial topology and convexity.

\bigskip

The Nash equilibria for two players can be formulated as a
\emph{linear complementarity problem}, the theory of which subsumes
both linear programming and two-player game theory (an introduction is
in \cite{cottleetalbook}).  Let $A$ and $B$ denote the $m \times n$
payoff matrices of the first and second players, respectively. By
definition, a pair $(\xx^*, \yy^*) \in \Delta_{m-1} \times
\Delta_{n-1}$ is a Nash equilibrium if and only if
\begin{equation}\label{eq:NELP}
  \xx^{*T} A\yy^* \geq\xx^T A\yy^* \quad \forall \xx \in \Delta_{m-1}
  \quad \textrm{ and } \quad\xx^{*T} B\yy^* \geq\xx^{*T} B\yy \quad
  \forall \yy \in \Delta_{n-1}.
\end{equation}
Here comes the linear complementarity formulation:

\begin{proposition}\label{p:lincomp}
  The pair $(\xx^*, \yy^*) \in \Delta_{m-1} \times \Delta_{n-1}$
  satisfies~\eqref{eq:NELP} if and only if there exist $\ve u^* \ge \ve
  0$, $\ve v^* \ge \ve 0$, $s \ge 0$, and $t \ge 0$ such that
  \[ A\yy^* + \ve u^*  =   s \ve 1, \quad 
  B^T\xx^* + \ve v^* = t \ve 1, \quad \text{and} \quad \xx^{*T}\ve u^* =
  \yy^{*T}\ve v^* = 0.\]
\end{proposition}

\noindent
Since all vectors are non-negative, the conditions on inner products
imply that the supports of $\xx^*$ and $\ve u^*$ are disjoint, and
similarly for $\yy^*$ and $\ve v^*$, hence the aforementioned
complementarity. The proof of Proposition~\ref{p:lincomp} goes as
follows. Start with a Nash equilibrium $(\xx^*, \yy^*)$ and let $s =
\xx^{*T} A\yy^*$ and $t = \xx^{*T} B\yy^*$; the values of $\ve u^*$
and $\ve v^*$ are forced and the complementarity conditions follow
from multiplying the first equation by $\xx^{*T}$ and the second
equation by $\yy^{*T}$. Conversely, the complementarity $\xx^{*T}\ve
u^* =0$ and $\xx^* \in \Delta_{m-1}$ yield that $s = \xx^{*T} A\yy^*$ (and
similarly $t = \xx^{*T} B\yy^*$); the positivity of $\xx^T \ve u^*$
for any $\xx \in \Delta_{m-1}$ implies that $\xx^*$ is a best response to
$\yy^*$; by a similar argument, $\yy^*$ is a best response to~$\xx^*$.

\bigskip

The linear complementarity formulation of Proposition~\ref{p:lincomp}
yields, after adequate rescaling, the standard method to compute
two-player Nash equilibria: find a non-trivial solution (i.e.,
other than $\xx = \yy = \ve 0$) to the linear complementarity system
\begin{equation}\label{eq:lincompspe}
  \pth{\begin{matrix} A & I_m\end{matrix}} \pth{\begin{array}{c} \yy
    \\ \ve u\end{array}} = \ve 1, \quad \pth{\begin{matrix} I_n &
    B^T \end{matrix}} \pth{\begin{array}{c} \ve v \\ \xx \end{array}}
  = \ve 1, \quad \text{and} \quad \xx^T\ve u = \yy^T\ve v = 0.
\end{equation}
The standard method to solve a linear complementarity system of this
form is the \emph{Lemke-Howson pivoting algorithm}~\cite{lemkehowson64} which
operates on feasible bases. A \emph{feasible basis} of a linear system
with non-negativity constraints is a set of indices of columns whose
induced sub-matrix has same rank as the system and defines a
non-negative solution. The trivial solution $\xx=\yy=\ve 0$ gives
feasible bases $I_1$ for $\pth{\begin{array}{c} \ve v
    \\ \xx \end{array}}$ and $J_1$ for $\pth{\begin{array}{c} \yy
    \\ \ve u\end{array}}$, and these bases are \emph{complementary}, i.e., have disjoint supports. Pick an (arbitrary) element $k_1
\notin I_1$. By Carath\'eodory's theorem (in the form of
Proposition~\ref{p:cara++}), $I_1 \cup \{k_1\}$ contains another
feasible basis $I_2$ for the system $B^T\xx + \ve v = \ve 1$.
Switching from $(I_1, J_1)$ to $(I_2,J_1)$ is our first pivot
step. Remark that $I_2$ and $J_1$ are no longer disjoint, as they
share $k_1$. To remedy this, note that the part of $I_2$ corresponding
to $\ve v$ lost some element $k_2$ which is also absent from the part
of $J_1$ corresponding to $\yy$. We can thus make $k_2$ enter $J_1$
without further degrading the complementarity; Carath\'eodory's theorem
ensures that $J_1 \cup \{k_2\}$ contains another feasible basis $J_2$,
and we pivot from $(I_2,J_1)$ to $(I_2,J_2)$. The part of $J_2$
corresponding to $\ve u$ lost some element $k_3$. If $k_3 = k_1$ then
$I_2$ and $J_2$ are complementary, and determine our non-trivial
solution, otherwise we continue pivoting by making $k_3$ enter $I_2$,
etc. until a pivot makes $k_1$ exit one of the bases; the pair at
hands is then complementary.

The Lemke-Howson algorithm is guaranteed to terminate under the
non-dege\-ne\-ra\-cy assumption that $\ve 1$ is not a positive
combination of less than $m$ columns of $\pth{\begin{matrix} A &
    I_m\end{matrix}}$ or less than $n$ columns of $\pth{\begin{matrix}
    I_n & B^T \end{matrix}}$. This follows from a non-degenerate
Carath\'eodory theorem (after Proposition~\ref{p:cara++}): any point
in the convex hull of $d+2$ points of $\R^d$ that is non-degenerate, i.e., is not in the convex hull of some $d$ of them, lies in the
convex hulls of exactly two $(d+1)$-element subsets. (Note that the
two non-degeneracy assumptions stated above are equivalent via the
convex/conic change of viewpoint.)  Now, for any pair~$(I, J)$ of
non-complementary feasible bases encountered by the algorithm, there
is exactly one index $k$ not in $I \cup J$. A pivot step makes $k$
enter either $I$ or $J$; in each case, the non-degenerate
Carath\'eodory theorem yields exactly one other pair of feasible
bases. Every pair~$(I, J)$ of non-complementary feasible bases
encountered by the algorithm thus has exactly two neighbors through
pivot steps. Since the trivial solution $(I_1,J_1)$ has exactly one
neighbor, there is no place where the walk can loop back. Remark that
this algorithm gives an alternate (constructive) proof of the
existence of a Nash equilibrium for two-player games.

The argument that proves termination also reveals that, from a
computational complexity point of view, linear complementarity systems
of the form of Equation~\eqref{eq:lincompspe} are in the PPAD
class. Let us point out that the Lemke-Howson algorithm can be
understood as a Sperner-type search for a fully-labeled simplex in a
pseudomanifold. Let $V = \{\pm1, \pm2,$ $\ldots, \pm(m+n)\}$, where
positive integers are understood as indices of columns of
$\pth{\begin{matrix} A & I_m\end{matrix}}$ and negative integers are
understood as (minus) indices of columns of $\pth{\begin{matrix} I_n &
    B^T \end{matrix}}$. Consider the simplicial complex $\K$ on~$V$
whose maximal simplices are the union of the complement of a feasible
basis of $\pth{\begin{matrix} A & I_m\end{matrix}}$ and the complement
of a feasible basis of $\pth{\begin{matrix} I_n &
    B^T \end{matrix}}$. The non-degenerate Carath\'eodory theorem
spelled out above ensures that $\K$ is a pseudomanifold without boundary. If every $i\in
V$ is labeled by its absolute value $|i|$, the fully-labeled simplices
correspond exactly to the complementary feasible bases.

\bigskip

As a byproduct of the linear complementarity formulation of
Proposition~\ref{p:lincomp}, we also get that the problem of computing
a Nash equilibrium for two players is well-posed from the point of
view of computational complexity: if the input involves only rational
data, there is an equilibrium that involves only rational data and has
encoding size polynomial in the input size (see for instance the
discussion in the survey of McKelvey and
McLennan~\cite{mckelvey1996computation}). This is in sharp contrast
with the case of three or more players: Nash's three-player poker
game~\cite{nashpaper2} shows that a three-player game with finitely
many strategies and rational payoff arrays may have a (unique) Nash
equilibrium with irrational coordinates.

The Lemke-Howson algorithm has exponential complexity in the
worst-case~\cite{savani-vonstengel} and solving general linear
complementarity problems is NP-hard~\cite{Chung1989}. 
The problem of computing a Nash equilibrium for two players is
PPAD-complete~\cite{chen+deng+teng2009} (see
also~\cite{daskalakis+goldberg+papadimitriu2006}). The intractability
for games with three or more players is even more stringent, as many
decision problems are $\exists \R$-complete, i.e., as difficult
as deciding the emptiness of a general semi-algebraic set; this
includes in particular deciding whether a $3$-player game has a Nash
equilibrium within $\ell_\infty$-distance $r$ from a given
distribution~$\ve x$~\cite[Corollary~5.5]{Schaefer2015} or the existence of
more than one equilibrium or of equilibria with payoff or support
conditions~\cite{gargetal2015}. Behind this $\exists \R$-completeness
lurks a more daunting fact: Datta's universality
theorem~\cite{Datta03} asserts that arbitrarily complicated
semi-algebraic sets can be encoded as sets of Nash equilibria
(formally: every real algebraic variety is isomorphic to the set of
mixed Nash equilibria of some $3$-player game). Whether the $\exists
\R$-completeness results stated above could follow from Datta's proof
is an interesting open question~\cite[Remark~5.6]{Schaefer2015} .

\begin{oproblem}
  Can Datta's universality theorem be improved to yield an efficient
  polynomial-time, reduction between any semi-algebraic set and the
  Nash equilibria of a game?
\end{oproblem}

\noindent
Another open question is whether both few players \emph{and} few
strategies per player already give rise to universality.

\begin{oproblem}
  Is there a universality result for the set of Nash equilibria of
  games with a constant number of players and a constant number of
  strategies?
\end{oproblem}

\subsubsection{Zero-sum games.}

The two-player games where what is won by a player is lost by the
other are called \emph{zero-sum} games; this is the case when the payoff matrices
satisfy $A=-B$ in the formulation~\eqref{eq:NELP}. In this special case, it is customary to
consider that one player aims at maximizing the payoff while the other
tries to minimize it. Nash's theorem then asserts that there exist $\xx^* \in
\Delta_{n-1}$ and $\yy^* \in \Delta_{m-1}$ such that
\begin{equation}\label{eq:nzs}
  \begin{array}{ll} \forall \yy \in \Delta_{m-1}, & {\xx^*}^T A \yy \ge {\xx^*}^T A \yy^* \\
     \forall \xx \in \Delta_{n-1}, & {\xx}^T A \yy^* \le {\xx^*}^T A  \yy^*.
    \end{array}
\end{equation}

\noindent
This readily implies
\begin{equation}\label{eq:minmax}
  \begin{aligned}
  {\xx^*}^T A \yy^* \ge & \max_{\xx\in\Delta_{n-1}} {\xx}^T A \yy^*  \ge \min_{\yy\in\Delta_{m-1}} \max_{\xx\in\Delta_{n-1}}  {\xx}^T A \yy \\
  & \ge \max_{\xx\in\Delta_{n-1}}  \min_{\yy\in\Delta_{m-1}} {\xx}^T A \yy   \ge \min_{\yy\in\Delta_{m-1}} {\xx^*}^T A \yy \ge {\xx^*}^T A \yy^*.
\end{aligned}
\end{equation}
The only inequality that does not follow from Equation~\eqref{eq:nzs},
the central one, holds in fact for arbitrary bivariate functions (see Section~\ref{s:lpkkt}). 
Altogether, this yields the min-max theorem of von Neumann.

\begin{theorem}\label{thm:VonNeumann}
  For any $A \in \R^{m\times n}$, 
   $$\max_{\xx\in\Delta_{n-1}}\min_{\yy\in\Delta_{m-1}}\xx^TA\yy=
  \min_{\yy\in\Delta_{m-1}}\max_{\xx\in\Delta_{n-1}}\xx^TA\yy.$$
\end{theorem}

\noindent
In game-theoretic language, the real number $\xx^{*T}A\yy^{*}$ is the
\emph{value} of the game. Von Neumann's theorem has a nice ``asynchronous''
interpretation: for any choice of a strategy by the minimizing player,
the maximizing player can respond so as to ensure a payoff at least
the value of the game. Moreover, if the maximizing player cares only
about achieving the value of the game, the strategy $\xx^*$ will work
regardless of what the opponent plays. (Of course, symmetric
statements hold for the minimizing player.) In zero-sum games, in each
of the (possibly many) Nash equilibria, every player gets the same
payoff. This is specific to zero-sum games and fails already for
broader types of two-player games. A classical result of Dantzig
\cite{dantzig-minmax} says that the minmax identity of von Neumann's
theorem can be proved, without help of combinatorial topology, via
linear programming duality and is thus polynomially solvable. We will
discuss more about this in Section~\ref{s:optimization}. More
generally, if the rank of $A+B$ is constant then the problem is
polynomial~\cite{kannan+theobald}.

\bigskip

Formulation~\eqref{eq:NELP} suggests an approximate relaxation and
the approximate Cara\-th\'e\-odory theorem~\ref{thm:caratheodory-p}
provides a positive complexity result. We say a mixed strategy pair
$(\xx,\yy)$ is {\em an $\varepsilon$-Nash equilibrium} if
$$\xx^T A\yy \geq\ee_i^T A\yy - \varepsilon \quad \forall \ i \in [n]
\qquad \textrm{ and }\qquad \xx^T B\yy\geq\xx^T B\ee_j - \varepsilon
\quad \forall \ j \in [n].$$ Intuitively, a mixed strategy pair is an
$\varepsilon$-Nash equilibrium if no player can benefit more than
$\varepsilon$, in expectation, by a unilateral deviation.  

The case when $A+B=0$ is precisely the case of zero-sum games, for
which we know efficient algorithms exist. Barman, using the
approximate Carath\'eodory theorem presented in Theorem
\ref{thm:caratheodory-p}, provided an extension in \cite{Barman}.

\begin{theorem}
  \label{thm:nash}
  Suppose that all entries of the payoff matrices $A, B$ lie in
  $[-1,1]$.  If the number of non-zero entries in $A+B$ is at most
  $s$, then an $\varepsilon$-Nash equilibrium of $(A,B)$ can be
  computed in time $n^{ O \left( \frac{\log
      (\max{(s,4)})}{\varepsilon^2} \right)}$.
\end{theorem}

\noindent
This, in particular, gives a polynomial-time approximation scheme for
Nash equilibrium in games with fixed column sparsity $s$. Moreover,
for arbitrary bi-matrix games -- since $s$ can be at most $n$ -- the
running time of this algorithm matches the best-known upper bound,
which was obtained by Lipton, Markakis, and Mehta \cite{liptonetal}.

\subsection{Two fair-division problems: cakes and necklaces} \label{cakes+necklaces}

In various situations players are eager to divide goods in a ``fair
way''. There are several examples of such fair-division problems where
our five theorems (or their variations) play a key role. We review
some famous examples, all with combinatorial-topological
proofs. Before we start we remark there are other interesting
mathematical challenges arising in distributing resources, some we
will not cover here, such as \emph{gerrymandering}, which is the
practice of drawing political maps to gain an advantage. See e.g.,
\cite{soberon-gerrymanderingetc} for connections to the five theorems
in this survey.

\subsubsection*{Cake cutting} 

A cake is a metaphor for a heterogeneous, divisible good, such as a
piece of land or an inheritance. We consider the problem of dividing a
cake between $r$ players in such a way that each player prefers his or
her part to any other part. We call this \emph{envy-free}. Let us
point out that the literature about fair division, including this and
other types of cake-cutting problems, is both old and huge; see
e.g.,~\cite{brams+taylor,bramsetal1,
  bramsetal2,roberson+webb,rothe-surveys, Steinhaus}. For example,
 one of the first envy-free division results was shown by
Dubins and Spanier \cite{dubins+spanier}.

One setting where pieces are connected is the following:
The cake is identified with the unit interval $\left[ 0, 1
  \right]$ and a {\em division} of the cake into $r$ pieces is an
$r$-tuple $\boldsymbol{x}=\left( x_{1}, \ldots , x_{r} \right) $, with
$x_{j}\geq 0$ for all $j \in[r]$ and $\sum_{j=1}^{r} x_{j} = 1$; in
other words, a division is a point $\xx$ in the $(r-1)$-dimensional
simplex $\Delta_{r-1}$. Here, $x_{j}$ represents the size of the
$j$-th piece, when ordered from left to right. The preferences of
player $i$ are modeled by a function $P^i$ mapping each division $\xx
\in \Delta_{r-1}$ to a nonempty subset of $[r]$ (indexing the pieces
that he or she prefers). A division $\xx$ is \emph{envy-free} if there
exists a choice of pairwise distinct indices, one from each
$P^i(\xx)$. 

It is natural to assume that the set $P^i(\xx)$ of preferences of
player $i$ never contains the index of a piece of size zero, (i.e.,
all players are hungry) and it is common to suppose that the
$P^i$'s are \emph{closed}, that is if $\lim_{n \to \infty} \xx_n =
\xx$ and $j \in P^i(\xx_n)$ for every $n$, then $j \in P^i(\xx)$. 
Stromquist~\cite{stromquist1} and Woodall \cite{woodall} proved the
following result independently.

\begin{theorem}\label{thm:division}
  Under the assumptions that all $r$ players are hungry and the
  preference functions $P^i$ are closed, there exists an envy-free
  division with connected pieces. 
\end{theorem}

\noindent
The original proof relies on the KKM theorem we saw in Section
\ref{s:topology}.  An unpublished proof due to Forest Simmons was
improved and adapted by Su~\cite{Su99}.  The idea is to refine
the usual derivation of Brouwer's theorem from Sperner's lemma.  Take
a sequence $(\T_n)$ of triangulations of $\Delta_{r-1}$ whose
edge-length goes to $0$. Assign every vertex to a player in a way that
every full-dimensional simplex of $\T_n$ has a vertex assigned to each
player; this may not be possible for any sequence of triangulations,
but taking iterated barycentric subdivisions does the job.  We label
or color every vertex $\xx$ assigned to player $i$ by some (arbitrary)
element from $P^i(\xx)$.  The assumption that players are hungry
ensures that this is a Sperner labeling. The limit of a converging
subsequence of fully-labeled simplices is, by the closedness
assumption, an envy-free division.

It may be disappointing in practice that one only gets an iterative
infinite process converging to an envy-free division, but it has been
shown in \cite{stromquist2} that there exists no finite procedure, if
you require each person to get a connected piece (i.e., the cake is
cut by a minimal number of cuts).  In fact, Aziz and Mackenzie
\cite{Aziz+Mackenzie2016} showed that there is in fact a bounded
finite procedure for $r$-person cake cutting bounded by a huge number of
steps, but this would involving breaking the cake into a ridiculous
number of pieces.  Thus for now one cannot get minimal number of cuts
in a finite procedure, but if you allow a lot of cuts, you get a
division that is impractical as it destroys the cake. This is why an
infinite process converging to an envy-free solution with a minimal
number of cuts makes sense. Let us comment that the difficulty of the
process is not so surprising perhaps. It is known that, in the
polynomial-time function model, where the utility functions are given
explicitly by polynomial-time algorithms, the envy-free cake-cutting
problem has the same complexity as finding a Brouwer's fixed point,
or, more formally, it is PPAD-complete~\cite{DeQiSa12}.

\medskip

The polytopal version of Sperner's lemma (Proposition
\ref{p:polytopalsperner} in Section \ref{s:convexgeo}) has many
interesting game theoretic applications.  Su \cite{su-hex} recently
gave a simple elegant proof that Hex does not end in a draw using it.
Cloutier, Nyman, and Su \cite{multicake} applied the polytopal
Sperner's lemma to the \emph{multi-cake multi-player fair-division}
problems. In this type of problems the players have several cakes to
choose from, but choices from one cake influence each other, e.g., for
a player the amount of vanilla cake may influence how much chocolate
cake to order, or after some vanilla cake the player may not want any
chocolate cake. Cloutier et al. asked whether there exists an integer
$r(q, m)$, independent of the preferences, such that there exists an
envy-free division of the $m$ cakes not requiring to divide each cake
into more than $r(q, m)$ pieces, some of which are assigned to each of
the $q$ players (some of the pieces can remain unassigned). Note that
Theorem \ref{thm:division} for a single cake asserts that $r(q, 1) =
q$. They also used the polytopal version of Sperner's lemma, and they
proved the existence of $r(2, 2)$ and $r(2, 3)$ and that $r(2, 2) = 3$
and $r(2, 3)\leq 4$. This means that two cakes can be divided into
three pieces each in such a way that two players receive the pieces
and everyone is satisfied with the fairness of division. Moreover they
asked whether $r(2, m)\leq m + 1$. Recently Lebert, Meunier, and
Carbonneaux \cite{lebertetal} have shown $r(2,m)$ exists for any $m
\geq 2$ and its value is at most $m(m-1)+1$. They used again the
polytopal version of Sperner's lemma and an inequality between the
matching number and the fractional matching number in $m$-partite
hypergraphs. Similarly, they showed $r(3, 2)$ exists and $r(3, 2) \leq
5$.  Several interesting open questions remain, consider for example:

\begin{oproblem}
  Can the bound on $r(2,m)$ be improved? Can one assure the existence
  of $r(q,m)$ for all values of $q,m$?
\end{oproblem}

Finally, there are other surprising variations of the cake-division
theorem.  Consider one where there are no cuts involved and one
divides objects that are not physically divisible. Consider a house
with $n$ rooms and a total rent amount to be divided among $n$
roommates.  Assume that for each possible division of the rent amount
each roommate can point to one or more rooms as preferred. Then, the
theorem proved by Su in \cite{Su99} states that there exists a
division of the rent and an assignment of rooms to each participant,
such that each player receives one of his/her preferred rooms.
Similar assumptions on the preference function that we made for
cake-cutting before Theorem \ref{thm:division} must hold again to make
this happen, for instance, now it is assumed that each roommate
prefers a free room over paying rent.  Once more the proof of this
theorem is grounded in Sperner's lemma. In most results on
fair-division it is assumed that no player is happy with an empty
piece of cake, but imagine a part of the cake is undesirable (burnt
cake anyone?) another recent variation in
\cite{Meunier+Zerbib,Segal-Halevi} considers the possibility players
may prefer an empty piece.

 It is well-known that Gale's colorful KKM theorem (see Section
 \ref{subsec:continuousversionspernertucker}) has interesting
 applications in economics, e.g., for the existence of economics
 equilibria. Now, Asada et al. \cite{asadaetal} used Gale's colorful
 KKM theorem to prove an extension, by Woodall \cite{woodall}, of Theorem \ref{thm:division}  
 and a similar extension of the rental-harmony result of Su. It turns
 out that there are envy-free cake divisions for any number of
 players, even if the preferences of one person remain secret! Say
 the situation is one where one of the cake-cutters (maybe the person
 celebrating a birthday) is not providing preference information, but
 still the cutting of the cake can be made without anyone being
 envious. Similarly, for deciding what the rent should be it suffices
 to consider the information of all but one of the roommates and still
 none of them will be jealous. The authors provided a rather nice
 existence proof of such fair divisions.  Recently Frick,
 Houston-Edwards, and Meunier gave an iterative approximation algorithm for the solution
 \cite{frickhoustonmeunier}.

\subsubsection*{Necklace splitting}

Another fair division problem asks for the fair splitting, between two
thieves, of an open necklace with beads threaded along the string.
Here, fair means that, for each type of bead, the number of beads
assigned to each of the thieves differ by at most one (say because the thieves
are unaware of the value of the beads). Perhaps surprisingly, this can
be achieved using only a few cuts (which turns out to be convenient
should the string material be precious). Contrary to the classical statement, with this notion of fairness we
do not need to add any conditions on divisibility.

\begin{theorem}[Necklace theorem]\label{thm:necklace}
  There exists a fair splitting of a necklace that uses no more cuts than there are bead types. 
\end{theorem}

\noindent

The result is optimal as the number of bead types many cuts are sometimes necessary
(e.g., if the beads of the same type are consecutive). Theorem~\ref{thm:necklace} was first proved by Goldberg and
West~\cite{GW85} and a simpler proof using the Borsuk-Ulam theorem was proposed by Alon and West~\cite{AW86}, 
who also came up with the above popular formulation.  More generally the challenge of finding a  division of an object into two 
portions so that each of $n$ people believes the portions are equal is called \emph{consensus-halving}. 
Note that necklace-splitting is a special case because different preferences are represented by different beads. 
In \cite{simmons+su} Simmons and Su showed how a non-constructive existence result on consensus-halving can be obtained from 
Borsuk-Ulam. They also showed, by a direct application of Tucker's lemma, how one can construct an approximate consensus-halving 
(up to a pre-specified tolerance level).

Later on, a combinatorial proof due to P\'alv\"olgyi~\cite{Pal09} used
 the \emph{octahedral} Tucker lemma, instead of the Borsuk-Ulam theorem, 
 for necklace-splitting. Here is a sketch of that idea. Let $n$ and $t$ denote, 
 respectively, the numbers of beads and bead types, let $a_j$ denote the 
 number of beads of type $j$, and assume, \emph{ad absurdum}, that 
 any fair splitting uses more than $t$ cuts. Every vector $\xx \in \{+,-,0\}^n$ defines a partial
assignment of the beads to the two thieves (identified with $+$ and
$-$ respectively). If every extension of $\xx$ into a complete
assignment has at most $t$ cuts, then no such extension can be a fair
splitting and there is some $j$ such that more than $\frac{a_j}2$
beads of type $j$ are assigned to the same thief by $\xx$; define
$\lambda(\xx)$ to be the smallest such index $j$, signed by the thief
who gets more than $\frac{a_j}2$ beads of that type. If there is an
extension of $\xx$ into a complete assignment with more than $t$ cuts,
define $\lambda(\xx)$ to be the maximum number of cuts achieved by a
completion, signed by the first component of that completion (that
this sign is well-defined is straightforward). This map $\lambda$
satisfies the condition of the octahedral Tucker lemma with $m=n-1$ and
therefore cannot exist, contradicting the initial assumption.

\bigskip

The necklace splitting problem naturally generalizes from $2$ to any
number $q$ of thieves. Alon~\cite{Alon:1987ur} showed that $(q-1)t$
cuts suffice (note that they are sometimes necessary). His proof first
assumes $q$ to be prime and replaces (in a non-trivial way) the
stronger form of the Borsuk-Ulam theorem due to B\'ar\'any, Shlosman, and Sz\H{u}cs (see Lemma 4.1 in \cite{Alon:1987ur} or Statement B' in \cite{Barany:1981vh}).
The case of general $q$ follows easily by a recursive argument. (The original proof assumed $a_j$ to
be divisible by $q$ but this was subsequently relaxed by Alon et
al.~\cite{alon2006algorithmic}, who obtained ``fair roundings'' via
integrality properties of flows.) It is not known, however, whether
P\'alv\"olgyi's proof can be adapted to the case with more than two
thieves; perhaps an ingredient for this could be
$\Z_q$-generalizations of the octahedral Tucker lemma~\cite{meunier2011chromatic,Z-GK}. Another open
question~\cite{Pal09} relates to the rounding:  For those $a_j$'s not divisible by $q$,
can one choose the thieves who get the additional beads in the fair
splitting? This is easily seen to be true for two thieves and it is also
true for three~\cite{alishahi2017fair}; it is open for $q \ge 4$.

\begin{oproblem}
  Is it possible to choose for each type $j$ the thieves who get $\lceil
  a_j/q\rceil$ and those who get $\lfloor a_j/q\rfloor$ in the fair
  splitting?
\end{oproblem}

\bigskip

For two thieves, a fair splitting can be computed in linear time for
$t=2$ and in $O(n^{t-2})$ time for $t \ge 3$~\cite{GW85}. Let us
mention again consensus-halving, the problem of dividing an object
into two portions so that each of $n$ players believes the portions
are equal appears in many contexts. In \cite{simmons+su} the
Borsuk-Ulam theorem and Tucker's lemma were used for this purpose. A
well-known challenge of computational complexity was to decide whether
the computation of a fair splitting is a PPA-complete problem
\cite{ppad-original}.  This was just recently settled by
Filos-Ratsikas and Goldberg
\cite{Filos-Ratsikas:2018,2018Filos-Ratsikas}.  The fact that
splitting necklace is PPA-complete implies that the algorithmic
version of the octahedral Tucker is PPA-complete as well, but there is
also a paper directly proving that the algorithmic Octahedral Tucker's
lemma is PPA-complete \cite{dengetal2017} (see Section
\ref{fixed-point-computing}).

The topic of fair-division is very active and once more we can only
point to a few additional types of results.  One may, for instance,
explore necklace splittings with the added constraint that adjacent
pieces of the necklace cannot be claimed by certain pairs of thieves;
for example, Asada et al. \cite{asadaetal} prove that four thieves on
a circle can share the beads of the necklace, with the restriction
that the two pairs of non-adjacent thieves will not receive adjacent
pieces of the necklace. There are also several nice high-dimensional
generalizations of (convex) splitting of booty; see
\cite{plave+pablo,mark+rade} and the references therein.

\section{Graphs}
\label{s-graphs}

Graphs are often used to model problems where pairwise interactions
are prominent.  This includes situations where graph-like structures
are apparent, for instance road or train networks, or situations as in
Euler's famous problem on the bridges of K\"onigsberg (although the
curious reader may check that Euler's original article does not use
any graph-related notion, but argues purely in terms of words coding
paths). In other situations, graphs are not evident, but exist
implicitly; for instance, to describe time dependencies between tasks
in scheduling problems.

Graph theory developed in many independent directions, driven both by
applications (e.g., finding graph matchings to resolve
assignments \cite{Sch03}) and deep structural questions (e.g.,
the graph minor theory~\cite{lovasz2006graph}). Its interaction with
(combinatorial) topology started in the mid-1970s with the proof by
Lov\'asz~\cite[Theorem~6.1]{Bjorner-survey} of the conjecture of Frank
and Maurer that any $k$-connected graph $G=(V,E)$ can be partitioned
into $k$ subsets that induce connected subgraphs, have prescribed size
(summing to $|V|$), and each contains a prescribed element. (This
result was independently given a non-topological proof by
Gy\H{o}ry~\cite{gyori1976division}.)  Lov\'asz~\cite{Lovasz:1978wu}
followed up shortly after with an astonishing solution to the Kneser
conjecture based on the Borsuk-Ulam theorem.

\subsection{Chromatic number of graphs}
\label{s-chromatic}

A \emph{coloring} of a graph by $k$ \emph{colors} is a map from its
vertex set into $[k]$; a coloring is \emph{proper} if any adjacent
vertices have different colors. The \emph{chromatic number} of a graph
is the smallest integer $k$ such that a proper coloring by $k$ colors
exists. Graph colorings arise in applications such as frequence
assignment~\cite{Aardal2007} or scheduling~\cite{de1997combinatorics}.

\begin{figure}[ht]
\centering
 \includegraphics[page=22]{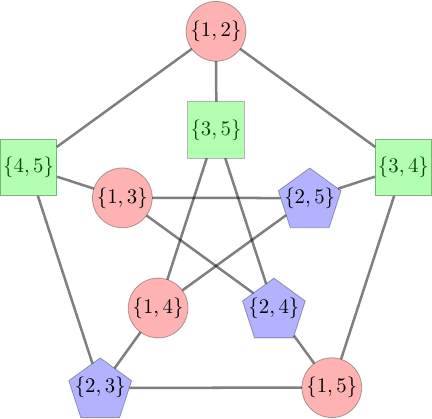}
\caption{A proper coloring with three colors: 1-\emph{red}, 2-\emph{blue}, 3-\emph{yellow} on the Kneser graph $\KG(5,2)$.} 
 \label{fig:kg52}
\end{figure}

Proving good lower bounds on chromatic numbers is usually a difficult
task, as one needs to 
show that \emph{all} colorings with fewer colors are improper. (In contrast, proving an upper bound on the chromatic
number of a given graph only requires to exhibit \emph{one} proper
coloring.) In some cases, such as the perfect graphs discussed in
Section~\ref{s:kernels}, sharp lower bounds can be obtained from the
existence of large cliques. Kneser graphs are archetypes where this
clique criterion fails dramatically. The {\em Kneser graph}
$\KG(n,k)$, where $n$ and $k$ are two integers, is the graph with
vertex set $\binom{[n]}{k}$ -- the $k$-element subsets of $[n]$ -- and
where two subsets form an edge if they are disjoint. When $n \ge
2k-1$, a natural coloring of $\KG(n,k)$ assigns to every $k$-element
subset that intersects $[n-2k+1]$ its minimal element and $n-2k+2$ to
all remaining subsets (see Figure~\ref{fig:kg52}: $\KG(5,2)$ is also
known as Petersen's graph). Kneser conjectured this to be optimal in
1955.

\bigskip

Lov\'asz approached the Kneser conjecture via a topological invariant.
Given a graph $G=(V,E)$, consider the simplicial complex $N(G)$ that
encodes its neighborhoods: $N(G)$ has  the same vertex set as $G$, and a
subset of the vertices forms a simplex whenever they have a common
neighbor in $G$. The key invariant is the (homotopy) connectivity of
$N(G)$: if $N(G)$ is $(k-1)$-connected, then $G$ is not
$(k+1)$-colorable. In the case of Kneser graphs, this yields a lower
bound that proves the conjecture. The interested reader may refer to
the book of Matou\v{s}ek~\cite[Section~3.3]{matousek2003using} and the
surveys of Bj\"orner~\cite[Section~6]{Bjorner-survey} and
B\'ar\'any~\cite[Section~5]{barany1993geometric} for details.

The idea of associating a simplicial complex to a graph and to relate
the chromatic number of the latter to topological properties of the
former has been especially fruitful and there are now various
complexes that can be used to obtain lower bounds for the chromatic
number of a graph, see \cite{Matousek:2002wl}
and~\cite[Chapter~19]{kozlov2007combinatorial} for surveys on that
approach. Recently Frick, used Tverberg-type results to show bounds for chromatic numbers of (generalized) Kneser 
graphs and hypergraphs \cite{FrickNerves,FrickKneser}.

\bigskip

The octahedral Tucker lemma emerged from a purely
combinatorial proof of the Kneser conjecture due to
Matou\v{s}ek~\cite{Matousek:2004hm}. Ziegler~\cite{Z-GK} then showed that his method can
combinatorialize other topological arguments for chromatic
numbers. Let us illustrate this method on a bound due to Dol'nikov
which deals with hypergraphs.

Recall that a \emph{hypergraph} is a pair $\HH = (V,E)$ where $V$ is a
finite set (the vertices) and $E$ a set of subsets of $V$ (the edges);
in particular, hypergraphs whose edges all have size two are
graphs. Given a hypergraph $\HH = (V,E)$ and a subset $S \subseteq V$
of its vertices, the hypergraph $\HH[S]$ induced on $S$ has vertex set
$S$ and edge set $\{e\in E \colon e\subseteq S\}$. A hypergraph $\HH$
is \emph{$2$-colorable} if $V$ can be colored so that no edge is
monochromatic. To any hypergraph $\HH= (V,E)$ we associate the {\em
  generalized Kneser graph} $\KG(\HH) = (V',E')$ where
\[ V' = E \quad \hbox{and} \quad E' = \{ef \colon e,f \in E \mbox{ and
} e \cap f = \varnothing\}.\]
In particular, the generalized Kneser graph of the $k$-uniform
hypergraph with vertex set $[n]$ is the usual Kneser graph
$\KG(n,k)$. The {\em colorability defect} $\cd(\HH)$ of a hypergraph
$\HH$ is the minimum number of vertices to be removed from $\HH$ to
ensure that the remaining vertices can be $2$-colored so that no edge
of $\HH$ is monochromatic (edges with at least one removed vertex are
discarded).

\begin{theorem}[Dol'nikov~\cite{dol1981transversals}]\label{t:Dol}
  For any hypergraph $\HH$, the chromatic number of $\KG(\HH)$ is at
  least $\cd(\HH)$.
\end{theorem}

\noindent
The combinatorial proof of Theorem~\ref{t:Dol} goes essentially as
follows. Consider a hypergraph $\HH$ with vertex set (identified with)
$[n]$. Given a proper coloring $c$ of $\KG(\HH)$ by $[t]$, we define a
map
$\lambda:\{+,-,0\}^n\setminus\{\zero\}\rightarrow\{\pm 1,\ldots,\pm m\}$
as follows. Consider $\xx \in \{+,-,0\}^n$. If at least one edge of
$\HH$ is entirely contained in $\xx^+$ or in $\xx^-$, then we choose such an edge with smallest color $a$ and we set
\[ \lambda(\xx)  = \varepsilon(n-\cd(\HH)+a),\]
where  $\varepsilon \in
\{-,+\}$ records which of $\xx^+$ or $\xx^-$ contains the edge. (The sign
$\varepsilon$ is unambiguously defined because the coloring is
proper.) If neither $\xx^+$ nor $\xx^-$ contains an edge of $\HH$,
then we define $\lambda(\ve x)$ to be $\varepsilon(|\ve x^+|+|\ve
x^-|)$, where the sign $\varepsilon$ is the first nonzero entry of
$\ve x$. In the latter case, the edges of $\HH$ contained in $\xx^+
\cup \xx^-$ induce a subgraph of $\HH$ that is $2$-colorable, so
$|\xx^+ \cup \xx^-| \le n-\cd(\HH)$. It then follows that $m =
n-\cd(\HH)+t$ suffices and that in either case, $\xx$ is mapped to
disjoint sets of labels, which helps checking that $\lambda$ satisfies
the condition of the octahedral Tucker lemma. As a consequence,
$n-\cd(\HH) + t\ge n$ and the announced inequality follows.

Since every graph is isomorphic to some (actually, many) generalized
Kneser graph, Theorem~\ref{t:Dol} provides a lower bound on the
chromatic number of any graph. In the case of Kneser graphs, this
bound is sharp. There exist refinements of the colorability defect that
yield better combinatorial bounds~\cite{alishahi2015chromatic}. Cases
of equality for Theorem~\ref{t:Dol} are remarkable for reasons related
to circular chromatic numbers~\cite{alishahi2016strengthening}
(see~\cite{zhu2006recent} for an introduction to circular chromatic
numbers). Deciding whether the chromatic number of $\KG(\HH)$ 
  equals  $\cd(\HH)$ is a natural problem asked in \cite{alishahi2016strengthening}.
 Very recently this has been proved to be NP-hard by Meunier and Mizrahi (personal communication).

\bigskip

Some generalizations of the Kneser conjecture are still open. A
$k$-element subset $A$ of $[n]$ is \emph{$s$-stable} if for any $i,j
\in A$ we have $s \le |i-j|\le n-s$ (or, equivalently, if $i$ and $j$
are distance at least $s$ apart in $\Z_n$). Let
$\KG_{s\text{-stab}}(n,k)$ denote the graph with vertices the
$s$-stable $k$-element subsets of $[n]$, and where two vertices span
an edge if they are disjoint. The graph $\KG_{2\text{-stab}}(n,k)$,
known as Schrijver's graph~\cite{schgraph}, has the same chromatic
number as the Kneser graph; since $\KG_{2\text{-stab}}(n,k)$ is a
subgraph of $\KG(n,k)$, this strengthens Lov\'asz's result. As a
special case of a conjecture on hypergraphs,
Meunier~\cite{meunier2011chromatic} proposed:

\begin{conjecture}
  For any $s \ge 2$ and $n \ge ks \ge 1$, the chromatic number
  of $\KG_{s\text{-stab}}(n,k)$ is $n-s(k-1)$.
\end{conjecture}

\noindent
Besides the case $s=2$, the conjecture is known for all even
$s$~\cite{chen2015multichromatic} and for $s\ge 4$ and $n$
sufficiently large~\cite{jonsson2012chromatic}. Some progress 
has been made by Chen \cite{chen2017}. See also \cite{FrickKneser} 
for questions related Kneser hypergraphs.
\medskip

A more systematic viewpoint recasts graph colorings as special cases
of graph homomorphisms~\cite{hahn1997graph}. A \emph{homomorphism}
from a graph $G=(V,E)$ to a graph $H=(W,F)$ is a map $f:V \to W$ such
that for every edge $vv' \in E$ the image $f(v)f(v')$ is an edge of
$H$, that is $f(v)f(v') \in F$. A proper coloring of $G$ with $k$ colors
corresponds to a homomorphism from $G$ to the complete graph with $k$
vertices. More generally, associating with every vertex of $G$ a
$k$-element subset of $[n]$ so that adjacent vertices have disjoint
subsets amounts to finding a homomorphism from $G$ to $\KG(n,k)$. The
structure of homomorphisms of Kneser graphs remains to be elucidated, as
is shown by the following, broadly open, conjecture.

\begin{conjecture}[Stahl~\cite{stahl1976n}]
  Let $n,k,k'$ be integers and let $q$ and $r$ be such that $k' =
  qk-r$ with $0 \le r <k$. There is a homomorphism $\KG(n,k) \to
  \KG(n',k')$ if and only if $n' \ge qn - 2r$.
\end{conjecture}

\noindent
For more details on partial progress on Stahl's conjecture,
see~\cite[$\mathsection 3.4$]{hahn1997graph}.

\subsection{Colorful independent sets}\label{subsubsec:hall}

A subset of vertices of a graph is \emph{independent} if there is no
edge between any pair of them. Independent sets are sometimes called \emph{stable sets}. 
This notion is central in graph theory. For instance, a proper coloring of a graph can be understood as a partition
of its vertices into independent sets, namely the sets of vertices
with same color. The search for independent sets is also natural in
applications such as the design of error-correcting
codes~\cite[$\mathsection 29$]{matouvsek2010thirty}.

\bigskip

Methods from combinatorial topology were particularly effective in
finding independent sets with color constraints, in the spirit of the
colorful theorems in combinatorial convexity. The following example
was first stated explicitly by Aharoni et al.~\cite{ABZ07}, who traced
it back to the proof of a result of Haxell~\cite[Theorem~3]{Ha95} on
hypergraph matchings.

\begin{theorem}\label{thm:haxell}
  Let $G$ be a colored graph with maximum degree $\Delta$. There
  exists an independent set of $G$ that intersects every color class
  of size at least $2\Delta$.
\end{theorem}

\noindent
(Note that the coloring of $G$ is not required to be proper.) It
suffices to prove the statement in the case where every color class
has size at least $2\Delta$ because deleting every vertex in a color
class of size less than $2\Delta$ preserves independence and does not
augment the maximum degree. The gist of the proof of Aharoni et
al.~\cite{ABZ07} is to apply Meshulam's lemma
(Proposition~\ref{p:meshulamlemma}) to the \emph{independence complex}
$\K$ of $G$, the simplicial complex consisting of its independent
sets. The connectivity of $\K$ can be controlled via a variety of
domination graph parameters; this principle, which underlies some
proofs of Meshulam~\cite{meshulam2003domination}, was made explicit by
Aharoni et al.~\cite[Theorem~2.3]{ABZ07} and given a detailed proof by
Adamaszek and Barmak~\cite{adamaszek2011lower}. In particular, given a
coloring of $V(G)$ by $[k]$, a subset $I \subseteq [k]$ and an integer
$j$, if no $2j+3$ vertices of $G$ dominate the vertices with colors in
$I$, then $\widetilde{H}_j\pth{\K[\lambda^{-1}(I)],\Z_2}$ is
trivial. (Recall that $X$ \emph{dominates} $Y$ if every vertex of $Y$
has a neighbor in $X$.)  Here, the condition holds for $j = |I|-2$
because dominating $2\Delta|I|$ vertices requires at least $2|I|$
vertices when the maximum degree is $\Delta$. By Meshulam's lemma (Proposition~\ref{p:meshulamlemma}), 
the independence complex therefore contains a colorful simplex; this is an
independent set that intersects every color class.

Theorem~\ref{thm:haxell} can be improved for graphs with more
structure as the following example shows~\cite{alishahi2017fair}.

\begin{theorem}\label{thm:fair_rep}
  In any coloring of a path, it is possible to delete a vertex of each
  color so that the remaining vertices can be partitioned into two
  independent sets $A$ and $B$ such that $-1 \le |A \cap U| - |B \cap
  U| \le 1$ for every color class $U$.
\end{theorem}

\noindent
(Again, the coloring does not need to be proper.) We sketch here a direct
proof based on the octahedral Tucker lemma in a way reminiscent of
P\'alvolgyi's proof of the necklace theorem (Theorem~\ref{thm:necklace}). 
We identify the vertex set of the path with $[n]$ and denote the color classes by $U_1,\ldots,U_t$. 
The existence of the two disjoint independent sets will be ensured via the notion of 
alternating subsequences, which has been useful in other similar contexts (e.g.,). 
A sequence of elements in $\{+,-,0\}^n$ is {\em alternating} if all terms are nonzero 
and any two consecutive terms are different. Given an $\ve x=(x_1,\ldots,x_n)\in\{+,-,0\}^n$, 
we denote by $\alt(\ve x)$ the maximum length of an alternating subsequence of $x_1,\ldots,x_n$.

The definition of the map $\lambda$ to which we will apply the octahedral 
Tucker lemma requires the quantity $s=\max\big\{\alt(\ve x)\colon \ve x\in\{+,-,0\}^n\;\,\mbox{s.t.}\;I(\ve x)=\varnothing\big\}$, where \[\begin{aligned}
I(\ve x)=\big\{i\in[t]\colon & |\ve x^+\cap U_i|=|\ve x^-\cap U_i|=|U_i|/2\\ & \mbox{or}\quad \max(|\ve x^+\cap U_i|,|\ve x^-\cap U_i|)>|U_i|/2\big\}.\end{aligned}\] Note that $s\geq 0$.

Consider a nonzero vector $\ve x\in\{+,-,0\}^n$. We distinguish 
two cases. In the case where $I(\ve x)\neq\varnothing$, we
set $\lambda(\ve x)=\pm (s+i')$, where $i'$ is the maximum element in
$I(\ve x)$ and where the sign is defined as follows. When $|\ve
x^+\cap U_{i'}|=|\ve x^-\cap U_{i'}|=|U_{i'}|/2$, the sign is $+$ if
$\min(\ve x^+\cap U_{i'})<\min(\ve x^-\cap U_{i'})$ and $-$
otherwise. When $\max(|\ve x^+\cap U_{i'}|,|\ve x^-\cap
U_{i'}|)>|U_{i'}|/2$, the sign is $+$ if $|\ve x^+\cap
U_{i'}|>|U_{i'}|/2$, and $-$ otherwise. In the case where $I(\ve
x)=\varnothing$, we set $\lambda(\ve x)=\pm\alt(\ve x)$, where the
sign is the first nonzero element of $\ve x$. Similarly to the
proofs of Theorems~\ref{thm:necklace} and~\ref{thm:fair_rep}, it can be checked
that the map $\lambda$ satisfies the condition of the octahedral Tucker lemma 
with $m=s+t$. We have thus $s+t\geq n$ and there exists $\ve
z'\in\{+,-,0\}^n$ such that $I(\ve z')=\varnothing$ and $\alt(\ve
z')\geq n-t$. It implies that there exists $\ve z\in\{+,-,0\}^n$ such
that $I(\ve z)=\varnothing$ and $\alt(\ve z)=|\ve z^+|+|\ve
z^-|=n-t$. Let $A=\ve z^+$ and $B=\ve z^-$. They are both independent
sets. Since $I(\ve z)=\varnothing$, we have $|A\cap U_i|+|B\cap
U_i|\leq |U_i|-1$ for all $i$. The fact that $|A|+|B|=n-t$ leads then
to $|A\cap U_i|+|B\cap U_i|=|U_i|-1$ for all $i$ and the conclusion
follows.

\bigskip

Many statements about independent sets have analogues in terms of
\emph{matchings}, where a matching in a graph is a set of disjoint
edges. This is natural since the matchings of a graph $G$ are the
independent sets of its \emph{line graph}, that is the graph in which
vertices are the edges of $G$, and where edges with a common endpoint
are connected. Matchings are important for theory and applications 
(see \cite{lovasz+plummer2009matching} for an excellent
book about matchings). For instance, many resource management
problems,  take after the following example: given a set
of tasks, a set of workers, and for each task the list of compatible
workers, assign to each task a different worker or report that no such
assignment exists. The worker/tasks compatibilities can be modeled by
a bipartite graph, so the question is whether there exists a matching
that covers the vertex class modeling the tasks. 

Colorful matchings still raise many questions, for instance the following well-known
conjecture due to Brualdi~\cite{B91} and Stein~\cite{S75}.

\begin{conjecture}
   If the edge set of $K_{n,n}$ (complete bipartite graph with $n$
   vertices in each side) is partitioned into sets $E_1,\ldots,E_n$ of
   size $n$ each, then there exists a matching in $K_{n,n}$ consisting
   of one edge from all but possibly one $E_i$.
\end{conjecture}

\noindent
A famous conjecture of Ryser about Latin squares~\cite{R67} asserts
that if $n$ is odd, then under the same condition as the Brualdi-Stein
conjecture, there exists a perfect matching intersecting each $E_i$
once.

\subsection{Kernels in graphs} \label{s:kernels}

A \emph{kernel} in a directed graph is a subset $K$ of the vertices
that is \emph{independent} (no two vertices of $K$ are joined by an
arc) and \emph{absorbing} (every vertex $v \notin K$ has an outgoing
arc $v \to u$ to a vertex $u \in K$). Kernels naturally arise in
certain combinatorial games, where they model the set of winning
positions~\cite{VoNMo44}, or in stable matchings, as the stable
matchings of a graph with preferences are the kernels of the
associated line graph~\cite[$\mathsection 33$]{TheBook}. Kernels
proved effective in revisiting classical questions and are for
instance at the heart of the proof by Galvin of Dinitz's conjecture on
list colorings~\cite[$\mathsection 33$]{TheBook}.

Not every directed graph has a kernel (consider a directed cycle of
length three); this is in sharp contrast with the non-directed case,
where the independent absorbing sets are the inclusion-maximal
independent sets. As shown by a series of works by Richardson, 
Duchet, Meyniel, Galeana-S\'anchez and Neumann-Lara~\cite{Duc80,DuMe83,GaNe84,Ric53}, 
a sufficient condition for the existence of a
kernel is that each odd directed cycle has two chords whose heads are
two consecutive vertices of the cycle. In particular, any acyclic
directed graph has a kernel; this situation is actually what
motivated, in the context of combinatorial games~\cite{VoNMo44}, the
introduction of kernels. In general, however, deciding if a directed
graph has a kernel is NP-complete~\cite{Chv73}.

\bigskip

Sperner's lemma comes up in the following relation between kernels and
perfect graphs. A graph is \emph{perfect} if for all its induced
subgraphs, including itself, the chromatic number is equal to the clique number. 
These are precisely the graphs for which the trivial lower bound on the chromatic
number is sharp for them and all their induced subgraphs. 
The relation between kernels and perfect graphs is a
special case of a conjecture of Berge and Duchet~\cite{BeDu83} proved
by Boros and Gurvich~\cite{BoGu96}.

\begin{theorem}\label{thm:BG}
  Any orientation of a perfect graph with no directed cycle of
  length three has a kernel.
\end{theorem}

\noindent
The original proof translates any directed graph into a coalitional
game, where the players are the cliques, and the outcomes are the
stable sets. Under the theorem's assumptions, the game has stability
properties that ensure, via results from coalitional game theory, the
existence of a ``non-rejecting'' outcome, which translates into the
desired kernel. A simpler and much more direct approach based on
Scarf's lemma was proposed by Aharoni and Holzman~\cite{AhHo98} and
further simplified by Kir\'aly and Pap~\cite{KiPa09} using Sperner's
lemma.

%

The proof of Kir\'aly and Pap goes as follows. Consider an orientation
$D=(V,A)$ of a perfect graph with no directed cycle of length $3$. Let
$P \subseteq \R^V$ denote the polyhedron of (possibly negative) vertex
weights summing to at most $1$ on every clique. This polyhedron has at
least one extreme point (assigning $1$ to any maximal independent set
of $D$ and $0$ to the remaining vertices does the job), so it is
pointed. Moreover, it has exactly $n$ independent extreme directions
(the $-\ve e_v$ where $\ve e_v$ is the unit vector associated to vertex
$v$). Every facet of $P$ corresponds to a clique on which the weights
sum to exactly~$1$. Since there is no directed cycle of length three,
every clique has a {\em source}, i.e., a vertex that is absorbed by all
other vertices of the clique. Label each facet by the source of the
corresponding clique. Note that a facet containing $-\ve e_v$ is not
labeled by $v$. By Corollary~\ref{cor:KP}, the polyhedron has an
extreme point $\ve\omega$ that is incident to facets of each label. Now,
consider the weights on $V$ defined by $\ve\omega$. If we could find an
independent set $K$ of $D$ intersecting every clique of weight~$1$,
this set $K$ would also be absorbing: indeed, every vertex $v$ labels
a facet incident to $\ve\omega$, so it is a source of a clique of weight
one which intersects $K$. The existence of $K$ follows from a classical
lemma of Lov\'asz: in a perfect graph, there exists an independent set
intersecting every clique of maximum weight. (This commonly used lemma
is perhaps difficult to find spelled out in this form; a standard way
to prove it is to use perfectness and coloring for cliques of maximum
cardinality, then apply the vertex replication lemma of perfect graphs~\cite{lovasz1972}, to allow
rational weights, then generalize to real weights using linear programming.)

%
%
%
\bigskip

Theorem~\ref{thm:BG} ensures the existence of a kernel, but the proof
does not give any efficient method to compute one.

\begin{oproblem}
  What is the complexity of computing a kernel in an orientation of a
  perfect graph with no directed cycle of length three?
\end{oproblem}

\section{Optimization}
\label{s:optimization}

Broadly speaking, mathematical optimization develops mathematical
tools for solving optimization problems. We illustrate in this section
how the theorems of Carath\'eodory, Sperner and Helly and their
variations provide original viewpoints on different aspects of this
field.

\subsection{Linear programming}

A \emph{linear program} (LP) asks for the minimum of a linear function
under a set of linear constraints and is usually written
\begin{equation}\label{eq:lp}
  \begin{aligned}
  \min \quad & \ve c^T \ve x \\
  \text{s.t.}\quad & A \ve x = \ve b\\
  & \ve x \geq \ve 0.
  \end{aligned}
\end{equation} 
Here, $A$ is a $m \times n$ matrix, $\ve x$ a vector of $n$
indeterminates, $\ve b$ and $\ve c$ vectors in $\R^m$ and~$\R^n$,
respectively, and $\xx \geq \ve 0$ means that each row of $\ve x$ is
non-negative. Linear programs may come in different presentations,
with $\max$ in place of $\min$ or possibly inequalities in place of
equalities; these presentations are essentially
equivalent~\cite[$\mathsection 4$]{matousek2007understanding}. Linear
programming is by now a central tool in operations research as it
allows to model a variety of resource management
problems~\cite[$\mathsection 2$]{matousek2007understanding} and can be
solved fairly effectively in practice. The theory of linear
programming builds on the study of systems of linear inequalities.
While this seems to be just small variation from linear algebra, 
linear programming was only systematized in the late 1940's.

\subsubsection{The simplex algorithm}\label{s:simplex}

Carath\'eodory's theorem underlies the \emph{simplex algorithm}, 
arguably the standard method to solve linear optimization problems.

On the one hand, Carath\'eodory's theorem gives a way to
discretize an a priori continuous problem. Indeed, the cone
version of Carath\'eodory's theorem ensures that if the system $A\ve x
= \ve b$ with $\ve x \ge \ve 0$ admits a solution, then it admits a
solution with support of size the rank of $A$. Such a support,
understood as a set of indices of columns of $A$, is a \emph{feasible
basis}. A closer inspection of the proof of Proposition~\ref{p:cara++}
reveals that the optimum of~\eqref{eq:lp}, if one exists, is attained
on a solution supported by a feasible basis. Since a feasible basis
determines a unique solution of $A\ve x = \ve b$, the optimum can be
found in finite (but possibly long) time by enumerating feasible bases
which are combinatorially described by their support.

On the other hand, Carath\'eodory's theorem, in the form of
Proposition~\ref{p:cara++}, also explains the pivoting mechanics of
the \emph{simplex algorithm}. Suppose there exists an optimum, and
that we have a feasible basis~$B$ determining a solution $\ve x^*$. It
turns out that if $\xx^*$ is not optimal, there exists $i \notin B$
such that increasing $x^*_i$ improves (i.e., decreases) the
objective $\ve c^T \ve x$. The set $B \cup \{i\}$ contains another
feasible basis, and it cannot define a worse solution than $\xx^*$
(again, a consequence of the proof of Proposition~\ref{p:cara++}). 
Switching to that new basis is a \emph{pivot step}. It is a non-trivial result 
of Bland that there exist rules for choosing non-cycling sequences of pivot steps,
see~\cite[$\mathsection$11.3]{Sch86} and~\cite[$\mathsection$5.8]{matousek2007understanding}. Broadly speaking, the simplex
algorithm starts by computing a feasible basis, and then it performs such
pivot steps until no entry outside the basis can be used to improve the
objective; the final basis then determines an optimal solution.

\subsection{Integer programming}

An \emph{integer program} (IP) adds integrality constraints to linear
programs by restricting all of the variables to take their
values over $\Z$ rather than $\R$, for instance
\begin{equation}\label{eq:ip}
  \begin{aligned}
  \min \quad & \ve c^T \ve x \\
  \text{s.t.}\quad & A \ve x = \ve b\\
  & \ve x \geq \ve 0, \ve x \in \Z^n.
  \end{aligned}
\end{equation} 
This variation arises naturally in the management of undivisible
resources or yes/no decision making; the emblematic example is
the \emph{knapsack problem} which asks, given a set of objects with
weights and values, for the subset of maximal value whose weight does not
exceed a given threshold.

The \emph{relaxation} of an integer linear program is the linear program
obtained by forgetting the integrality conditions, as is~\eqref{eq:lp}
for~\eqref{eq:ip}. In general, the solution to the relaxed linear
program provides a bound on the solution to the integer program, a
lower bound in the case of~\eqref{eq:lp} and~\eqref{eq:ip}. Linear
programming and relaxation play a fundamental role in combinatorial
algorithms; we refer the reader to the books~\cite{V01, WS11} for more
detail.

What we just said applies also on \emph{mixed} integer programs in which only some of the variables are required to be integer.

\subsubsection{Sparsity of integer solutions}

Carath\'eodory's theorem readily measures the sparsity of optimal
solutions. For example, Theorem~\ref{Bound_via_siegel} provides the
following bound on the support of an optimal solution.

\begin{corollary}
 Let $A \in \Z^{m \times d}$, $\ve b \in \Z^m$ and $\ve
c \in \Z^d$. The integer point of the polyhedron $\{\ve
x \in \R^d \colon A\ve x= \ve b, \ \ve x\geq \ve 0\}$ that minimizes
$\ve c^T \ve x$ has at most $$2(m+1) \log(2\sqrt{(m+1)}M)$$ non-zero
components, where $M$ is the largest of the entries of $A, \ve c$ in
absolute value.
\end{corollary}

\noindent
Similar results have been used,  for instance, for solving of bin-packing
problems, see e.g., \cite{eisenbrandshmonin-caratheodory, goemans+rothvoss}. 
See \cite{Alievetal2018} for an application to the sparsity of optimal solutions and tighter bounds
for special cases such as knapsack problems.

\subsubsection{Graver bases} \label{s:hilbertbasisinIP1}

Another example of the influence of Carath\'eodory's theorem is the
use of \emph{Hilbert bases} by \emph{Graver's optimization methods}. 
Although we present these ideas for integer programs, they apply more 
broadly, for instance to convex integer optimization problems, with respect to 
a convex objective function composed with linear functions, or convex separable
functions, see~\cite[$\mathsection 3$ and $\mathsection
4$]{alggeo4optimization} and~\cite{shmuelbook}. Consider the integer
program~\ref{eq:ip} and assume, as is usually the case in practice, that
$A \in \Q^{m \times n}$. We can decompose the polyhedron $A \ve x
= \ve 0$ into $2^n$ cones
\[\{A \ve x = \ve 0\} = \bigcup_{\ve \varepsilon \in \{-1,1\}^n} \{A \ve x
= \ve 0, \varepsilon_1x_1 \ge 0, \varepsilon_2x_2 \ge
0, \ldots, \varepsilon_nx_n \ge 0\}.\]
Each of these cones is pointed and rational, so it has a Hilbert
basis~\cite[Corollary 2.6.4]{alggeo4optimization}. The union of these
$2^d$ Hilbert bases is called the \emph{Graver basis} of $A$. Note
that Seb\H{o}'s integer Carath\'eodory theorem
(Theorem~\ref{thm:sebo}) ensures that any integer point in the
polyhedron $A \ve x = \ve b$ can be written as a non-negative integer linear 
combination of at most $2n-2$ vectors from the Graver basis of $A$. Moreover,
Graver~\cite{Graver1975} established the following remarkable
augmentation property: any non-optimal feasible solution of the
integer program~\eqref{eq:ip} can be improved by adding some suitable
vector from the Graver basis of $A$. Hence, any integer program of the
form~\eqref{eq:ip} can be solved by first computing a Graver basis for
$A\xx = \ve b$, then computing a feasible solution, and finally
improving this solution by a greedy walk on the set of integer
solutions, the candidate steps being provided by the vectors of the
Graver basis.

\begin{figure}[ht]
\centering
\includegraphics[page=23]{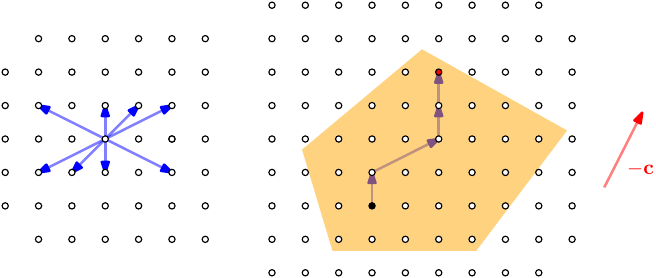}
 \caption{Illustration, in a planar projection, of Graver basis
   methods for $A \ve x = \ve b$ for $A = (1 \ 2 \ 1)$. The Graver
   basis of $A\ve x = \ve 0$ consists of
   $(2,-1,0),(0,-1,2),(1,0,-1),(1,-1,1)$, and their opposites. On the
   left, the neighbors of an integer point through the Graver
   basis. On the right, the cone $\ve x \ge 0$ (shaded) and a walk
   from an arbitrary (black) feasible integer point to the (red)
   integer point optimal for the given direction $\ve c$.}
 \label{fig:graveraction}
\end{figure}

The main obstacle to the practical use of Graver bases is their
potentially exponential size. For general matrices, deciding if a
given set of vectors is a Hilbert basis is already
coNP-complete~\cite{durandetal}. The good news, from the last decade of work, are that for
highly structured matrices, such as those with regular block
decompositions, Graver bases can be computed efficiently and are actually manageable for
optimization (see details in Chapters 3 and 4 of~\cite{alggeo4optimization} and the extensive presentation in \cite{shmuelbook}).

\subsection{LP duality}
\label{s:lpduality}

An important idea for the study of linear programs is the notion of
\emph{LP duality}. This idea naturally arises from the question of
certifying the quality of a solution to a linear program. For example,
the objective function value of an optimal solution to
\begin{equation}\label{eq:lp2}
  \begin{aligned}
    \min \quad & 5 x_1 + 3 x_2 \\ \text{s.t.}\quad & \pth{\begin{matrix}
        2 & 3 \\ 1 & 2 \end{matrix}}
    \pth{\begin{matrix}x_1\\x_2\end{matrix}} \ge \pth{\begin{matrix} 2
        \\ 5\end{matrix}}\\ & \ve x \geq \ve 0
  \end{aligned}
\end{equation} 
can be seen to be at least $2$ by looking at the first constraint and
at least $\frac{15}2$ by scaling the second constraint by
$\frac32$. More generally, the solution of a linear
program~\eqref{primal} can be bounded from below by the solution of a
linear program~\eqref{dual}, where

\noindent
\hfill
\begin{minipage}{5cm}
  \begin{align*}\tag{P} \label{primal}
    \min \quad & \ve c^T \ve x\\
    \text{s.t.}\quad & A \ve x \geq \ve b\\
    & \ve x \geq \ve 0
  \end{align*} 
\end{minipage}
\quad
and
\quad
\begin{minipage}{5cm}
  \begin{align*}\tag{D} \label{dual}
    \max \quad & \ve b^T \ve y\\
    \text{s.t.}\quad & A^T \ve y \leq \ve c\\
    & \ve y\geq \ve 0.
  \end{align*}
\end{minipage}
\hfill

\noindent
The linear programs~\eqref{primal} and~\eqref{dual} are said to be
\emph{dual} to one another; the variable $\ve y$ of the dual program can be
interpreted as the weights of a linear combination of the constraints
of the primal program (and conversely). This relation, called
\emph{weak linear programming duality}, can be strengthened.

\begin{theorem}[Strong Duality] \label{strongLPduality}
  Given two dual linear programs, if at least one is feasible, then
  they have the same optimal value.
\end{theorem}

\noindent

The duality theory of linear optimization has many applications such
as fast certification of solutions or primal-dual algorithms and in proving combinatorial theorems~\cite{Sch86}. 
But it also plays a role in discrete geometry, for example in the proof of $(p,q)$ theorems \cite{Alon92pq}. 
In the following subsection we are going to prove Theorem \ref{strongLPduality} in an atypical way.

\subsubsection{LP duality from the MinMax theorem}

LP duality is the classical favorite approach to prove von Neumann's MinMax
theorem (Theorem~\ref{thm:VonNeumann}) for two-player zero-sum games,
as we mentioned in Section~\ref{s:nash-complexity}. Going in the other
direction, Dantzig~\cite{dantzig-minmax} proposed a deduction of the
strong duality theorem from the minimax theorem; as we explain below,
Dantzig's proof required a detour, in some cases, via Farkas' lemma,
another result equivalent to Theorem~\ref{strongLPduality}. The
impression of equivalence between minimax theorem, Nash equilibria for
zero-sum two-player games and strong duality theorem nevertheless
lingered and became a folklore theorem. It is only recently that
Adler~\cite{adler-minmaxvsduality} filled-in the missing case to give
a genuine direct equivalence between these three cornerstones.

Dantzig's approach proceeds as follows. The weak duality already
proves one inequality. The other inequality reduces to
finding a solution $(\ve x, \ve y) \ge 0$ of the system~\eqref{combinedPD}:
\[\tag{P+D}\label{combinedPD} \left\{\begin{array}{rr}
 A\ve x -\ve b&\geq 0,\\
 -A^T\ve y+\ve c&\geq 0,\\
\ve b^T \ve y-\ve c^T \ve x&\geq 0,\\
\end{array}\right.\]
This system rewrites
\begin{align*}
\underbrace{\pth{\begin{matrix}
  0 & A & -\ve b \\ 
  -A^T & 0 & \ve c \\ 
  \ve b^T & -\ve c^T & 0
\end{matrix}}}_{M} \pth{\begin{matrix} \ve y \\ \ve x \\ 1 \end{matrix}} \ge \ve 0,
\end{align*}
so our task is to find a vector $\ve z$ with positive last component
and such that $M\ve z \ge \ve 0$. Consider the zero-sum game with
payoff matrix $-M$ and let $(\ve s^*,\ve t^*)$ be a Nash
equilibrium. The matrix $-M$ has dimension $n \times n$, so one may pit
each strategy against itself to see that the value $v$ of the game satisfies
\[ {\ve s^*}^T (-M) \ve s^* \ge v \ge {\ve t^*}^T (-M) \ve t^*.\]
Since $M$ is antisymmetric, this implies that $v=0$. Moreover, for any
$\ve z \in \Delta_{n-1}$ we have ${\ve s^*}^T (-M) \ve z \ge 0$, so
$M {\ve s^*} \ge \ve 0$. Since $\ve s^* \in \Delta_{n-1}$, writing
$\ve s^* = (\begin{matrix} \widetilde{\ve y} & \widetilde{\ve x} &
  \widetilde{u}\end{matrix})^T$ leads to the desired solution whenever
$\widetilde{u} \neq 0$. When $\widetilde{u} = 0$, Dantzig concluded by
a separate use of Farkas' lemma (the ``incompleteness'' in his
derivation). Adler was able to complete this missing case, without appealing to Farkas.

\subsubsection{Totally Dual Integral polyhedra} \label{s:hilbertbasisinIP2}

Applied optimization models typically involve \emph{rational
  polyhedra}, which are expressed as systems of linear inequalities
with rational coefficients. An important question for computation is
whether a rational polyhedron is \emph{integral}, that is whether all
its vertices have integer coordinates. Indeed, for integral polyhedra,
integer optimization (which is typically very hard) becomes linear
optimization (which is considered tractable). Let us see how Hilbert
bases help when looking for rational polyhedra that are integral.  In
what follows, we consider a rational polyhedron $P=\{\xx \colon A\xx
\le \ve b\}$ with $A$ and $\ve b$ rational.

Checking whether $P$ is integral is a finite process, as one can
simply list all the vertices. The following structural result allows
us to bypass this tedious enumeration in various situations. Observe
that if $P$ is integral, then for every integral vector $\ve w$, the
value $\max \{\ve w^T\xx : \ve x \in P\}$ is an integer (indeed, it is
the inner product of two integral vectors). Surprisingly, a rational
polyhedron is integral if and only if it satisfies this condition.
This equivalence, due to Edmonds and Giles~\cite{edmondsgiles}, is
still not a practical way to detect integral polyhedra (the set of
candidate vectors $\ve w$ is infinite) but it suggests to look at
things via duality. Indeed, the strong LP duality (Theorem
\ref{strongLPduality}) states
$$\max\{\ve w^T \ve x \colon A \xx \leq \ve b\}=\min\{\ve y^T \ve b \colon A^T \ve y = \ve
w, \ve y\geq 0\},$$
so $P$ is integral if the vector $\ve b$ is integral and the
right-hand side minimization problem has an integral optimal
\emph{value} for every integral vector $\ve w$. (Note that in general
integrality properties are not preserved through linear programming duality.) A system
of inequalities $A \ve x \leq \ve b$ is \emph{totally dual integral
}(TDI) if the right-hand side minimization problem above has an
integral solution for every integral vector $\ve w$ (for which the
optimum is finite). A rational polyhedron $P$ that can be represented
by a TDI system where $\ve b$ is integral is thus integral. The
converse is true and any integral polyhedron can be represented by a
TDI system of inequalities (but which, in practice, may not be easy to
find). Let us stress that TDIness is a property of the system of
inequalities, not of the underlying polyhedron: Giles and
Pulleyblank~\cite{gilespulleyblanktdi} proved that for every rational
system of inequalities $A \ve x \leq \ve b$, there is a rational
number $\alpha$ such that $\alpha A \ve x \leq \alpha \ve b$ is
TDI. They also proved, using Carath\'eodory-style properties, that a
system of inequalities $A \ve x \leq \ve b$ is TDI if and only if for
every face $F$ of $P$, the rows of $A$ which are
active in $F$ form a Hilbert basis for that supporting cone. This
makes checking TDIness a finite process, but still not a practical one
as checking whether a system of vectors forms a Hilbert bases is not
efficient. See~\cite{durandetal} and references therein for computational issues.

TDIness and related notions such as box-TDIness often shed new light on results in combinatorial optimization, for instance on the matroid 
intersection theorem. Consider two matroids $M_1$ and $M_2$ over the
same ground set $S$, understood as their sets of independent sets. Any matroid
has an associated \emph{matroid polytope}, obtained by taking the
convex hull of its (indicator vectors of) independent sets. It turns out that the
convex hull of the independent sets in $M_1 \cap M_2$ coincides with
the intersection of the matroid polytopes of both matroids. This is
remarkable, as in general $\conv (A \cap B)$ is different from
$\conv(A) \cap \conv(B)$. See~\cite[$\mathsection$41]{Sch03} for more
on this topic.

A special case of TDIness allows for linear programming proofs of combinatorial results. A matrix $A$ is \emph{totally unimodular} (TU) if every square submatrix has determinant in $\{0,-1,+1\}$. For such a matrix, the polyhedron $\{\ve x: A\ve x \leq \ve b, \ve x\geq \zero\}$ is integral for every integral vector $\ve b$. TU matrices give rise to TDI systems. They are completely characterized and
are very important in combinatorics and optimization (see
\cite{seymourTU,truemper+walter}). Since the transpose of a TU matrix is TU again, the strong duality theorem in linear programming (Theorem~\ref{strongLPduality}) provides alternative proofs of the K\H{o}nig theorem
(the maximum cardinality of a matching in
a bipartite graph is equal to the minimum cardinality of a set of
nodes intersecting each edge), the
K\H{o}nig-Rado theorem (the maximum cardinality of a stable set in a bipartite graph without isolated vertices is equal
to the minimum number of edges needed to cover all nodes), and the integrality of the Max-Flow-Min-Cut theorem. In all these three cases, the matrix $A$ is the
node-arc incidence matrix of a directed graph, and it is an easy exercise to check that such a matrix is TU.

\subsection{Convex optimization}

Using linearization techniques, one may apply ideas from the theory of
linear programming, and its duality, to more general optimization problems of the form
\begin{equation}\tag{P'} \label{program}
  \begin{array}{rlll}\min & f(\ve x) \\
    \mbox{s.t.} & h_j(\ve x)\leq 0 & j=1,\ldots,q 
  \end{array}
\end{equation}
where $f$ and $h_j$'s are differentiable functions
$\R^n\rightarrow\R\cup\{-\infty,+\infty\}$ (not just linear as
before).

\subsubsection{The KKT conditions from LP duality.}
 A milestone in mathematical programming is the following
necessary optimality condition, due to Karush, Kuhn, and
Tucker~\cite{boyd2004convex}.

\begin{theorem}[KKT condition]\label{kkt}
  Let $\ve x^*$ be a feasible solution of the problem~\eqref{program}
  such that the constraints are qualified at $\ve x^*$.  If $\ve x^*$
  is a local optimum, then there are nonnegative real numbers $\mu_1,
  \mu_2, \ldots, \mu_q$ such that
  \[\displaystyle{\nabla f(\ve x^*)+\sum_{j=1}^q\mu_j\nabla h_j(\ve
    x^*)=\zero}\qquad\mbox{and}\qquad \forall j, \quad \mu_jh_j(\ve
    x^*)=0.\]
\end{theorem}

\noindent
The requirement that ``the constraints are qualified'' is a regularity
condition on the feasible domain $F = \{\ve x \colon h_j(\ve x) \le 0,
j \in [q]\}$ near $\ve x^*$. We do not spell out this rather technical
condition but give a sufficient requirement. Call a direction $\ve d$
\emph{feasible at} $\ve x^*$ if $F$ contains a segment of positive
length with endpoint $\ve x^*$ and direction $\ve d$. (This is where
we are approximating: the adequate notion of feasibility is somewhat
more flexible.) The constraints are qualified at $\xx^*$ if the
closure of the cone of feasible directions coincides with the tangent
cone at $\ve x^*$, that is $\{\ve d \colon \nabla h_j(\ve x^*)\cdot\ve
d\leq 0 \text{ for all } j \in [q] \hbox{ s.t. } h_j(\ve x^*) =0\}$.
 Even this coarser requirement may prove tedious to check, and
several simpler sufficient conditions were investigated; the above
criterion readily yields that affine constraints are qualified in any
feasible point; another important case is that of convex
differentiable constraints, which are qualified at any feasible point
provided there exists a point satisfying strictly every
constraint~\cite[$\mathsection$5.5]{boyd2004convex}.

The factors $\mu_i$ in Theorem~\ref{kkt} are the \emph{Lagrange
  multipliers}. The strong duality theorem is classically equivalent
to the following lemma of Farkas~\cite[$\mathsection$ 7.3]{Sch86}.

\begin{lemma}[Farkas' lemma] \label{farkas}
  Let $A$ be a real matrix and let $\bb$ be a vector. There exists
  $\ve x \geq\zero$ such that $A\ve x=\bb$ if and only if $\ve
  y\cdot\bb\geq 0$ for every $\ve y$ such that $A^T\yy\geq\zero$.
\end{lemma}

\noindent
Theorem~\ref{kkt} can be deduced from Farkas' lemma via the following
 linearization argument. Let $\ve x^*$ be a feasible solution
 of \eqref{program} and write
\[  A=- \pth{\pth{\nabla h_j(\ve
    x^*)}_{j\in J}} \quad \hbox{where} \quad J = \{j\in[q] \colon
h_j(\ve x^*)=0\}.\]
Let $\ve x^*$ be a local optimum of~\eqref{program}. Since $f$ is
differentiable, $\nabla f(\ve x^*)\cdot\ve d\geq 0$ for every
direction~$\ve d$ feasible at $\xx^*$. Moreover, since the constraints
are qualified at $\ve x^*$, any direction $\ve d$ satisfying $A^T \ve d
\ge \ve 0$ is in the closure of the cone of directions feasible at
$\xx^*$, hence
\[ \forall \ve d \in \R^n \hbox{ s.t. } A^T \ve d \ge \ve 0, \quad \nabla f(\ve x^*)\cdot\ve d\geq 0.\]
By Farkas' lemma, this is equivalent to the existence of a vector $\ve
\mu'$ in $\R_+^{J}$ such that $A \ve \mu'=\nabla f(\ve
x^*)$. Completing $\ve \mu'$ into $\ve \mu$ by zeroes yields Theorem
\ref{kkt}.

\subsubsection{Strong duality in convex programming from the KKT conditions.} 
\label{s:lpkkt}

The KKT condition (Theorem~\ref{kkt}) can, in turn, be used to prove a
strong duality theorem for convex programming, generalizing
Theorem~\ref{strongLPduality}. Let us introduce the Lagrangian function
\[ \L(\ve x,\ve \mu)=f(\ve x)+\sum_{j=1}^q\mu_j h_j(\ve
    x).\]
Since
\[ \sup_{\ve \mu \in \R_+^{q}} \L(\ve x,\ve \mu) = \left\{\begin{array}{ll} f(\xx) & \hbox{ if } \mu_jh_j(\ve
    x) \le 0\\ +\infty & \hbox{ otherwise.} \end{array}\right.\]
\eqref{program} is equivalent to 
\[
  \begin{array}{rl}\min & \sup \{\L(\ve x, \ve \mu) \colon \ve \mu \in \R_+^q\} \\
    \mbox{s.t.} & \ve x\in\R^n.
  \end{array}
\]
The same argument as in Equation~\eqref{eq:minmax} yields
\begin{equation}\label{eq:dualitefaibleconvex}
\displaystyle{  \inf_{\ve x\in\R^n} \  \sup \{\L(\ve x, \ve \mu) \colon \ve \mu \in \R_+^q\} \ge \sup_{\mu\in\R_+^q} \ \inf \{\L(\ve x, \ve \mu) \colon \ve \xx \in \R^n\}}.
\end{equation}
Finding the right-hand side term consists in solving the following \emph{dual} program:
\begin{equation}\tag{D'}\label{dualprogram}
  \begin{array}{rlll}\max & g(\ve \mu) \qquad \hbox{where} \quad g(\ve \mu) = \inf\{ \L(\ve x, \ve \mu) \colon \ve \xx \in \R^n\}\\
    \mbox{s.t.} & \ve \mu \ge \ve 0.
  \end{array}
\end{equation}
This dual program always asks to maximize a concave function. In the
case where Problem~\eqref{program} is a linear program, this notion of
duality coincides with the LP duality introduced in
Section~\ref{s:lpduality}.

\begin{proposition}[Strong duality for convex optimization]
  Suppose that in~\eqref{program}, $f$ and $h_1, \ldots, h_q$ are
  convex functions. Suppose moreover that the constraints are
  qualified at every feasible solution. If \eqref{program} has an
  optimal solution, then the dual program has one too and the optimal
  values of \eqref{program} and \eqref{dualprogram} coincide.
\end{proposition}

\noindent
The proof goes as follows. Let $\ve x^*$ be an optimal solution of
\eqref{program}. By Theorem~\ref{kkt}, there exists $\ve
\mu^*\in\R_+^q$ such that
 \[\displaystyle{\nabla f(\ve x^*)+\sum_{j=1}^q\mu_j^*\nabla h_j(\ve
    x^*)=\zero}\qquad\mbox{and}\qquad \forall j, \quad \mu_j^*h_j(\ve
 x^*)=0.\]
 On the one hand, we have
 \[ \L(\xx^*,\ve \mu^*)=f(\ve x^*)+\sum_{j=1}^q\mu_j^* h_j(\ve    x^*) = f(\xx^*).\]
 On the other hand, we have
 \[ \nabla_{\xx}\L(\xx,\ve \mu^*) = \nabla f(\ve x)+\sum_{j=1}^q\mu_j^*\nabla h_j(\ve
 x), \quad \hbox{so} \quad  \nabla_{\xx}\L(\xx^*,\ve \mu^*)=\zero. \]
 The map $\xx \mapsto \L(\xx,\ve \mu^*)$ is convex because $\ve \mu^*
 \ge \zero$, so $\xx^*$ is a global minimum and $g(\ve \mu^*)=
 \L(\xx^*,\ve \mu^*) = f(\xx^*)$. Together with the weak duality of
 Equation~\eqref{eq:dualitefaibleconvex}, this ensures that $\ve \mu^*$ is
 an optimal solution for \eqref{dualprogram}.

\subsection{Sampling approaches}

Let us now consider optimization problems of
the form
\begin{equation}\label{eq:minwit}
  \begin{aligned}
    \min \quad & f(\ve x)\\
    \text{s.t.}\quad & \ve x \in C_1 \cap C_2 \cap \cdots \cap C_m
      \end{aligned}
\end{equation} 
where one minimizes a function over an intersection of subsets $C_i$,
the \emph{constraints}. Such problems include linear programming (when
$f$ linear and the $C_i$ are halfspaces), convex programming (when $f$
is convex and the $C_i$ are convex sets) or their integral or mixed
analogues (via restrictions to $\Z^k \times \R^d$).

\subsubsection{Witness sets of constraints}

A first use of Helly's theorem concerns the removal of redundant
constraints defining an optimal solution. To begin, consider the linear program
\begin{equation}\label{eq:lpwit}
  \begin{aligned}
  \min \quad & \ve c^T \ve x \\
  \text{s.t.}\quad & A \ve x \le \ve b\\
  & \xx \in \R^d,
  \end{aligned}
\end{equation} 
where $A \in \R^{m \times n}$ (this form differs from those seen so
far but is equivalent~\cite[$\mathsection
  4$]{matousek2007understanding}). Assume that the problem is feasible
and let $t$ denote its solution. The set $A \ve x \le \ve b$ of
feasible solutions is an intersection of $m$ halfspaces, and their
common intersection with $\ve c^T \ve x < t$, another halfspace, is
empty. By Helly's theorem, some $d+1$ of these $m+1$ halfspaces must
have empty intersection, and $\ve c^T \ve x < t$ must be one of
them. It follows that we may drop from~\eqref{eq:lpwit} all but some
(carefully chosen) $d$ constraints without changing the
solution. Recall that Helly's theorem for halfspaces is dual to
Carath\'eodory's theorem; in the dual, this argument yields that an
optimal solution can be realized by a feasible bases (as introduced in
Section~\ref{s:simplex}).

More generally, given a problem of the form~\eqref{eq:minwit}, a
subset $W \subseteq \{C_1, C_2, \ldots C_m\}$ is a \emph{witness set
of constraints} if
\[  \begin{aligned}
  & \min \quad f(\ve x)  & \quad = \qquad & \min \quad f(\ve x)\\
  & \text{s.t.}\quad \ve x \in \bigcap_{i=1}^m C_i & & \text{s.t.} \quad  \ve x \in \bigcap_{C \in W} C, \\
  \end{aligned}\]
and $W$ is inclusion-minimal for that property; in other words, a
witness set is a non-redundant set of constraints that defines the
same optimum as the entire problem. The above argument gives, \emph{mutatis
mutandis}, that the witness sets of any feasible mixed-convex
programming over $\Z^k \times \R^d$ have size at most $2^k(d+1)-1$;
this bound increases by $1$ for unfeasible programs.

The relation between Helly's theorem and witness sets extends beyond
convexity, as observed by Amenta~\cite{amenta94}. We say that a family
$\F$ admits a \emph{Helly-type theorem} with constant $h$ if the
non-empty intersection of every $h$-element subset of $\F$ implies the
non-empty intersection of $\F$. Consider a minimization problem of the
form~\eqref{eq:minwit}, where $f$ is a function from a space $X$ to a
space $V$. We assume that $V$ is equipped with a total order $\prec$, so that
the minimization question makes sense. For any $v \in V$ we let $L_v
= \{x \in X\colon f(x) \prec v\}$.

\begin{proposition}\label{prop:Amenta}
Let $h \in \N$. If every family $\{C_1,C_2, \ldots, C_m, L_v\}$ admits a Helly-type
theorem with constant $h$, then any witness set of constraints
of~\eqref{eq:minwit} has size at most $h$ (actually, $h-1$ if the
problem is feasible).
\end{proposition}

\noindent
Note that we make no assumption on $f$ (not even continuity!), $X$ or
$V$. The proof for the feasible case goes as follows. Let $\F
= \{C_1,C_2, \ldots, C_m\}$ and let $s$ denote the solution
to~\eqref{eq:minwit}. For $G \subseteq \{C_1, C_2, \ldots, C_m\}$ define
$s(G)$ as the minimum of $f$ over $\bigcap_{C \in G}C$ and put
\[ s' = \max\{ s(G)\colon G \subseteq \F \hbox{ and } |G| = h-1\}.\]
On the one hand, $s \ge s(G)$ for every $G \subseteq \F$, we have
$s \ge s'$. On the other hand, every $h$ elements of $\F$ intersect
(because~\eqref{eq:minwit} is feasible) and $L_{s'}$ intersects any
$h-1$ elements of $\F$ (by definition of $s'$); thus, every $h$
elements in $\F \cup \{L_{s'}\}$ intersect, and the Helly-type theorem
on $\F \cup \{L_{s'}\}$ ensures that $s \le s'$. By minimality, any
witness set thus has size at most $h-1$.

\subsubsection{Combinatorial algorithms for linear programming}

Devising algorithms for linear programming with provably good
complexity has been a major challenge for the past 70 years. The
interior point method of Karmarkar~\cite{Karmarkar1984} and the
analysis of the ellipsoid method by
Khachyian~\cite{khachiyan1980polynomial} only showed that the
complexity of LP is polynomial in the number $m$ of constraints, the
number $d$ of variables, and the bit complexity $L$ of the entries of
the matrix $A$ and vectors $\ve b$ and $\ve c$. A major question thus
remains:

\begin{oproblem}
Is there an algorithm that solves linear programming in time
polynomial in the number of constraints and the number of
variables, assuming that arithmetic operations on input numbers have
unit cost?
\end{oproblem}

\noindent
(This is problem number nine in Smale's list of open
problems~\cite{smale1998mathematical}.) An algorithm with complexity
polynomial in $m$ and $d$ in the unit cost model is called
\emph{strongly polynomial}. Although the simplex algorithm proves
effective in practice, no choice of pivoting rule is known to ensure a
number of step polynomial in the number $m$ of constraints and the
number $d$ of variables; in fact, for every pivot rule whose
worst-case complexity is established, that complexity is at least
exponential in $m$ and $d$~(see \cite{avis+friedmann,cunninghampivot,friedmann,klee+minty} 
and the references therein).  Although no strongly polynomial
time algorithm is known, partial progress was made through the 1980's
and 1990's via combinatorial random sampling algorithms; this approach
hinges on the fact that the bounded size of witness sets allows to
throw away redundant constraints quickly.

Let us illustrate the basic idea of combinatorial random sampling
algorithms in its simplest form, due to Seidel~\cite{seidel1991small}
(see also~\cite[$\mathsection 4$]{MMMC}). Consider
\begin{equation}\label{eq:lpsei}
  \begin{aligned}
  \min \quad & \ve c^T \ve x \\
  \text{s.t.}\quad &\xx \in \R^d,\\
  & \ve x \in H_1 \cap H_2 \cap \cdots \cap H_m,\\
  \end{aligned}
\end{equation} 
where each $H_i$ is a halfspace in $\R^d$. Pick $t \in [m]$ uniformly
at random and let $\ve s_t$ denote the solution to the linear program
with the constraint $H_t$ removed. The idea is to compute $\ve s_t$
recursively, then deduce $\ve s$ from $\ve s_t$: we check in $O(d)$
time whether $\ve s_t$ belongs to $H_t$. If $\ve s_t \in H_t$, then
$\ve s = \ve s_t$, and if $\ve s_t \notin H_t$, then $\ve s$ must
belong to the hyperplane bounding $H_t$ and be the solution to the
linear program
\begin{equation}\label{eq:lpsei2}
  \begin{aligned}
  \min \quad & \ve c^T \ve x \\
  \text{s.t.}\quad &\xx \in \partial H_t,\\
  & \ve x \in \bigcap_{i \in [m]\setminus \{t\}} (H_i \cap \partial H_t).
  \end{aligned}
\end{equation} 
This new linear program has $m-1$ constraints in $d-1$ variables and
can be obtained from~\eqref{eq:lpsei} in time $O(dm)$. Altogether, the
expected time $T(m,d)$ to compute $\ve s$ writes
\[ T(m,d) = T(m-1,d) + O(d) + \mathbb{1}_{\ve s_t \notin H_t} T(m-1,d-1).\]
Observe that $\ve s_t \notin H_t$ if and only if $H_t$ belongs to
\emph{every} witness set for~\eqref{eq:lpsei}. Since the size of
witness sets is at most $d$, the event $\ve s_t \notin H_t$ occurs
with probability at most $\frac{d}{m}$ when $t$ is chosen uniformly at
random. Altogether, this recursion solves to $T(m,d) = O(d!m)$, which is
the running time of Seidel's algorithm. In fact, Seidel's algorithm builds on an
idea by Clarkson, which we describe next.

The \emph{iterated reweighting} method of Clarkson~\cite{clarkson-lv}
consists of assigning a \emph{weight} $w_i$, initially set to $1$, to
every constraint $H_i$ and iterating a simple process: Sample $O(d^2)$
constraints with probabilities proportional to their current weights
and solve the problem on these constraints.  If the solution is
feasible, we are done; otherwise double the weights of all violated
constraints and reiterate. It remains to be shown that, almost surely,
the algorithm terminates. This can be seen by comparing the growth
rates of the total weight of the system, and the weights of some
witness set $W$. On the one hand, every unsuccessful iteration must
double the weight of at least one constraint in $W$. On the other
hand, as the constraints are chosen in each iteration with probability
proportional to their current weight, the expected total weight of the
constraints violated at any given iteration is $O\pth{\frac{\sum_i
    w_i}{d}}$~\cite{C16}. Thus, after $k$ iterations, $W$ has weight
$\Omega(d 2^{k/d})$ but the total weight of all the constraints in the
system is $O\pth{m\pth{1+\frac{O(1)}d}^k}$. Putting these two bounds
together implies that the algorithm terminates, with high probability,
within $O(d \log m)$ iterations. At each iteration one has to solve a
linear program on $O\left(d^2\right)$ constraints, which can be done
in time $O(d)^{d/2+O(1)}$, say by the simplex method, and compute the
set of violating constraints, which takes time $O(md)$.  This implies
an overall running time of $O(d^2 m \log m + O(d)^{d/2+O(1)} \log m)$.

Clarkson's approach was later improved by Matou\v{s}ek et
al.~\cite{msw-sblp-92} to achieve an expected time complexity
of~$O(d^2m+e^{O(\sqrt{d \log d})})$. A similar bound was obtained
independently by Kalai~\cite{kalai1992subexponential} via a randomized
pivot rule for the simplex algorithm. Clarkson's algorithm was subsequently 
derandomized; see~\cite{C16} and the references therein for the latest developments.

\subsubsection{Combinatorial abstractions of LP} 

Clarkson's algorithm uses very little structure from linear
programming, namely the abilities to solve a small-size problem and to
decide if a given solution violates a given constraint.  Surprisingly,
there are many computational problems for which these two operations
can be performed effectively to find a solution and are not
``linear''.  A simple example is the computation of the smallest
enclosing circle of a finite point set $P\subset \R^2$. Here, the
constraints are the points of $P$, the candidate solutions are the
circles, and a circle violates a point if it does not enclose
it. Observe that subsets of points that minimally define their
enclosing centers have size at most three, so in this case witness
sets again have size at most three. It turns out that Clarkson's
algorithm readily applies to this problem. This is not an LP in
disguise: a generic instance may have a witness sets of size two
\emph{or} three.

Various combinatorial abstractions of LP were studied in order to
understand precisely what class of problems can be solved with the
randomized approach we described before.  Consider an abstract set of
constraints numbered from $1$ to $m$, and an objective function that
associates to any set $S \subseteq [m]$ the value $f(S) \in \R$ of the
optimum when only the constraints in $S$ are considered. A natural
black-box model allows to compute $f(S)$ when $S$ has bounded size
(independently of $m$), or to decide violations asking whether $f(S
\cup \{i\}) \stackrel{?}= f(S)$.  It turns out that Clarkson's
algorithm can compute effectively $f(S)$ in this abstract model under
three assumptions~\cite{msw-sblp-92}: (i) that $f$ be decreasing under
inclusion, that is $f(S) \le f(T)$ whenever $S \subseteq T \subseteq
[m]$, (ii) that $f$ be \emph{local} in the sense that
\[\begin{aligned}
  f(S) \neq f(S \cup \{i\}) & \Leftrightarrow f(T)  \neq  f(T \cup \{i\})\\
  & \forall S  \subseteq T \subseteq [m] \hbox{ such that } f(S)=f(T), \hbox{ and } \forall i \in [m],
  \end{aligned}\]
and (iii) that witness sets have bounded size, where a witness set $S$
is a minimal subset of $[m]$ with $f(S)=f([m])$. Functions $f$
satisfying~(i) and~(ii) are called \emph{LP-type problems} or
\emph{generalized linear programming} problems.

Any generic problem of the form~\eqref{eq:minwit} is LP-type; here by
generic we mean that for every subset of constraints, $f$ achieves its
minimum over the intersection of those constraints in only a bounded
number of points. As noted by Proposition~\ref{prop:Amenta},
controlling the size of witness set for such problems, and thus the
effectiveness of Clarkson's approach is a matter of Helly-type
theorems.  Later a generalization, called \emph{violator spaces}~\cite{ViolatorSpaces2008}, 
 was shown to give the precise family of problems solved by Clarkson's approach.

\subsubsection{Chance-constrained optimization}\label{s:CCP}

  Consider the problem of computing, given $n$ points in the plane, the
smallest disk containing a given proportion of these points (say
$70\%$). More generally, given a probability measure $\mu$ in the
plane, and a positive number $\varepsilon$, one may consider the optimization problem
\begin{equation*}
  \begin{split}
    \min \quad & r \\
    \text{s.t.} \quad &  \Pr\left[\|\xx-\yy\|_2 \le r \right] \geq 1-\varepsilon, \\
      & \xx \in \R^2, 
  \end{split}
\end{equation*}
where $\yy$ is a random point chosen from the probability distribution
$\mu$. This quantitative variation of the smallest enclosing circle
problem, discussed above, is no longer LP-type but can still be solved
effectively. The technique relies, again, on the fact that witness
sets have bounded size and applies more generally to \emph{chance
  constrained problem} (CCP). A CCP asks to optimize a function of a
variable $\xx \in \R^d$ under constraints depending on a parameter
$\ve w \in \Delta$, of the form:
\begin{equation*}
  \begin{split}
    CCP(\varepsilon) =\min \quad & g(\xx) \\ \text{s.t.} \quad &
    \Pr\pth{f(\ve x,\ve w) \leq 0} \geq 1-\varepsilon,\; \xx \in K.
  \end{split}
\end{equation*}
Here, $g$ is a convex function, the probability is taken relative to a
measure $\mu$ on the space $\Delta$ of parameters, $f(\xx,\ve w)$ is
measurable with respect to $\ve w$, $f(\cdot, \ve w)$ is convex for
every $\ve w$, and $K$ is a convex set. This type of optimization
problem naturally arises when modeling with uncertain
constraints~\cite{shapirodentchevarusz}.

An approach to solve CCP, initiated by Calafiore and Campi
\cite{calafiorecampi2005,calafiorecampi2006}, is to sample $\ve w^1,
\ve w^2, \ldots, \ve w^N$ from $\mu$ and solve the deterministic
convex program
\begin{align*}
   SCP(N)= \min \quad & g(\xx) \\
    \text{s.t.} \quad & f(\xx, \ve w^i) \leq 0 ,\quad i=1,2,\ldots, N,\; \xx \in K.
\end{align*}
For any $\delta \in (0,1)$, if $N \geq \frac{2d}{\varepsilon}
\ln\frac{1}{\varepsilon} + \frac{2}{\varepsilon} \ln\frac{1}{\delta}+2d$, 
then the solution to $SCP(N)$ is a solution to $CCP(\varepsilon)$ with
probability at least
$1-\delta$~\cite{calafiorecampi2005,calafiorecampi2006}.

The proof in~\cite{DLetal-ccopt} goes as follows. For $\ve x\in \R^d$, let
$V (\ve x)= \Pr\pth{f(x,w)>0}$ so that we are interested in ensuring
$V(\xx) < \varepsilon$. Each of the $N$ constraints $f(\cdot, \ve w^i)
\leq 0$ is convex so any witness set has size at most $d$. Now, for
every $I \in \binom{[N]}d$ we define
\[ \Gamma_N^I=\left\{(\ve w^1,\dots, \ve w^N) \in \Delta^N\colon (\ve w^i)_{i \in I}
\text{ is a witness set} \right\}.\]
Note that the $\Gamma_N^I$ decompose $\Delta^N$ according to which are
the indices of the $d$ witness constraints. Let $\xx^*$ denote the
optimal solution of $SCP(N)$ and $\xx^I$ the solution to the convex
program defined by the constraints $\{\ve w^i \colon i\in I\}$
alone. The probability of failure $\Pr \pth{V(\xx^*)\ge\varepsilon}$ is less than or equal to
\[
   \sum_{I \in \binom{[N]}d} \Pr\pth{ \left\{\left(\ve w^1,\dots,
     \ve w^N\right) \in \Gamma_N^I \colon V(\ve x_I) \geq
     \varepsilon\right\}}.
  \]
The summand corresponding to index set $I$ rewrites as
\[\begin{aligned}
\Pr&\left[ \left\{\left(\ve w^i\right)_{i\in I}\colon  V(\xx_I)  \geq {\varepsilon} \right\}\right] \\
& \prod_{j \notin I} \Pr\left[ \left\{\left(\ve w^i\right)_{i \in I}:  f\left(\xx_I, \ve w^j\right) \leq 0\right\} \Bigm| \left\{\left(\ve w^i\right)_{i\in I}\colon  V(\xx_I) \geq {\varepsilon} \right\} \right].
\end{aligned}\]
The first factor is at most $1$ and each of the following $N-d$
factors is at most $\varepsilon$. Altogether, we get $ \Pr
\pth{V(\xx^*)} \le \binom{N}d \pth{1-\varepsilon}^{N-d}$ and the
announced bound follows.
  
\subsubsection{$S$-optimization} \label{centerpointsinoptima}

There are many situations where one wants to optimize under convex
constraints while restricting the solutions to belong to some set $S$;
these are called \emph{$S$-optimization problems}. This allows to
model complicated constraints like mixed-integer constraints ($S=\R^d
\times \Z^k$), sparsity constraints (e.g., compressed sensing), or
complementarity constraints.

Several of the techniques described above generalize if the
intersections of $S$ with convex sets of the ambient space admits a
Helly-type theorem. We denote its Helly number by $h(S)$, for example
$h(\R^d\times \Z^k) = 2^k(d+1)$. Proposition~\ref{prop:Amenta} yields
that witness sets have size at most $h(S)-1$, and the reader can check
that the analysis of chance constraint programs in Section~\ref{s:CCP}
applies: the solution to $SCP(N)$ is a solution to $CCP(\varepsilon)$
with probability at least $1-\delta$ if $N \geq
\frac{2h(S)-2}{\varepsilon} \ln\frac{1}{\varepsilon} + \frac{2}{\varepsilon}
\ln\frac{1}{\delta}+2h(S)-2$. In fact, the proof presented above
differs from the initial argument~\cite{calafiorecampi2005,calafiorecampi2006} and was found
when generalizing CCP to the $S$-optimization
setup~\cite{DLetal-ccopt}. Helly-type theorems have been obtained for
various sets $S$~\cite{Averkov:2013uo,deloeraetal2015Helly,garber,queyranne+tardella}.

Let us now turn our attention to the problem of minimizing a
  convex function~$g$ over an arbitrary subset $S \subseteq \R^d$. We
  assume that $g$ is given by a first-order evaluation oracle and that
  $S$ is nonempty and closed. The {\em cutting plane method}
  that we now present allows to approximate the solution. To allow us
  to control the quality of the approximation, we fix a finite measure
  $\mu$ supported on $S$. The algorithm starts with a convex set $E_0$
  that contains the solution in its interior and is such that
  $\mu(\inter E_0) >0$. It then builds a sequence $\{E_i\}$ where each
  $E_i$ is also a convex set that contains the solution in its
  interior. Given $E_{i-1}$, we select a point $\xx_i \in (\inter
  E_{i-1})\cap S$ and compute $g(\xx_i)$ and a subgradient $\ve h_i
  \in \partial g(\xx_i)$. We set $\xx^\star:=\textrm{argmin}_{\xx \in
    \{\xx_1, \ldots, \xx_i\}}g(\xx)$ and define $E_i$ in a way that
  ensures that
  \[ E_{i}\;  \supseteq \;\{ \xx \in E_0 \colon   g(\xx^\star)-g(\xx_j)\ge \ve h_j^T (\xx-\xx_j),\; \forall j\in[i] \},\]
  and that $\mu(\inter E_i)$ is non-increasing. We stop when
  $\mu(\inter E_i)$ is smaller than the desired error and return
  $\xx^*$. This approach leaves many details unspecified, in
  particular the precise definition of $E_i$ and the way to choose the
  points $\xx_i$. When $S=\R^d$ and $\mu$ is the Lebesgue measure, a
  possible implementation is the classical \emph{ellipsoid
    method}~\cite{khachiyan1980polynomial}. When $S = \Z^d$ and $\mu$ is the counting measure for
  $\Z^d$, we obtain cutting plane algorithms for convex integer optimization problems. 
  Another variant of this method which uses random sampling was explored by Bertsimas and
  Vempala~\cite{Bertsimas:2004fp}.

  The choice of the points $\xx_i$ in the cutting plane method is
  important. A particularly good choice are the Tukey centers. Given a
  vector $\ve u \in \s^{d-1}$ and a point $\xx \in \R^d$ we let
  $H^+(\ve u,\xx)$ denote the halfspace $\{\ve y\in\R^d \colon \ve u^T
  (\ve y-\xx) \ge 0\}$. Consider the function
  \begin{equation}\label{eq:max-value}
    \F(S,\mu):=\max_{\xx\in S}\inf_{\ve u\in \s^{n-1}}\mu(H^+(\ve u,\xx)).
  \end{equation}
  A {\em Tukey    center} is a point that attains the maximum value of
  $\F(S,\mu)$. Lower bounds on $\F(S,\mu)$, and therefore on the depth
  of Tukey centers (see Section~\ref{depthsec}), can often be obtained
  from Helly-type theorems. For instance, if $S = \R^d$ and $\mu$ is
  the counting measure of a finite subset of $\R^d$, then the
  centerpoint theorem (\ref{thm:centerpointthm}) ensures that
  $\F(S,\mu) \ge \frac1{d+1}|P|$. The proof of the centerpoint theorem
  from Helly's theorem extends to the setting of $S$-convexity: if
  $S\subseteq \R^d$ is nonempty and closed and $\mu$ is finite and
  supported on $S$, then $\F(S,\mu)\ge
  \frac{1}{h(S)}\mu(\R^d)$~\cite{basu+oertel}. For instance, Doignon's theorem
  ensures that if $S = \Z^d$ and $\mu_C$ counts integer points inside
  a compact set $C$, then $\F(S,\mu_C) \ge \frac{|C \cap \Z^d|}{2^d}$.
It turns out that choosing $\xx_i$ among the points maximizing
  $\F(\mu_i,S)$, where $\mu_i$ is the restriction of $\mu$ to
  $\inter(E  _{i-1})$, gives the best running times among cutting plane
  algorithms for convex minimization over~$S$~\cite{basu+oertel}.
  
  More notions generalize to $S$-optimization. For instance there
  exists an analogue of the strong duality theorem for $S$-convex
  optimization under some natural
  conditions~\cite{basuetal:S-optimality}. We will meet Tukey centers
  again in Section \ref{depthsec}, and we conclude this section with a challenge:

 \begin{oproblem} 
   What is the complexity of computing Tukey centers for given $S$ and
   $\mu$? E.g., can one compute an exact Tukey center of the integer points of
   a convex polytope in polynomial time in the input size?
\end{oproblem}

\section{Data point sets}
\label{s:datapoints}

In this section we consider some computer science results that are
either applications of the main theorems or strongly related to
them. There are simply far too many results for us to do justice to
even a small number of them, thus we restrict ourselves to a few
central themes, including classification, geometric shape analysis,
and partitioning of $n$ points in $d$-dimensional Euclidean
spaces. These will involve an interplay between affine geometric and
topological techniques, offering the usual mix of advantages and
drawbacks of the two: affine tools -- Radon's lemma, Helly's theorem,
linear programming duality, simplicial decompositions -- will imply
fast algorithms, though apply to a restricted group of geometric
objects. Topological tools -- the Borsuk-Ulam theorem, Tucker's
lemma---yield broader structural statements, though, at this moment,
settling the algorithmic feasibility of these methods remains a major
open problem.

\subsection{Equipartitioning: ham sandwich theorem and its relatives}

We say that a hyperplane $h$ \emph{bisects} a set $P$ of points if the
two open halfspaces defined by $h$ contain at most $ \frac{|P|}{2} $
points of $P$. Note that if $|P|$ is odd, then a point of $P$ must
necessarily lie on $h$.  The famous \emph{ham sandwich theorem} is the
starting point of a large number of results concerning
equipartitioning of geometric objects with other geometric objects.

\begin{theorem}   
  Let $P_1, \ldots, P_d$ be finite point sets in $\R^d$. Then there
  exists a hyperplane $h$ that simultaneously bisects each $P_i$, $i =
  1, \ldots, d$.
\label{thm:hamsandwich}
\end{theorem}

The theorem holds more generally for finite Borel measures that
evaluate to zero on every affine hyperplane.  All known proofs of
Theorem~\ref{thm:hamsandwich} are essentially topological in nature. A
classical proof follows from the Borsuk-Ulam theorem by identifying
points of $\mathbb{S}^d$ with hyperplanes in $\R^d$, and where the
function $f: \mathbb{S}^d \to \R^d$ encodes the ``unbalance'' of the
$d$ point sets (more generally, measures) with respect to that
hyperplane~\cite{matousek2003using}.  We now outline another proof of
the ham sandwich theorem using Tucker's lemma.  The proof we present
was found independently by Holmsen and by the third author (both
unpublished).

For simplicity we will assume that the given point sets $P_1, \ldots,
P_d$ are in general position, and let $\bigcup_{i=1}^d P_i = \{\ve
p_1, \ldots, \ve p_n\}$.  Each pair
$(\aaa,b)\in(\R^d\times\R_+)\setminus\{\zero,0\}$ induces a {\em sign
  pattern} $\xx\in\{+,-,0\}^n$ with $x_j$ being the sign of $\aaa\cdot
p_j+b$. We will apply Tucker's lemma (Proposition~\ref{tuckerlemma})
on an abstract simplicial complex induced by these sign patterns to
show that there is a pair $(\aaa,b)$ whose sign pattern $\xx$ is such
that $\{p_j\colon j\in \xx^+\}$ and $\{p_j\colon j\in \xx^-\}$ each
contain at most half of each $P_i$. Then the hyperplane
$\{\yy\in\R^d\colon\aaa\cdot\yy+b=0\}$ bisects each $P_i$.

Denote by $\PP$ the partially ordered set of all achievable sign
patterns endowed with the partial order $\preceq$ (see
Section~\ref{s:topology} before Tucker's octahedral lemma). It is a
well-known result from oriented matroid theory that the order complex
$\T$ of $\PP$ is a triangulation of $\B^d$ that is symmetric on its
boundary. Suppose, for contradiction, that there is no hyperplane
bisecting each $P_i$. Given an $\xx$ in $\PP$, we define
$\lambda(\xx)$ to be $\varepsilon i$, where $i$ is the smallest index
such that either $\{\ve p_j\colon j\in \xx^+\}\cap P_i$ or $\{\ve
p_j\colon j\in \xx^-\}\cap P_i$ contains more than half of the points
of $P_i$, and where $\varepsilon$ is $+$ if it is the first set and is
$-$ if it is the second. This map $\lambda$ is clearly antipodal on
the boundary of $\T$ and labels the vertices of $\T$ with the elements
of $\{\pm 1,\ldots,\pm d\}$. According to the Tucker lemma, there
exists an edge $\uu\vv$ of $\T$ with
$\lambda(\uu)+\lambda(\vv)=0$. Without loss of generality, we assume
that $\uu\preceq\vv$ and that $-\lambda(\uu)=\lambda(\vv)=k$ for some
$k\in[d]$. By definition of $\lambda$, we have then $|\{\ve p_j\colon
j\in \uu^-\}\cap P_k|>|P_k|/2$ and $|\{\ve p_j\colon j\in \vv^+\}\cap
P_k|>|P_k|/2$. Combined with $\uu^-\subseteq\vv^-$, it implies that
$|\{\ve p_j\colon j\in \vv^-\cup\vv^+\}\cap P_k|>|P_k|$, a
contradiction.

We will see further applications of the ham sandwich theorem later on, 
but for now we point out that it gives another proof of Theorem~\ref{thm:necklace}:
given an open necklace with $t$ types of beads to be  divided equally between two thieves, embed the beads
of the necklace along a moment curve in $\R^t$, and use a hyperplane $h$ guaranteed
by Theorem~\ref{thm:hamsandwich} to bisect each type of bead.
As any hyperplane intersects 
a  moment curve at $t$ points, $h$ splits the open necklace
into $t+1$ pieces that can then be divided among the two thieves.

Note that from the above discussion, we have the following
``computational hierarchy'': computing a solution
to this variation of the octahedral Tucker problem is harder than computing
a ham sandwich cut (note also that this implies
that the latter is in PPA), which is harder than computing
a solution to the fair splitting necklace problem.
In particular, the Filos-Ratsikas and Goldberg's paper \cite{2018Filos-Ratsikas} proves that 
computing the ham sandwich cut is PPA-complete (see Section \ref{cakes+necklaces}).

As far as the efficiency aspects of Theorem~\ref{thm:hamsandwich} are concerned, 
a line bisecting two given point sets of total $n$ points  in $\R^2$ can be computed
in $O(n)$ time. In $\R^d$, the best algorithm
to computing a ham sandwich cut for $d$ point sets in $\R^d$
runs in time $O(n^{d-1})$~\cite{LMS94}; in
fact the algorithm proposed presents a new proof of Theorem~\ref{thm:hamsandwich}
that proceeds by induction on the dimension and thus is
more amenable to efficient algorithm design.

By now there are dozens of variants of the ham sandwich theorem, and
generalizations to other types of bisecting and bisected objects, We
now present a few nice examples.  For the ham sandwich
theorem, see also \cite{matousek2003using,Zivaljevic:Handbook} and
references therein.

One variation is called the \emph{center transversal theorem}: given
$s+1$ point sets in $\R^d$ where $s \in \{0, \ldots, d-1\}$, there
exists an $s$-dimensional affine subspace $h$ of $\R^d$ such that any
hyperplane containing $h$ has at least $\frac{1}{(d-s+1)}$-th fraction
of the points of $P_i$ on each side, for $i= 0, \ldots, s+1$. In fact,
it has been conjectured that the constant $\frac{1}{(d-s+1)}$ can be
replaced by $\frac{(s+1)}{(d+s+1)}$, see
\cite{BMN08}. Theorem~\ref{thm:hamsandwich} is the case $s = d-1$, and
the case $s=0$ is the important centerpoint theorem that we
encountered in the introduction and will visit again later.  We refer
the reader to the book~\cite{matousek2003using} for many variants of
this and related theorems in $\R^2$ and higher dimensions.  Here are
now two famous conjectures.

\begin{conjecture}
Let $P$ be a set of points in $\R^4$. Then there exists a
set $\HH$ of four hyperplanes such that each of the resulting $16$
open regions
of $\R^4 \setminus \HH$ contains at most $\frac{|P|}{16}$ points of $P$.
\end{conjecture}

In the papers \cite{pfag2015,pfag2016}, the authors showed that in $\R^5$ it is 
indeed possible to find four hyperplanes that divide the set into 16 equal parts.
See also \cite{SimonS}. We should mention here a related theorem of 
Yao and Yao~\cite{YaoYao85} (see also Theorem~\ref{thm:simplicialpartitions}): 
Given a set $P$ of $n$ points in $\R^d$, one can partition $\R^d$ into $2^d$ regions 
such that the interior of each region contains at most $\frac{n}{2^d}$ 
points of $P$, and any hyperplane intersects the interior of at most $2^d-1$ regions.

Next we consider another conjecture by Tverberg and
Vre{\'c}ica~\cite{TverbergVrecica}.
 
\begin{conjecture}
  Let $0 \le k \le d-1$ be a given parameter, and let $P_0, P_1, \ldots, P_k$ be finite
  point sets in $\R^d$. If $|P_i| = (r_i-1)(d-k+1)+1$ for
  $i=0,1,\ldots, k$, then each $P_i$ can be partitioned into $r_i$
  parts $P_{i,1}, P_{i,2}, \ldots, P_{i,r_i}$ such that the sets
  $\{\conv(P_{i,j}) \colon 0 \le i \le k, 1 \le j \le r_i\}$ can be
  intersected by a $k$-dimensional affine space.
\end{conjecture}

Note that Tverberg's theorem is the case
when $k=0$ in the above statement.

If one wants to partition more than $d$ point sets in $\R^d$,
then hyperplanes are often insufficient; however the following
important variant of the ham sandwich theorem,   due to Stone and Tukey~\cite{ST42}, 
shows that then polynomials of a sufficiently high degree can be used to do the partition.

We say that a $d$-variate   
polynomial $f \in \R[x_1, \ldots, x_d]$
\emph{bisects} a point set $P \subseteq \R^d$ if it evaluates
to negative on at most $\frac{|P|}{2}$ points of $P$
and likewise evaluates to positive on at most  $\frac{|P|}{2}$ points of $P$.
Note that for the case of polynomials
of degree one, this coincides
with our earlier definition of bisection for hyperplanes.

\begin{theorem}
Let $P_1, \ldots, P_s$ be finite point sets in $\R^d$. Then there exists
a $d$-variate polynomial $f  \in \R[x_1, \ldots, x_d]$ of degree  
$O(s^{\frac{1}{d}})$ such that $f$ bisects each $P_i$, $i= 1, \ldots, s$.
\label{thm:tukeystone}
\end{theorem}

The idea is to reduce the above problem
to the usual ham sandwich theorem in a suitably high dimension.
As a $d$-variate polynomial of degree $D$ has $d' = \binom{D+d}{d}-1$ monomials (aside
from the constant term), identify each  such monomial  with a distinct dimension of $\R^{d'}$.
Then each $d$-variate polynomial $f$ of degree $D$ can be identified 
with a hyperplane in $\R^{d'}$,  where the coefficients defining the hyperplane
  in $\R^{d'}$ (i.e., the $d'$ coordinates of the normal vector of the hyperplane)
correspond to the coefficients of the corresponding monomials of $f$.
This also  gives a mapping -- called \emph{Veronese mapping} -- of the points in $P_1 \cup \cdots \cup P_s$
to $\R^{d'}$, where the $i$-th coordinate of a point $p \in \R^d$ 
is the value of the corresponding monomial on $p$.
One can  now use Theorem~\ref{thm:hamsandwich} on
the $d$ sets corresponding to the lifted points of
$P_1, \ldots, P_s$ to get a hyperplane $h$ in $\R^{d'}$  that bisects
each of the $s$ lifted sets. Note that to use the ham sandwich theorem, we require
$s \leq d' = \binom{D+d}{d}-1$ and thus need to satisfy the constraint $D = \Omega( s^{\frac{1}{d}} )$.
The ham sandwich hyperplane $h$ corresponds to the required $d$-variate
polynomial in   $\R^{d}$, of degree $O( s^{\frac{1}{d}} )$.

\subsection{Parametrized partitioning of data via geometric methods}
\label{s:geompart}

So far the partitioning statements have been of the type where the
input geometric configuration precisely fixes the output type -- e.g.,
given $d$ point sets in $\R^d$ fixes the output of
Theorem~\ref{thm:hamsandwich} to be a hyperplane. Or given $s$ point
sets in $\R^d$ in Theorem~\ref{thm:tukeystone} fixes the output to be
a polynomial of degree $O(s^{\frac{1}{d}})$.  Now we consider
statements where, besides the input geometric configuration, one is
also given an independent parameter $r$ and the complexity of the
output is a function of both the geometric configuration as well as
the value of $r$. Thus one gets a hierarchy of output structures
(varying with $r$), and one is free to choose the value of $r$
depending on the precise problem at hand. This turns out to be very
useful for designing hierarchical data structures where one can pick
the value of $r$ to maximize computational efficiency.
  
These kinds of partitioning statements -- which we call
\emph{parameterized spatial partitioning} -- have been a key theme in
discrete and computational geometry for both algorithmic and proof
purposes. Consider, for example, \emph{Hopcroft's problem} studied in
the early 1980s: given a set of lines $L$ and a set of points $P$, an
\emph{incidence} is a pair $(\ve p, l)$, where $\ve p \in P$, $l \in
L$, and the point $\ve p$ lies on the line $l$.  Then given $L$ and
$P$, is there an efficient method to determine if there exists at
least one incidence between them?  It is not difficult to see that a
spatial partition that either partitions the points or partitions the
line (in a suitable sense) is useful for decomposing the original
problem into several problems of smaller size; the current best
algorithm with a running time of $2^{O(\log^* n)} \cdot
n^{\frac{4}{3}}$ is based on such techniques.  We refer the reader
to~\cite{E96} for details on this and other related results on
Hopcroft's problem.

For a more mathematical application in a similar setting, consider the
following question posed by Erd\H os~\cite{E46}: what is the maximum
number of incidences between any $n$ points and any $n$ lines in the
plane?  Erd\H os observed that as the bipartite incidence graph
between points and lines is $K_{2,2}$-free, this is upper-bounded by
$O\big(n^{\frac{3}{2}}\big)$. More generally, the maximum number of
incidences between $m$ lines and $n$ points is $O( n \sqrt{m} )$. On
the other hand, a set of points in a `grid-like' configuration
exhibits $\Omega\big(n^{\frac{4}{3}}\big)$ incidences, and Erd\H os
conjectured this to be, asymptotically, the right bound.

This question was resolved affirmatively by Szemer\'edi and Trotter by
a complicated combinatorial argument~\cite{ST83}. We sketch here a
beautiful and simple proof of this theorem by Clarkson et
al. \cite{CEGSW90} that showcases the use of spatial partitioning for
proving combinatorial bounds. Given a set $L$ of $n$ lines in the
plane, suppose that for any parameter $r>1$, there exists a partition
of the plane into $t = O(r^2)$ (possibly unbounded) interior-disjoint
triangles $\Pi = \{ \triangle_1, \ldots, \triangle_t \}$ such that
each triangle $\triangle_i$, for $i= 1 \ldots t$, intersects at most
$\frac{n}{r}$ lines of $L$ in its interior.  Now one can partition the
incidences between $P$ and $L$ into those for which the points lie on
boundary of some triangle, and those for which the points lie in the
interior of some triangle of $\Pi$.  It is not difficult to see that
the former can be only $O(nr)$; on the other hand, a triangle
$\triangle_i$ intersecting $m_i$ lines of $L$ in its interior and
containing $n_i$ points of $P$ can contain $O( n_i \sqrt{m_i} )$
incidences in its interior (via the graph-theoretic bound). This gives
the overall number of incidences lying in the interior of triangles to
be $\sum_i O( n_i \sqrt{m_i} ) = O(n \sqrt{ \frac{n}{r} })$. Thus for
any $r>1$, the total number of incidences is bounded by $O\left(nr +
\frac{n^{\frac{3}{2}}}{\sqrt{r}}\right)$, and setting $r=
\Theta(n^{\frac{1}{3}})$ gives the desired bound!

We now elaborate on this structural partitioning problem and its
variations.  The key behind the proof is the partition of $\R^2$ into
triangles, each of which intersects ``proportionally few'' lines of
$L$. More generally, in any dimension one can show the existence of a
similar partition of hyperplanes~\cite{C93}.

\begin{theorem}
  Given a set $\HH$ of hyperplanes in $\R^d$ and a parameter $r \geq
  1$, there exists a partition of $\R^d$ into $\Theta( r^d )$
  interior-disjoint simplices such that the interior of any simplex
  intersects at most $\frac{|\HH|}{r}$ hyperplanes of $\HH$.
\label{thm:cuttings}
\end{theorem}

Such a partition is called a $ \frac{1}{r} $-cutting of $\HH$.  Note
that the bound of $\Theta(r^d)$ cannot be improved: each simplex can
contain at most $(\frac{|\HH|}{r})^d$ vertices induced by $d$-tuples
of $\HH$ in its interior, and so there must be $\Omega(r^d)$ simplices
to account for all the $\Theta( |\HH|^d )$ vertices induced by $\HH$.

The intuition behind Theorem~\ref{thm:cuttings} is as follows.  Pick
each hyperplane of $\HH$ into a random sample $R$ independently with
probability $p$ (to be set later).  Then $\Ex [|R| ]= |\HH| \cdot p$,
and so $R$ partitions $\R^d$ into $O\big( (|\HH|p)^d \big)$ induced
cells, each of which can then be further partitioned into simplices,
to get a partition of $\R^d$ into an expected total of $O\big(
(|\HH|p)^d \big)$ simplices.  Furthermore note that each such simplex
intersects, in expectation, $\frac{1}{p}$ hyperplanes of $\HH$ in its
interior. Setting $p = \frac{r}{|\HH|}$ gives the required statement.
This argument was done `in expectation', and it is non-trivial to
convert it with the same asymptotic bounds to where each simplex is
guaranteed to intersect no more than $\frac{|\HH|}{r}$ hyperplanes of
$\HH$.
 
The proof of Theorem~\ref{thm:cuttings} is usually presented in the
more general abstract framework of the theory of $\eps$-nets, whose
setting we briefly describe now.  Given a base set of elements $X$, a
set system $\RC$ on $X$, and a parameter $0 \leq \eps \leq 1$, call a
set $N \subseteq X$ an $\eps$-net for $\RC$ if $N$ contains at least
one element from each $R \in \RC$ with $|R| \geq \eps \cdot |X|$.  In
the case of $\frac{1}{r}$-cuttings, the set $X$ is the set of
hyperplanes $\HH$, and $R \in \RC$ if and only if there exists a
simplex $\Delta$ with $R = \mathrm{int}(\Delta) \cap \HH$. Then a
$\frac{1}{r}$-cutting can be constructed by taking an $\eps$-net $N$
for $\RC$ with $\eps = \frac{1}{r}$ and partitioning $\R^d$ into $O(
|N|^d )$ simplices using $N$. As the interior of each simplex induced
by $N$ does not intersect any hyperplane of $N$, it can only intersect
less than $\eps \cdot |\HH| = \frac{|\HH|}{r}$ hyperplanes of $\HH$,
as desired.

Bounds on $\eps$-nets have been extensively studied for set systems
satisfying a combinatorial condition, called the \emph{VC
  dimension}~\cite{VC71}: given a set system $(X, \RC)$, define the
projection of $\RC$ onto a set $Y \subseteq X$, denoted $\RC|_Y$, as
the set system $\RC|_Y = \{ Y \cap R \colon R \in \RC\}$.  We say that
a subset $Y$ is \emph{shattered} by $\RC$ if all $2^{|Y|}$ subsets of
$Y$ can be realized by intersection with some set of $\RC$, i.e., if
$|\RC|_Y| = 2^{|Y|}$.  Then the VC dimension of $\RC$ is defined to be
the size of the largest set that is shattered by $\RC$.

The VC dimension plays an important role in the theory of
set systems derived from geometric configurations due to the fact
that the VC dimension of such systems is usually quite small. 
For example, consider the set system where $X$ is a finite set of $n$
points in $\R^d$, and the subsets in $\RC$ are derived
from intersection with halfspaces; here $R \in \RC$
if and only if there exists a halfspace $h$ such that
$R = h \cap X$. The VC dimension
of this set system is $d+1$; in other words, given any 
set $X$ of $d+2$ points in $\R^d$, one cannot `separate' all subsets
of $X$ by intersection with halfspaces. This is an immediate consequence
of Radon's lemma (recall that Radon's lemma
is the case of Tverberg's theorem for two parts, $r = 2$, namely that
any set $P$ of $d+2$ points in $\R^d$ can be partitioned
into two subsets $P_1, P_2$ 
such that $\conv(P_1) \cap \conv(P_2) \neq \varnothing$): the convex hulls of the Radon partitions 
intersect, and thus cannot be separated by a hyperplane. On
the other hand, any set of $d+1$ points in general position 
can be shattered by halfspaces in $\R^d$. By Veronese maps, 
this implies more generally bounded VC dimension for  set systems induced
by geometric objects of bounded algebraic complexity (see~\cite[Chapter 10]{Mbook}).

Returning to $\eps$-nets, building on the work of Vapnik and Chervonenkis~\cite{VC71}, 
Haussler and Welzl~\cite{Haussler:1987fr} showed the existence of small $\eps$-nets 
as a function of the VC dimension of a set system.

\begin{theorem}
Let $(X, \RC)$ be a finite set system with VC dimension at most $d\geq 1$. 
Then for any real parameter $\eps > 0$, there exists
an $\eps$-net for $\RC$ of size $O(\frac{d}{\eps} \log \frac{1}{\eps} )$.
\label{thm:hwepsnetsthm}
\end{theorem}

The power of this theorem derives from the fact that the size
is independent of the number of elements in $X$ and the number
of sets in $\RC$. 
Combined with the VC dimension bound of $d+1$ on set systems
induced on points in $\R^d$ by halfspaces, Theorem~\ref{thm:hwepsnetsthm}
implies the existence of $\eps$-nets of size 
$O \left( \frac{d}{\eps} \log \frac{1}{\eps} \right)$, 
which has been shown to be optimal~\cite{KMP16}. 
The bounds of Theorem~\ref{thm:hwepsnetsthm} can be further improved
for many geometric set systems, and recent
work presents a unified framework for these bounds~\cite{CGKS12, MDG16,V10}.
We refer the reader to the books~\cite[Chapter 15]{PA95},~\cite[Chapter 10]{Mbook},~\cite[Chapter 4]{C00}
and~\cite[Chapter 5]{M99} for a more detailed exposition on $\eps$-nets and their many applications.
 
The theory of VC dimension fails to 
help in the construction of $\eps$-nets when
it is unbounded.  A basic case is the set system induced by convex
objects in $\R^d$; namely given a set $P$ of $n$ points in $\R^d$ and 
a parameter $\eps>0$, one would like to show the existence
of a small set $Q \subset \R^d$,  called a \emph{weak $\eps$-net}
for $P$ induced by convex objects, such that any convex object 
containing at least $\eps n$ points of $P$ must
contain at least one point of $Q$. 
Note here that $Q$ can be
any set of points in $\R^d$, and is not just limited
to being a subset of $P$ -- hence the term ``weak''.

An initial bound of $O(\frac{1}{\eps^{d+1}})$ on the size of $Q$
was shown by Alon et al.~\cite{ABFK92} (their proof uses the colorful Carath\'eodory's theorem together
with Tverberg's theorem; see the discussion about simplicial depth in Section \ref{depthsec}), 
and this was improved to  $\tilde{O}(\frac{1}{\eps^{d}})$  by Chazelle et al.~\cite{CEGGSW93}. 
This was  improved even further, by logarithmic factors, by Matou{\v{s}}ek
and Wagner~\cite{MW02} whose elegant proof
we outline now. Partition $P$ into $t$ equal sized subsets
$\PP = \{P_1, \ldots, P_t\}$,
for a parameter $t$ that is chosen optimally, such
that any hyperplane intersects the convex hulls
of $O(t^{1-1/d})$ subsets of $\PP$. Let $C$ be a convex object
containing $\eps n$ points of $P$. When $C$ intersects ``few''
sets of $\PP$,  the proportion of points of $C$
contained in some $P_i \in \PP$ is higher than $\eps$, hence
$C$ can be hit by a set constructed inductively for $P_i$.
Otherwise, $C$ intersects many sets of $\PP$. In this case,
pick one arbitrary point from each set of $\PP$ intersecting $C$.
Let $\ve q$ be  the centerpoint of those points. Then $\ve q$ must be contained in $C$, as otherwise the hyperplane
separating $\ve q$ from $C$ must intersect the convex hulls
of many sets in $\PP$, a contradiction to the definition of $\PP$.  

For the case $d=2$, Rubin~\cite{rubin18} proved
the bound of $O( \frac{1}{\eps^{1.5+\delta}} )$, where $\delta > 0$
is an arbitrarily small constant.

Finding asymptotically optimal
bounds for  weak $\eps$-nets induced by convex objects in $\R^d$ is a
tantalizing open problem. The best known lower bounds  
of $\Omega( \frac{1}{\eps} (\log \frac{1}{\eps})^d )$~\cite{BMV11}
are quite far from the upper bounds. On the other hand,
there are partial results that indicate the upper bounds that
can be improved; e.g., it is known~\cite{MR08} that one can
construct weak $\eps$-nets from $\tilde{O}(\frac{1}{\eps})$ points of $P$
in $\R^d$: pick a set $Q'$ of
$\tilde{O}(\frac{1}{\eps})$ points that form an $\eps$-net for the set
system induced by the intersection of $d$ halfspaces. Then 
adding points lying in Tverberg partitions
of carefully chosen subsets of $Q'$ results in a weak $\eps$-net,
though of size $O(\frac{1}{\eps^{d+2}})$.
For a different formulation, one can also fix an integer parameter 
$k > 0$, and then ask for the minimization problem. Find $\varepsilon = \varepsilon(k)$ such that
for any set $P$ of points in $\R^d$, there exists a set $Q \subset \R^d$ of $k$ 
points such that all convex objects containing at least $\varepsilon \cdot |P|$ points of 
$P$ are hit by $Q$ (see~\cite{MR09b}).

\begin{oproblem}
  What is the asymptotically best bound for the size
  of the smallest weak $\eps$-net for the set system 
  induced by convex objects in $\R^d$?
\end{oproblem}

We next turn to the algorithmic aspects of spatial partitioning.
There has been substantial work on improving the constants in the bounds
on $\eps$-nets, as they are directly linked to the approximation
ratios of   algorithms for the geometric hitting set problem~\cite{BG95, ERS05}.
Given a set $P$ of points in $\R^d$ and a set $\mathcal{O}$
of geometric objects, the goal is to compute a minimum 
subset of $P$ that hits all the objects in $\mathcal{O}$. 
This can be written as an integer program, which is then
approximated as follows: (a) solve the linear relaxation
of the integer program (i.e., 
the linear program 
obtained by replacing the integral constraints with real ones), and (b) assign  weights 
to the points of $P$ according to the linear program, and
finally (c) compute a $\frac{1}{W}$-net $N$ for the set system
induced by $\mathcal{O}$, where $W$ is the value of the   linear programming
relaxation. As the weight
of the points contained in each object of $\mathcal{O}$ is at
least $1$ by the linear program, 
$N$ is a hitting set for $\mathcal{O}$. The size
of $N$ is then bounded by the size of a $\frac{1}{W}$-net; e.g.,
if an $\eps$-net exists of size $\frac{c}{\eps}$ for some constant $c$,
then $N$ has size at most $c \cdot W$ -- i.e., at most $c$ times
the optimal solution. We refer the reader
to the recent surveys~\cite{CMR16,MV16} for precise bounds on $\eps$-nets
and cuttings for various geometric set systems as well as the current-best
algorithms for computing such nets.

Cuttings, together with linear programming duality (or alternatively, Farkas' lemma), 
can be used to derive another important partitioning tool
we have already encountered in the construction of weak $\eps$-nets -- the 
simplicial partition theorem~\cite{M92}.
\begin{theorem}
Let $P$ be a set of points in $\R^d$, and $t > 1$ a given integer parameter.
Then there exists a partition of $P$ into $t$ sets $\PP = \{P_1, \ldots, P_t\}$
such that $(a)$ $|P_i| \leq \frac{2 |P|}{t}$
for all $i = 1, \ldots, t$ and $(b)$ any hyperplane
intersects the convex-hull of $O\big( t^{1-\frac{1}{d}} \big)$
sets of $\PP$.
\label{thm:simplicialpartitions}
\end{theorem}
We outline the general idea of the proof. Given $P$, one
can first ``discretize'' the space of all possible hyperplanes
in $\R^d$
into a finite set $\HH$ of hyperplanes, so that one only
needs to  construct
a partition $\PP$ such that each hyperplane of $\HH$ intersects the convex hulls
of $O(t^{1-\frac{1}{d}})$ sets of $\PP$. Construct a $\Theta(\frac{1}{t^{1/d}})$-cutting
for $\HH$, consisting of at most $t$ simplices, and let $\PP'$ be
the collection of points of $P$ lying in each simplex of the cutting.
Then, on average, $|P'| = \Theta(\frac{|P|}{t})$ for each $P' \in \PP'$, and furthermore
the convex hull of each $P'$ is intersected by at most $\frac{|\HH|}{t^{1/d}}$ hyperplanes
of $\HH$. This is, in a suitable sense, the dual of the statement that we want, where
each hyperplane should intersect few cells in the cutting (and thus few convex hulls
of sets in $\P'$). However we are not far off -- the total number of intersecting pairs of hyperplanes of $\HH$ with 
convex-hulls of sets in $\PP'$ is $O(t \cdot \frac{|\HH|}{t^{1/d}}) = O( t^{1-\frac{1}{d}} \cdot \HH)$. 
In other words, the ``average'' hyperplane in $\HH$ intersects
the convex hulls of $O(t^{1-\frac{1}{d}})$ sets of $\PP'$.  Now LP duality~\cite{HP09} (or
Farkas' lemma~\cite{GKM14}) shows the existence of the required partition.

\subsection{Parametrized partitioning of data via topological methods}

Now  we move to more
recent approaches to parameterized spatial partitioning. These methods use
equipartitioning results such as the ham sandwich
theorem and Theorem~\ref{thm:tukeystone} to show the existence
of polynomials that induce parametrized partitions of $\R^d$. 
The resulting statements  are similar in spirit and use to
those we saw earlier for affine objects. The main advantage
of these newer approaches is that partitioning space with polynomials
circumvents several difficult technical issues
that arise when dealing with piece-wise linear objects. On the other hand, new computational
challenges arise in the topological approaches, as they often do not lend
themselves to efficient algorithm design.
  
We present now a ``polynomial version'' of Theorem~\ref{thm:simplicialpartitions},
discovered by Guth and  Katz~\cite{GK15}. For more on the impact of this theorem
see the book \cite{GuthBook}. In what follows given a polynomial $f$, we denote by  $Z(f)$ the zero set of $f$.

\begin{theorem}
Let $P$ be a set of $n$ points in $\R^d$, and let $r>1$
be a given parameter. Then there exists
a $d$-variate polynomial $f$ of degree $O(r^{\frac{1}{d}})$ such that
each connected component of $\R^d \setminus Z(f)$ contains at most $\frac{|P|}{r}$ points.
\label{thm:polynomialpartitions}
\end{theorem}

Here the points lying
in the components of $\R^d \setminus Z(f)$ play the role of the sets
in Theorem~\ref{thm:simplicialpartitions}.
Observe that any hyperplane $h$ in $\R^d$ intersects
$O\big( (\deg(f))^{d-1} \big) = O(r^{1-\frac{1}{d}})$ different
components of $\R^d \setminus Z(f)$, quantitatively the same
bound as in Theorem~\ref{thm:simplicialpartitions}.
In contrast with Theorem~\ref{thm:simplicialpartitions}, these components are interior disjoint,
though they will, in general,  not be convex.

We sketch a proof:
partition $P$ into two equal sized
sets $P_1, P_2$ by a polynomial, say $f_0$,
of degree $O(1)$ (using Theorem~\ref{thm:hamsandwich}, for example).  
Then partition these two point sets 
$P_1$ and $P_2$ into four equal disjoint subsets $P_3, P_4, P_5, P_6$ using
a polynomial, say $f_1$, of degree $O(2^{\frac{1}{d}})$ 
via Theorem~\ref{thm:tukeystone}. Continuing, let $f_i$
be a polynomial of degree at most $O((2^i)^{\frac{1}{d}})$ that equipartitions $2^{i}$ equal-sized disjoint subsets of $P$,
for $i = 0, \ldots, \log r$. Note here
that as long as $2^i \leq d$, a hyperplane suffices
for $f_i$.
The key observation now is that the polynomial $f$  formed
by taking the product of all these polynomials -- namely
$f = \prod_{i=0}^{\log r} f_i$ -- is
the required polynomial: as the zero set
of $f$ is simply the union of the zero sets of all the $f_i$'s,  
  each connected region of $\R^d \setminus Z(f)$ can contain 
at most $\frac{|P|}{2^{\log r}} = \frac{|P|}{r}$ points of $P$.
The degree of $f$ can be bounded as
$$\deg(f) \leq \sum_{i=0}^{\log r} \deg(f_i) \leq \sum_{i=0}^{\log r} O\left( (2^i)^{\frac{1}{d}} \right) 
= O\left( r^{\frac{1}{d}} \right).$$
 Theorem~\ref{thm:polynomialpartitions} gives
another proof of the Szemer\'edi-Trotter  theorem to
bound the number of incidences between a set $P$ of $n$
points and a set $L$ of $n$ lines in the plane, this time by
partitioning the points instead of the lines: apply Theorem~\ref{thm:polynomialpartitions}
on $P$ with $r = n^{\frac{2}{3}}$ to get a polynomial $f$, with $\deg(f) = O(n^{\frac{1}{3}})$. 
Note that each line of $L$ intersects $O(\sqrt{r}) = O( n^{\frac{1}{3}})$ components of $\R^2 \setminus Z(f)$
and  each of the $O(r) = O(n^{\frac{2}{3}})$ components 
of $\R^2 \setminus Z(f)$ contains at most $O(\frac{n}{r})= O(n^{\frac{1}{3}})$ points of $P$.
A simple calculation by summing up
the incidences across each component induced by $f$ shows that the overall
number of incidences is bounded by $O(n^{\frac{4}{3}})$.

On the computational side, an efficient algorithm to compute the partition 
guaranteed by Theorem~\ref{thm:polynomialpartitions} was discovered
by Agarwal et al.~\cite{AMS13}, who then used it to construct
efficient data-structures for answering range queries
with constant-complexity semi-algebraic sets as ranges, in time
close to $O(n^{1-\frac{1}{d}})$. 

We now move on to our last topic in this section,
a recent beautiful theorem of Guth~\cite{G15} that gives a very
general theorem that implies many statements that were previously
regarded as unconnected.

\begin{theorem}
Let $\Gamma$ be a set of $k$-dimensional varieties in $\R^d$,
each defined by at most $m$ polynomial equations
of degree at most $D$. For any parameter $r \geq 1$, there exists
a $d$-variate polynomial $f$ of degree $O(r^{\frac{1}{d-k}})$, so that each
component of $\R^d \setminus Z(f)$ intersects
at most $\frac{|\Gamma|}{r}$ varieties of $\Gamma$. The constant
in the asymptotic notation depends on $D, m$ and $d$.
\label{thm:polynomialvaritiespartitioning}
\end{theorem}

Note that this theorem implies Theorem~\ref{thm:polynomialpartitions}
by setting $k = 0$ and implies Theorem~\ref{thm:cuttings} (strictly speaking, a polynomial
version of it where $\R^d$ is partitioned into the components
of $\R^d \setminus Z(f)$ instead of simplices)
by setting $k= d-1$. See~\cite{BBZ16} for an exciting new
alternate proof of Theorem~\ref{thm:polynomialvaritiespartitioning} that extends it to a setting 
with several algebraic varieties.

We present a brief sketch of the proof of Theorem~\ref{thm:polynomialvaritiespartitioning},
which builds upon the proof of Theorem~\ref{thm:polynomialpartitions}
in a natural way. The goal, as before, is to find
polynomials $f_0, \ldots, f_s$, where $s \leq \log r^{\frac{d}{d-k}}$
and $\deg(f_i) \leq 2^{\frac{i}{d}}$ for each $i=0, \ldots, s$. 
Then  the required polynomial will be 
$$ f = \prod_{i=0}^{ \log r^{\frac{d}{d-k}}} f_i, \quad \text{with} \quad \deg(f) \leq \sum_{i=0}^{\log r^{\frac{d}{d-k}}} \deg(f_i) \leq \sum_{i=0}^{\log r^{\frac{d}{d-k}}} O\left( 2^{\frac{i}{d}} \right) 
= O\left( r^{\frac{1}{d-k}} \right),$$
as required. To see the idea behind the proof of Theorem~\ref{thm:polynomialvaritiespartitioning},
trace  the 
proof of Theorem~\ref{thm:polynomialpartitions} backwards: the polynomials $f_i$'s
were constructed using Theorem~\ref{thm:tukeystone},
whose proof used the ham sandwich theorem in a suitably
high dimension, whose proof identified the coefficients
of the polynomial with points of the sphere $\myS^{d'}$ in a suitable
dimension $d'$ and then applied the Borsuk-Ulam theorem. 
Note that the $f_i$'s were constructed
\emph{independently} via an iterative argument, where
first $f_0$ was used to partition the given point set $P$
into two sets; these two sets were then equipartitioned with $f_1$ and so on.
This approach relied crucially on the fact that a point, outside the zero set
of any polynomial $f_i$, lies in precisely one cell induced by their product. 
This fact fails for $k$-dimensional varieties when $k > 0$:
assume one has constructed polynomials $f_0, \ldots, f_{s'-1}$, $s' < s$, 
such that each component
induced by $f' = \Pi_{i=0}^{s'-1} f_i$ has the same number of incidences with
the varieties in $\Gamma$. Then, even if
the next polynomial $f_{s'}$ equipartitions the incidences of each 
component of $f'$, that does not imply that each 
component induced by the polynomial $f' \cdot f_{s'}$ will have the same number
of incidences with $\Gamma$. 

The key idea here is to compute $f_0, \ldots, f_s$ \emph{simultaneously}. Thus we will identify, as before, the coefficients
of each $f_i$ with some $\myS^{d_i}$, but now instead of applying
the Borsuk-Ulam theorem separately for each $i$, we will consider the
product of these spheres, and apply the following natural variant
of the Borsuk-Ulam theorem to get the required polynomials
$f_0, \ldots, f_s$ in one step.
For an integer $s\geq 1$, let $X_s = \prod_{i=1}^s {\myS}^{2^{i-1}}$, and
 for each $i \in[s]$, define the functions
\begin{align*}
Fl_i(\ve x_1, \ldots, \ve x_{i-1}, \ve x_i, \ve x_{i+1}, \ldots, \ve x_s) = (\ve x_1, \ldots, \ve x_{i-1}, -\ve x_i, \ve x_{i+1}, \ldots, \ve x_s).
\end{align*}

\noindent For each $v \in \mathbb{Z} \setminus \{0\}$, let $f_v: X_s \rightarrow \R$
be a continuous function  with the property that $f_v(Fl_i(\ve x)) = (-1)^{v_i}f_v(\ve x)$.
Then there exists a point $\ve x \in X_s$ where $f_v(\ve x) = 0$ for all $v \in \mathbb{Z} \setminus \{0\}$.
 
We conclude this section with an open problem -- the  affine  version of Theorem~\ref{thm:polynomialvaritiespartitioning}.

\begin{conjecture}
For any set $\HH_1$ of $n$ $k_1$-dimensional flats,
any set $\HH_2$ of $m$ $k_2$-dimensional flats  in $\R^d$,
and any integer $r$, there exists a partition of $\R^d$ into
$O(r^d)$ simplices such that
(a) each simplex intersects $O\big(\frac{n}{r^{d-k_1}}\big)$ flats of $\HH_1$, and
(b) each flat in $\HH_2$ intersects $O\big(r^{k_2}\big)$ simplices.
\end{conjecture}

\subsection{Depth of point sets} \label{depthsec}

Data science aims to understand the features of data sets.
The goal of data depth measures is to generalize the idea 
of the statistical median of a set of reals to higher dimensions:
the data consists of a finite set $P$ of  points in $\RR^d$,
and  the goal  is to compute a point $\ve q \in \RR^d$ that is a
  ``combinatorial center'' of the data $P$.
 As we will see, there are several natural ways to measure data depth, 
and they are related to each other in sometimes surprising ways. 

\begin{figure}
\centering
  \begin{minipage}{0.5\textwidth}
  \centering
  \includegraphics[width=50mm,page=24]{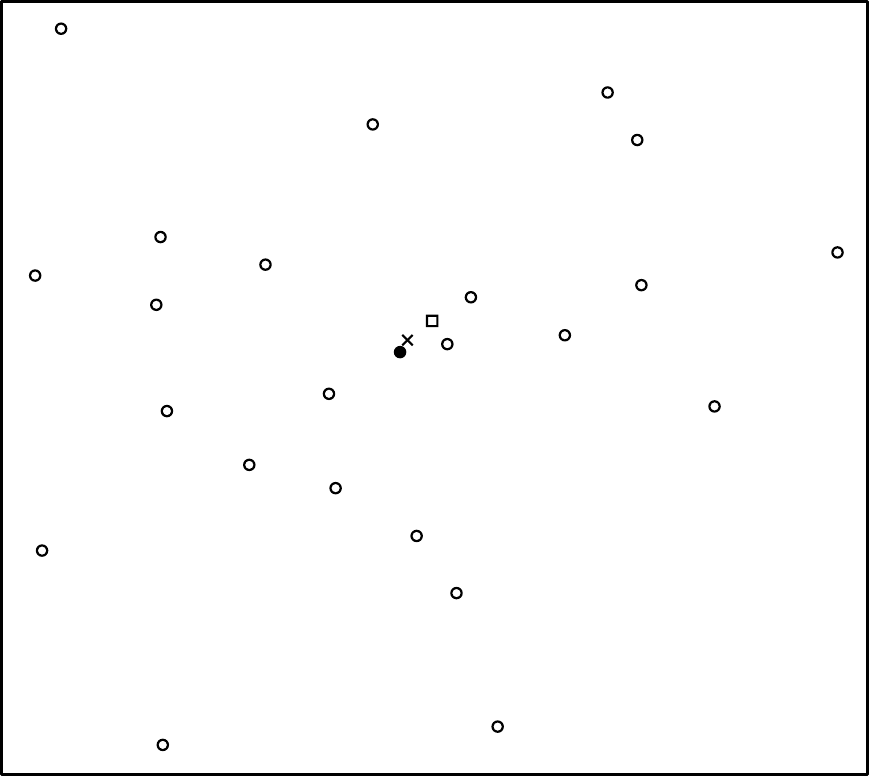}
  \end{minipage}%
  \begin{minipage}{0.5\textwidth}
  \centering
  \includegraphics[width=50mm,page=25]{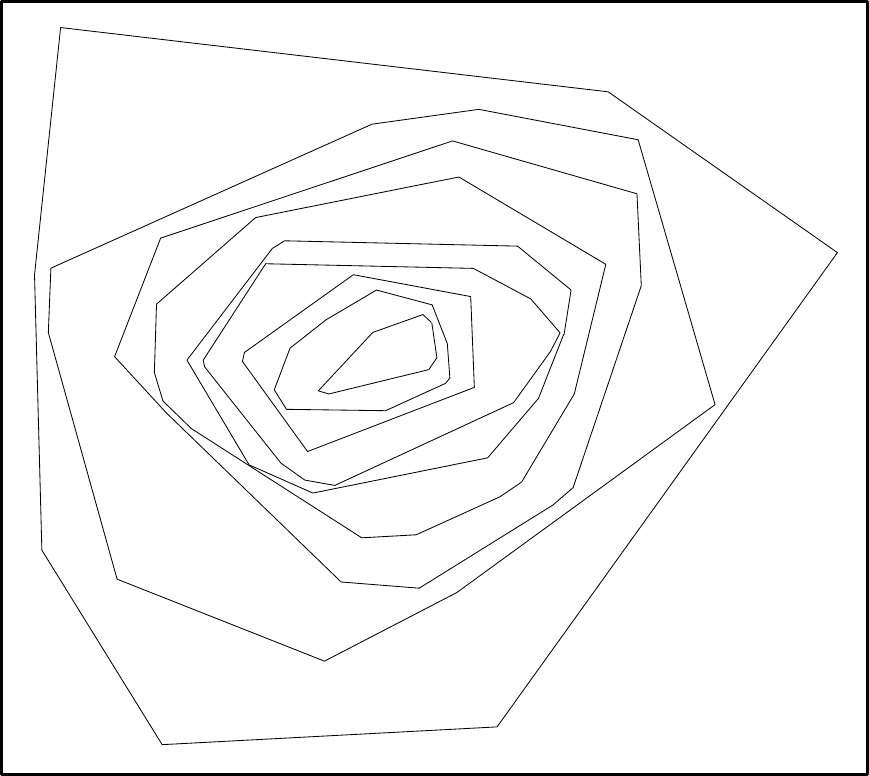}
  \end{minipage}
\caption{(left) A set of points in $\R^2$ together with three centers
under  halfspace   (cross),
simplicial   (box) and Oja   (disk) depth measures, and (right)
the $\beta$-deep regions under halfspace depth measure.}
\label{fig:datadepth}
\end{figure}

Figure~\ref{fig:datadepth} shows a set of points in $\RR^2$ (circles), with
``combinatorial centers'' under three different measures: halfspace depth (cross),
simplicial depth (box) and Oja depth (shaded disk). As the figure shows,
the points for these three measures are geometrically close.
Given integers $d$ and $n$, let
$\P^d_n$ be the point set in $\RR^d$ of size $n$,
with at least $\lfloor \frac{n}{d+1} \rfloor$ points placed
at each of the vertices of the standard simplex.
Slightly perturb each point so that all $n$ points are distinct,
and in general position.
This point set will be very useful for the remainder
of this section, as it essentially captures both the intuition
as well as the worst-case behavior with respect to many
depth measures.

\subsubsection*{Halfspace depth.} \label{half-space depth:subsec}
Given a set $P$ of $n$ points in $\RR^d$, the
\emph{halfspace depth} of a point $\ve q \in \RR^d$  with respect to $P$ 
is the minimum number of points of $P$ 
in any closed halfspace containing $q$:
$$  \Halfspacedepth(\ve q, P) = \min_{\textrm{halfspace } H, \ \ve q \in H} |H \cap P|.$$ 
Define the halfspace depth of $P$ as the maximum
halfspace depth of any point in $\RR^d$ (this has also been
called \emph{Tukey depth}~\cite{T75}).
The separation theorem implies
that any point outside $\conv(P)$ has halfspace depth zero.
It is a non-trivial fact, first shown by Rado in 1947, that points of high halfspace depth   exist
for every point set.

\begin{theorem}[Centerpoint theorem~\cite{R47}]
Any set  of $n$ points in $\RR^d$ has half\-space
depth at least $\lceil \frac{n}{d+1} \rceil$. 
\label{thm:centerpointthm}
\end{theorem}

Recall such a point is called a \emph{centerpoint} of $P$ and we saw its importance and generalizations
in various places along this survey (e.g., in Section \ref{centerpointsinoptima}). It turns out the centerpoint 
theorem is optimal, in the sense that the bound $\lceil \frac{n}{d+1} \rceil$ cannot be improved
$\P^d_n$ is an example of  a point set where it is not possible to do better.
By now there are several proofs of the centerpoint theorem:
using Brouwer's fixed-point theorem~\cite{C69}, using Helly's theorem~\cite{Mbook}, 
following from Tverberg's theorem, and an elementary extremal argument
by induction on the dimension $d$~\cite{MR09b}. 
Perhaps the following proof is the simplest: 
observe that any point $q \in \RR^d$ hitting all convex  objects
containing greater than $\frac{d}{d+1}n$ points of $P$ is a centerpoint, 
whose existence now follows from Helly's theorem.

The  centerpoint theorem and its generalizations have found several applications
in combinatorial geometry, statistics, geometric algorithms, meshing, and 
related areas. A beautiful example is by Miller and Thurston~\cite{MT90}, who showed
that given a set $\D$ of $n$ disjoint disks in the plane, there exists
another disk $B \subset \RR^2$ intersecting $O(\sqrt{n})$ disks of $\D$,
and with at least $\frac{n}{4}$ disks of $\D$ lying completely
in the two connected components of $\RR^2$ induced by $B$. To see this,
use an inverse stereographic projection to lift the centers of the disks of $\D$
to a set $P$ of points lying on a carefully chosen sphere
in $\RR^3$; then, with high probability,
 the image of the intersection of a random hyperplane through
the centerpoint of $P$ with the sphere is the required disk $D$!
  
A point of highest halfspace depth with respect to $P$ is called a
\emph{Tukey median} of $P$. It may not be unique.
In general, the set of points of halfspace depth at least $\beta n$,
for   $0 \leq \beta \leq 1$, forms a convex region called the $\beta$-deep region of $P$. 
It is the intersection of all halfspaces containing more than $(1-\beta)n$ 
points of $P$. Each facet of this region is supported by a hyperplane that passes 
through $d$ points of $P$. Figure~\ref{fig:datadepth} shows that set of all such
regions for the earlier point set. Mart{\'{\i}}nez-Sandoval and {Tamam}, 
gave a generalization of Tukey depth in \cite{leo+tamam} which connects depth
to fractional Helly theorems. 

\textit{Algorithms.} There has been considerable work on the algorithmic question
of computing points of large halfspace depth. We first discuss the 
case in $\R^2$, which is by now settled. A centerpoint 
of $n$ points in $\RR^2$ can be computed in $O(n)$ time~\cite{JM93},
the key  tool being the linear-time
algorithm for computing ham sandwich partitions of two point sets in the plane. 
Chan~\cite{C04} gave an $O(n \log n)$ time randomized algorithm for computing a
point of the highest halfspace depth, i.e., a Tukey median for a set of points in the plane.
The set of all depth contours of $n$ points in $\RR^2$
can be computed in time $O(n^2)$~\cite{MRRSSSS01}.
A real-time GPU-based algorithm for computing the set of all deep regions
of a two-dimensional point set was given in~\cite{KMV02}.
Turning to $\RR^d$, $d \geq 3$, the current-best algorithms for both
computing any centerpoint and the highest depth point take $O(n^{d-1})$ time~\cite{C04}. Clarkson et al.~\cite{CEMST96} presented an iterative method
to compute approximate centerpoints:
the algorithm constructs a $(d+2)$-ary tree $T$, where the $n$ leaves
of $T$ are the input points, and each internal node represents
the Radon point (namely, the unique
point lying in the common intersection
of the   convex-hulls of the two Tverberg partitions of these $d+2$ points)
of its $d+2$ children. This method was
improved to the current best algorithm~\cite{MW13} which computes
a point of halfspace depth at least $\frac{n}{4(d+1)^3}$ in time $d^{O(\log d)}n$.
In fact, this method computes an approximate Tverberg partition; namely,
a partition of $P$ into
$\lceil \frac{n}{4(d+1)^3} \rceil$ sets whose convex-hulls
have a common intersection.

\begin{oproblem} 
Can a centerpoint of  $n$ points in $\RR^3$ be computed
in   $\tilde{O} (n)$ time?
\end{oproblem}

\subsubsection*{Simplicial Depth.} 
A straightforward implication
of Proposition~\ref{p:cara++} is 
that given a set $P$ of $n$ points in $\R^{d}$, 
any point $\ve q \in \conv(P)$ lies in the convex hull
of at least $n - d$ tuples of $\binom{P}{d+1}$.
In fact, any centerpoint must be
contained in many more simplices spanned by $(d+1)$-tuples
of points in $P$. Not surprisingly,
the number of $(d+1)$-tuples of $P$ whose convex hull
contains $\ve q$ is positively correlated 
to the halfspace depth of $\ve q$. 
This leads to the related   depth measure
of \emph{simplicial depth}, first defined by Liu~\cite{L88}, 
which is the number of simplices spanned by points of $P$ containing a given point:
$$\Simplicialdepth(\ve q, P) =  \left|  \left\{ Q \in \binom{P}{d+1} \colon \ve q \in \conv(Q)  \right\}  \right|.$$ 
The simplicial depth of $P$ is the highest simplicial depth of any point $\ve q \in \RR^d$. 
As mentioned earlier, there is a close relation between  halfspace depth
and simplicial depth; the current best bound~\cite{W03} 
shows that a point of halfspace depth $\tau n$ has simplicial depth at least
\begin{equation*} 
\frac{(d+1)\tau^d - 2d\tau^{d+1}}{(d+1)!} \cdot n^{d+1}  - O(n^d).
\end{equation*}

B\'{a}r\'{a}ny~\cite{baranys-caratheodory} showed that the colorful Carath\'eodory's theorem together
with Tverberg's theorem implies that there always exists a point contained in at least 
$\frac{1}{d! (d+1)^{d+1}} \cdot n^{d+1}$ simplices spanned by $P$.
Let $c_d$ be a constant such that any set $P$ of $n$ points  has simplicial depth at least 
$c_d \cdot n^{d+1}$. The optimal dependency on  $c_d$  is a long-standing open problem.
Bukh, Matou\v{s}ek, and Nivash~\cite{BMN08} constructed $n$ points in $\RR^d$ so 
that no point in $\RR^d$ is contained in, up to lower-order terms, more than $(\frac{n}{d+1})^{d+1}$ 
simplices defined by $P$. Furthermore, they conjectured that this is the optimal bound.

\begin{conjecture} 
Any set  of $n$ points in $\RR^d$ has simplicial depth at least $(\lceil \frac{n}{(d+1)}\rceil)^{d+1}$.
\end{conjecture}
 
For $d=2$, a positive answer to the above conjecture
was  known already in 1984 by Boros and F\"{u}redi~\cite{BF82}. 
Bukh~\cite{B05}  gave a beautiful short proof:   
the required point set is the common intersection point of three lines in $\RR^2$, 
having a common intersection and where each of the six induced cones contain at least 
$\frac{n}{6}$ points of $P$; the existence of such three lines follows via 
an elementary topological argument. For $d=3$, an elementary argument shows that $c_3 \geq 0.0023$~\cite{BMRR10b}.
Using algebraic topology machinery, Gromov~\cite{G10} improved the bound to the value $c_3 \geq 0.0026 $. 
This bound for $\RR^3$ has since been improved even further by Matou\v{s}ek and Wagner~\cite{MW14} to $0.00263$, 
and then by Kr\'al et al.~\cite{KMS12} to $0.0031$.  In fact, Gromov proves the bound for general $d$, 
showing that 
$$  c_d \geq \frac{2d}{(d+1)!^2(d+1)}.$$

His proof has since been simplified by Karasev~\cite{K12}. Using the concepts of Gale diagrams and secondary polytopes, 
one can observe that the concept of simplicial depth is   equivalent to a different maximization problem: given a set of points 
on the sphere $\myS^d$, what is the triangulation that uses those points (but possibly not all) as vertices that has the largest 
number of $d$-dimensional simplices? (see~\cite[Chapter 5]{triangbook} and the references therein).

We conclude this discussion of simplicial depth
with \emph{colorful simplicial depth}, which was introduced by Deza et al.~\cite{Dezaetal2006}.
Consider a set of points $P$ in $\R^d$ partitioned into $(d+1)$
color classes  $P = P_0 \cup \cdots \cup P_d$.
Suppose that $P$ has the property that the origin $\ve o$ is in the relative interior of each $\conv(P_i)$, for $0 \le i \le d$. 
Recall from Section~\ref{caratheodory-thms} that a colorful simplex is a simplex where each of the vertices comes from a different $P_i$.
While the colorful Carath\'eodory theorem asserts the existence of at least
one colorful simplex containing $\ve o$,  one can further ask about the  number of distinct
colorful simplices containing $\ve o$ that must always exist.
Define the \emph{colorful simplicial depth} of  $P$, denoted   $\csd(P)$, as the number of colorful simplices in $P$ containing $\ve o$.
Deza et al.~\cite{Dezaetal2006} proposed some lower bounds on the colorful simplicial depth, and 
they conjectured that if $|P_i| = d+1$ for  $0 \leq i \leq d$,   then $\csd(P) \geq d^2 + 1$. This was proven by Sarrabezolles~\cite{S15}. The bound is optimal by the work of 
Deza et al.~\cite{Dezaetal2006}. They also conjectured  the following upper bound which was 
shown by Adiprasito \etal~\cite{adiprasito2016colorful}:
let $P = P_0 \cup \ldots \cup P_d$ be a point set in $\R^d$ with $|P_i| \ge 2$ for all $0 \le i \le d$. 
If no colorful simplex $S$ spanned by $P$ of dimension $d-1$ contains the origin $\ve o$ in its convex hull, then $ \csd(P) \ \le \ 1 + \prod_{i=0}^d (|P_i| - 1)$.

\textit{Algorithms.} For the case of the plane, computing the simplicial depth of a query point
can be done in time $O(n \log n)$~\cite{GSW92},
which is optimal. Aloupis et al.~\cite{ALST03} presented
an algorithm to compute a point of highest simplicial
depth in $\RR^2$ in time $O(n^4)$. Using the fact that
finding the highest simplicial depth is Gale dual to the problem
of finding a maximum triangulation of points in the sphere, then
one can set up an integer program to find a point of largest
simplicial depth for any point configuration (see \cite[Chapter 8]{triangbook}).
A GPU-based algorithm for computing simplicial depth and colorful simplicial depth
of point sets in the plane was given in~\cite{KMV06}. 

\subsubsection*{Ray-shooting depth.}
It turns out that the previous two measures -- halfspace depth
and simplicial depth -- are further related to each other
via an even more general depth measure, called \emph{ray-shooting depth}. 
Given a set $P$ of $n$ points in $\RR^d$, 
let $E_P$ be the set of all ${n \choose d}$ $(d-1)$-simplices spanned
by points of $P$. Given a point $\ve q \in \RR^d$ and a direction
$\ve u\in \mathbb{S}^{d-1}$, let $r(\ve q, \ve u)$ be the half-infinite ray
from $\ve q$ in direction $\ve u$. Then the \emph{ray-shooting depth}
of a point $\ve q \in  \RR^d$ is defined as
$$\Rayshootingdepth(\ve q, P) = \min_{\ve u\in \mathbb{S}^{d-1}} \big| \big\{ e \in E_P \colon  r(\ve q,\ve u) \cap e \neq \varnothing \big\} \big|.$$

The ray-shooting depth of $P$ is the maximum ray-shooting depth of any point in $\RR^d$.
The notion of ray-shooting depth was first introduced in Fox et al.~\cite{FGLNP11}, who
proved the following using Brouwer's fixed point theorem.

\begin{theorem} 
Any set  of $n$ points in $\RR^2$ has ray-shooting depth at least $\frac{n^2}{9}$.
\end{theorem}

Note that any point realizing the maximum ray-shooting depth must have halfspace depth at least
$\frac{n}{3}$ \emph{and} simplicial depth at least $\frac{n^3}{27}$: let $\ve q$ be a point with
ray-shooting depth at least $\frac{n^2}{9}$. Then any line through $\ve q$ must intersect
at least $\frac{2n^2}{9}$ segments in $E_P$, so both halfspaces defined by it
must contain at least $\frac{n}{3}$ points. For simplicial depth, 
consider, for each point $\ve p \in P$, the ray from $\ve q$ in the direction $\overrightarrow{\ve p \ve q}$. 
Then for every edge $\{\ve p_i, \ve p_j\}$ that intersects this ray, the triangle defined
by $\{\ve p, \ve p_i, \ve p_j\}$ must contain $\ve q$. Summing up these triangles
over all points, each triangle can be counted three times, and so
$\ve q$ lies in at least $\frac{n^3}{27}$ distinct triangles spanned by $P$.

The problem of showing the existence
of a point with large ray-shooting depth is open in higher dimensions.

\begin{conjecture}
Any set of $n$ points in $\RR^d$ has ray-shooting depth at least $\big( \lceil \frac{n}{d+1} \rceil\big)^d$.
\end{conjecture}

Other notions of ray-shooting depth for convex sets, instead of points sets, were studied in \cite{Fulek:2009vz}.

\textit{Algorithms.}
The proof in~\cite{FGLNP11} is topological and does not give a method to
compute such a point. A combinatorial proof,
together with efficient algorithms were obtained in Mustafa et al.~\cite{MRS11},
who showed how compute a point of ray-shooting depth 
at least $\frac{n^2}{9}$ in time $\tilde{O}(n^2)$.

\subsubsection*{Oja depth.}
It turns out that ray-shooting depth
is related to another older 
depth measure, the \emph{Oja depth} of a point set,
first defined by Oja~\cite{O83}.
Assume, without loss of generality, that $\vol(\conv(P)) = 1$. Then define
the Oja depth of a point $\ve q$ with respect to $P$ as
$$\Ojadepth(\ve q, P) =  \sum_{\substack{P' \subseteq P\\|P'|=d}} \vol \big(\conv(P' \cup \{\ve q\})\big).$$

The Oja depth of $P$ is the minimum Oja depth over all $\ve q \in \RR^d$. 
It is easy to see that the Oja depth of $\P^d_n$ is at least $(\frac{n}{d+1})^d$.
The conjecture~\cite{CDILM10} is that the lower bound given by $\P^d_n$ is essentially tight.

\begin{conjecture} 
The Oja depth of any set  of $n$ points in $\RR^d$ is
at most $\big(\frac{n}{d+1}\big)^d.$
\end{conjecture}

The conjecture has been  resolved only for the case $d=2$
by Mustafa \etal~\cite{MTW14}.
For general $d$, it can be shown that the center of mass of $P$ has Oja
depth at most ${n \choose d}/(d+1)$~\cite{CDILM10}.
This estimate can be improved via ray-shooting depth, as the Oja depth
of any point $\ve q$ which has ray-shooting depth at least $\frac{n^2}{9}$ is
at least $\frac{n^2}{7.2}$. The reason is  the number of triangles  spanned by pairs of points in $P$
and the point $\ve q$,  containing any point $\ve p\in \RR^2$,
is at most the number of edges spanned by $P$
intersecting the ray $\overrightarrow{\ve q \ve p}$, which is at most $\frac{n^2}{4} - \frac{n^2}{9} = \frac{n^2}{7.2}$. 
Integrating over all $\ve p \in \RR^2$ gives the required bound.
A  calculation in $\RR^d$ gives the current-best bound~\cite{MTW14}:

\begin{theorem} 
Every set  of $n$ points in $\RR^d$, $d \geq 3$, has Oja depth at most 
$$ \frac{2n^{d}}{2^{d}d!} - \frac{2d}{(d+1)^2(d+1)!} {n \choose d} + O(n^{d-1}). $$
\end{theorem}

\textit{Algorithms.}
For the case $d=2$, Rousseeuw and Ruts~\cite{RR96} presented
a  $O(n^5 \log n)$ time algorithm for computing the lowest depth point, which 
was then improved to the current-best algorithm with running time $O(n \log^3 n)$~\cite{ALST03}. A 
point of Oja depth at most $\frac{n^2}{9}$ 
can be computed in $O(n \log n)$ time~\cite{MTW14}.
For general $d$, various heuristics for computing points with low Oja depth were given by 
Ronkainen, Oja and Orponen~\cite{ROO01}. 

\subsubsection*{Regression depth.}
The next depth measure, unlike
earlier measures, is a combinatorial
analogue of fitting a hyperplane through a set of points.
Therefore it will be more convenient to state
it in the dual. Given a point $\ve p \in \RR^d$, let $p^*$
be its dual hyperplane, and for a set
of points $P$, let  $P^* = \{p^* \colon \ve p \in P\}$.
Then define the \emph{regression depth} of a point as
$$\Regressiondepth(\ve q, P) = \min_{\ve u\in \mathbb{S}^{d-1}} \big| \big\{  H \in P^* \colon r(\ve q, \ve u) \cap H \neq \varnothing \big\} \big|.$$
The regression depth of a set $P$ of points in $\RR^d$
is the maximum regression depth of any point $\ve q \in \RR^d$.
It was introduced  by Rousseeuw and Hubert~\cite{RH99},
who showed that any set $P$ of $n$ points in $\RR^2$ has
regression depth at least $\lceil \frac{n}{3} \rceil$. Their proof
is elegant: given the set $P$ of $n$ points, let $P_1, P_2, P_3$
be a partition of $P$ by consecutive $x$-coordinate values,
and where $|P_i| \leq \lceil \frac{n}{3} \rceil$ for $i=1,2,3$.
Then the required line is the ham sandwich cut 
of the two sets $P_1 \cup P_2$ and $P_2 \cup P_3$.
The optimal bound for general $d$ was discovered later.

\begin{theorem}[\cite{ABET00, K08, Mi02}]
Any set  of $n$ points in $\RR^d$ has regression depth at least $\lceil \frac{n}{d+1} \rceil$.
\end{theorem}

\noindent
Given a set $X \subseteq \R^d$ and a point $q \in \R^d$, the closest
point in $X$ to $q$ (if it exists) is denoted by $\closestpoint(q,
X)$. The proof in~\cite{K08} deduces it from the centerpoint theorem:
define the function $f(\ve q)$ that maps $\ve q \in \RR^d$ to a
centerpoint of the set $\{ \closestpoint(\ve q, p^*) \colon \ve p \in
P\}$.  This can be done so that $f(\cdot)$ is continuous and maps a
sufficiently large ball to itself; then observe that the dual of any fixed
point of $f(\cdot)$ is the required hyperplane.

\textit{Algorithms.} The method of~\cite{RH99} gives a 
linear-time algorithm for computing a point of regression
depth at least $\lceil \frac{n}{3} \rceil$ immediately, as it uses
only ham sandwich cuts. A point of maximum regression depth in $\RR^2$ can
be computed in time $O(n \log n)$~\cite{MMRSSS08}, 
improving upon an earlier $O(n\log^2 n)$  time algorithm~\cite{LS03}.
For $d \geq 3$, the best algorithm takes time $O(n^d)$~\cite{MMRSSS08}
to compute a point of maximum regression depth.

\subsubsection*{The $k$-centerpoint conjectures.}
It turns out that  many of the discussed depth measures
are special cases of the following more general conjecture,
first proposed by Mustafa et al.~\cite{MRS15}. See also the related paper \cite{Karasev:2011jv}.

\begin{conjecture} 
For any set $P$ of $n$ points in $\RR^d$ and any integer 
$0 \leq k \leq d$, there exists a point
$\ve q \in \RR^d$ such that any $(d-k)$-half flat through $\ve q$ 
intersects at least $\big(\lceil \frac{n}{d+1} \rceil \big)^{k+1}$ of the  $k$-simplices spanned by $P$.
\end{conjecture}

The case $k=0$ corresponds to halfspace depth, 
$k=d$ to simplicial depth, and $k=d-1$ to ray-shooting depth.

 It is not hard to show that given a set $P$ of $n$ points in $\RR^d$, and an integer $0 \leq k \leq d-1$, there exists a point 
$\ve q \in \RR^d$ such that any $(d-k)$-half flat through $\ve q$ intersects at least 
$$ \max \left\{  {\frac{n}{d+1} \choose k+1}, \ \ \ \frac{2d}{(d+1)(d+1)! {n \choose d-k}} \cdot {n \choose d+1}  \right\} $$
$k$-simplices spanned by $P$. For simplicity assume that $|P|$ is a multiple of $(d+1)$.
The proof follows from the use of Tverberg's theorem to partition $P$ into $t = \frac{n}{(d+1)}$ sets $P_1, \ldots, P_t$
such that there exists a point $\ve q$ with $\ve q \in \conv(P_i)$ for all $i$. Consider any $(d-k)$-dimensional
half-flat $\mathcal{F}$ through $\ve q$, where $\partial \mathcal{F}$ is a $(d-k-1)$-dimensional flat containing $\ve q$. Project
$\mathcal{F}$ onto a $(k+1)$-dimensional subspace $\mathcal{H}$ orthogonal to $\partial \mathcal{F}$ 
such that the projection of $\mathcal{F}$ is a ray $\ve r$ in $\mathcal{H}$, and $\partial \mathcal{F}$ and $\ve q$ are 
projected to the point $\ve q'$. And let $P'_1, \ldots, P'_t$ be the projected sets whose convex-hulls now contain the point $\ve q'$. 
Then note that the $k$-dimensional simplex spanned by $(k+1)$ points $Q' \subset P'$ intersects the ray $\ve r$ if and only if
the $k$-dimensional simplex defined by the corresponding set $Q$ in $\RR^d$ intersects the flat $\mathcal{F}$. Now
apply the single-point version (i.e., given any point $\ve s\in \RR^d$ and $d$ sets $P_1, \ldots, P_d$
in $\RR^d$ such that each $\conv(P_i)$ contains the origin, there exists a $d$-simplex spanned 
by $\ve s$ and one point from each $P_i$ 
which also contains the origin) of colorful Carath\'{e}odory's theorem  to every $(k+1)$-tuple
of sets, say $P'_1, \ldots, P'_{k+1}$,  together with the point $\ve s$ at infinity in the direction 
antipodal to the direction of $\ve r$ to get a 
`colorful' simplex, defined by $\ve s$ and one point from each $P'_i$, and containing $\ve q'$. Then
the ray $\ve r$ must intersect the $k$-simplex defined by the $(k+1)$ points of $P'$, and 
so the corresponding points of $P$ in $\RR^d$ span a $(k+1)$-simplex intersecting $\mathcal{F}$. 
In total, we get $n/(d+1) \choose k+1$ of the $k$-simplices intersecting $\mathcal{F}$.
Another way is to use the simplicial depth bound of Gromov, that
given any set $P$ of $n$ points in $\RR^d$, there exists a point $\ve q$
lying in $2d/((d+1)(d+1)!) \cdot {n \choose d+1}$ $d$-simplices. Now take any
$(d-k)$-half flat through $\ve q$. It must intersect at least one $k$-simplex
of each $d$-simplex containing it, where each $k$-simplex 
is counted at most $n \choose d-k$ times. This implies the stated bound.

In the plane ($d=2$) the centerpoint theorem, can be re-stated as follows:
   Given a set P of $n$ points in the plane, there exists a point $q$ such that if you take any line L passing through $q$, and move 
   it continuously until it arrives outside the convex hull of $P$, then along this motion, the line will intersect at least $n/3$ points of $P$.
The following more general statement has also been conjectured in~\cite{MRS15}:

\begin{conjecture}\label{lastconjecture}
Given a set $P$ of $n$ points in $\R^d$, and an integer $0 \leq k \leq d$, 
there exists a point $q \in \R^d$ such that the following holds:
Let $F_q, F_{o}$ be two $\left(d-k-1\right)$-flats, such that $q \in F_q$ and
$F_{o}$ does not intersect the convex-hull of $P$. 
Then any continuous motion family of $\left(d-k-1\right)$-flats, starting at $F_q$ and ending at $F_o$, must intersect
at least $\left( \lceil \frac{n}{d+1} \rceil \right)^{k+1}$ $k$-simplices spanned by  $P$.
\end{conjecture}

In Conjecture \ref{lastconjecture}, the case $k=d$ gives a `-$1$'-flat moving
to infinity, which can be treated as a stationary point.
The validity of these conjectures for the planar case $d=2$ follows from  the work of Gromov~\cite{G10}.

\section*{Acknowledgments} 
J. A. De Loera  was partially supported by the LabEx
Bezout (ANR-10-LABX-58). He is grateful to the labex BEZOUT and the CERMICS
research center at \'Ecole National des Ponts et Chauss\'ees for the support
received, and the enjoyable and welcoming environment in which the
topics in this paper were discussed and built over several visits. 
De Loera was also partially supported by NSF grant DMS-1522158.
Mustafa was supported by the grant ANR SAGA (JCJC-14-CE25-0016-01).
Goaoc was partially supported by Institut Universitaire de France.

The authors are truly grateful to the two anonymous referees who
provided a very large and very detailed set of comments, insights, and
corrections. Their reading of this work was truly invaluable and it
took a lot of effort and time.  Very special thanks to Peter Landweber
who gave us numerous corrections too. We are also grateful to the
following colleagues who gave us support and many comments and references: Ian Agol,
Imre B\'ar\'any, Pavle Blagojevi\'c, Boris Bukh, Sabrina Enriquez, Florian Frick, Alexei Garber, 
Tommy Hogan, Andreas Holmsen, Wolfgang Mulzer, Oleg Musin, Deborah Oliveros, J\'anos
Pach, D\H{o}m\H{o}t\H{o}r P\'alv\H{o}lgyi, Axel Parmentier, Lily Silverstein, Pablo Sober\'on, 
Francis Su, Uli Wagner, and G\"unter M. Ziegler.  We are grateful for their help.  \bibliographystyle{acm}
\bibliography{references}

\def\cprime{$'$}
\begin{thebibliography}{100}

\bibitem{Aardal2007}
{\sc Aardal, K.~I., van Hoesel, S. P.~M., Koster, A. M. C.~A., Mannino, C., and
  Sassano, A.}
\newblock Models and solution techniques for frequency assignment problems.
\newblock {\em Ann. Oper. Res. 153}, 1 (2007), 79--129.

\bibitem{adamaszek2011lower}
{\sc Adamaszek, M., and Barmak, J.~A.}
\newblock On a lower bound for the connectivity of the independence complex of
  a graph.
\newblock {\em Discrete Math. 311}, 21 (2011), 2566--2569.

\bibitem{adiprasito+barany+mustafa}
{\sc {Adiprasito}, K., {B{\'a}r{\'a}ny}, I., and {Mustafa}, N.}
\newblock Theorems of {C}arath\'eodory, {H}elly, and {T}verberg without
  dimension.
\newblock {\em \rm ArXiv:1806.08725\/} (2018).

\bibitem{adiprasito2016colorful}
{\sc Adiprasito, K.~A., Brinkmann, P., Padrol, A., Paták, P., Patáková, Z.,
  and Sanyal, R.}
\newblock Colorful simplicial depth, {M}inkowski sums, and generalized {G}ale
  transforms.
\newblock {\em Int. Math. Res. Not. IMRN\/} (2017), rnx184.

\bibitem{adler-minmaxvsduality}
{\sc Adler, I.}
\newblock The equivalence of linear programs and zero-sum games.
\newblock {\em Internat. J. Game Theory 42}, 1 (2013), 165--177.

\bibitem{AMS13}
{\sc Agarwal, P.~K., Matou{\v{s}}ek, J., and Sharir, M.}
\newblock On range searching with semialgebraic sets. {II}.
\newblock {\em SIAM J. Comput. 42}, 6 (2013), 2039--2062.

\bibitem{ABZ07}
{\sc Aharoni, R., Berger, E., and Ziv, R.}
\newblock Independent systems of representatives in weighted graphs.
\newblock {\em Combinatorica 27\/} (2007), 253--267.

\bibitem{AhHo98}
{\sc Aharoni, R., and Holzman, R.}
\newblock Fractional kernels in digraphs.
\newblock {\em J. Combin. Theory Ser. {B} 73\/} (1998), 1--6.

\bibitem{Maya+Yo+Raymond}
{\sc Ahmed, M., De~Loera, J., and Hemmecke, R.}
\newblock Polyhedral cones of magic cubes and squares.
\newblock In {\em Discrete and computational geometry}, vol.~25 of {\em
  Algorithms Combin.} Springer, Berlin, 2003, pp.~25--41.

\bibitem{TheBook}
{\sc Aigner, M., and Ziegler, G.~M.}
\newblock {\em Proofs from THE BOOK}, 5th~ed.
\newblock Springer, 2014.

\bibitem{aisenbergetal2015}
{\sc Aisenberg, J., Bonet, M.~L., and Buss, S.}
\newblock 2-d {T}ucker is {PPA}-complete.
\newblock {\em Electronic Colloquium on Computational Complexity {(ECCC)} 22\/}
  (2015), 163.

\bibitem{alievetal}
{\sc Aliev, I., Bassett, R., De~Loera, J.~A., and Louveaux, Q.}
\newblock A quantitative {D}oignon-{B}ell-{S}carf theorem.
\newblock {\em Combinatorica 37}, 3 (2017), 313--332.

\bibitem{Alievetal2018}
{\sc Aliev, I., De~Loera, J., Eisenbrand, F., Oertel, T., and Weismantel, R.}
\newblock The support of integer optimal solutions.
\newblock {\em SIAM J. Optim. 28}, 3 (2018), 2152--2157.

\bibitem{sparseintcaratheo}
{\sc Aliev, I., {De Loera}, J.~A., Oertel, T., and O'Neill., C.}
\newblock {Sparse solutions of linear {D}iophantine equations.}
\newblock {\em SIAM J. Appl. Algebra Geom. 1\/} (2017), 239--253.

\bibitem{alishahi2017}
{\sc Alishahi, M.}
\newblock Colorful subhypergraphs in uniform hypergraphs.
\newblock {\em Electron. J. Combin. 24}, 1 (2017), Research Paper P1.23.

\bibitem{alishahi2015chromatic}
{\sc Alishahi, M., and Hajiabolhassan, H.}
\newblock On the chromatic number of general {K}neser hypergraphs.
\newblock {\em J. Combin. Theory Ser. B 115\/} (2015), 186--209.

\bibitem{alishahi2016strengthening}
{\sc Alishahi, M., Hajiabolhassan, H., and Meunier, F.}
\newblock Strengthening topological colorful results for graphs.
\newblock {\em European J. Combin. 64\/} (2017), 27 -- 44.

\bibitem{alishahi2017fair}
{\sc Alishahi, M., and Meunier, F.}
\newblock Fair splitting of colored paths.
\newblock {\em Electr. J. Comb. 24}, 3 (2017), P3.41.

\bibitem{Alon:1987ur}
{\sc Alon, N.}
\newblock {Splitting necklaces}.
\newblock {\em Adv. Math. 63}, 3 (1987), 247--253.

\bibitem{ABFK92}
{\sc Alon, N., B{\'a}r{\'a}ny, I., F{\"u}redi, Z., and Kleitman, D.~J.}
\newblock Point selections and weak $\epsilon$-nets for convex hulls.
\newblock {\em Combin. Probab. Comput. 1\/} (1992), 189--200.

\bibitem{AKMM}
{\sc Alon, N., Kalai, G., Matou{\v{s}}ek, J., and Meshulam, R.}
\newblock Transversal numbers for hypergraphs arising in geometry.
\newblock {\em Adv. in Appl. Math. 29}, 1 (2002), 79 -- 101.

\bibitem{Alon92pq}
{\sc Alon, N., and Kleitman, D.~J.}
\newblock {Piercing convex sets and the Hadwiger-Debrunner $(p,q)$-problem}.
\newblock {\em Adv. Math. 96}, 1 (1992), 103--112.

\bibitem{alon2006algorithmic}
{\sc Alon, N., Moshkovitz, D., and Safra, S.}
\newblock Algorithmic construction of sets for $k$-restrictions.
\newblock {\em ACM Trans. Algorithm. 2}, 2 (2006), 153--177.

\bibitem{AW86}
{\sc Alon, N., and West, D.}
\newblock The {B}orsuk-{U}lam theorem and bissection of necklaces.
\newblock {\em Proc. Amer. Math. Soc. 98\/} (1986), 623--628.

\bibitem{ALST03}
{\sc Aloupis, G., Langerman, S., Soss, M., and Toussaint, G.}
\newblock Algorithms for bivariate medians and a {F}ermat-{T}orricelli problem
  for lines.
\newblock {\em Comput. Geom. 26}, 1 (2003), 69--79.

\bibitem{amenta94}
{\sc Amenta, N.}
\newblock Helly theorems and generalized linear programming.
\newblock {\em Discrete Comput. Geom. 12\/} (1994), 241---261.

\bibitem{ABET00}
{\sc Amenta, N., Bern, M., Eppstein, D., and Teng, S.-H.}
\newblock Regression depth and center points.
\newblock {\em Discrete Comput. Geom. 23}, 3 (2000), 305--323.

\bibitem{amentaetal2017helly}
{\sc Amenta, N., De~Loera, J.~A., and Sober{\'o}n, P.}
\newblock Helly's theorem: new variations and applications.
\newblock In {\em Algebraic and Geometric Methods in Discrete Math.}, vol.~685
  of {\em Contemporary Mathematics}. American Mathematical Soc., 2017,
  pp.~55--95.

\bibitem{ABBFM09}
{\sc Arocha, J.~L., B\'ar\'any, I., Bracho, J., Fabila, R., and Montejano, L.}
\newblock Very colorful theorems.
\newblock {\em Discrete Comput. Geom. 42}, 2 (2009), 142--154.

\bibitem{Arora+Barak:complexity}
{\sc {Arora}, S., and {Barak}, B.}
\newblock {\em {Computational complexity. A modern approach.}}
\newblock Cambridge University Press, 2009.

\bibitem{asadaetal2}
{\sc Asada, M., Chen, R.~Frick, F., Huang, F., Polevy, M., Stoner, D., Tsang,
  L., and Wellner, Z.}
\newblock On {R}eay's relaxed {T}verberg conjecture and generalizations of
  {C}onway's thrackle conjecture.
\newblock {\em Electronic J. Combin. 25}, (3) (2018), P3.16.

\bibitem{asadaetal}
{\sc Asada, M., Frick, F., Pisharody, V., Polevy, M., Stoner, D., Tsang, L.,
  and Wellner, Z.}
\newblock Fair division and generalizations of {S}perner- and {KKM}-type
  results.
\newblock {\em SIAM J. Discrete Math. 32}, (1) (2018), 591--610.

\bibitem{aumann}
{\sc Aumann, R.}
\newblock Subjectivity and correlation in randomized strategies.
\newblock {\em J. Math. Econom. 1}, 1 (1974), 67--96.

\bibitem{Averkov:2013uo}
{\sc Averkov, G.}
\newblock {On maximal S-free sets and the {H}elly number for the family of
  S-convex sets}.
\newblock {\em SIAM J. Discrete Math.\/} (2013).

\bibitem{Averkovetal-tightbounds}
{\sc {Averkov}, G., {Gonz{\'a}lez Merino}, B., {Schymura}, M., {Paschke}, I.,
  and {Weltge}, S.}
\newblock {Tight bounds on discrete quantitative {H}elly numbers}.
\newblock {\em Adv. Appl. Math. 89\/} (2017), 76 -- 101.

\bibitem{AW2012}
{\sc Averkov, G., and Weismantel, R.}
\newblock Transversal numbers over subsets of linear spaces.
\newblock {\em Adv. Geom. 12}, 1 (2012), 19--28.

\bibitem{avis+friedmann}
{\sc Avis, D., and Friedmann, O.}
\newblock An exponential lower bound for {C}unningham's rule.
\newblock {\em Math. Program. 161}, 1-2 (2017), 271--305.

\bibitem{Aziz+Mackenzie2016}
{\sc Aziz, H., and Mackenzie, S.}
\newblock A discrete and bounded envy-free cake cutting protocol for any number
  of agents.
\newblock In {\em 57th {A}nnual {IEEE} {S}ymposium on {F}oundations of
  {C}omputer {S}cience---{FOCS} 2016}. IEEE Computer Soc., Los Alamitos, CA,
  2016, pp.~416--427.

\bibitem{Bajmoczy:1979bj}
{\sc Bajm{\'o}czy, E.~G., and B{\'a}r{\'a}ny, I.}
\newblock {On a common generalization of {B}orsuk's and {R}adon's theorem}.
\newblock {\em Acta Math. Acad. Sci. Hungar. 34}, 3-4 (1979), 347--350.

\bibitem{baranys-caratheodory}
{\sc B{\'a}r{\'a}ny, I.}
\newblock A generalization of {C}arath\'eodory's theorem.
\newblock {\em Discrete Math. 40}, 2-3 (1982), 141--152.

\bibitem{barany1993geometric}
{\sc B{\'a}r{\'a}ny, I.}
\newblock Geometric and combinatorial applications of {Borsuk's} theorem.
\newblock In {\em New trends in discrete and computational geometry}. Springer,
  1993, pp.~235--249.

\bibitem{50tverbergsurvey}
{\sc B{\'a}r{\'a}ny, I., Blagojevi{\'c}, P.~V.~M., and Ziegler, G.~M.}
\newblock Tverberg's theorem at 50: extensions and counterexamples.
\newblock {\em Notices Amer. Math. Soc. 63}, 7 (2016), 732--739.

\bibitem{baranykatchalskipach}
{\sc B{\'a}r{\'a}ny, I., Katchalski, M., and Pach, J.}
\newblock Quantitative {H}elly-type theorems.
\newblock {\em Proc. Amer. Math. Soc. 86}, 1 (1982), 109--114.

\bibitem{Barany:1992tx}
{\sc B{\'a}r{\'a}ny, I., and Larman, D.~G.}
\newblock {A colored version of {T}verberg's theorem}.
\newblock {\em J. Lond. Math. Soc. 2}, 2 (1992), 314--320.

\bibitem{baranymatousek}
{\sc B{\'a}r{\'a}ny, I., and Matou{\v{s}}ek, J.}
\newblock A fractional {H}elly theorem for convex lattice sets.
\newblock {\em Adv. Math. 174}, 2 (2003), 227--235.

\bibitem{barany1995caratheodory}
{\sc B{\'a}r{\'a}ny, I., and Onn, S.}
\newblock {Carath{\'e}odory{\textquoteright}s theorem, colourful and
  applicable}.
\newblock {\em {\rm} in Intuitive Geometry 6\/} (1995), 11--21.

\bibitem{baranyonn-colorfulLP}
{\sc B{\'a}r{\'a}ny, I., and Onn, S.}
\newblock Colourful linear programming and its relatives.
\newblock {\em Math. Oper. Res. 22}, 3 (1997), 550--567.

\bibitem{Barany:1981vh}
{\sc B{\'a}r{\'a}ny, I., Shlosman, S.~B., and Sz{\"u}cs, A.}
\newblock {On a topological generalization of a theorem of {T}verberg}.
\newblock {\em J. Lond. Math. Soc. 2}, 1 (1981), 158--164.

\bibitem{Tverbergplusminus}
{\sc {B{\'a}r{\'a}ny}, I., and {Sober{\'o}n}, P.}
\newblock {Tverberg plus minus}.
\newblock {\em Discrete Comput. Geom.\/} (2018).

\bibitem{barany+soberonsurvey}
{\sc B\'ar\'any, I., and Sober\'on, P.}
\newblock Tverberg's theorem is 50 years old: A survey.
\newblock {\em Bull. Amer. Math. Soc.\/} (2018).

\bibitem{Barman}
{\sc Barman, S.}
\newblock Approximating {N}ash equilibria and dense bipartite subgraphs via an
  approximate version of {Carath\'eodory's} theorem.
\newblock In {\em Proceedings {ACM} on Symposium on Theory of Computing,
  {(STOC)}\/} (2015), pp.~361--369.

\bibitem{Bar2002}
{\sc Barvinok, A.}
\newblock {\em A course in Convexity}, vol.~54 of {\em Grad. Stud. Math.}
\newblock American Mathematical Society, Providence, RI, 2002.

\bibitem{BMRR10b}
{\sc Basit, A., Mustafa, N.~H., Ray, S., and Raza, S.}
\newblock Hitting simplices with points in $\mathbb{R}^3$.
\newblock {\em Discrete Comput. Geom. 44}, 3 (2010), 637--644.

\bibitem{basuetal:S-optimality}
{\sc {Basu}, A., {Conforti}, M., {Cornu{\'e}jols}, G., {Weismantel}, R., and
  {Weltge}, S.}
\newblock {Optimality certificates for convex minimization and {H}elly
  numbers}.
\newblock {\em Oper. Res. Lett. 45}, (6) (2017), 671--674.

\bibitem{basu+oertel}
{\sc Basu, A., and Oertel, T.}
\newblock Centerpoints: {A} link between optimization and convex geometry.
\newblock In {\em Proceedings Integer Programming and Combinatorial
  Optimization, {(IPCO)}\/} (2016), pp.~14--25.

\bibitem{Bell:1977tm}
{\sc Bell, D.}
\newblock {A theorem concerning the integer lattice}.
\newblock {\em Studies in App. Math. 56}, 2 (1977), 187--188.

\bibitem{BeDu83}
{\sc Berge, C., and Duchet, P.}
\newblock Seminar, {MSH}i ({M}aison des {S}ciences de l'{H}omme), 1983.

\bibitem{bertsekasbook}
{\sc Bertsekas, D.~P.}
\newblock {\em Convex optimization algorithms}.
\newblock Athena Scientific, Nashua, NH, 2015.

\bibitem{Bertsimas:2004fp}
{\sc Bertsimas, D., and Vempala, S.}
\newblock {Solving convex programs by random walks}.
\newblock {\em J. ACM 51}, 4 (July 2004), 540--556.

\bibitem{Bezdek2003}
{\sc Bezdek, K., and Blokhuis, A.}
\newblock The {R}adon number of the three-dimensional integer lattice.
\newblock {\em Discrete Comput. Geom. 30}, 2 (2003), 181--184.

\bibitem{Birch:1959ii}
{\sc Birch, B.}
\newblock {On $3N$ points in a plane}.
\newblock {\em Math. Proc. Cambridge Philos. Soc. 55}, 04 (Oct. 1959),
  289--293.

\bibitem{Bjorner-survey}
{\sc Bj\"{o}rner, A.}
\newblock Topological methods.
\newblock In {\em Handbook of Combinatorics (Vol. 2)}. MIT Press, 1995,
  pp.~1819--1872.

\bibitem{BBZ16}
{\sc Blagojevi{\'{c}}, P. V.~M., Blagojevi{\'{c}}, A. S.~D., and Ziegler,
  G.~M.}
\newblock Polynomial partitioning for several sets of varieties.
\newblock {\em J. Fixed Point Theory Appl. 19}, 3 (2017), 1653--1660.

\bibitem{pfag2015}
{\sc Blagojevi{\'c}, P.~V.~M., Frick, F., Haase, A., and Ziegler, G.~M.}
\newblock Topology of the {G}r{\"u}nbaum-{H}adwiger-{R}amos hyperplane mass
  partition problem.
\newblock {\em Trans. Amer. Math. Soc. (to appear)\/} (2015).

\bibitem{pfag2016}
{\sc Blagojevi{\'c}, P.~V.~M., Frick, F., Haase, A., and Ziegler, G.~M.}
\newblock Hyperplane mass partitions via relative equivariant obstruction
  theory.
\newblock {\em Doc. Math. 21\/} (2016), 735--771.

\bibitem{Blagojevic:2014ul}
{\sc Blagojevi{\'c}, P. V.~M., Frick, F., and Ziegler, G.~M.}
\newblock {Tverberg plus constraints}.
\newblock {\em Bull. Lond. Math. Soc. 46}, (5) (2014), 953--967.

\bibitem{blafrickzie}
{\sc Blagojevi\'c, P.~V.~M., Frick, F., and Ziegler, G.~M.}
\newblock Barycenters of polytope skeleta and counterexamples to the
  topological {T}verberg conjecture, via constraints.
\newblock {\em J. Europ. Math. Soc. (to appear)\/} (2017).

\bibitem{blamaszie}
{\sc Blagojevi{\'c}, P.~V.~M., Matschke, B., and Ziegler, G.~M.}
\newblock Optimal bounds for a colorful {T}verberg-{V}re{\'c}ica type problem.
\newblock {\em Adv. Math. 226}, 6 (2011), 5198--5215.

\bibitem{Blagojevicetal2015}
{\sc Blagojevi{\'c}, P.~V.~M., Matschke, B., and Ziegler, G.~M.}
\newblock Optimal bounds for the colored {T}verberg problem.
\newblock {\em J. Eur. Math. Soc. (JEMS) 017}, 4 (2015), 739--754.

\bibitem{plave+pablo}
{\sc Blagojevi\'c, P. V.~M., and Sober\'on, P.}
\newblock Thieves can make sandwiches.
\newblock {\em Bull. Lond. Math. Soc. 50}, 1 (2018), 108--123.

\bibitem{Blagojevic+Zieglersurvey}
{\sc Blagojevi{\'{c}}, P. V.~M., and Ziegler, G.~M.}
\newblock {\em Beyond the Borsuk--Ulam Theorem: The Topological {T}verberg
  Story}.
\newblock Springer International Publishing, Cham, 2017, pp.~273--341.
\newblock chapter of ``A Journey Through Discrete Math.: A Tribute to
  Ji{\v{r}}{\'i} Matou{\v{s}}ek'', edited by M. Loebl, J. Ne{\v{s}}et{\v{r}}il,
  and R. Thomas.

\bibitem{BHPR16}
{\sc Blum, A., Har{-}Peled, S., and Raichel, B.}
\newblock Sparse approximation via generating point sets.
\newblock In {\em Proceedings {ACM-SIAM} Symposium on Discrete Algorithms,
  {(SODA)}\/} (2016), pp.~548--557.

\bibitem{BF82}
{\sc Boros, E., and F\"{u}redi, Z.}
\newblock The number of triangles covering the center of an $n$-set.
\newblock {\em Geom. Dedicata 17}, 1 (1984), 69--77.

\bibitem{BoGu96}
{\sc Boros, E., and Gurvich, V.}
\newblock Perfect graphs are kernel solvable.
\newblock {\em Discrete Math. 159\/} (1996), 35--55.

\bibitem{nerve}
{\sc Borsuk, K.}
\newblock On the imbedding of systems of compacta in simplicial complexes.
\newblock {\em Fund. Math. 35}, 1 (1948), 217--234.

\bibitem{boyd2004convex}
{\sc Boyd, S., and Vandenberghe, L.}
\newblock {\em Convex optimization}.
\newblock Cambridge university press, 2004.

\bibitem{brams+taylor}
{\sc Brams, S., and Taylor, A.}
\newblock {\em Fair Division: From cake-cutting to dispute resolution}.
\newblock Cambridge University Press, 1996.

\bibitem{bramsetal1}
{\sc Brams, S.~J., Jones, M.~A., and Klamler, C.}
\newblock {$N$}-person cake-cutting: there may be no perfect division.
\newblock {\em Amer. Math. Monthly 120}, 1 (2013), 35--47.

\bibitem{bramsetal2}
{\sc Brams, S.~J., Kilgour, D.~M., and Klamler, C.}
\newblock How to divide things fairly.
\newblock {\em Math. Mag. 88}, 5 (2015), 338--348.

\bibitem{Brazitikos2017}
{\sc Brazitikos, S.}
\newblock Quantitative {H}elly-type theorem for the diameter of convex sets.
\newblock {\em Discrete Comput. Geom. 57}, 2 (2017), 494--505.

\bibitem{BG95}
{\sc Br\"{o}nnimann, H., and Goodrich, M.}
\newblock Almost optimal set covers in finite {VC}-dimension.
\newblock {\em Discrete Comput. Geom. 14}, 4 (1995), 463--479.

\bibitem{B91}
{\sc Brualdi, R.~A., and Ryser, H.~J.}
\newblock {\em Combinatorial Matrix Theory}.
\newblock Cambridge University Press, 1991.

\bibitem{brunsetal}
{\sc Bruns, W., Gubeladze, J., Henk, M., Martin, A., and Weismantel, R.}
\newblock A counterexample to an integer analogue of {C}arath\'eodory's
  theorem.
\newblock {\em J. Reine Angew. Math. 510\/} (1999), 179--185.

\bibitem{B05}
{\sc Bukh, B.}
\newblock A point in many triangles.
\newblock {\em Electron. J. Combin. 13}, 1 (2006), Note 10, 3.

\bibitem{BMV11}
{\sc Bukh, B., Matou{\v{s}}ek, J., and Nivasch, G.}
\newblock Lower bounds for weak epsilon-nets and stair-convexity.
\newblock {\em Israel J. Math. 182}, 1 (2011), 199--228.

\bibitem{BMN08}
{\sc Bukh, B., Matou{\v{s}}ek, J.~r., and Nivasch, G.}
\newblock Stabbing simplices by points and flats.
\newblock {\em Discrete Comput. Geom. 43}, 2 (2010), 321--338.

\bibitem{calafiorecampi2005}
{\sc Calafiore, G., and Campi, M.~C.}
\newblock Uncertain convex programs: randomized solutions and confidence
  levels.
\newblock {\em Math. Program. 102}, 1, Ser. A (2005), 25--46.

\bibitem{calafiorecampi2006}
{\sc Calafiore, G.~C., and Campi, M.~C.}
\newblock The scenario approach to robust control design.
\newblock {\em IEEE Trans. Automat. Control 51}, 5 (2006), 742--753.

\bibitem{originalCaratheodory}
{\sc Carath{\'e}odory, C.}
\newblock \"{U}ber den {V}ariabilit\"atsbereich der {K}oeffizienten von
  {P}otenzreihen, die gegebene {W}erte nicht annehmen.
\newblock {\em Math. Ann. 64}, 1 (1907), 95--115.

\bibitem{carl1985inequalities}
{\sc Carl, B.}
\newblock Inequalities of {Bernstein-Jackson}-type and the degree of
  compactness of operators in {B}anach spaces.
\newblock {\em Ann. Inst. Fourier 35}, 3 (1985), 79--118.

\bibitem{C69}
{\sc Chakerian, G.~D.}
\newblock Some intersection properties of convex bodies.
\newblock {\em Proc. Amer. Math. Soc. 18}, 1 (1967), 109--112.

\bibitem{Chakerian-HellyBrouwer}
{\sc Chakerian, G.~D.}
\newblock Intersection and covering properties of convex sets.
\newblock {\em Amer. Math. Monthly 76}, 7 (1969), 753--766.

\bibitem{C04}
{\sc Chan, T.~M.}
\newblock An optimal randomized algorithm for maximum {T}ukey depth.
\newblock In {\em Proceedings {ACM-SIAM} Symposium on Discrete Algorithms
  ({(SODA)})\/} (2004), pp.~430--436.

\bibitem{C16}
{\sc Chan, T.~M.}
\newblock Improved deterministic algorithms for linear programming in low
  dimensions.
\newblock In {\em Proceedings {ACM-SIAM} Symposium on Discrete Algorithms,
  {SODA}\/} (2016), pp.~1213--1219.

\bibitem{CGKS12}
{\sc Chan, T.~M., Grant, E., K{\"o}nemann, J., and Sharpe, M.}
\newblock Weighted capacitated, priority, and geometric set cover via improved
  quasi-uniform sampling.
\newblock In {\em Proceedings {ACM-SIAM} Symposium on Discrete Algorithms
  ({(SODA)})\/} (2012), pp.~1576--1585.

\bibitem{C93}
{\sc Chazelle, B.}
\newblock Cutting hyperplanes for divide-and-conquer.
\newblock {\em Discrete Comput. Geom. 9\/} (1993), 145--158.

\bibitem{C00}
{\sc Chazelle, B.}
\newblock {\em {The Discrepancy Method: Randomness and Complexity}}.
\newblock Cambridge University Press, 2000.

\bibitem{CEGGSW93}
{\sc Chazelle, B., Edelsbrunner, H., Grigni, M., Guibas, L., Sharir, M., and
  Welzl, E.}
\newblock Improved bounds on weak $\epsilon$-nets for convex sets.
\newblock In {\em Proceedings {ACM} symposium on Theory of computing
  ({STOC})\/} (1993), pp.~495--504.

\bibitem{CDILM10}
{\sc Chen, D., Devillers, O., Iacono, J., Langerman, S., and Morin, P.}
\newblock {O}ja centers and centers of gravity.
\newblock {\em Comput. Geom. 46}, 2 (2013), 140--147.

\bibitem{chen2015multichromatic}
{\sc Chen, P.-A.}
\newblock On the multichromatic number of $s$-stable {K}neser graphs.
\newblock {\em J. Graph Theory 79}, 3 (2015), 233--248.

\bibitem{chen2017}
{\sc {Chen}, P.-A.}
\newblock {On the chromatic number of almost $s$-stable Kneser graphs}.
\newblock {\em \rm ArXiv:1711.06621\/} (2017).

\bibitem{chen2009complexity}
{\sc Chen, X., and Deng, X.}
\newblock On the complexity of 2-{D} discrete fixed point problem.
\newblock {\em Theoret. Comput. Sci. 410}, 44 (2009), 4448--4456.

\bibitem{chen+deng+teng2009}
{\sc Chen, X., Deng, X., and Teng, S.-H.}
\newblock Settling the complexity of computing two-player {N}ash equilibria.
\newblock {\em J. ACM 56}, 3 (2009), Art. 14, 57.

\bibitem{CMR16}
{\sc Cheong, O., Mulmuley, K., and Ramos, E.}
\newblock Randomization and derandomization.
\newblock In {\em Handbook of Discrete Comput. Geom.}, J.~E. Goodman,
  J.~O'Rourke, and C.~D. T\'oth, Eds. CRC Press LLC, 2016, to appear.

\bibitem{chestnutetal}
{\sc {Chestnut}, S.~R., {Hildebrand}, R., and {Zenklusen}, R.}
\newblock {Sublinear bounds for a quantitative {D}oignon-{B}ell-{S}carf
  Theorem}.
\newblock {\em SIAM J. Discrete Math. 32}, (1) (2018), 352--371.

\bibitem{Chung1989}
{\sc Chung, S.~J.}
\newblock {NP}-completeness of the linear complementarity problem.
\newblock {\em J. Optim. Theory Appl. 60}, 3 (1989), 393--399.

\bibitem{Chv73}
{\sc Chv\'atal, V.}
\newblock On the computational complexity of finding a kernel.
\newblock Tech. rep., Centre de Recherches Math\'ematiques, Universit\'e de
  Montr\'eal, 1973.

\bibitem{clarkson-rs}
{\sc Clarkson, K.~L.}
\newblock New applications of random sampling in computational geometry.
\newblock {\em Discrete Comput. Geom. 2}, 2 (1987), 195--222.

\bibitem{clarkson-lv}
{\sc Clarkson, K.~L.}
\newblock {L}as {V}egas algorithms for linear and integer programming when the
  dimension is small.
\newblock {\em J. ACM 42}, 2 (1995), 488--499.

\bibitem{CEGSW90}
{\sc Clarkson, K.~L., Edelsbrunner, H., Guibas, L., Sharir, M., and Welzl, E.}
\newblock Combinatorial complexity bounds for arrangement of curves and
  spheres.
\newblock {\em Discrete Comput. Geom. 5\/} (1990), 99--160.

\bibitem{CEMST96}
{\sc Clarkson, K.~L., Eppstein, D., Miller, G., Sturtivant, C., and Teng, S.}
\newblock Approximating center points with iterative {R}adon points.
\newblock {\em Internat. J. Comput. Geom. Appl. 6}, 3 (1996), 357--377.

\bibitem{clarkson-shor}
{\sc Clarkson, K.~L., and Shor, P.~W.}
\newblock Applications of random sampling in computational geometry, {II}.
\newblock {\em Discrete Comput. Geom. 4}, 5 (1989), 387--421.

\bibitem{multicake}
{\sc Cloutier, J., Nyman, K.~L., and Su, F.~E.}
\newblock Two-player envy-free multi-cake division.
\newblock {\em Math. Social Sci. 59}, 1 (2010), 26--37.

\bibitem{cook+fonlupt+schrijver}
{\sc Cook, W., Fonlupt, J., and Schrijver, A.}
\newblock An integer analogue of {C}arath\'eodory's theorem.
\newblock {\em J. Combin. Theory Ser. {B} 40}, 1 (1986), 63--70.

\bibitem{cottleetalbook}
{\sc Cottle, R., Pang, J., and Stone, R.}
\newblock {\em The Linear Complementarity Problem}.
\newblock Society for Industrial and Applied Mathematics, 2009.

\bibitem{collectionmathdecisionsgames}
{\sc Crisman, K.-D., and Jones, M.~A.}, Eds.
\newblock {\em The mathematics of decisions, elections, and games\/} (2014),
  vol.~624 of {\em Contemporary Mathematics}, American Mathematical Society,
  Providence, RI.

\bibitem{cunninghampivot}
{\sc Cunningham, W.~H.}
\newblock Theoretical properties of the network simplex method.
\newblock {\em Math. Oper. Res. 4}, 2 (1979), 196--208.

\bibitem{dantzig-minmax}
{\sc Dantzig, G.~B.}
\newblock Constructive proof of the {M}in-{M}ax theorem.
\newblock {\em Pacific J. Math. 6\/} (1956), 25--33.

\bibitem{DGKsurvey63}
{\sc Danzer, L., Gr{\"u}nbaum, B., and Klee, V.}
\newblock Helly's theorem and its relatives.
\newblock In {\em Proc. {S}ympos. {P}ure {M}ath., {V}ol. {VII}}. Amer. Math.
  Soc., Providence, R.I., 1963, pp.~101--180.

\bibitem{daskalakis+goldberg+papadimitriu2006}
{\sc Daskalakis, C., Goldberg, P.~W., and Papadimitriou, C.~H.}
\newblock The complexity of computing a {N}ash equilibrium.
\newblock In {\em S{TOC}'06: {P}roceedings of the 38th {A}nnual {ACM}
  {S}ymposium on {T}heory of {C}omputing}. ACM, New York, 2006, pp.~71--78.

\bibitem{Datta03}
{\sc Datta, R.~S.}
\newblock Universality of {N}ash equilibria.
\newblock {\em Math. Oper. Res. 28}, 3 (2003), 424--432.

\bibitem{MMMC}
{\sc de~Berg, M., Cheong, O., van Kreveld, M., and Overmars, M.}
\newblock {\em {Computational Geometry, Algorithms and Application}}.
\newblock Springer, 2008.
\newblock Third edition.

\bibitem{alggeo4optimization}
{\sc De~Loera, J.~A., Hemmecke, R., and K{\"{o}}ppe, M.}
\newblock {\em Algebraic and Geometric Ideas in the Theory of Discrete
  Optimization}, vol.~14 of {\em {MOS-SIAM} Series on Optimization}.
\newblock {SIAM}, 2013.

\bibitem{lowdimTverberg}
{\sc {De Loera}, J.~A., {Hogan}, T., {Meunier}, F., and {Mustafa}, N.}
\newblock Integer and mixed integer {T}verberg numbers.
\newblock {\em Proceedings European Congress of Computational Geometry\/}
  (2018).

\bibitem{DeLoeraetal-tv-2018}
{\sc {De Loera}, J.~A., {Hogan}, T.~A., {Oliveros}, D., and {Yang}, D.}
\newblock {Tverberg-Type Theorems with Trees and Cycles as (Nerve) Intersection
  Patterns}.
\newblock {\em \rm ArXiv:1808.00551\/} (2018).

\bibitem{DeLoeraetalquant}
{\sc De~Loera, J.~A., La~Haye, R., Rolnick, D., and Sober{\'o}n, P.}
\newblock Quantitative combinatorial geometry for continuous parameters.
\newblock {\em Discrete Comp. Geom. 57}, 2 (2017), 318--334.

\bibitem{DLetal-ccopt}
{\sc De~Loera, J.~A., La~Haye, R.~N., Oliveros, D., and Rold\'an-Pensado, E.}
\newblock Beyond chance-constrained convex mixed-integer optimization: Two
  sampling algorithms within {$S$}-optimization.
\newblock {\em J. Convex Anal. 25}, 1 (2015), 201--218.

\bibitem{deloeraetal2015Helly}
{\sc De~Loera, J.~A., La~Haye, R.~N., Oliveros, D., and Rold{\'a}n-Pensado, E.}
\newblock Helly numbers of algebraic subsets of $\mathbb{R}^d$.
\newblock {\em Advances in Geometry 17}, 4 (2017), 473--482.

\bibitem{integertverberg2017}
{\sc De~Loera, J.~A., La~Haye, R.~N., Rolnick, D., and Sober{\'o}n, P.}
\newblock Quantitative {T}verberg theorems over lattices and other discrete
  sets.
\newblock {\em Discrete Comput. Geom. 58}, (2) (2017), 435--448.

\bibitem{DeLoera:2002hj}
{\sc De~Loera, J.~A., Peterson, E., and Su, F.}
\newblock {A Polytopal Generalization of Sperner's Lemma}.
\newblock {\em J. Combin. Theory Ser. {A} 100}, 1 (Oct. 2002), 1--26.

\bibitem{triangbook}
{\sc De~Loera, J.~A., Rambau, J., and Santos, F.}
\newblock {\em Triangulations}, vol.~25 of {\em Series Algorithms and
  Computation in Mathematics}.
\newblock Springer-Verlag, Berlin, 2010.
\newblock Structures for algorithms and applications.

\bibitem{Longueville-book}
{\sc de~Longueville, M.}
\newblock {\em A Course in Topological Combinatorics}.
\newblock Universitext. Springer-Verlag New York, 2013.

\bibitem{mark+rade}
{\sc de~Longueville, M., and \v{Z}ivaljevi\'c, R.~T.}
\newblock Splitting multidimensional necklaces.
\newblock {\em Adv. Math. 218}, 3 (2008), 926--939.

\bibitem{acyclic}
{\sc De~Verdi{\`e}re, {\'E}.~C., Ginot, G., and Goaoc, X.}
\newblock Helly numbers of acyclic families.
\newblock {\em Adv. Math. 253\/} (2014), 163--193.

\bibitem{de1997combinatorics}
{\sc de~Werra, D.}
\newblock The combinatorics of timetabling.
\newblock {\em European J. Oper. Res. 96}, 3 (1997), 504--513.

\bibitem{dengetal2017}
{\sc Deng, X., Feng, Z., and Kulkarni, R.}
\newblock Octahedral tucker is {PPA}-complete.
\newblock In {\em Electronic Colloquium on Computational Complexity {(ECCC)}\/}
  (2017), vol.~24, p.~118.

\bibitem{DeQiSa12}
{\sc Deng, X., Qi, Q., and Saberi, A.}
\newblock Algorithmic solutions for envy-free cake cutting.
\newblock {\em Oper. Res. 60\/} (2012), 1461--1476.

\bibitem{Dezaetal2006}
{\sc Deza, A., Huang, S., Stephen, T., and Terlaky, T.}
\newblock Colourful simplicial depth.
\newblock {\em Discrete Comput. Geom. 35}, 4 (2006), 597--615.

\bibitem{diestel-graphtheory}
{\sc Diestel, R.}
\newblock {\em Graph Theory}, 3rd~ed., vol.~173 of {\em Grad. Texts in Math.}
\newblock Springer-Verlag Berlin Heidelberg, 2005.

\bibitem{dobbins2015point}
{\sc Dobbins, M.~G.}
\newblock A point in a $nd$-polytope is the barycenter of $n$ points in its
  $d$-faces.
\newblock {\em Invent. Math. 199}, 1 (2015), 287--292.

\bibitem{MR1430097}
{\sc Dodson, C. T.~J., and Parker, P.~E.}
\newblock {\em A user's guide to algebraic topology}, vol.~387 of {\em
  Mathematics and its Applications}.
\newblock Kluwer Academic Publishers Group, Dordrecht, 1997.

\bibitem{Doi1973}
{\sc Doignon, J.-P.}
\newblock {Convexity in cristallographical lattices}.
\newblock {\em J. Geom. 3}, 1 (Mar. 1973), 71--85.

\bibitem{Dold:1983wr}
{\sc Dold, A.}
\newblock {Simple proofs of some Borsuk-Ulam results}.
\newblock {\em Contemp. Math. 19\/} (1983), 65--69.

\bibitem{dol1981transversals}
{\sc Dol'nikov, V.}
\newblock Transversals of families of sets.
\newblock {\em Studies in the theory of functions of several real variables
  (Russian)\/} (1981), 30--36.

\bibitem{dubins+spanier}
{\sc Dubins, L.~E., and Spanier, E.~H.}
\newblock How to cut a cake fairly.
\newblock {\em Amer. Math. Monthly 68\/} (1961), 1--17.

\bibitem{Duc80}
{\sc Duchet, P.}
\newblock Graphes noyau-parfaits.
\newblock {\em Ann. Discrete Math. 9\/} (1980), 93--101.

\bibitem{Duchetsurvey1987}
{\sc Duchet, P.}
\newblock Convexity in combinatorial structures.
\newblock In {\em Proceedings of the 14th Winter School on Abstract Analysis\/}
  (1987), Circolo Matematico di Palermo, pp.~261--293.

\bibitem{DuMe83}
{\sc Duchet, P., and Meyniel, H.}
\newblock Une g\'en\'eralisation du th\'eor\`eme de {R}ichardson sur
  l'existence de noyaux dans les graphes orient\'es.
\newblock {\em Discrete Math. 43\/} (1983), 21--27.

\bibitem{durandetal}
{\sc Durand, A., Hermann, M., and Juban, L.}
\newblock On the complexity of recognizing the {H}ilbert basis of a linear
  diophantine system.
\newblock {\em Theoret. Comput. Sci. 270}, 1 (2002), 625 -- 642.

\bibitem{eckhoff1985upper}
{\sc Eckhoff, J.}
\newblock An upper-bound theorem for families of convex sets.
\newblock {\em Geom. Dedicata 19}, 2 (1985), 217--227.

\bibitem{Eckhoff:1993survey}
{\sc Eckhoff, J.}
\newblock {Helly, {R}adon, and {C}arath\'eodory type theorems}.
\newblock In {\em Handbook of convex geometry, Vol.\ A, B}. North-Holland,
  Amsterdam, 1993, pp.~389--448.

\bibitem{Eckhoff:2000jw}
{\sc Eckhoff, J.}
\newblock {The partition conjecture}.
\newblock {\em Discrete Math. 221}, 1-3 (2000), 61--78.

\bibitem{edmondsgiles}
{\sc Edmonds, J., and Giles, R.}
\newblock A min-max relation for submodular functions on graphs.
\newblock In {\em Studies in Integer Programming}, B.~K. P.L.~Hammer,
  E.L.~Johnson and G.~Nemhauser, Eds., vol.~1 of {\em Annals of Discrete Math.}
  Elsevier, 1977, pp.~185 -- 204.

\bibitem{eisenbrandshmonin-caratheodory}
{\sc Eisenbrand, F., and Shmonin, G.}
\newblock Carath\'eodory bounds for integer cones.
\newblock {\em Oper. Res. Lett. 34}, 5 (2006), 564--568.

\bibitem{E46}
{\sc Erd\H{o}s, P.}
\newblock On sets of distances of $n$ points.
\newblock {\em Amer. Math. Monthly 53}, 5 (1946), 248--250.

\bibitem{E96}
{\sc Erickson, J.}
\newblock New lower bounds for {H}opcroft's problem.
\newblock {\em Discrete Comput. Geom. 16}, 4 (1996), 389--418.

\bibitem{etessami+yannakakis}
{\sc Etessami, K., and Yannakakis, M.}
\newblock On the complexity of {N}ash equilibria and other fixed points.
\newblock {\em SIAM J. Comput. 39}, 6 (2010), 2531--2597.

\bibitem{ERS05}
{\sc Even, G., Rawitz, D., and Shahar, S.}
\newblock Hitting sets when the {VC}-dimension is small.
\newblock {\em Inform. Process. Lett. 95\/} (2005), 358--362.

\bibitem{kyfan1952}
{\sc Fan, K.}
\newblock A generalization of {T}ucker's combinatorial lemma with topological
  applications.
\newblock {\em Ann. of Math. (2) 56\/} (1952), 431--437.

\bibitem{Fan68}
{\sc Fan, K.}
\newblock Simplicial maps from an orientable $n$-pseudomanifold into ${S}^m$
  with the octahedral triangulation.
\newblock {\em J. Combin. Theory 2}, 4 (1967), 588 -- 602.

\bibitem{Farb}
{\sc Farb, B.}
\newblock Group actions and {H}elly's theorem.
\newblock {\em Adv. Math. 222}, 5 (2009), 1574 -- 1588.

\bibitem{Filos-Ratsikas:2018}
{\sc Filos-Ratsikas, A., and Goldberg, P.~W.}
\newblock Consensus halving is {PPA}-complete.
\newblock In {\em Proceedings {ACM} Symposium on Theory of Computing,
  {(STOC)}\/} (2018), pp.~51--64.

\bibitem{2018Filos-Ratsikas}
{\sc {Filos-Ratsikas}, A., and {Goldberg}, P.~W.}
\newblock {The Complexity of Splitting Necklaces and Bisecting Ham Sandwiches}.
\newblock {\em \rm ArXiv:1805.12559\/} (2018).

\bibitem{Flores:NichtEinbettbar-1933}
{\sc Flores, A.~I.}
\newblock {\"Uber die {E}xistenz $n$-dimensionaler {K}omplexe, die nicht in den
  $\R^{2n}$ topo\-logisch einbettbar sind.}
\newblock {\em Ergeb. Math. Kolloqu. 5\/} (1933), 17--24.

\bibitem{FGLNP11}
{\sc Fox, J., Gromov, M., Lafforgue, V., Naor, A., and Pach, J.}
\newblock Overlap properties of geometric expanders.
\newblock In {\em Proceedings {ACM-SIAM} Symposium on Discrete Algorithms
  ({SODA})\/} (2011), pp.~1188--1197.

\bibitem{freundtodd}
{\sc Freund, R.~M., and Todd, M.~J.}
\newblock A constructive proof of {T}ucker's combinatorial lemma.
\newblock {\em J. Combin. Theory Ser. {A} 30}, 3 (1981), 321 -- 325.

\bibitem{frick15}
{\sc Frick, F.}
\newblock {Counterexamples to the topological {T}verberg conjecture}.
\newblock {\em Oberwolfach Reports 12}, (1) (2015), 318--321.

\bibitem{FrickNerves}
{\sc Frick, F.}
\newblock Intersection patterns of finite sets and of convex sets.
\newblock {\em Proc. Amer. Math. Soc. 145}, 7 (2017), 2827--2842.

\bibitem{FrickKneser}
{\sc Frick, F.}
\newblock Chromatic numbers of stable {K}neser hypergraphs via topological
  tverberg-type theorems.
\newblock {\em Int. Math. Res. Not. IMRN\/} (2018), rny135.

\bibitem{frickhoustonmeunier}
{\sc {Frick}, F., {Houston-Edwards}, K., and {Meunier}, F.}
\newblock {Achieving rental harmony with a secretive roommate}.
\newblock {\em Amer. Math. Monthly, to appear\/} (2017).

\bibitem{frick+zerbib}
{\sc Frick, F., and Zerbib, S.}
\newblock Colorful coverings of polytopes and piercing numbers of colorful
  $d$-intervals.
\newblock {\em Combinatorica (to appear) 21\/} (2018).

\bibitem{friedmann}
{\sc Friedmann, O.}
\newblock A subexponential lower bound for {Z}adeh's pivoting rule for solving
  linear programs and games.
\newblock In {\em Proceedings Integer Programming and Combinatorial
  Optimization, {(IPCO)}\/} (2011), pp.~192--206.

\bibitem{Fulek:2009vz}
{\sc Fulek, R., Holmsen, A.~F., and Pach, J.}
\newblock {Intersecting convex sets by rays}.
\newblock {\em Discrete Comput Geom 42}, 3 (2009), 343--358.

\bibitem{fulek2010}
{\sc Fulek, R., and Pach, J.}
\newblock A computational approach to {Conway's} thrackle conjecture.
\newblock {\em Comput. Geom. 44}, 6 (2011), 345 -- 355.

\bibitem{gale+hex}
{\sc Gale, D.}
\newblock The game of {H}ex and the {B}rouwer fixed-point theorem.
\newblock {\em Amer. Math. Monthly 86}, 10 (1979), 818--827.

\bibitem{Gale-colorfulkkm}
{\sc Gale, D.}
\newblock Equilibrium in a discrete exchange economy with money.
\newblock {\em Internat. J. Game Theory 13}, 1 (1984), 61--64.

\bibitem{GaNe84}
{\sc Galeana-S\'anchez, H., and Neumann-Lara, V.}
\newblock On kernels and semikernels of digraphs.
\newblock {\em Discrete Math. 48\/} (1984), 67--76.

\bibitem{garber}
{\sc {Garber}, A.}
\newblock {On Helly number for crystals and cut-and-project sets}.
\newblock {\em \rm ArXiv:1605.07881\/} (2016).

\bibitem{nataliaetal2017}
{\sc Garc{\'i}a-Col{\'i}n, N., Raggi, M., and Rold{\'a}n-Pensado, E.}
\newblock A note on the tolerant {T}verberg theorem.
\newblock {\em Discrete Comput. Geom. 58}, 3 (Oct 2017), 746--754.

\bibitem{gargetal2015}
{\sc Garg, J., Mehta, R., Vazirani, V.~V., and Yazdanbod, S.}
\newblock {ETR}-completeness for decision versions of multi-player (symmetric)
  {N}ash equilibria.
\newblock In {\em Proceedings International Colloquium Automata, Languages, and
  Programming, {(ICALP)}\/} (2015), pp.~554--566.

\bibitem{ViolatorSpaces2008}
{\sc G{\"a}rtner, B., Matou\v{s}ek, J., R{\"u}st, L., and {\v{S}}kovro{\v{n}},
  P.}
\newblock Violator spaces: Structure and algorithms.
\newblock {\em Discrete Appl. Math. 156}, 11 (2008), 2124--2141.

\bibitem{GKM14}
{\sc Giannopoulos, P., Konzack, M., and Mulzer, W.}
\newblock Low-crossing spanning trees: an alternative proof and experiments.
\newblock In {\em Proceedings European Workshop on Computational Geometry\/}
  (2014).

\bibitem{gijswijt+regts}
{\sc Gijswijt, D., and Regts, G.}
\newblock Polyhedra with the integer {C}arath\'eodory property.
\newblock {\em J. Combin. Theory Ser. B 102}, 1 (2012), 62--70.

\bibitem{GSW92}
{\sc Gil, J., Steiger, W.~L., and Wigderson, A.}
\newblock Geometric medians.
\newblock {\em Discrete Math. 108}, 3 (1992), 37--51.

\bibitem{gilboa+zemel}
{\sc Gilboa, I., and Zemel, E.}
\newblock Nash and correlated equilibria: Some complexity considerations.
\newblock {\em Games and Econom. Behav. 1}, 1 (1989), 80--93.

\bibitem{gilespulleyblanktdi}
{\sc Giles, F., and Pulleyblank, W.}
\newblock Total dual integrality and integer polyhedra.
\newblock {\em Linear Algebra and its Applications 25\/} (1979), 191 -- 196.

\bibitem{BH-journal}
{\sc Goaoc, X., Pat\'ak, P., Pat\'agov\'a, Z., Tancer, M., and Wagner, U.}
\newblock {\em Bounding {H}elly Numbers via {B}etti Numbers}.
\newblock Springer International Publishing, Cham, 2017, pp.~407--447.
\newblock chapter of ``A Journey Through Discrete Math.: A Tribute to
  Ji{\v{r}}{\'i} Matou{\v{s}}ek'', edited by M. Loebl, J. Ne{\v{s}}et{\v{r}}il,
  and R. Thomas.

\bibitem{goemans+rothvoss}
{\sc Goemans, M.~X., and Rothvo{\ss}, T.}
\newblock Polynomiality for bin packing with a constant number of item types.
\newblock In {\em Proceedings {ACM-SIAM} Symposium on Discrete Algorithms,
  {(SODA)}\/} (2014), pp.~830--839.

\bibitem{GW85}
{\sc Goldberg, C., and West, D.}
\newblock Bisection of circle colorings.
\newblock {\em SIAM J. Algeb. Discrete Meth. 6\/} (1985), 93--106.

\bibitem{GoldbergPapa17}
{\sc Goldberg, P.~W., and Papadimitriou, C.~H.}
\newblock {TFNP:} an update.
\newblock In {\em Proceedings International Conference in Algorithms and
  Complexity, {(CIAC)}\/} (2017), pp.~3--9.

\bibitem{Graver1975}
{\sc Graver, J.~E.}
\newblock {On the foundations of linear and integer linear programming {I}}.
\newblock {\em Math. Program. 9}, 1 (1975), 207--226.

\bibitem{G10}
{\sc Gromov, M.}
\newblock Singularities, expanders and topology of maps. {P}art 2: From
  combinatorics to topology via algebraic isoperimetry.
\newblock {\em Geom. Func. Anal. 20\/} (2010), 416--526.

\bibitem{handbookofconvexgeo}
{\sc Gruber, P., and Wills, J.}, Eds.
\newblock {\em Handbook of Convex Geometry}.
\newblock North-Holland, Amsterdam, 1993.

\bibitem{gruber2007convex}
{\sc Gruber, P.~M.}
\newblock {\em Convex and discrete geometry}.
\newblock Springer-Verlag, Berlin, 2007.

\bibitem{G15}
{\sc Guth, L.}
\newblock Polynomial partitioning for a set of varieties.
\newblock {\em Math. Proc. Cambridge Philos. Soc. 159\/} (2015), 459--469.

\bibitem{GuthBook}
{\sc Guth, L.}
\newblock {\em {Polynomial methods in combinatorics}}.
\newblock American Mathematical Society, 2016.

\bibitem{GK15}
{\sc Guth, L., and Katz, N.}
\newblock On the {E}rd{\H o}s distinct distances problem in the plane.
\newblock {\em Ann. of Math. 181}, 1 (2015), 155--190.

\bibitem{gyori1976division}
{\sc Gy\H{o}ri, E.}
\newblock On division of graphs to connected subgraphs.
\newblock In {\em Combinatorics (Proc. Fifth Hungarian Colloq., Keszthely,
  1976)\/} (1976), vol.~1, pp.~485--494.

\bibitem{hahn1997graph}
{\sc Hahn, G., and Tardif, C.}
\newblock Graph homomorphisms: structure and symmetry.
\newblock In {\em Graph symmetry}. Springer, 1997, pp.~107--166.

\bibitem{HaSaScZi09}
{\sc Hanke, B., Sanyal, R., Schultz, C., and Ziegler, G.~M.}
\newblock Combinatorial {S}tokes formulas via minimal resolutions.
\newblock {\em J. Combin. Theory Ser. {A} 116\/} (2009), 404--420.

\bibitem{HP09}
{\sc {Har-Peled}, S.}
\newblock {Approximating Spanning Trees with Low Crossing Number}.
\newblock {\em \rm ArXiv:0907.1131\/} (2009).

\bibitem{Haussler:1987fr}
{\sc Haussler, D., and Welzl, E.}
\newblock {$\epsilon$-nets and simplex range queries}.
\newblock {\em Discrete Comput Geom 2}, 1 (Dec. 1987), 127--151.

\bibitem{Ha95}
{\sc Haxell, P.~E.}
\newblock A condition for matchability of hypergraphs.
\newblock {\em Graphs Combin. 11\/} (1995), 245--248.

\bibitem{HT06}
{\sc Hell, S.}
\newblock {\em Tverberg-type theorems, and the Fractional {H}elly property}.
\newblock PhD thesis, T.U. Berlin, 2006.

\bibitem{hell1}
{\sc Hell, S.}
\newblock On the number of {T}verberg partitions in the prime power case.
\newblock {\em European J. Combin. 28}, 1 (2007), 347--355.

\bibitem{Hell2008}
{\sc Hell, S.}
\newblock On the number of {B}irch partitions.
\newblock {\em Discrete Comput. Geom. 40}, 4 (2008), 586.

\bibitem{Hoffman:1979ix}
{\sc Hoffman, A.~J.}
\newblock {Binding constraints and {H}elly numbers}.
\newblock {\em Ann. New York Acad. Sci. 319}, 1 (1979), 284--288.

\bibitem{Holmsen-colorful}
{\sc Holmsen, A.}
\newblock The intersection of a matroid and an oriented matroid.
\newblock {\em Adv. Math. 290\/} (2016), 1 -- 14.

\bibitem{HPT08}
{\sc Holmsen, A., Pach, J., and Tverberg, H.}
\newblock Points surrounding the origin.
\newblock {\em Combinatorica 28}, 6 (2008), 633--644.

\bibitem{holmsen+wenger}
{\sc Holmsen, A., and Wenger, R.}
\newblock Helly-type theorems and geometric transversals.
\newblock In {\em Handbook of discrete and computational geometry, third
  edition}, CRC Press Ser. Discrete Math. Appl. CRC, Boca Raton, FL, 2016.

\bibitem{Holmsen+Karasev}
{\sc Holmsen, A.~F., and Karasev, R.}
\newblock Colorful theorems for strong convexity.
\newblock {\em Proc. Amer. Math. Soc. 145}, 6 (2017), 2713--2726.

\bibitem{NervesMinors}
{\sc {Holmsen}, A.~F., {Kim}, M., and {Lee}, S.}
\newblock {Nerves, minors, and piercing numbers}.
\newblock {\em \rm ArXiv:1706.05181\/} (2017).

\bibitem{PBS-sperner}
{\sc Houston-Edwards, K.}
\newblock Splitting rent with triangles.
\newblock {\em Infinite series, PBS digital studios February 23\/} (2017).
\newblock \url{https://www.youtube.com/watch?v=48oBEvpdYSE}.

\bibitem{JM93}
{\sc Jadhav, S., and Mukhopadhyay, A.}
\newblock Computing a centerpoint of a finite planar set of points in linear
  time.
\newblock {\em Discrete Comput. Geom. 12\/} (1994), 291--312.

\bibitem{jonsson2012chromatic}
{\sc Jonsson, J.}
\newblock On the chromatic number of generalized stable {K}neser graphs, 2012.
\newblock Retrieved from
  \url{http://www.math.kth.se/~jakobj/doc/submitted/stablekneser.pdf}.

\bibitem{Kalai:1984bg}
{\sc Kalai, G.}
\newblock {Intersection patterns of convex sets}.
\newblock {\em Israel J. Math. 48}, 2-3 (June 1984), 161--174.

\bibitem{kalai1992subexponential}
{\sc Kalai, G.}
\newblock A subexponential randomized simplex algorithm.
\newblock In {\em Proceedings {ACM} Symposium on Theory of computing
  {(STOC)}\/} (1992), pp.~475--482.

\bibitem{kalai1995survey}
{\sc Kalai, G.}
\newblock Combinatorics and convexity.
\newblock In {\em Proceedings International Congress of Mathematicians: August
  3--11, 1994 Z{\"u}rich, Switzerland\/} (Basel, 1995), S.~D. Chatterji, Ed.,
  Birkh{\"a}user Basel, pp.~1363--1374.

\bibitem{Kalai01algebraicshifting}
{\sc Kalai, G.}
\newblock Algebraic shifting.
\newblock In {\em Computational Commutative Algebra and Combinatorics\/}
  (2002), T.~Hibi, Ed., vol.~33 of {\em Advanced Studies in Pure Mathematics},
  Mathematical Soc. of Japan, pp.~121--163.

\bibitem{kalai_conjectures}
{\sc Kalai, G.}
\newblock Combinatorial expectations from commutative algebra.
\newblock In {\em Combinatorial Commutative Algebra}, I.~Peeva and V.~Welker,
  Eds., vol.~1(3). Oberwolfach Reports, 2004, pp.~1729--1734.

\bibitem{Kalai:2005cm}
{\sc Kalai, G., and Meshulam, R.}
\newblock A topological colorful {H}elly theorem.
\newblock {\em Adv. Math. 191}, 2 (2005), 305 -- 311.

\bibitem{kannan+theobald}
{\sc Kannan, R., and Theobald, T.}
\newblock Games of fixed rank: a hierarchy of bimatrix games.
\newblock In {\em Proceedings {ACM-SIAM} Symposium on Discrete Algorithms,
  {SODA}\/} (2007), pp.~1124--1132.

\bibitem{K08}
{\sc Karasev, R.~N.}
\newblock Dual theorems on central points and their generalizations.
\newblock {\em Sb. Math. 199}, 10 (2008), 1459--1479.

\bibitem{Karasev-survey}
{\sc Karasev, R.~N.}
\newblock Topological methods in combinatorial geometry.
\newblock {\em Russian Math. Surveys 63\/} (2008), 1031--1078.

\bibitem{Karasev:2011jv}
{\sc Karasev, R.~N.}
\newblock {Tverberg-Type Theorems for Intersecting by Rays}.
\newblock {\em Discrete Comput. Geom. 45}, 2 (Mar. 2011), 340--347.

\bibitem{K12}
{\sc Karasev, R.~N.}
\newblock A simpler proof of the
  {B}oros-{F}\"uredi-{B}\'ar\'any-{P}ach-{G}romov theorem.
\newblock {\em Discrete Comput. Geom. 47}, 3 (2012), 492--495.

\bibitem{Karmarkar1984}
{\sc Karmarkar, N.}
\newblock A new polynomial-time algorithm for linear programming.
\newblock {\em Combinatorica 4}, 4 (1984), 373--395.

\bibitem{Kat79frac}
{\sc Katchalski, M., and Liu, A.}
\newblock {A problem of geometry in $\R^n$}.
\newblock {\em Proc. Amer. Math. Soc. 75}, 2 (1979), 284--288.

\bibitem{Kay:1971uf}
{\sc Kay, D., and Womble, E.}
\newblock {Axiomatic convexity theory and relationships between the
  Carath{\'e}odory, {H}elly, and {R}adon numbers}.
\newblock {\em Pacific J. Math. 38}, 2 (1971), 471--485.

\bibitem{khachiyan1980polynomial}
{\sc Khachiyan, L.~G.}
\newblock Polynomial algorithms in linear programming.
\newblock {\em USSR Computational Mathematics and Mathematical Physics 20}, 1
  (1980), 53--72.

\bibitem{KiPa09}
{\sc Kir\'aly, T., and Pap, J.}
\newblock A note on kernels and {S}perner's lemma.
\newblock {\em Discrete Appl. Math. 157\/} (2009), 3327--3331.

\bibitem{klee+minty}
{\sc Klee, V., and Minty, G.}
\newblock How good is the simplex algorithm?
\newblock In {\em Inequalities III, Proc. Third Sympos., Univ. California, Los
  Angeles, CA 1969; dedicated to the memory of T.S. Motzkin\/} (1972), Academic
  Press, pp.~159--175.

\bibitem{kozlov2007combinatorial}
{\sc Kozlov, D.}
\newblock {\em Combinatorial algebraic topology}, vol.~21 of {\em Algorithms
  and Computation in Mathematics}.
\newblock Springer Science \& Business Media, 2007.

\bibitem{KMS12}
{\sc Kr{\'{a}}l, D., Mach, L., and Sereni, J.}
\newblock A new lower bound based on {G}romov's method of selecting heavily
  covered points.
\newblock {\em Discrete Comput. Geom. 48}, 2 (2012), 487--498.

\bibitem{krano}
{\sc {Krasnosel{\cprime}ski{\u\i}, M. A.}}
\newblock On a proof of {H}elly's theorem on sets of convex bodies with common
  points.
\newblock {\em Vorone\v z. Gos. Univ. Trudy. Fiz.-Mat. Sb. 33\/} (1954),
  19--20.

\bibitem{KMV02}
{\sc Krishnan, S., Mustafa, N.~H., and Venkatasubramanian, S.}
\newblock Hardware-assisted computation of depth contours.
\newblock In {\em Proceedings {ACM-SIAM} Symposium on Discrete Algorithms
  ({(SODA)})\/} (2002), pp.~558--567.

\bibitem{KMV06}
{\sc Krishnan, S., Mustafa, N.~H., and Venkatasubramanian, S.}
\newblock Statistical data depth and the graphics hardware.
\newblock In {\em Data Depth: Robust Multivariate Analysis, Computational
  Geometry and Applications}, D.~Souvaine, Ed. DIMACS Series in Discrete Math.
  and Theoret. Comput. Sci., 2006, pp.~223--250.

\bibitem{KMP16}
{\sc Kupavskii, A., Mustafa, N.~H., and Pach, J.}
\newblock New lower bounds for epsilon-nets.
\newblock In {\em Proceedings Symposium on Computational Geometry, {(SoCG)}\/}
  (2016), pp.~54:1--54:16.

\bibitem{LS03}
{\sc Langerman, S., and Steiger, W.~L.}
\newblock The complexity of hyperplane depth in the plane.
\newblock {\em Discrete Comput. Geom. 30}, 2 (2003), 299--309.

\bibitem{lebertetal}
{\sc Lebert, N., Meunier, F., and Carbonneaux, Q.}
\newblock Envy-free two-player \emph{m}-cake and three-player two-cake
  divisions.
\newblock {\em Oper. Res. Lett. 41}, 6 (2013), 607--610.

\bibitem{lemkehowson64}
{\sc Lemke, C.~E., and Howson, Jr., J.~T.}
\newblock Equilibrium points of bimatrix games.
\newblock {\em J. Soc. Indust. Appl. Math. 12\/} (1964), 413--423.

\bibitem{liptonetal}
{\sc Lipton, R.~J., Markakis, E., and Mehta, A.}
\newblock Playing large games using simple strategies.
\newblock In {\em Proceedings {ACM} Conference on Electronic Commerce,
  {(EC)}\/} (2003), pp.~36--41.

\bibitem{L88}
{\sc Liu, R.}
\newblock On a notion of simplicial depth.
\newblock {\em Proc. Natl. Acad. Sci. USA 85\/} (1988), 1732--1734.

\bibitem{LMS94}
{\sc Lo, C., Matou{\v{s}}ek, J., and Steiger, W.}
\newblock Algorithms for ham-sandwich cuts.
\newblock {\em Discrete Comput. Geom. 11}, 4 (1994), 433--452.

\bibitem{lovasz1972}
{\sc Lov{\'a}sz, L.}
\newblock Normal hypergraphs and the perfect graph conjecture.
\newblock {\em Discrete Math.}, 2 (1972), 253--267.

\bibitem{Lovasz:1978wu}
{\sc Lov{\'a}sz, L.}
\newblock {Kneser's conjecture, chromatic number, and homotopy}.
\newblock {\em J. Combin. Theory Ser. A 25}, 3 (1978), 319--324.

\bibitem{lovasz2006graph}
{\sc Lov{\'a}sz, L.}
\newblock Graph minor theory.
\newblock {\em Bull. Amer. Math. Soc. 43}, 1 (2006), 75--86.

\bibitem{lovasz+plummer2009matching}
{\sc Lov{\'a}sz, L., and Plummer, M.~D.}
\newblock {\em Matching theory}.
\newblock AMS Chelsea Publishing, Providence, RI, 2009.

\bibitem{LStheorem}
{\sc Lyusternik, L., and Schnirel'mann, L.}
\newblock M\'ethodes topologiques dans les probl\`emes variationnels.
\newblock {\em Gosudarstvennoe Izdat.\/} (1930).
\newblock French transl., Actualit\'es Sci. Indust., no. 118, Hermann, Paris,
  1934.

\bibitem{originalkreinmilman}
{\sc M.~Krein, D.~M.}
\newblock On extreme points of regular convex sets.
\newblock {\em Studia Math. 9\/} (1940), 133--138.

\bibitem{Mabillard2014}
{\sc Mabillard, I., and Wagner, U.}
\newblock {Eliminating {T}verberg Points, I. An Analogue of the {W}hitney
  Trick}.
\newblock In {\em Proceedings Symposium on Computational Geometry (SoCG)\/}
  (2014), pp.~171--180.

\bibitem{Martinez-Sandovaletal18}
{\sc Mart{\'{\i}}nez{-}Sandoval, L., Rold{\'{a}}n{-}Pensado, E., and Rubin, N.}
\newblock Further consequences of the colorful {H}elly hypothesis.
\newblock In {\em Proceedings Symposium on Computational Geometry, {(SoCG)}\/}
  (2018), pp.~59:1--59:14.

\bibitem{leo+tamam}
{\sc {Mart{\'{\i}}nez-Sandoval}, L., and {Tamam}, R.}
\newblock {Depth with respect to a family of convex sets}.
\newblock {\em \rm ArXiv:1612.03435\/} (2016).

\bibitem{M92}
{\sc Matou{\v{s}}ek, J.}
\newblock Efficient partition trees.
\newblock {\em Discrete Comput. Geom. 8\/} (1992), 315--334.

\bibitem{M99}
{\sc Matou{\v{s}}ek, J.}
\newblock {\em {Geometric Discrepancy: An Illustrated Guide}}.
\newblock Springer, 1999.

\bibitem{Mbook}
{\sc Matou{\v{s}}ek, J.}
\newblock {\em Lectures on Discrete Geometry}, vol.~212 of {\em Grad. Texts in
  Math.}
\newblock Springer-Verlag, New York, 2002.

\bibitem{matousek2003using}
{\sc Matou{\v{s}}ek, J.}
\newblock {\em {Using the {Borsuk-Ulam} theorem: lectures on topological
  methods in combinatorics and geometry}}.
\newblock Universitext. Springer-Verlag, Berlin, 2003.
\newblock Written in cooperation with A. Bj\"orner and G. M. Ziegler.

\bibitem{Matousek:2004hm}
{\sc Matou{\v{s}}ek, J.}
\newblock {A Combinatorial Proof of {K}neser{\textquoteright}s Conjecture}.
\newblock {\em Combinatorica 24}, 1 (Jan. 2004), 163--170.

\bibitem{Matousek:2004cs}
{\sc Matou{\v{s}}ek, J.}
\newblock {Bounded VC-Dimension Implies a Fractional {H}elly Theorem}.
\newblock {\em Discrete Comput. Geom. 31}, 2 (Feb. 2004), 251--255.

\bibitem{matouvsek2010thirty}
{\sc Matou{\v{s}}ek, J.}
\newblock {\em Thirty-three miniatures: Mathematical and Algorithmic
  applications of Linear Algebra}, vol.~53 of {\em Student Mathematical
  Library}.
\newblock Amer. Math. Soc., 2010.

\bibitem{matousek2007understanding}
{\sc Matou{\v{s}}ek, J., and G{\"a}rtner, B.}
\newblock {\em Understanding and using linear programming}.
\newblock Springer Science \& Business Media, 2007.

\bibitem{Matousek:2002wl}
{\sc Matou{\v{s}}ek, J., and Ziegler, G.~M.}
\newblock {Topological lower bounds for the chromatic number: A hierarchy}.
\newblock {\em Jahresber. Dtsch. Math.-Ver. 106\/} (Aug. 2004), 71--90.

\bibitem{msw-sblp-92}
{\sc Matou\v{s}ek, J., Sharir, M., and Welzl, E.}
\newblock A subexponential bound for linear programming.
\newblock {\em Algorithmica 16}, 4-5 (1996), 498--516.

\bibitem{MW02}
{\sc Matou\v{s}ek, J., and Wagner, U.}
\newblock New constructions of weak epsilon-nets.
\newblock {\em Discrete Comput. Geom. 32}, 2 (2004), 195--206.

\bibitem{MW14}
{\sc Matou\v{s}ek, J., and Wagner, U.}
\newblock On {G}romov's method of selecting heavily covered points.
\newblock {\em Discrete Comput. Geom. 52}, 1 (2014), 1--33.

\bibitem{mckelvey1996computation}
{\sc McKelvey, R.~D., and McLennan, A.}
\newblock Computation of equilibria in finite games.
\newblock {\em Handbook of computational economics 1\/} (1996), 87--142.

\bibitem{mclennan2005}
{\sc McLennan, A.}
\newblock The expected number of {N}ash equilibria of a normal form game.
\newblock {\em Econometrica 73}, 1 (2005), 141--174.

\bibitem{McLennan+Tourky2008}
{\sc McLennan, A., and Tourky, R.}
\newblock Using volume to prove {S}perner's lemma.
\newblock {\em Econom. Theory 35}, 3 (2008), 593--597.

\bibitem{Mes01}
{\sc Meshulam, R.}
\newblock The clique complex and hypergraph matching.
\newblock {\em Combinatorica 21\/} (2001), 89--94.

\bibitem{meshulam2003domination}
{\sc Meshulam, R.}
\newblock Domination numbers and homology.
\newblock {\em J. Combin. Theory Ser. A 102}, 2 (2003), 321--330.

\bibitem{Me05}
{\sc Meunier, F.}
\newblock A $\mathbb{Z}_q$-{F}an theorem.
\newblock Tech. rep., Laboratoire Leibniz-IMAG, Grenoble, 2005.

\bibitem{meunier-sperner}
{\sc Meunier, F.}
\newblock Sperner labellings: a combinatorial approach.
\newblock {\em J. Combin. Theory Ser. {A} 113}, 7 (2006), 1462--1475.

\bibitem{meunier2011chromatic}
{\sc Meunier, F.}
\newblock The chromatic number of almost stable {K}neser hypergraphs.
\newblock {\em J. Combin. Theory Ser. {A} 118}, 6 (2011), 1820--1828.

\bibitem{deza+meunier:cara}
{\sc Meunier, F., and Deza, A.}
\newblock A further generalization of the colourful {C}arath{\'e}odory theorem.
\newblock In {\em Discrete Geometry and Optimization}, K.~Bezdek, A.~Deza, and
  Y.~Ye, Eds. Springer International Publishing, Heidelberg, 2013,
  pp.~179--190.

\bibitem{PPAD-ColorfulCaratheodory}
{\sc Meunier, F., Mulzer, W., Sarrabezolles, P., and Stein, Y.}
\newblock The rainbow at the end of the line---a {{PPAD}} formulation of the
  colorful {C}arath\'eodory theorem with applications.
\newblock In {\em Proceedings {ACM}-{SIAM} {S}ymposium on {D}iscrete
  {A}lgorithms, {(SODA)}\/} (2017), pp.~1342--1351.

\bibitem{meunier-sarrabezolles2}
{\sc Meunier, F., and Sarrabezolles, P.}
\newblock Colorful linear programming, {N}ash equilibrium, and pivots.
\newblock {\em Discrete Appl. Math. 240\/} (2018), 78--91.

\bibitem{Meunier+Zerbib}
{\sc {Meunier}, F., and {Zerbib}, S.}
\newblock {Envy-free divisions of a partially burnt cake}.
\newblock {\em \rm ArXiv:1804.00449\/} (2018).

\bibitem{MT90}
{\sc Miller, G.~L., and Thurston, W.~P.}
\newblock Separators in two and three dimensions.
\newblock In {\em Proceedings {ACM} Symposium on Theory of Computing
  ({STOC})\/} (1990), pp.~300--309.

\bibitem{MRRSSSS01}
{\sc Miller, K., Ramaswami, S., Rousseeuw, P., Sellar{\`{e}}s, J., Souvaine,
  D.~L., Streinu, I., and Struyf, A.}
\newblock Fast implementation of depth contours using topological sweep.
\newblock In {\em Proceedings {ACM-SIAM} Symposium on Discrete Algorithms
  ({SODA})\/} (2001), pp.~690--699.

\bibitem{mirza+vondrak}
{\sc Mirzakhani, M., and Vondr{\'{a}}k, J.}
\newblock Sperner's colorings, hypergraph labeling problems and fair division.
\newblock In {\em Proceedings {ACM-SIAM} Symposium on Discrete
  Algorithms,{(SODA)}\/} (2015), pp.~873--886.

\bibitem{MR3726616}
{\sc Mirzakhani, M., and Vondr\'ak, J.}
\newblock {\em Sperner's colorings and optimal partitioning of the simplex}.
\newblock Springer International Publishing, Cham, 2017, pp.~615--631.
\newblock chapter of ``A Journey Through Discrete Math.: A Tribute to
  Ji{\v{r}}{\'i} Matou{\v{s}}ek'', edited by M. Loebl, J. Ne{\v{s}}et{\v{r}}il,
  and R. Thomas.

\bibitem{Mi02}
{\sc Mizera, I.}
\newblock On depth and deep points: A calculus.
\newblock {\em Ann. Statist. 30}, 6 (2002), 1681--1736.

\bibitem{motzkin-double}
{\sc Motzkin, T.~S., Raiffa, H., Thompson, G.~L., and Thrall, R.}
\newblock The double description method.
\newblock In {\em Contributions to the theory of games II\/} (1953), H.~Kuhn
  and A.~Tucker, Eds., vol.~28 of {\em Annals of Mathematical Studies},
  Princeton University Press.

\bibitem{Mulzer:2013je}
{\sc Mulzer, W., and Stein, Y.}
\newblock {Algorithms for tolerated {T}verberg partitions}.
\newblock In {\em Algorithms and Computation}. Springer Berlin Heidelberg,
  Berlin, Heidelberg, Jan. 2013, pp.~295--305.

\bibitem{mulzerstein2015}
{\sc Mulzer, W., and Stein, Y.}
\newblock {Computational Aspects of the Colorful Carath{\'{e}}odory Theorem}.
\newblock In {\em Proceedings Symposium on Computational Geometry, {(SoCG)}\/}
  (2015), pp.~44--58.

\bibitem{MW13}
{\sc Mulzer, W., and Werner, D.}
\newblock Approximating {T}verberg points in linear time for any fixed
  dimension.
\newblock {\em Discrete Comput. Geom. 50}, 2 (2013), 520--535.

\bibitem{munkres1984elements}
{\sc Munkres, J.~R.}
\newblock {\em Elements of algebraic topology}.
\newblock Addison-Wesley Menlo Park, 1984.

\bibitem{Musin:2012tq}
{\sc Musin, O.~R.}
\newblock {Extensions of {S}perner and {T}ucker's lemma for manifolds}.
\newblock {\em J. Combin. Theory Ser. A 132\/} (2015), 172--187.

\bibitem{Musin:2016}
{\sc Musin, O.~R.}
\newblock Homotopy invariants of covers and {KKM}-type lemmas.
\newblock {\em Algebr. Geom. Topol. 16}, 3 (2016), 1799--1812.

\bibitem{Musin:2017}
{\sc Musin, O.~R.}
\newblock K{KM} type theorems with boundary conditions.
\newblock {\em J. Fixed Point Theory Appl. 19}, 3 (2017), 2037--2049.

\bibitem{MDG16}
{\sc Mustafa, N.~H., Dutta, K., and Ghosh, A.}
\newblock {A Simple Proof of Optimal Epsilon-nets}.
\newblock {\em Combinatorica\/} (2017).
\newblock \url{https://doi.org/10.1007/s00493-017-3564-5}.

\bibitem{MR08}
{\sc Mustafa, N.~H., and Ray, S.}
\newblock Weak $\eps$-nets have a basis of size $o(1/\eps \log 1/\eps)$.
\newblock {\em Comput. Geom. 40}, 1 (2008), 84--91.

\bibitem{MR09b}
{\sc Mustafa, N.~H., and Ray, S.}
\newblock An optimal extension of the centerpoint theorem.
\newblock {\em {Comput. Geom.} 42}, 6 (2009), 505--510.

\bibitem{MR16}
{\sc Mustafa, N.~H., and Ray, S.}
\newblock An optimal generalization of the colorful {C}arath{\'{e}}odory
  theorem.
\newblock {\em Discrete Math. 339}, 4 (2016), 1300--1305.

\bibitem{MRS11}
{\sc Mustafa, N.~H., Ray, S., and Shabbir, M.}
\newblock Ray-shooting depth: Computing statistical data depth of point sets in
  the plane.
\newblock In {\em Proceedings European Symposium on Algorithms ({ESA})\/}
  (2011), pp.~506--517.

\bibitem{MRS15}
{\sc Mustafa, N.~H., Ray, S., and Shabbir, M.}
\newblock $k$-centerpoints conjectures for pointsets in $\mathbb{R}^d$.
\newblock {\em Internat. J. Comput. Geom. Appl. 25}, 3 (2015), 163--186.

\bibitem{MTW14}
{\sc Mustafa, N.~H., Tiwary, H.~R., and Werner, D.}
\newblock A proof of the {O}ja depth conjecture in the plane.
\newblock {\em Comput. Geom. 47}, 6 (2014), 668--674.

\bibitem{MV16}
{\sc Mustafa, N.~H., and Varadarajan, K.}
\newblock Epsilon-approximations and epsilon-nets.
\newblock In {\em Handbook of Discrete Comput. Geom.}, J.~E. Goodman,
  J.~O'Rourke, and C.~D. T\'oth, Eds. CRC Press LLC, 2017.

\bibitem{nashpaper}
{\sc Nash, Jr., J.~F.}
\newblock Equilibrium points in {$n$}-person games.
\newblock {\em Proc. Natl. Acad. Sci. USA 36\/} (1950), 48--49.

\bibitem{nashpaper2}
{\sc Nash, Jr., J.~F.}
\newblock Non-cooperative games.
\newblock {\em Ann. of Math. 54\/} (1951), 286--295.

\bibitem{nashhex}
{\sc Nash, Jr., J.~F.}
\newblock Some games and machines for playing them.
\newblock Tech. Rep. D-1164, Rand Corporation, 1952.

\bibitem{Naszodi2016}
{\sc Nasz\'{o}di, M.}
\newblock Proof of a conjecture of {B}\'{a}r\'{a}ny, {K}atchalski and {P}ach.
\newblock {\em Discrete Comput. Geom. 55}, 1 (2016), 243--248.

\bibitem{Nisanetal-2007}
{\sc Nisan, N., Roughgarden, T., Tardos, E., and Vazirani, V.~V.}
\newblock {\em Algorithmic Game Theory}.
\newblock Cambridge University Press, New York, NY, USA, 2007.

\bibitem{N62}
{\sc Novikoff, A.~B.}
\newblock On convergence proofs on perceptrons.
\newblock In {\em Proceedings of the Symposium on the Mathematical Theory of
  Automata\/} (1962), pp.~615--622.

\bibitem{Nyman+Su2013}
{\sc Nyman, K.~L., and Su, F.~E.}
\newblock A {B}orsuk-{U}lam equivalent that directly implies {S}perner's lemma.
\newblock {\em Amer. Math. Monthly 120}, 4 (2013), 346--354.

\bibitem{O83}
{\sc Oja, H.}
\newblock Descriptive statistics for multivariate distributions.
\newblock {\em Statist. Probab. Lett. 1\/} (1983), 327--332.

\bibitem{onn+radon}
{\sc Onn, S.}
\newblock On the geometry and computational complexity of {R}adon partitions in
  the integer lattice.
\newblock {\em SIAM J. Discrete Math. 4}, 3 (1991), 436--446.

\bibitem{shmuelbook}
{\sc Onn, S.}
\newblock {\em Nonlinear discrete optimization, an algorithmic theory}.
\newblock Zurich Lectures in Advanced Mathematics. European Mathematical
  Society (EMS), Z\"urich, 2010.

\bibitem{Oza87}
{\sc {\"O}zaydin, M.}
\newblock {Equivariant maps for the symmetric group}.
\newblock unpublished preprint, University of Winsconsin-Madison, 17 pages,
  available at http://digital.library.wisc.edu/1793/63829, 1987.

\bibitem{PA95}
{\sc Pach, J., and Agarwal, P.~K.}
\newblock {\em Combinatorial Geometry}.
\newblock John Wiley \& Sons, 1995.

\bibitem{palvolgyi20092d}
{\sc P{\'a}lv{\"o}lgyi, D.}
\newblock $2d$-{T}ucker is {PPAD}-complete.
\newblock In {\em International Workshop on Internet and Network Economics\/}
  (2009), Springer, pp.~569--574.

\bibitem{Pal09}
{\sc P{\'a}lv{\"o}lgyi, D.}
\newblock Combinatorial necklace splitting.
\newblock {\em Electron. J. Combin. 16}, 1 (2009), Research Paper 79, 8.

\bibitem{ppad-original}
{\sc Papadimitriou, C.~H.}
\newblock On the complexity of the parity argument and other inefficient proofs
  of existence.
\newblock {\em J. Comput. Syst. Sci. 48}, 3 (1994), 498--532.

\bibitem{pinchasi-fracHelly}
{\sc Pinchasi, R.}
\newblock A note on smaller fractional {H}elly numbers.
\newblock {\em Discrete Comput. Geom. 54}, 3 (2015), 663--668.

\bibitem{pisier1980remarques}
{\sc Pisier, G.}
\newblock Remarques sur un r{\'e}sultat non publi{\'e} de {B.} {M}aurey.
\newblock {\em S{\'e}minaire Analyse fonctionnelle\/} (1980), 1--12.

\bibitem{Por-2018}
{\sc {Por}, A.}
\newblock {Universality of vector sequences and universality of Tverberg
  partitions}.
\newblock {\em \rm ArXiv:1805.07197\/} (2018).

\bibitem{prescott+su}
{\sc Prescott, T., and Su, F.~E.}
\newblock A constructive proof of {K}y {F}an's generalization of {T}ucker's
  lemma.
\newblock {\em J. Combin. Theory Ser. {A} 111}, 2 (2005), 257--265.

\bibitem{queyranne+tardella}
{\sc Queyranne, M., and Tardella, F.}
\newblock {C}arath\'eodory, {H}elly, and {R}adon numbers for sublattice and
  related convexities.
\newblock {\em Math. Oper. Res. 42}, 2 (2017), 495--516.

\bibitem{R47}
{\sc Rado, R.}
\newblock A theorem on general measure.
\newblock {\em J. Lond. Math. Soc. 21\/} (1947), 291--300.

\bibitem{originalRadon}
{\sc Radon, J.}
\newblock Mengen konvexer {K}\"orper, die einen gemeinsamen {P}unkt enthalten.
\newblock {\em Math. Ann. 83}, 1-2 (1921), 113--115.

\bibitem{Reay-mem}
{\sc Reay, J.}
\newblock Generalizations of a theorem of {C}arath\'eodory.
\newblock {\em Memoirs of the American Mathematical Society 54\/} (1965), 50
  pages.

\bibitem{Ric53}
{\sc Richardson, M.}
\newblock Solutions of irreflexive relations.
\newblock {\em Ann. of Math. 58\/} (1953), 573--590.

\bibitem{roberson+webb}
{\sc Robertson, J., and Webb, W.}
\newblock {\em Cake-cutting Algorithms: Be fair if you can}.
\newblock AK Peters/CRC Press, 1998.

\bibitem{rolnick+soberon2}
{\sc Rolnick, D., and Sober{\'{o}}n, P.}
\newblock Algorithmic aspects of {T}verberg's theorem.
\newblock {\em \rm ArXiv:1601.03083\/} (2016).

\bibitem{Rolnick+Soberon}
{\sc Rolnick, D., and Sober\'on, P.}
\newblock Quantitative $(p,q)$-theorems in combinatorial geometry.
\newblock {\em Discrete Math. 340}, 10 (2017), 2516 -- 2527.

\bibitem{ROO01}
{\sc Ronkainen, T., Oja, H., and Orponen, P.}
\newblock Computation of the multivariate {O}ja median.
\newblock In {\em R. Dutter and P. Filzmoser}. International Conference on
  Robust Statistics, 2003.

\bibitem{rothe-surveys}
{\sc Rothe, J.}, Ed.
\newblock {\em Economics and computation}.
\newblock Springer Texts in Business and Economics. Springer, Heidelberg, 2016.
\newblock An introduction to algorithmic game theory, computational social
  choice, and fair division, With illustrations by Irene Rothe.

\bibitem{Roudneff-ejc}
{\sc Roudneff, J.-P.}
\newblock Partitions of points into simplices with {$k$}-dimensional
  intersection. {I}. {T}he conic {T}verberg's theorem.
\newblock {\em European J. Combin. 22}, 5 (2001), 733--743.
\newblock Combinatorial geometries (Luminy, 1999).

\bibitem{RR96}
{\sc Rousseeuw, P., and Ruts, I.}
\newblock Algorithm {AS} 307: Bivariate location depth.
\newblock {\em J. R. Stat. Soc. Ser. C. Appl. Stat. 45\/} (1996), 516--526.

\bibitem{RH99}
{\sc Rousseeuw, P.~J., and Hubert, M.}
\newblock Regression depth.
\newblock {\em J. Amer. Statist. Assoc. 94}, 446 (1999), 388--402.

\bibitem{rubin18}
{\sc Rubin, N.}
\newblock An improved bound for weak epsilon-nets in the plane.
\newblock In {\em Proceedings {IEEE} Symposium on Foundations of Computer
  Science {(FOCS})\/} (2018).

\bibitem{R67}
{\sc Ryser, H.~J.}
\newblock Neuere {P}robleme der {K}ombinatorik.
\newblock In {\em Vortr{\"a}ge {\"u}ber Kombinatorik\/} (July 1967),
  Oberwolfach, Mathematisches Forschunginstitute.

\bibitem{Sarkaria:1992vt}
{\sc Sarkaria, K.~S.}
\newblock {Tverberg{\textquoteright}s theorem via number fields}.
\newblock {\em Israel J. Math. 79}, 2 (1992), 317--320.

\bibitem{Sarkaria-TT}
{\sc Sarkaria, K.~S.}
\newblock {Tverberg partitions and Borsuk-Ulam theorems}.
\newblock {\em Pacific J. Math. 196}, 1 (2000), 231--241.

\bibitem{S15}
{\sc Sarrabezolles, P.}
\newblock The colourful simplicial depth conjecture.
\newblock {\em J. Comb. Theory, Ser. {A} 130\/} (2015), 119--128.

\bibitem{savani-vonstengel}
{\sc Savani, R., and von Stengel, B.}
\newblock Exponentially many steps for finding a {N}ash equilibrium in a
  bimatrix game.
\newblock In {\em Proceedings {IEEE} Symposium on Foundations of Computer
  Science {(FOCS)}\/} (2004), pp.~258--267.

\bibitem{scarf1967}
{\sc Scarf, H.~E.}
\newblock The approximation of fixed points of a continuous mapping.
\newblock {\em SIAM J. Appl. Math. 15\/} (1967), 1328--1343.

\bibitem{game-scarf}
{\sc Scarf, H.~E.}
\newblock The core of an {$N$} person game.
\newblock {\em Econometrica 35\/} (1967), 50--69.

\bibitem{Sca1977}
{\sc Scarf, H.~E.}
\newblock An observation on the structure of production sets with
  indivisibilities.
\newblock {\em Proc. Natl. Acad. Sci. USA 74}, 9 (1977), 3637--3641.

\bibitem{Schaefer2015}
{\sc Schaefer, M., and {\v{S}}tefankovi{\v{c}}, D.}
\newblock Fixed points, {N}ash equilibria, and the existential theory of the
  reals.
\newblock {\em Theory Comput. Syst.\/} (2015), 1--22.

\bibitem{topotverbergwinding}
{\sc Sch{\"{o}}neborn, T., and Ziegler, G.~M.}
\newblock The topological {T}verberg theorem and winding numbers.
\newblock {\em J. Comb. Theory, Ser. {A} 112}, 1 (2005), 82--104.

\bibitem{schgraph}
{\sc Schrijver, A.}
\newblock Vertex-critical subgraphs of {K}neser graphs.
\newblock {\em Nieuw Arch. Wiskd, III. Ser. 26\/} (1978), 454--461.

\bibitem{Sch86}
{\sc Schrijver, A.}
\newblock {\em Theory of linear and integer programming}.
\newblock Wiley, 1986.

\bibitem{Sch03}
{\sc Schrijver, A.}
\newblock {\em Combinatorial Optimization: polyhedra and efficiency}.
\newblock Springer, 2003.

\bibitem{simplotopes}
{\sc Seacrest, T., and Su, F.}
\newblock A lower bound technique for triangulations of simplotopes.
\newblock {\em SIAM J. Discrete Math. 32}, 1 (2018), 1--28.

\bibitem{sebocaratheodory}
{\sc Seb\H{o}, A.}
\newblock Hilbert bases, {C}arath\'eodory's theorem and combinatorial
  optimization.
\newblock In {\em Integer programming and combinatorial optimization (Waterloo,
  1990)}. Univ. of Waterloo Press, 1990, pp.~431--455.

\bibitem{Segal-Halevi}
{\sc {Segal-Halevi}, E.}
\newblock Fairly dividing a cake after some parts were burnt in the oven.
\newblock {\em \rm ArXiv:1704.00726\/} (2017).

\bibitem{seidel1991small}
{\sc Seidel, R.}
\newblock Small-dimensional linear programming and convex hulls made easy.
\newblock {\em Discrete Comput. Geom. 6}, 1 (1991), 423--434.

\bibitem{seymourTU}
{\sc Seymour, P.~D.}
\newblock Decomposition of regular matroids.
\newblock {\em J. Combin. Theory Ser. B 28}, 3 (1980), 305--359.

\bibitem{shapirodentchevarusz}
{\sc Shapiro, A., Dentcheva, D., and Ruszczy{\'n}ski, A.}
\newblock {\em Lectures on stochastic programming. Modeling and theory},
  2nd~ed.
\newblock MOS-SIAM Series on Optimization. Society for Industrial and Applied
  Mathematics (SIAM), Philadelphia, PA; Mathematical Optimization Society,
  Philadelphia, PA, 2014.

\bibitem{siegel}
{\sc Siegel, C.~L.}
\newblock { \"Uber einige Anwendungen diophantischer Approximationen}.
\newblock {\em Abh. der Preus. Akad. der Wissenschaften. Phys.--math.}, 1
  (1929), 209--266.

\bibitem{simmons+su}
{\sc Simmons, F., and Su, F.~E.}
\newblock Consensus-halving via theorems of {B}orsuk-{U}lam and {T}ucker.
\newblock {\em Math. Social Sci. 45}, 1 (2003), 15--25.

\bibitem{SimonS}
{\sc {Simon}, S.}
\newblock {Hyperplane Equipartitions Plus Constraints}.
\newblock {\em \rm ArXiv:1708.00527\/} (2017).

\bibitem{smale1998mathematical}
{\sc Smale, S.}
\newblock Mathematical problems for the next century.
\newblock {\em Math. Intelligencer 20}, 2 (1998), 7--15.

\bibitem{soberon-gerrymanderingetc}
{\sc Sober\'on, P.}
\newblock Gerrymandering, sandwiches, and topology.
\newblock {\em Notices Amer. Math. Soc. 64}, 9 (2017), 1010--1013.

\bibitem{pablo2018}
{\sc Sober\'on, P.}
\newblock Robust {T}verberg and colourful {C}arath\'eodory results via random
  choice.
\newblock {\em Combin. Probab. Comput. 27}, 3 (2018), 427--440.

\bibitem{Soberon:2012er}
{\sc Sober{\'o}n, P., and Strausz, R.}
\newblock {A generalisation of {T}verberg's theorem}.
\newblock {\em Discrete Comput. Geom. 47\/} (2012), 455--460.

\bibitem{spencer+su}
{\sc Spencer, G., and Su, F.~E.}
\newblock The {LSB} theorem implies the {KKM} lemma.
\newblock {\em Amer. Math. Monthly 114}, 2 (2007), 156--159.

\bibitem{stahl1976n}
{\sc Stahl, S.}
\newblock $n$-tuple colorings and associated graphs.
\newblock {\em J. Combin. Theory Ser. B 20}, 2 (1976), 185--203.

\bibitem{Hex-Tucker}
{\sc Stehlik, M.}
\newblock Personal communication.

\bibitem{S75}
{\sc Stein, S.~K.}
\newblock Transversals of {L}atin squares and their generalizations.
\newblock {\em Pacific J. Math. 59\/} (1975), 567--575.

\bibitem{Steinhaus}
{\sc Steinhaus, H.}
\newblock Sur la division pragmatique.
\newblock {\em Econometrica 17}, (Supplement) (1949), 315--319.

\bibitem{originalsteinitz}
{\sc Steinitz, E.}
\newblock Bedingt konvergente {R}eihen und konvexe {S}ysteme. i-ii-iii.
\newblock {\em J. Reine Angew. Math. 143\/} (1913), 128--175.

\bibitem{ST42}
{\sc Stone, A., and Tukey, J.}
\newblock Generalized ``sandwich'' theorems.
\newblock {\em Duke Math. J. 9\/} (1942), 356--359.

\bibitem{stromquist1}
{\sc Stromquist, W.}
\newblock How to cut a cake fairly.
\newblock {\em Amer. Math. Monthly 87}, 8 (1980), 640--644.

\bibitem{stromquist2}
{\sc Stromquist, W.}
\newblock Envy-free cake divisions cannot be found by finite protocols.
\newblock {\em Electron. J. Combin. 15}, 1 (2008), Research Paper 11, 10.

\bibitem{Su99}
{\sc Su, F.}
\newblock Rental harmony: {S}perner's lemma in fair division.
\newblock {\em Amer. Math. Monthly 106\/} (1999), 930--942.

\bibitem{su-hex}
{\sc Su, F.}
\newblock A rectangular {S}perner’s lemma implies the {H}ex theorem, 2018.
\newblock working manuscript.

\bibitem{NYT-sperner}
{\sc Sun, A.}
\newblock To divide the rent, start with a triangle.
\newblock {\em New York Times April 29\/} (2014).
\newblock
  \url{https://www.nytimes.com/2014/04/29/science/to-divide-the-rent-start-with-a-triangle.html?_r=1}.

\bibitem{ST83}
{\sc Szemer{\'e}di, E., and Trotter, W.}
\newblock Extremal problems in discrete geometry.
\newblock {\em Combinatorica 3}, 3 (1983), 381--392.

\bibitem{tancer2013intersection}
{\sc Tancer, M.}
\newblock Intersection patterns of convex sets via simplicial complexes: A
  survey.
\newblock In {\em Thirty essays on geometric graph theory}. Springer, 2013,
  pp.~521--540.

\bibitem{todd1977number}
{\sc Todd, M.~J.}
\newblock The number of necessary constraints in an integer program: a new
  proof of {S}carf's theorem.
\newblock Tech. Rep. 355, Cornell University, 1977.

\bibitem{todd1993new}
{\sc Todd, M.~J., and Tun{\c{c}}el, L.}
\newblock A new triangulation for simplicial algorithms.
\newblock {\em SIAM J. Discrete Math. 6}, 1 (1993), 167--180.

\bibitem{T75}
{\sc Tukey, J.}
\newblock Mathematics and the picturing of data.
\newblock In {\em Proceedings {I}nternational {C}ongress of {M}athematicians\/}
  (1975), pp.~523--531.

\bibitem{Tverberg:1966tb}
{\sc Tverberg, H.}
\newblock {A generalization of {R}adon{\textquoteright}s theorem}.
\newblock {\em J. Lond. Math. Soc. 41}, 1 (1966), 123--128.

\bibitem{Tverberg1981}
{\sc Tverberg, H.}
\newblock A generalization of {R}adon's theorem. {II}.
\newblock {\em Bull. Austral. Math. Soc. 24}, 3 (1981), 321--325.

\bibitem{TverbergVrecica}
{\sc Tverberg, H., and Vre{\'c}ica, S.}
\newblock On generalizations of {R}adon's theorem and the ham sandwich theorem.
\newblock {\em European J. Combin. 14}, 3 (1993), 259 -- 264.

\bibitem{vaalerbest}
{\sc Vaaler, J.}
\newblock { The best constant in Siegel's lemma}.
\newblock {\em Monatsh. Math. 140\/} (2003), 71--89.

\bibitem{convexityspaces-vandevel}
{\sc van~de Vel, M.}
\newblock {\em Theory of convex structures}, vol.~50 of {\em North-Holland
  Mathematical Library}.
\newblock North-Holland Publishing Co., Amsterdam, 1993.

\bibitem{vanKampen:KomplexeInEuklidischenRaeumen-1932}
{\sc van Kampen, E.~R.}
\newblock {Komplexe in {E}uklidischen R\"aumen}.
\newblock {\em Abh. Math. Sem. Univ. Hamburg 9\/} (1932), 72--78.

\bibitem{MMRSSS08}
{\sc van Kreveld, M.~J., Mitchell, J. S.~B., Rousseeuw, P., Sharir, M.,
  Snoeyink, J., and Speckmann, B.}
\newblock Efficient algorithms for maximum regression depth.
\newblock {\em Discrete Comput. Geom. 39}, 4 (2008), 656--677.

\bibitem{VC71}
{\sc Vapnik, V.~N., and Chervonenkis, A.~Y.}
\newblock On the uniform convergence of relative frequencies of events to their
  probabilities.
\newblock {\em Theory Probab. Appl. 16}, 2 (1971), 264--280.

\bibitem{V10}
{\sc Varadarajan, K.}
\newblock Weighted geometric set cover via quasi uniform sampling.
\newblock In {\em Proceedings {ACM} Symposium on Theory of Computing
  ({STOC})\/} (2010), pp.~641--648.

\bibitem{V01}
{\sc Vazirani, V.}
\newblock {\em Approximation algorithms}.
\newblock Springer, Berlin, New York, Paris, 2001.

\bibitem{Volovikov-TT}
{\sc Volovikov, A.}
\newblock On a topological generalization of the {T}verberg theorem.
\newblock {\em Math. Notes 59}, 3 (1996), 324--326.

\bibitem{VoNMo44}
{\sc Von~Neumann, J., and Morgenstern, O.}
\newblock {\em Theory of Games and Economic Behavior}.
\newblock Princeton University Press, 1944.

\bibitem{vucicrade}
{\sc Vu{\v{c}}i{\'{c}}, A., and {\v{Z}}ivaljevi{\'{c}}, R.~T.}
\newblock Note on a conjecture of {S}ierksma.
\newblock {\em Discrete Comput. Geom. 9}, 4 (1993), 339--349.

\bibitem{W03}
{\sc Wagner, U.}
\newblock {\em On $k$-Sets and Applications}.
\newblock PhD thesis, ETH Zurich, 2003.

\bibitem{truemper+walter}
{\sc Walter, M., and Truemper, K.}
\newblock Implementation of a unimodularity test.
\newblock {\em Math. Program. Comput. 5}, 1 (2013), 57--73.

\bibitem{WS11}
{\sc Williamson, D.~P., and Shmoys, D.~B.}
\newblock {\em {The Design of Approximation Algorithms}}.
\newblock Cambridge University Press, 2011.

\bibitem{woodall}
{\sc Woodall, D.~R.}
\newblock Dividing a cake fairly.
\newblock {\em J. Math. Anal. Appl. 78}, 1 (1980), 233--247.

\bibitem{YaoYao85}
{\sc Yao, A.~C., and Yao, F.~F.}
\newblock A general approach to $d$-dimensional geometric queries.
\newblock In {\em Proceedings {ACM} Symposium on Theory of
  Computing,{(STOC)}\/} (1985), pp.~163--168.

\bibitem{zhu2006recent}
{\sc Zhu, X.}
\newblock Recent developments in circular colouring of graphs.
\newblock In {\em Topics in discrete mathematics}. Springer, 2006,
  pp.~497--550.

\bibitem{Z-GK}
{\sc Ziegler, G.~M.}
\newblock Generalized {K}neser coloring theorems with combinatorial proofs.
\newblock {\em Invent. Math. 147}, 3 (2002), 671--691.

\bibitem{3nziegler}
{\sc Ziegler, G.~M.}
\newblock 3{N} colored points in a plane.
\newblock {\em Notices Amer. Math. Soc. 58}, 4 (2011), 550--557.

\bibitem{MR2728499}
{\sc {\v{Z}}ivaljevi\'c, R.~T.}
\newblock Oriented matroids and {K}y {F}an's theorem.
\newblock {\em Combinatorica 30}, 4 (2010), 471--484.

\bibitem{Zivaljevic:Handbook}
{\sc {\v Z}ivaljevi{\'c}, R.~T.}
\newblock {\em Topological methods in discrete geometry, Chapter 21, in}.
\newblock CRC Press LLC, Boca Raton FL, 2017, pp.~551--580.
\newblock edited by J.E. Goodman, J. O'Rourke, and C. D. T\'oth.

\bibitem{colortv-ZivaljevicV92}
{\sc {\v Z}ivaljevi{\'c}, R.~T., and Vre\'{c}ica, S.~T.}
\newblock The colored {T}verberg's problem and complexes of injective
  functions.
\newblock {\em J. Comb. Theory, Ser. {A} 61}, 2 (1992), 309--318.

\end{thebibliography}

\end{document}